\documentclass[leqno,10pt]{amsart}

\PassOptionsToPackage{dvipsnames}{xcolor}

\usepackage[utf8]{inputenc}
\usepackage[english]{babel}
\usepackage{epigraph}
\usepackage{mathtools}  
\usepackage{amssymb, amsthm, mathrsfs}
\usepackage{braket}     
\usepackage{esint}      
\usepackage{xfrac, nicefrac}

\usepackage{graphicx}
\usepackage{color}

\usepackage{enumitem}   
\usepackage[toc]{appendix}
\usepackage{csquotes}   
\usepackage[colorlinks, pdfpagelabels, pdfstartview=FitH, bookmarksopen=true, bookmarksnumbered = true, linkcolor = NavyBlue, plainpages = false, hypertexnames = false, citecolor = OliveGreen] {hyperref}   
\usepackage{orcidlink} 

\usepackage[a4paper, margin=2.5cm]{geometry}
\allowdisplaybreaks

\numberwithin{equation}{section}

\newcommand{\R}{\mathbb{R}} \newcommand{\N}{\mathbb{N}} 
\newcommand{\dx}{\,dx}

\newcommand{\vs}{\medskip}

\def\eps{\varepsilon}
\newcommand{\x}{\textnormal{\texttt{x}}}
\newcommand{\bI}{{\mathbf I}}
\newcommand{\bi}{\mathfrak{b}}

\newcommand{\ratio}{L/\nu}

\newcommand{\smallzero}{{\scalebox{.8}{$\scriptscriptstyle 0$}}}
\newcommand{\smallone}{{\scalebox{.8}{$\scriptscriptstyle 1$}}}
\newcommand{\smalltwo}{{\scalebox{.8}{$\scriptscriptstyle 2$}}}
\newcommand{\smallthree}{{\scalebox{.8}{$\scriptscriptstyle 3$}}}

\newcommand{\smallxi}{{\scalebox{.8}{$\scriptscriptstyle \xi$}}}
\newcommand{\smallH}{{\scalebox{.8}{$\scriptscriptstyle H$}}}

\DeclareMathOperator{\divergence}{div}
\DeclareMathOperator*{\osc}{osc}

\def\mean#1{\mathchoice%
          {\mathop{\kern 0.2em\vrule width 0.6em height 0.69678ex depth -0.58065ex
                  \kern -0.8em \intop}\nolimits_{\kern -0.4em#1}}%
          {\mathop{\kern 0.1em\vrule width 0.5em height 0.69678ex depth -0.60387ex
                  \kern -0.6em \intop}\nolimits_{#1}}%
          {\mathop{\kern 0.1em\vrule width 0.5em height 0.69678ex
              depth -0.60387ex
                  \kern -0.6em \intop}\nolimits_{#1}}%
          {\mathop{\kern 0.1em\vrule width 0.5em height 0.69678ex depth -0.60387ex
                  \kern -0.6em \intop}\nolimits_{#1}}}
                  
\newcommand{\aveint}[2]{\mathchoice
  {\mathop{\kern 0.2em\vrule width 0.6em height 0.7ex depth -0.6ex \kern -0.8em \intop}\nolimits_{\kern -0.45em#1}^{#2}}
  {\mathop{\kern 0.1em\vrule width 0.5em height 0.7ex depth -0.6ex \kern -0.6em \intop}\nolimits_{#1}^{#2}}
  {\mathop{\kern 0.1em\vrule width 0.5em height 0.7ex depth -0.6ex \kern -0.6em \intop}\nolimits_{#1}^{#2}}
  {\mathop{\kern 0.1em\vrule width 0.5em height 0.7ex depth -0.6ex \kern -0.6em \intop}\nolimits_{#1}^{#2}}}

\makeatletter
\DeclareFontFamily{OMX}{MnSymbolE}{}
\DeclareSymbolFont{MnLargeSymbols}{OMX}{MnSymbolE}{m}{n}
\SetSymbolFont{MnLargeSymbols}{bold}{OMX}{MnSymbolE}{b}{n}
\DeclareFontShape{OMX}{MnSymbolE}{m}{n}{
    <-6>  MnSymbolE5
   <6-7>  MnSymbolE6
   <7-8>  MnSymbolE7
   <8-9>  MnSymbolE8
   <9-10> MnSymbolE9
  <10-12> MnSymbolE10
  <12->   MnSymbolE12
}{}
\DeclareFontShape{OMX}{MnSymbolE}{b}{n}{
    <-6>  MnSymbolE-Bold5
   <6-7>  MnSymbolE-Bold6
   <7-8>  MnSymbolE-Bold7
   <8-9>  MnSymbolE-Bold8
   <9-10> MnSymbolE-Bold9
  <10-12> MnSymbolE-Bold10
  <12->   MnSymbolE-Bold12
}{}

\let\llangle\@undefined
\let\rrangle\@undefined
\DeclareMathDelimiter{\llangle}{\mathopen}%
                     {MnLargeSymbols}{'164}{MnLargeSymbols}{'164}
\DeclareMathDelimiter{\rrangle}{\mathclose}%
                     {MnLargeSymbols}{'171}{MnLargeSymbols}{'171}
\makeatother

\newcommand{\linethrough}{\mathpalette\@thickbar}
\makeatletter
\newcommand{\@thickbar}[2]{{#1\mkern0mu\vbox{
  \sbox\z@{$#1#2\mkern-0.5mu$}%
  \dimen@=\dimexpr\ht\tw@-\ht\z@+2\p@\relax
  \hrule\@height0.5\p@
  \vskip\dimen@ \box\z@}}}
\makeatother

\newcommand{\data}{\textnormal{\texttt{data}}}
\newcommand{\gdata}{\textnormal{\texttt{gdata}}}
\newcommand{\esp}{\textnormal{\texttt{t}}}
\newcommand{\es}{\textnormal{\texttt{s}}}
\newcommand{\ti}{\mathfrak{t}}
\newcommand{\el}{\mathfrak{l}}
\newcommand{\benergy}{\texttt{En}_b}

\makeatletter
\newcommand{\pushright}[1]{\ifmeasuring@#1\else\omit\hfill$\displaystyle#1$\fi\ignorespaces}
\newcommand{\pushleft}[1]{\ifmeasuring@#1\else\omit$\displaystyle#1$\hfill\fi\ignorespaces}
\makeatother

\newtheorem{theorem}{Theorem}[section]
\newtheorem*{theorem*}{Theorem}
\newtheorem{lemma}[theorem]{Lemma}
\newtheorem{proposition}[theorem]{Proposition}
\newtheorem*{proposition*}{Proposition}
\newtheorem{corollary}[theorem]{Corollary}
\newtheorem{definition}[theorem]{Definition}
\newtheorem{remark}[theorem]{Remark}
\newtheorem*{remark*}{Remark}

\title[{Potential theory for double phase equations with measure data}]{Potential theory for double phase equations\\ with measure data}

\author[Baroni]{Paolo Baroni \orcidlink{0000-0002-2683-4025}}
\address{Paolo Baroni \\ Dipartimento di Scienze Matematiche, Fisiche e Informatiche, Universit\`a di Parma, Parco Area delle Scienze, 53/A, 43124 Parma, Italy\\ ORCID ID: 0000-0002-2683-4025}
\email{\url{paolo.baroni@unipr.it}}


\begin{document}

\begin{abstract}
We develop a potential theory for elliptic equations with double phase growth and signed measure data. We establish existence of solutions obtained by approximation, pointwise bounds for the flux in terms of linear Riesz potentials, and nonlinear potential estimates for the solution in terms of generalized Wolff potentials; a quantitative decay estimate for the flux excess is the main ingredient and also yields a potential criterion for gradient continuity.  The results cover both singular and degenerate equations, including the range $p>n$.
\end{abstract}

\maketitle

\tableofcontents

\epigraph{\emph{Je n'ai fait celle-ci plus longue que parce que\\ je n'ai pas eu le loisir de la faire plus courte.}}
{B.~Pascal, \emph{Les Provinciales}, XVI, 4 d\'ecembre 1656}

\vspace{5mm}

\section{Introduction}

In this paper we develop a measure-data theory for double phase equations, including linear and nonlinear potential estimates; these results complement the regularity theory for energy solutions and minimizers in \cite{BCMgen}. We consider solutions, in the sense specified in Definition~\ref{def.solution}, to the measure data problem
\begin{equation}\label{equation}
-\divergence A(x,Du)=\mu\qquad\text{in $\Omega\subset\mathbb{R}^n$, \quad $n\geq2$},
\end{equation}
where $\Omega$ is a bounded domain, $A\colon\Omega\times\mathbb{R}^n\longrightarrow\mathbb{R}^n$ is a Carath\'eodory vector field displaying double phase growth and $\mu$ is a signed Radon measure with finite total variation $|\mu|(\Omega)<\infty$. The model case for $A(\,\cdot\,)$ is
\begin{equation}\label{model.A}
A_{\texttt d \texttt p}(x,z):=\bigl[|z|^{p-2}+a(x)|z|^{q-2}\bigr]z=\frac{h(x,z)}{|z|}\,z
\end{equation}
for every $x\in\Omega$ and every $z\in\mathbb{R}^n\setminus\{0\}$, extended by continuity at $z=0$, where
\begin{equation}\label{h.def}
h(x,z):=|z|^{p-1}+a(x)|z|^{q-1};
\end{equation}
we shall refer to the vector field $A(\,\cdot\,,Du)$ as the \emph{flux} and to the scalar $h(x,z)$, with a slight abuse of notation, as the \emph{flux magnitude}: indeed \eqref{a-h.weak} will only give $|A(x,z)|\approx h(x,z)$ and not equality. We are indeed going to assume that $A$ behaves like $A_{\texttt d \texttt p}$, that is,
\begin{equation}\label{approximate.A}
A(x,z)\approx A_{\texttt d \texttt p}(x,z)
\end{equation}
see \eqref{a-h.weak}: the relation in~\eqref{approximate.A} is meant at this stage only formally, and the precise growth, coercivity, and monotonicity assumptions on the vector field, which vary with the class of results considered, but always stay consistent with the double phase structure, will be stated below in Paragraphs~\ref{regular.A} and~\ref{rough.A}. The switching coefficient $a(\,\cdot\,)$ is such that
\begin{equation}\label{a.holder}
a\in C^{0,\alpha}(\Omega)\quad \text{for some $\alpha\in(0,1]$,}\qquad\qquad a(x)\geq0\quad \text{for all $x\in\Omega$}. 
\end{equation}
Denoting by $|z|^{p-2}z$ and $a(x)|z|^{q-2}z$ the $p$- and the weighted $q$-part, respectively, the model vector field $A_{\texttt d \texttt p}(\,\cdot\,,\,\cdot\,)$ in \eqref{model.A} coincides simply with their sum, with the non-negative weight $a(\,\cdot\,)$ acting as a modulating factor for the $q$-growth. Since the operator degenerates to a pure $p$-growth regime on the zero set $\{a=0\}$, $a(\,\cdot\,)$ is appropriately referred to as the \emph{switching coefficient}. Outside this locus, comparing the magnitudes of the two components yields the intrinsic ratio $a(x)|z|^{q-p}$ for $z \neq 0$, indicating that the $q$-part takes over in the asymptotic regime $|z|\to\infty$ precisely when $q>p$ (see \cite{DeFiMi23} for a detailed comparison involving standard ellipticity ratios). Whether regularity persists across scales relies on the interplay between $p,q$, and the H\"older exponent $\alpha$ across the boundary of $\{a=0\}$ - as highlighted by a series of illuminating counterexamples~\cite{BDS2020, BS2021, BalciDieningSurnachevMeyers, DeFilippisDeFilippisPiccinini}, which establish that regularity can fail if the growth gap $q-p$ is too large compared to the H\"older exponent of $a(\,\cdot\,)$; in this light, conditions \eqref{bounds} identify the natural structural window admissible within our measure-data framework.

\medskip

The present work addresses existence below the natural energy threshold as well as potential bounds for both solutions and the flux, encompassing singular, mixed, and degenerate regimes for both $p\leq n$ and $p>n$, under the corresponding structural assumptions. The technical machinery developed hereafter rests on several foundational blocks: Theorem~\ref{existence.SOLA} establishes existence and global flux estimates for $p\leq n$, Theorem~\ref{main.thm} provides linear potential bounds for the flux, Proposition~\ref{qualitative.flux.oscillation} and Corollary~\ref{cor:flux.continuity} yield flux and gradient continuity, and Theorem~\ref{main.thm.2} provides nonlinear potential estimates for the solution. Before stating these results, we specify the assumptions on the exponents and coefficients and define the notions of solution used throughout the paper.

\subsection{Exponents and balance conditions}\label{exponents.balance}

We assume that the exponents $p$ and $q$ satisfy the relation
\begin{equation}\label{2-1/n}
 2-\frac1n<p\leq q,
\end{equation}
the restriction from below being customary already in the case where $a(\,\cdot\,)\equiv0$, which gives back the well-known $p$-Laplacian operator, see~\cite{DuMiJFA, KMGuide, KMVectorial, MinMil} and also~\cite{SongYounZatorska}: indeed, the gradient of its fundamental solution belongs to $L^1_{\rm loc}$ if and only if $p>2-1/n$. In the double phase setting, under the hypotheses of Theorem~\ref{existence.SOLA}, the Dirichlet problem with zero boundary data admits a SOLA belonging to $W^{1,1}_0(\Omega)$; in particular, the lower bound ensures that the fundamental solution, i.e., the radial distributional solution to \eqref{equation} corresponding to the Dirac delta datum and the radial switching coefficient $a(x)=|x|^\alpha$, has gradient in $L^1_{\rm loc}$, see Paragraphs~\ref{radial.model} and~\ref{gradient.asymptotics}.

\medskip 

Moreover, we impose the following structural relation between the exponents $p,q$ and the H\"older exponent $\alpha$ we are considering:
\begin{subequations}\label{bounds}
\begin{align}
\frac{q-1}{p-1} &< 1 +\frac\alpha{n-1} &&\text{if $p\leq n$,}\label{subcritical.bound}\\[1mm]
\frac{q}{p} &\leq 1+\frac\alpha n &&\text{if $p>n$.}\label{bound.supercritical}
\end{align}
\end{subequations}
Equivalently, in terms of the ratio $(q-1)/(p-1)$,
\[
\frac{q-1}{p-1}<1+\frac\alpha n\cdot\frac n{n-1}\quad\text{if $p\leq n$},\qquad\qquad
\frac{q-1}{p-1} \leq  1 +\frac\alpha n\cdot \frac p{p-1}\quad\text{if $p>n$}.
\]
The balance conditions \eqref{bounds} quantify how the H\"older regularity of $a(\,\cdot\,)$ compensates for the gap between the growth exponents $p$ and $q$; from the technical viewpoint, they allow us to control the errors due to the spatial oscillation of $a(\,\cdot\,)$ in the comparison and freezing arguments. The first condition is natural in view of the behavior of the fundamental solution: it is the condition under which the $q$-term remains of lower order along it, see Paragraph~\ref{radial.model} and Section~\ref{sec.fundamental.solution}. While the relation among $p,q$ and $\alpha$ above will remain in force for all our regularity results, the existence theorem will hold under weaker hypotheses.

\medskip

There is an evident parallel between the formal dichotomy in \eqref{bounds} and the different functional frameworks dictated by measure data in the regimes $p\leq n$ and $p>n$: when $p\leq n$, solutions need not belong to the natural energy space, which leads us to work with solutions locally constructed via approximation - a setting where condition \eqref{subcritical.bound} was first introduced by Song and Youn~\cite{SongYoun} in the degenerate range $2\leq p\leq n$, and subsequently by Song, Youn, and Zatorska-Goldstein~\cite{SongYounZatorska} in the singular case $2-1/n<p<2$, to establish local Calder\'on-Zygmund gradient estimates for SOLAs under structural assumptions on the vector field analogous to \eqref{A.strong}.

When $p>n$, finite Radon measures belong to the dual of the energy space, and the higher integrability of homogeneous comparison maps allows us to cover the borderline case of equality in \eqref{bound.supercritical}. The same endpoint threshold appears in the regularity theory for minimizers of double phase functionals~\cite{BCMU}; for a non-perturbative treatment of general $(p,q)$-growth functionals and a discussion of borderline cases in the general framework, see~\cite{DeFiMi23, DeFiMi25} and the introductions therein. In the mixed regime $p<n<q$, the exponents $p$ and $q$ lie on opposite sides of the dimension $n$. On balls away from the zero set $\{a=0\}$, the $q$-growth places measure data in the dual of the natural energy space; near $\{a=0\}$, however, the coercivity falls back to the $p$-growth regime, where general measures no longer lie in the dual space and one must work in the setting of very weak solutions. Developing a unified approach that covers both regimes within a single proof is precisely what motivates our analysis of the case $p>n$ as well.

At its heart, our approach relies on a perturbative mechanism driven by the gain of a positive power of the radius in the comparison estimates: when $p\leq n$, this gain arises directly from the strict inequality in \eqref{subcritical.bound} (see Paragraph \ref{heuristic.proofs}), whereas for $p>n$ it is secured by the higher integrability of the comparison maps, even when equality holds in \eqref{bound.supercritical}.

\medskip

Double phase problems involving a mild phase transition and measure data have seen significant progress in recent years. Their prototypical vector field is given by
\[
A_{\log}(x,z)=|z|^{p-2}z+a(x)|z|^{p-2}\log(e+|z|)z,
\]
and we refer to \cite{ByunYoun2018,ByunChoYoun2021,ByunShinSong2022} for key contributions covering linear potential estimates and fractional differentiability. Indeed, the logarithmic size of the growth gap simplifies perturbation arguments similarly to what happens for variable-exponent problems, where localization produces a growth gap of the same logarithmic order. The double phase structure considered in \eqref{equation} under the assumption \eqref{approximate.A}, conversely, exhibits a fixed power-like gap between the two polynomial growths, and the coupling with the regularity of the coefficient requires a substantially more delicate analysis.

\subsection{The radial model and the scale of the estimates}\label{radial.model}

We consider the radial fundamental solution $\mathcal G_{p,q}$ to \eqref{equation} for $A(x,z)=A_{\texttt d \texttt p}(x,z)$, with $C^{0,\alpha}$ radial weight and Dirac datum
\begin{equation}\label{radial}
a(x)=|x|^\alpha,\qquad \mu=\delta_0.
\end{equation}
By radial symmetry, the flux is a constant multiple of the outward unit normal on each sphere $\partial B_r(0)$, and integrating the equation over $B_r(0)$ shows that this multiple is $-1/[n\omega_nr^{n-1}]$. The flux magnitude therefore agrees with that of the fundamental solution $\mathcal G_p$ of the $p$-Laplacian with the same Dirac datum:
\begin{equation}\label{fundamental.flux.magnitude}
h\bigl(x,D\mathcal G_{p,q}(x)\bigr)=\frac{1}{n\omega_n |x|^{n-1}}={|D\mathcal G_p(x)|}^{p-1}.
\end{equation}
Under the balance condition \eqref{bounds}, the gradient has the same behavior near the origin as in the $p$-Laplacian case; more precisely, as shown in Section~\ref{sec.fundamental.solution},
\[
|D\mathcal G_{p,q}(x)|\approx |x|^{-\frac{n-1}{p-1}}.
\]
It follows that
\begin{equation}\label{int.sharp.Lebesgue}
D\mathcal G_{p,q}\in L^s_{\mathrm{loc}}\qquad\text{if and only if}\qquad s<\frac{n(p-1)}{n-1},
\end{equation}
with a weak-$L^{n(p-1)/(n-1)}$ bound at the endpoint. In particular, the gradient is locally integrable precisely when $p>2-1/n$, showing that this lower bound is sharp. Integrating the radial profile also gives the usual $p$-Laplacian behavior for the solution itself: $\mathcal G_{p,q}$ has a power singularity for $p<n$, a logarithmic singularity for $p=n$, and is locally H\"older continuous for $p>n$.

\subsection{Solutions and existence}\label{Par:solutions}

Under \eqref{approximate.A}, the double phase energy density related to \eqref{equation} is
\begin{equation}\label{H.def}
H(x,z):=\frac{|z|^p}{p}+a(x)\frac{|z|^q}{q}
\end{equation}
for every $x\in\Omega$ and $z\in \mathbb{R}^n$; this naturally induces the energy functional
\[
w\longmapsto\int_\Omega H(x,Dw)\dx.
\]
In the model case \eqref{model.A}, we have $A_{\texttt d\texttt p}(x,z)=\partial_z H(x,z)$, meaning that the homogeneous equation is simply its Euler-Lagrange equation. Even for the broader class of operators treated in this paper, $H(\,\cdot\,,\,\cdot\,)$ still dictates the growth and coercivity conditions in \eqref{A.weak} and \eqref{A.strong}: that makes the Musielak-Orlicz-Sobolev space
\[
W^{1,H(\cdot)}(\Omega):=\biggl\{u\in W^{1,1}(\Omega):\rho_{1,H}(u):=\int_\Omega \bigl[H(x,u)+H(x,Du)\bigr]\dx<\infty\biggr\}
\]
the right functional setting; it is usually endowed with the Luxemburg norm, defined as the infimum of all $\lambda>0$ such that the Musielak-Orlicz-Sobolev modular $\rho_{1,H}(\,\cdot\,)$ of $u/\lambda$ is at most one, see \cite[Section~6.1]{HH}. We will also use the zero-trace space, defined by
\[
W^{1,H(\cdot)}_0(\Omega):=\overline{C_c^\infty(\Omega)}^{\,W^{1,H(\cdot)}(\Omega)},
\]
where the closure is meant in the sense of normed spaces. Since $H$ satisfies a uniform $\Delta_2$ condition, a function $u\in W^{1,H(\cdot)}(\Omega)$ belongs to $W^{1,H(\cdot)}_0(\Omega)$ if and only if one can find a sequence $\{\varphi_j\}\subset C_c^\infty(\Omega)$ such that
\[
\lim_{j\to\infty}\int_\Omega\bigl[H(x,|u-\varphi_j|)+H(x,|Du-D\varphi_j|)\bigr]\dx=0,
\]
compare with \cite[Lemma~6.4.5 and Proposition~7.2.2]{HH}.

\medskip

We shall use this energy formulation for the approximating problems when $p\leq n$ and for the original measure data problem when $p>n$:
\begin{definition}[Local energy solution]\label{def.local.energy.solution}
A function $u\in W^{1,H(\cdot)}_{\rm loc}(\Omega)$ is called a \emph{local energy solution} to \eqref{equation} if
\begin{equation}\label{en.solution}
\int_\Omega \langle A(x,Du),D\varphi\rangle\dx=\int_\Omega \varphi\,d\mu
\end{equation}
for every $\varphi\in C_c^\infty(\Omega)$.
\end{definition}

\subsubsection{The case $p\leq n$} 

In this case we will mainly work under the strict, measure data bound \eqref{subcritical.bound}. The peculiar feature of \eqref{equation} under \eqref{approximate.A} when $p\leq n$ is that one cannot expect to have, in general, local energy solutions in $W^{1,H(\cdot)}_{\rm loc}(\Omega)$; see for instance our discussion on the fundamental solution in Paragraph~\ref{radial.model} and the more detailed account in Section~\ref{sec.fundamental.solution}, in particular Paragraph~\ref{gradient.asymptotics}. As customary in the theory of measure data problems in the nonlinear setting, one is therefore led to consider distributional solutions enjoying the minimal integrability level ensuring that the same distributional formulation is meaningful, the so-called very weak solutions. In our case, under the assumption \eqref{approximate.A} that is going to be clarified either in \eqref{A.strong}, \eqref{A.weak} or \eqref{A.veryweak}, this will have the following meaning: 
\begin{definition}
A function $u\in W^{1,1}_{\rm loc}(\Omega)$ is called a \emph{very weak solution} to \eqref{equation} if $A(\,\cdot\,,Du)\in L^1_{\rm loc}(\Omega;\mathbb{R}^n)$ and the distributional formulation \eqref{en.solution} is satisfied for every $\varphi\in C^\infty_c(\Omega)$.
\end{definition}
By the coercivity and growth assumptions in \eqref{A.veryweak}, the flux integrability required in the preceding definition is equivalent to $h(\,\cdot\,,Du)\in L^1_{\rm loc}(\Omega)$; under the balance condition \eqref{subcritical.bound}, the local Sobolev embedding also gives $u\in L^{q-1}_{\rm loc}(\Omega)$ and hence $h(\,\cdot\,,u)\in L^1_{\rm loc}(\Omega)$. It is convenient to introduce the notation
\[
W^{1,h(\cdot)}(\Omega):=\biggl\{u\in W^{1,1}(\Omega):\int_\Omega \bigl[h(x,u)+h(x,Du)\bigr]\dx<\infty\biggr\}
\]
together with the corresponding local variant $W^{1,h(\cdot)}_{\rm loc}(\Omega)$, consisting of those functions belonging to $W^{1,h(\cdot)}(U)$ for every open set $U\Subset\Omega$. The requirement $u\in W^{1,1}(\Omega)$ (or $u\in W^{1,1}_{\rm loc}(\Omega)$) is explicitly assumed since, when $p<2$, the integrability of $h(x,Du)$ alone is not sufficient for the corresponding $L^1$ integrability of $Du$, and $z\mapsto h(x,z)$ is not a Young function.

As a measure datum does not, in general, belong to the dual of the energy space, standard monotonicity methods do not apply directly to \eqref{equation}, and the construction of very weak solutions requires suitable approximation and compactness procedures: in this paper we adopt the framework of SOLAs (Solutions Obtained as Limits of Approximations) in the following local form, adapted from \cite[Definition~2.1]{MingionePalatucci} to the double-phase setting, with an additional uniform local $L^1$ bound on $h(\,\cdot\,,Du_j)$.

\begin{definition}[Local SOLA]\label{def.SOLA}
A very weak solution $u\in W^{1,1}_{\rm loc}(\Omega)$ to \eqref{equation} is called a local SOLA if there exists a sequence of functions $\mu_j\in L^\infty_{\rm loc}(\Omega)$ with
\begin{equation}\label{SOLA.uniform.meas}
\mu_j\dx\stackrel{*}{\rightharpoonup}\mu\quad \text{in }\mathcal M_{\rm loc}(\Omega), \qquad \limsup_{j\to\infty}|\mu_j|(B)\leq|\mu|(\overline B) \quad \text{for every ball $B\Subset\Omega$}
\end{equation}
such that the local energy solutions $u_j\in W^{1,H(\cdot)}_{\rm loc}(\Omega)$, in the sense of Definition~\ref{def.local.energy.solution}, to
\[
-\divergence A(x,Du_j)=\mu_j\qquad\text{in }\Omega,
\]
satisfy
\begin{equation}\label{SOLA.convergence}
u_j\rightharpoonup u\qquad\text{weakly in }W^{1,1}_{\rm loc}(\Omega),\qquad \sup_j\int_{B}h(x,Du_j)\dx<\infty\quad\text{for every ball $B\Subset\Omega$}.
\end{equation}
We call any such sequence $\{(u_j,\mu_j)\}$ an admissible approximating sequence for $u$.
\end{definition}
Notice that the condition on the total variations in \eqref{SOLA.uniform.meas} is satisfied by standard mollifications of the measure; we shall discuss the approximation procedure and its convergence properties in Section~\ref{basic.measure.data.comparisons}.

\medskip

As usual, a concrete construction of local SOLAs is obtained by considering the Dirichlet problem (see \cite{BG89,BG92} and also \cite{CianchiMazya} for related results with Orlicz growth): given $g\in W^{1,H(\cdot)}(\Omega)$, we consider
\begin{equation}\label{CD.sola}
\begin{cases}
 -\divergence  A(x,Du)=\mu\qquad\text{in $\Omega$}\\[2mm]
u\in g+ W^{1,1}_0(\Omega)
\end{cases};
\end{equation}
one can consider energy solutions to the approximating problems and obtain a SOLA by a suitable compactness argument. In particular, this path has already been followed and the existence of SOLAs for homogeneous Dirichlet data under the measure-data balance condition \eqref{subcritical.bound} has been recently established in the aforementioned paper by Song and Youn; see \cite[Proposition 2.3]{SongYoun} for the full range $2-1/n<p\le n$. Since we here consider nonhomogeneous boundary data and we work under the forthcoming weak structural assumptions \eqref{A.veryweak}, we include the existence argument for SOLAs to \eqref{CD.sola}, together with the global flux estimate \eqref{mathcalD}. For the existence theorem, we consider Carath\'eodory vector fields $A\colon\Omega\times\mathbb{R}^n\longrightarrow\mathbb{R}^n$ satisfying the following coercivity, growth, and qualitative strict monotonicity conditions:
\begin{equation}\label{A.veryweak}
\begin{cases}
\displaystyle\langle A(x,z),z\rangle \geq \nu\, H(x,z)\\[2mm]
|A(x,z)|\leq L\,h(x,z)\\[2mm]
\langle A(x,z_1)-A(x,z_2),z_1-z_2\rangle>0
\end{cases}
\end{equation}
for some constants $0<\nu\leq L<\infty$ and for almost every $x\in\Omega$ and all $z,z_1,z_2\in\mathbb{R}^n$ with $z_1\neq z_2$, with the reference functions defined in \eqref{H.def} and \eqref{h.def}. These are the double phase counterpart of the classical assumptions considered by Boccardo and Gallou\"et for the study of equations with standard $p$-growth and measure data, see~\cite{BG89,BG92}.

\medskip

Under these assumptions, we obtain the following existence and global integrability result.

\begin{theorem}\label{existence.SOLA}
Consider the Dirichlet problem \eqref{CD.sola}, where the vector field $A(\,\cdot\,)$ is a Carath\'eodory function satisfying \eqref{A.veryweak} with $p\leq n$ fulfilling \eqref{2-1/n}, $a\in L^\infty(\Omega)$ and $a(\,\cdot\,)\geq0$ almost everywhere in $\Omega$, $g\in W^{1,H(\cdot)}(\Omega)$ and $\mu$ is a finite signed Radon measure on $\Omega$. If moreover
\[
\begin{cases}
\displaystyle\frac{q-1}{p-1}<\frac{n}{n-p}\qquad&\text{when }p<n\\[3mm]
q< +\infty \qquad&\text{when }p=n 
\end{cases},
\]
then \eqref{CD.sola} admits a local SOLA $u\in W^{1,h(\cdot)}(\Omega)$ in the sense of Definition~\ref{def.SOLA}, with $u-g\in W^{1,h(\cdot)}(\Omega)\cap W^{1,1}_0(\Omega)$. The global estimate
\begin{equation}\label{mathcalD}
\int_\Omega h(x,Du)\dx\leq c\biggl[\mathfrak D+\|a\|_{L^\infty(\Omega)}^{\frac1q}\mathfrak D^{\frac{p'}{q'}}\biggr]\qquad\text{with}\qquad \mathfrak D:=|\mu|(\Omega)+\biggl(\int_\Omega H(x,Dg)\dx\biggr)^{\frac1{p'}}
\end{equation}
holds for some constant $c\equiv c(n,p,q,\nu,L,|\Omega|)$ and we can further select an approximating sequence $\{u_j\}$ obeying the same bound uniformly in $j$. As for the superlevel sets, when $p<n$ we obtain
\[
\sup_{\varsigma>0}\varsigma^{\frac{n(p-1)}{n-p}}\bigl|\{x\in\Omega:|u(x)-g(x)|>\varsigma\}\bigr|\leq c\,\mathfrak D^{\frac n{n-p}},
\]
whereas in the borderline case $p=n$ if $\mathfrak D>0$, we have
\[
\bigl|\{x\in\Omega:|u(x)-g(x)|>\varsigma\}\bigr|\leq c|\Omega|\exp\biggl(-\frac{\varsigma}{c\,\mathfrak D^{\frac1{n-1}}}\biggr)\qquad\text{for all }\varsigma>0.
\]
\end{theorem}
The proof of Theorem~\ref{existence.SOLA} uses truncation estimates and compactness under \eqref{A.veryweak}. For the local regularity results when $p\leq n$, we derive estimates uniform along the approximating sequence under \eqref{a.holder} to \eqref{bounds} and the structural assumptions of the result in question. Lemmas~\ref{SOLA.flux.stability} and~\ref{SOLA.oscillation.stability} then give the corresponding estimates for local SOLAs.

\subsubsection{The case $p>n$}

In this case measures embed in the dual of the energy space, and we shall work under the bound \eqref{bound.supercritical} and consider local energy solutions in the sense of Definition \ref{def.local.energy.solution}. Note that under \eqref{a.holder} and \eqref{bound.supercritical}, the weak formulation \eqref{en.solution} extends to every $\varphi\in W^{1,H(\cdot)}(\Omega)$ with compact support in $\Omega$, by approximation in the energy space with smooth functions supported in a fixed compact subset of $\Omega$; see \cite[Lemma~6.4.5 and Propositions~7.2.2-7.2.3]{HH}. Since $p>n$, these approximations also converge uniformly to the continuous representative of $\varphi$, which justifies passage to the limit in the integral against $\mu$.

\medskip

Unless otherwise specified, in the general case we shall use the following notion:
\begin{definition}[Local solution]\label{def.solution}
A local solution (or simply a solution) to \eqref{equation} is a local SOLA in the sense of Definition~\ref{def.SOLA} when $p\leq n$ and a local energy solution in the sense of Definition~\ref{def.local.energy.solution} when $p>n$.
\end{definition}
Whenever we deal with a local energy solution for the full range of $p$, we will state it explicitly; in this case it will be implied that the datum $\mu$ is sufficiently regular if $p\leq n$, for instance having density in $L^\infty_{\rm loc}(\Omega)$.

\subsection{Regular vector fields}\label{regular.A}

We now turn to the regularity of solutions. The structural assumptions are tailored to the quantity being estimated: linear potential bounds for the flux require the regularity assumptions below, whereas nonlinear potential bounds for the solution itself will be proved under the weaker assumptions of Paragraph~\ref{rough.A}.

\medskip

Here we consider the case where the growth condition in \eqref{approximate.A} is formulated at the $C^1$ level. We assume that $A\colon\Omega\times\R^n\to\R^n$ is a Carath\'eodory vector field such that, for every $x\in\Omega$, $A(x,0)=0$ and $A(x,\,\cdot\,)\in C^1(\R^n\setminus\{0\};\R^n)$. Moreover, we assume that there exist constants $0<\nu\leq L<\infty$ such that
\begin{equation}\label{A.strong}
\begin{cases}
|A(x,z)|\,|z|+|\partial_z A(x,z)|\,|z|^2\leq L H(x,z)\\[2mm]
\big\langle\partial_z A(x,z)\xi,\xi\big\rangle\geq\nu\bigl[|z|^{p-2}+a(x)|z|^{q-2}\bigr]|\xi|^2\\[2mm]
|A(x_1,z)-A(x_2,z)|\leq L\omega(|x_1-x_2|)\bigl[h(x_1,z)+h(x_2,z)\bigr]+L|a(x_1)-a(x_2)|\,|z|^{q-1}
\end{cases},
\end{equation}
for all $x,x_1,x_2\in\Omega$, every $z\in\mathbb{R}^n\setminus\{0\}$ and every $\xi\in\mathbb{R}^n$. Here $\omega(\,\cdot\,)$ is a non-decreasing, concave modulus of continuity, continuous at zero and satisfying $\omega(0)=0$. The structural assumptions adopted here, inspired by \cite{BCMgen}, naturally include product-type vector fields of the form
\begin{equation}\label{model.A.strong}
A(x,z)=b_p(x)A_p(z)+a(x)b_q(x)A_q(z),
\end{equation}
where $A_p$ and $A_q$ satisfy standard autonomous $p$- and $q$-growth conditions, respectively. The coefficients $b_p(\,\cdot\,)$ and $b_q(\,\cdot\,)$ are bounded away from zero and infinity, and possess a sufficiently regular modulus of continuity. A basic model example is of course
\begin{equation}\label{basic.example.VF}
A(x,z)=b(x)\bigl[|z|^{p-2}+a(x)|z|^{q-2}\bigr]z, 
\end{equation}
with $b(\,\cdot\,)$ positive, bounded and quantitatively continuous.

\begin{remark}
{\rm The $C^1$ regularity of $A(x,\,\cdot\,)$ away from the origin and the first two inequalities in \eqref{A.strong} may also be assumed for almost every $x\in\Omega$. After choosing the representative of $A(\,\cdot\,)$ that is continuous in $x$, one can regularize the frozen fields in $z$ with uniform structural constants and pass to the limit in the corresponding estimates. We formulate these assumptions for every $x\in\Omega$ so that the autonomous regularity results can be applied directly to the frozen fields; see Remark~\ref{approx.vf}.}
\end{remark}

\medskip

The coercivity and growth assumptions in \eqref{A.veryweak} (or the forthcoming \eqref{A.weak}) imply $\nu q^{-1}h(x,z)\leq|A(x,z)|\leq Lh(x,z)$ for almost every $x\in\Omega$ and every $z\in\mathbb{R}^n$, but under the stronger assumptions \eqref{A.strong}, the upper bound improves to
\begin{equation}\label{a-h.weak}
\frac{\nu}{q}h(x,z)\leq|A(x,z)|\leq\frac{L}{p}h(x,z)\qquad\text{for almost every $x\in\Omega$ and every $z\in\mathbb{R}^n$}.
\end{equation}
We recall that a function is Dini continuous if it admits a modulus of continuity $\omega(\,\cdot\,)$ satisfying
\begin{equation}\label{Dini}
\int_0^R\omega(\varrho)\,\frac{d\varrho}{\varrho}<\infty\qquad\text{for some $R>0$.}
\end{equation}
Dini conditions involving fractional powers of the modulus already emerge naturally in the standard $p$-growth theory; for instance, when $p>2$, Duzaar and Mingione~\cite{DuMiAJM} assume that $\omega^{2/p}$ satisfies \eqref{Dini} in their Wolff potential estimate for the gradient. In our derivation of the linear potential estimate for the flux, fractional powers of $\omega$ likewise enter through the harmonic approximation argument, leading us to impose the stronger requirement that
\begin{equation}\label{Dini.powers}
\int_0^1[\omega(\varrho)]^s\,\frac{d\varrho}{\varrho}<\infty\qquad\text{for every $s>0$;}
\end{equation}
by virtue of \eqref{dyadic.Dini.smallness}, this ensures that for each fixed $s>0$ and $\delta\in(0,1)$,
\[
\lim_{R\searrow0}\sum_{j=0}^\infty[\omega(\delta^jR)]^s=0.
\]
Within the model framework \eqref{model.A.strong}, condition \eqref{Dini.powers} is required to hold for both the moduli of continuity of $b_p(\,\cdot\,)$ and $b_q(\,\cdot\,)$, whereas the oscillation of $a(\,\cdot\,)$ is controlled separately by \eqref{a.holder}.

Given any ball $B_R(x_0)\subset\Omega$, we denote by $\bI_1^{|\mu|}(x_0,R)$ the truncated Riesz potential of order one of the total variation measure $|\mu|$, defined as in \cite{KMGuide} by
\[
\bI_1^{|\mu|}(x_0,R):=\int_0^R\frac{|\mu|(B_\varrho(x_0))}{\varrho^{n-1}}\,\frac{d\varrho}{\varrho}.
\]
This potential precisely controls the sum of the normalized masses $|\mu|(B_\varrho(x_0))/\varrho^{n-1}$ arising in the comparison estimates on concentric balls, as detailed in \eqref{dyadic.potential.bound}. Pointwise gradient bounds of this flavor, despite being classical for the Poisson equation, have a relatively recent yet rich history: for autonomous uniformly elliptic equations with quadratic growth, Mingione~\cite[Theorem~1.1]{MingioneJEMS} bounded each $|D_i u(x_0)|$ in terms of the Riesz potential of $|\mu|$ and a local average of $|D_i u|$, whereas for scalar equations of $p$-Laplacian type, corresponding pointwise estimates for $|Du|^{p-1}$ were established in \cite[Theorem~1.1]{DuMiJFA} for $2-1/n<p\leq2$ and in \cite[Theorem~1.3]{KMI1} for the much more delicate degenerate autonomous case $p\geq2$, with subsequent extensions to classes of autonomous equations with degenerate Orlicz growth provided in \cite[Theorem~1.1]{B}. The main result of this paper, contained in Theorem~\ref{main.thm} below, is that linear potential estimates continue to hold not only for appropriate powers of the modulus of the gradient of solutions to \eqref{equation}, but also for the whole flux.

\medskip

To formulate the local estimates below, we need to fix an open subset
\begin{equation}\label{Omega0}
\Omega_0\Subset\Omega,\qquad\text{setting}\qquad d_{\Omega_0}:={\rm dist}(\Omega_0,\partial\Omega)>0,
\end{equation}
and we define an intermediate open set by
\begin{equation}\label{Omega1}
\Omega_1:=\bigl\{x\in\Omega:{\rm dist}(x,\Omega_0)<d_{\Omega_0}/4\bigr\},
\end{equation}
so that $\Omega_0\Subset\Omega_1\Subset\Omega$ and $B_{2R}(x_0)\subset\Omega_1$ whenever $x_0\in\Omega_0$ and $0<R\leq d_{\Omega_0}/8$. When $p\leq n$, the flux magnitude $h(\,\cdot\,,Du)$ is locally integrable by the very definition of a local SOLA, whereas for $p>n$ we work with local energy solutions, so that the integral of $H(\,\cdot\,,Du)$ over $\Omega_1$ is also finite; accordingly, we set
\begin{equation}\label{energy.M}
\mathcal M:=
\begin{cases}
\displaystyle\int_{\Omega_1}h(x,Du)\dx,&\text{if }p\leq n,\\[3mm]
\displaystyle\int_{\Omega_1}\bigl[H(x,Du)+h(x,Du)\bigr]\dx,&\text{if }p>n.
\end{cases}
\end{equation}
When $p\leq n$ and $u$ is the Dirichlet SOLA provided by Theorem~\ref{existence.SOLA}, the global estimate \eqref{mathcalD} bounds $\mathcal M$ in terms of the data of the Dirichlet problem.

\begin{theorem}[Linear potential estimates for the flux]\label{main.thm} 
Let $u$ be a local solution to \eqref{equation} in the sense of Definition~\ref{def.solution}; assume that \eqref{a.holder} to \eqref{bounds} hold and that $A(\,\cdot\,)$ satisfies \eqref{A.strong} with a modulus $\omega(\,\cdot\,)$ fulfilling \eqref{Dini.powers}.

Let $\Omega_0\Subset\Omega$ as in \eqref{Omega0}; then there exist a constant $c\geq1$, depending only on $\data$ and $\gdata$, and a threshold radius $R_0\in(0,\min\{1,d_{\Omega_0}/8\}]$, additionally depending on $\omega(\,\cdot\,)$, such that for almost every $x_0\in\Omega_0$ and every $R\in(0,R_0]$ it holds
\begin{equation}\label{potential.estimate.Du}
\bigl|A(x_0,Du(x_0))\bigr|\approx_{p,q,\nu,L}h\bigl(x_0,Du(x_0)\bigr)\leq c\mean{B_R(x_0)}h(x,Du)\dx+c\,\bI_1^{|\mu|}(x_0,R)+c.
\end{equation}
\end{theorem}

Details on how constants and radii depend on the parameters can be found in Section~\ref{notation}, where the precise role of the intermediate set $\Omega_1$ and of the quantity $\mathcal M$ - which are absorbed into $\gdata$ - is explicitly addressed. For the model case \eqref{model.A.strong}, Theorem~\ref{main.thm} should still hold if \eqref{Dini.powers} is replaced by the standard Dini condition \eqref{Dini} on $b_p(\,\cdot\,)$ and $b_q(\,\cdot\,)$. Indeed, since both coefficients are decoupled from $a(\,\cdot\,)$, one could likely estimate their oscillations separately, a point we leave aside in this paper.

\begin{proposition}[Qualitative oscillation estimate for the flux]\label{qualitative.flux.oscillation}
Under the assumptions of Theorem~\ref{main.thm}, for every $\varepsilon\in(0,1]$ there exist constants $C_\varepsilon\geq1$, $\tau_\varepsilon\in(0,1/4]$ and a radius $R_\varepsilon\in(0,R_0/2]$ with the following property. Let $x_0\in\Omega_0$ and $R\in(0,R_\varepsilon]$, and suppose that
\begin{equation}\label{Linfinity.bound}
{\bigl\|h(\,\cdot\,,Du)\bigr\|}_{L^\infty(B_{2R}(x_0))}\leq\lambda 
\end{equation}
for some $\lambda\geq1$. Then
\begin{equation}\label{qualitative.flux.oscillation.estimate}
\bigl|A(x,Du(x))-A(y,Du(y))\bigr|\leq\varepsilon\lambda+C_\varepsilon\bigl[\bI_1^{|\mu|}(x,R)+\bI_1^{|\mu|}(y,R)\bigr]
\end{equation}
whenever $x,y\in B_R(x_0)$ are Lebesgue points of $A(\,\cdot\,,Du)$ such that 
\[
|x-y|\leq\tau_\varepsilon R.
\]
The constants $C_\varepsilon,\tau_\varepsilon$ depend only on $\data,\gdata$ and $\varepsilon$, while $R_\varepsilon$ additionally depends on $\omega(\,\cdot\,)$.
\end{proposition}
The previous proposition has an immediate consequence: 
\begin{corollary}[Continuity of the flux]\label{cor:flux.continuity}
Under the assumptions of Theorem~\ref{main.thm}, suppose additionally that
\begin{equation}\label{uniform.vanishing.potential}
\lim_{R\searrow0}\sup_{x\in K}\bI_1^{|\mu|}(x,R)=0\qquad\text{for every compact set }K\Subset\Omega;
\end{equation}
then $A(\,\cdot\,,Du)$ admits a continuous representative in $\Omega$, and $u$ admits a representative in $C^1(\Omega)$. In particular, the same conclusions hold if $\mu=f\dx$ with $f\in L^{n,1}_{\rm loc}(\Omega)$.
\end{corollary}

\medskip

Our approach to gradient estimates is organized around the oscillation of the flux $A(\,\cdot\,,Du)$. For every ball $B\Subset\Omega$, we denote its ($L^1$-)excess by
\begin{equation}\label{def.flux.excess}
E(Du,B):=\mean{B}\bigl|A(x,Du)-(A(\,\cdot\,,Du))_B\bigr|\dx
\end{equation}
and call it the \emph{flux excess}. The decay estimate for this excess in Lemma~\ref{quantitative.excess.decay} is the main tool in the proofs of the potential and continuity estimates above.

In the super-dimensional case $p>n$, the restriction to compactly contained subsets is also needed to obtain uniform higher integrability estimates for the gradients of the homogeneous comparison maps, especially when dealing with the borderline balance $q/p=1+\alpha/n$; see Theorem~\ref{hi.correct.exponent}.

\subsection{Rough vector fields}\label{rough.A}

In this part of the paper we consider irregular vector fields and the corresponding results will be weaker. We assume that $A$ is a Carath\'eodory vector field, $A(x,0)=0$, and that
\begin{equation}\label{A.weak}
\begin{cases}
\displaystyle\langle A(x,z),z\rangle \geq \nu H(x,z)\\[2mm]
|A(x,z)|\leq Lh(x,z)\\[2mm]
\big\langle A(x,z_1)-A(x,z_2),z_1-z_2 \big\rangle \geq \nu\Bigl[{\big(|z_1|+|z_2|\big)}^{p-2}+a(x){\big(|z_1|+|z_2|\big)}^{q-2}\Bigr]{|z_1-z_2|}^2
\end{cases}
\end{equation}
for some constants $0<\nu\leq L<\infty$, almost every $x\in\Omega$, and all $z,z_1,z_2\in\mathbb{R}^n$; if $p<2$ or $q<2$, the right-hand side of the last inequality is taken to vanish at $z_1=z_2=0$. Such assumptions cover operators of the form \eqref{basic.example.VF} featuring measurable, bounded, and uniformly elliptic coefficients. Notice that \eqref{A.strong} implies \eqref{A.weak}, up to rescaling $\nu$ and $L$ by factors depending only on $p$ and $q$; the sole nontrivial step follows from the chain of bounds that can also be found in \eqref{difference.a.general}. The quantitative monotonicity in the third line of \eqref{A.weak} will come into play later, proving to be essential both in deriving our comparison lemmas and throughout the regularity proof.

\medskip

Nonlinear potential bounds for solutions trace back to the seminal paper of Kilpel\"ainen and Mal\'y \cite{KilpelainenMaly1994}, see also \cite{KMGuide} for an overview of later advances on measure-data problems. For the $p$-Laplace system with vector-valued measure data, Kuusi and Mingione~\cite[Theorems~1.1 and~1.4]{KMVectorial} proved Wolff potential bounds for $u$ when $p>2-\frac1n$ and Riesz potential bounds for $Du$ when $p>2$. Within the double phase framework, we recall the contributions of Skrypnik and Buryachenko, who derived in \cite{SkrypnikBuryachenko2017}  pointwise estimates for solutions, and also the subsequent developments that include the Orlicz-growth extensions in \cite{ChlebickaGiannettiZatorska2024}, as well as the generalized Wolff potential bounds proved by Benyaiche and Khlifi \cite{BenyaicheKhlifi2023} for non-negative supersolutions. 

The following theorem gives an upper potential bound for solutions with finite signed measure data under \eqref{bounds} and \eqref{A.weak}.

\begin{theorem}[Nonlinear potential estimates for the solution]\label{main.thm.2}
Let $u$ be a local solution to \eqref{equation} in the sense of Definition~\ref{def.solution}. Assume that \eqref{a.holder} to \eqref{bounds} hold and that $A$ satisfies \eqref{A.weak}.

For $\Omega_0\Subset\Omega$ as in \eqref{Omega0}, and $\Omega_1$ and $\mathcal M$ defined by \eqref{Omega1} and \eqref{energy.M}, there exist constants $c\geq1$ and $R_0\in(0,\min\{1,d_{\Omega_0}/8\}]$, depending only on $\data$ and $\gdata$, such that
\begin{equation}\label{potential.estimate.u}
|u(x_0)|\leq c\int_0^R h_{x_0}^{-1}\biggl(\frac{|\mu|(B_\rho(x_0))}{\rho^{n-1}}\biggr)\,d\rho+c\mean{B_R(x_0)}|u|\dx
\end{equation}
for every Lebesgue point $x_0\in\Omega_0$ of $u$ and every $R\in(0,R_0]$. Here $h_{x_0}(t):=t^{p-1}+a(x_0)t^{q-1}$ for $t\geq0$.
\end{theorem}

\begin{remark}\label{rem:potential.u.sup}
{\rm Under the assumptions of Theorem~\ref{main.thm.2}, the potential in \eqref{potential.estimate.u} admits an equivalent formulation using the supremum of $a(\,\cdot\,)$ at each scale, which might be useful in certain situations. More precisely, for 
\[
h_\rho(t):=t^{p-1}+\bigl[\sup_{B_\rho(x_0)}a(\,\cdot\,)\bigr]t^{q-1},
\]
we have
\[
\int_0^R h_\rho^{-1}\biggl(\frac{|\mu|(B_\rho(x_0))}{\rho^{n-1}}\biggr)\,d\rho\approx_{\data,\gdata}\int_0^R h_{x_0}^{-1}\biggl(\frac{|\mu|(B_\rho(x_0))}{\rho^{n-1}}\biggr)\,d\rho.
\]
}
\end{remark}

\medskip

\subsection{Outline of the proofs}\label{heuristic.proofs}

In this paper, the main technical tool we employ to establish the pointwise gradient estimates - which represent the most challenging part - is the powerful indirect linearization strategy of Kuusi and Mingione~\cite{KMI1, KMGuide}, adapted here to our framework as heuristically outlined below. The core idea, which consists in comparing the solution with homogeneous maps to isolate the scales where nondegenerate behavior prevails, allows one to sharpen the measure dependence in the comparison bounds and ultimately recover the linear Riesz potential available for the Poisson equation. Adapting this scheme to the double phase setting requires monitoring the interplay between the two growth regimes at every step. Rather than separating the phases throughout the iteration, we formulate all estimates in terms of the natural quantities dictated by the full operator - ultimately provided by the flux of the equation - while carefully accounting for the error terms generated by the spatial variation of the coefficients and by the phase transitions.

\medskip

Our quantitative potential estimates also rely on the approach developed by Avelin, Kuusi, and Mingione in~\cite{AKM}, particularly regarding how oscillation terms and explicit measure remainders enter the decay of the flux excess through a linearization scheme. Excess decay estimates for the flux have recently been established outside this setting: when dealing with obstacle problems with measure data - where potential estimates were studied by Byun, Song, and Youn~\cite{ByunSongYoun2023} - the same authors proved a decay bound for an $L^\kappa$-excess of $A(Du)$ ($0<\kappa<1$) combining a linear measure term with an obstacle contribution in the range $1<p\leq2-1/n$~\cite{ByunSongYoun2024}. Chlebicka, Song, Youn, and Zatorska-Goldstein~\cite{ChlebickaSongYounZatorska2025} obtained a decay estimate for the $L^1$-excess of the local part of the flux featuring a linear measure remainder combined with a nonlocal tail when considering mixed local-nonlocal equations with $p\geq2$. In the present double phase setting, Lemma~\ref{quantitative.excess.decay} provides the corresponding decay bound for the $L^1$-excess of the full flux $A(x,Du)$, cleanly separating the contributions of the coefficients $\mathcal{S}(R)\lambda$ from those of the measure $|\mu|(B_R)/R^{n-1}$.

\medskip

To make the underlying mechanism as transparent as possible, let us illustrate how it operates in the degenerate case $p\geq 2$ for model operators as in~\eqref{basic.example.VF}, under the hypotheses of Theorem~\ref{main.thm}: once we fix a contraction factor $\delta\in(0,1/32)$, depending only on $\data$ and $\gdata$, admissible in Lemma~\ref{quantitative.excess.decay} for $\varepsilon=\vartheta=1/2$, we can find a constant $\mathcal B\geq 1$, still depending solely on $\data$ and $\gdata$, such that - shrinking $R_0$ if necessary - on every ball $B_R(x_0)$ with $x_0\in\Omega_0$ and $0<R\leq R_0$, if
\begin{equation}\label{test.assumption}
\frac{\lambda}{2}\leq\mean{B_{\delta R}}h(x,Du)\dx,\qquad\mean{B_R}h(x,Du)\dx\leq\lambda,\qquad\frac{|\mu|(B_R)}{R^{n-1}}\leq\frac{\lambda}{\mathcal B}
\end{equation}
for some $\lambda\geq1$, then
\begin{multline}\label{decay.heuristic}
\mean{B_{\delta^\smalltwo R}}\bigl|A(x,Du)-\bigl(A(\,\cdot\,,Du)\bigr)_{B_{\delta^\smalltwo R}}\bigr|\dx\\
\leq\frac12\mean{B_{\delta R}}\bigl|A(x,Du)-\bigl(A(\,\cdot\,,Du)\bigr)_{B_{\delta R}}\bigr|\dx+c\,\mathcal S(R)\lambda+c\,\frac{|\mu|(B_R)}{R^{n-1}},
\end{multline}
where $c$ depends only on $\data,\gdata$ and $\mathcal S$ is defined in \eqref{small.S}. The radius $R_0$ additionally depends on $\omega(\,\cdot\,)$, and all these choices are independent of $\lambda$. In the iteration of Paragraph~\ref{main.proof}, this estimate controls the accumulated oscillations of the flux: the measure contributions are bounded by the truncated linear Riesz potential, while \eqref{Dini.powers} ensures the summability of the coefficient errors over geometric scales; see Paragraph~\ref{dyadic.potentials}.

\medskip

More precisely, let us describe the indirect linearization procedure along with some intermediate estimates of independent interest in the measure data range $p\leq n$. For the purposes of this heuristic discussion, we restrict ourselves to the model case \eqref{model.A} with $2\leq p\leq n$, temporarily setting aside the intermediate comparison maps and the fixed changes of radius required in the rigorous proof. For this a priori comparison argument, $u$ is a local energy solution to \eqref{equation}, as in the approximating problems of Definition~\ref{def.SOLA}. Given a ball $B_R=B_R(x_0)\Subset\Omega$, we define $v$ and $w$ by the Dirichlet problems
\[
\begin{cases}
\divergence A(x,Dv)=0&\text{in $B_R$},\\[2mm]
v=u&\text{on $\partial B_R$},
\end{cases}
\qquad
\begin{cases}
\divergence A(x_m,Dw)=0&\text{in $B_{R/2}$},\\[2mm]
w=v&\text{on $\partial B_{R/2}$},
\end{cases}
\]
where $x_m$ is a minimizer of $x\mapsto a(x)$ in $\overline{B_{R/2}}$; consequently, the frozen vector field $A(x_m,z)$ exhibits growth governed by the frozen flux magnitude $|z|^{p-1}+a(x_m)|z|^{q-1}$. In the specific model framework \eqref{model.A}, one might be tempted to seek a $p$-harmonic approximation of $v$ on balls fulfilling the $p$-phase condition \eqref{p.phase}, thereby treating the $q$-term purely as a perturbation. However, for a general operator satisfying \eqref{A.strong}, the frozen operator may fail to reduce to the $p$-Laplacian even when $a(x_m)=0$; we therefore employ the same frozen comparison problem for $w$ across both phases.

The integral non-degeneracy and boundedness bounds in \eqref{test.assumption}, merged with the rough comparison estimates we derive as a first step of the proof, see Proposition~\ref{rough.comparison.DuDw}, schematically translate into integral bounds on the more regular solution $w$:
\[
\frac{\lambda}{2}\lesssim\mean{B_{\delta R}}h(x_m,Dw)\dx,\qquad\mean{B_{R/2}}h(x_m,Dw)\dx\lesssim\lambda.
\]
The regularity theory for the frozen problem - specifically the quantitative gradient oscillation estimates - now comes strongly into play, allowing us to pass from integral bounds to pointwise ones, provided $\delta$ and the comparison error are chosen sufficiently small:
\begin{equation}\label{comparability}
\frac{\lambda}{2}\lesssim\inf_{B_{\delta R}}h(x_m,Dw),\qquad\sup_{B_{R/4}}h(x_m,Dw)\lesssim\lambda;
\end{equation}
see \eqref{non.degeneracy} and \eqref{pointwise.1}-\eqref{non.degeneracy.Dw1}.

\medskip

This, in some sense, allows us to linearize the equation. We evaluate the scalar coefficient at $Dw$ and define the resulting flux error by
\[
\mathcal E:=A(x,Du)-\bigl[|Dw|^{p-2}+a(x)|Dw|^{q-2}\bigr]Du;
\]
the equation can then be rewritten as
\begin{equation}\label{heu.frozen}
-\divergence\Bigl(\bigl[|Dw|^{p-2}+a(x)|Dw|^{q-2}\bigr]Du\Bigr)=\mu+\divergence\mathcal E. 
\end{equation}
At this point, the perturbative character of the vector field comes into play. In $B_{\delta R}$, we write
\begin{align}\label{heu.1}
-\divergence\Bigl(\frac{h(x_m,Dw)}{|Dw|}Du\Bigr)&=-\divergence\Bigl(\bigl[|Dw|^{p-2}+a(x)|Dw|^{q-2}\bigr]Du\Bigr)+\divergence\Bigl([a(x)-a(x_m)]|Dw|^{q-2}Du\Bigr)\notag\\
&=\mu+\divergence\mathcal E+\divergence\Bigl([a(x)-a(x_m)]|Dw|^{q-2}Du\Bigr).
\end{align}
On the left-hand side, \eqref{comparability} gives
\[
\frac{h(x_m,Dw)}{|Dw|}\approx\frac{\lambda}{h_{x_m}^{-1}(\lambda)}=:\kappa_\lambda\qquad\text{in }B_{\delta R}:
\]
$\kappa_\lambda$ describes the size of this coefficient, which still varies with $x$ through $Dw$. On the right-hand side, employing the H\"older continuity of $a(\,\cdot\,)$ and again \eqref{comparability}, we have
\[
|a(x)-a(x_m)|\,|Dw|^{q-2}\lesssim R^\alpha|Dw|^{q-2}\lesssim R^\alpha\lambda^{\frac{q-2}{p-1}}.
\]
We now invoke another essential ingredient of the theory developed here, namely the flux magnitude decay established in Proposition~\ref{decay.measure.data},
\[
\mean{B_R}h(x,Du)\dx\lesssim R^{-(n-1)};
\]
see Corollary~\ref{flux.decay}. This shows that, at small scales, the averaged flux magnitude exhibits at most the same growth as the flux magnitude of the fundamental solution (compare with \eqref{fundamental.flux.magnitude}). Applying this bound at scale $\delta R$ and appealing to the lower bound in \eqref{test.assumption}, we deduce the basic coupling between levels and radii
\[
\lambda\lesssim R^{-(n-1)},
\]
with an implicit constant depending on the fixed parameter $\delta$. Consequently, we can normalize the size of the perturbative term by the natural scaling of the main $p$-part times a positive power of $R$:
\[
R^\alpha\lambda^{\frac{q-2}{p-1}}=R^\alpha\lambda^{\frac{q-p}{p-1}}\lambda^{\frac{p-2}{p-1}}\lesssim R^{2s_0}\lambda^{\frac{p-2}{p-1}},\qquad 2s_0:=\alpha-\frac{q-p}{p-1}(n-1),
\]
where the crucial point is that $s_0>0$ by virtue of \eqref{subcritical.bound}. Since $h_{x_m}(t)\geq t^{p-1}$, we also have
\[
\kappa_\lambda =\frac{\lambda}{h_{x_m}^{-1}(\lambda)} \geq\lambda^{\frac{p-2}{p-1}},
\]
and this ensures that the perturbation term in \eqref{heu.1} is controlled by $R^{2s_0}\kappa_\lambda$. Locally, for $R$ sufficiently small, this perturbation can be reabsorbed into the leading coefficients on $B_{\delta R}$. Rearranging the terms in \eqref{heu.1} and dividing by $\kappa_\lambda$, we obtain
\[
-\divergence\left(\frac{|Dw|^{p-2}+a(x)|Dw|^{q-2}}{\kappa_\lambda}Du\right)=\frac{\mu}{\kappa_\lambda}+\divergence\Bigl(\frac{\mathcal E}{\kappa_\lambda}\Bigr),\quad\text{with}\quad \frac{|Dw|^{p-2}+a(x)|Dw|^{q-2}}{\kappa_\lambda}\approx1,
\]
in $B_{\delta R}$, with implicit constants independent of $\lambda$. Treating $Dw$ as a frozen weight, the operator on the left-hand side is linear in $u$ and uniformly elliptic on $B_{\delta R}$, whereas the right-hand side still depends on $Du$ through $\mathcal E$. Under condition \eqref{test.assumption}, the equation thus behaves as a linear elliptic problem with ellipticity $\kappa_\lambda$. On $B_R$, the datum $\mu$ shifts $Du$ by a term of order $|\mu|(B_R)/[R^{n-1}\kappa_\lambda]$; multiplying back by $\kappa_\lambda$ to estimate the flux $A(x,Du)$ cancels the ellipticity weight and leaves a contribution of size $|\mu|(B_R)/R^{n-1}$. After averaging, the remaining flux error is controlled by the same measure term and by a positive power of $R$ times the averaged flux magnitude. We obtain \eqref{decay.heuristic} by coupling the excess decay of the homogeneous frozen equation (Proposition~\ref{varphi-harmonic.are.regular}) with the explicit error structure; a complete and rigorous proof can be found in Lemma~\ref{quantitative.excess.decay}.

\section{The fundamental solution}\label{sec.fundamental.solution}

In this section we compute the radial fundamental solution corresponding to \eqref{radial} and study its behavior near the origin and its local integrability.

\subsection{The radial solution and its conserved flux}\label{distr.solution}

For each $r>0$, the map $h_{\texttt r}(r,t):=t^{p-1}+r^\alpha t^{q-1}$ is continuous and strictly increasing from $[0,+\infty)$ onto itself. This allows us to define the positive radial gradient profile
\begin{equation}\label{fundamental.gradient.inverse}
f_{p,q}(r):=h_{\texttt r}(r,\,\cdot\,)^{-1}\Bigl(\frac{1}{n\omega_n r^{n-1}}\Bigr),\qquad r>0,
\end{equation}
which by definition satisfies
\begin{equation}\label{flux.fundamental}
r^{n-1}\bigl[f_{p,q}(r)^{p-1}+r^\alpha f_{p,q}(r)^{q-1}\bigr]=\frac{1}{n\omega_n},\qquad r>0.
\end{equation}
Fixing the additive constant so that
\begin{equation}\label{fundamental.profile.integral}
F_{p,q}(r):=1-\int_1^r f_{p,q}(s)\,ds,\qquad r>0,
\end{equation}
we set $\mathcal G_{p,q}(x):=F_{p,q}(|x|)$ for $x\neq0$. Clearly we have $F_{p,q}(1)=1$ and $F_{p,q}'=-f_{p,q}<0$, so we can write
\[
D\mathcal G_{p,q}(x)=-f_{p,q}(|x|)\frac{x}{|x|},
\]
and the corresponding flux vector field takes the form
\[
A_{\mathtt r}(x,D\mathcal G_{p,q}):=|D\mathcal G_{p,q}|^{p-2}D\mathcal G_{p,q}+|x|^\alpha|D\mathcal G_{p,q}|^{q-2}D\mathcal G_{p,q}=-\frac{x}{n\omega_n|x|^n}.
\]
In particular, it belongs to $L^1_{\rm loc}(\R^n;\R^n)$ and is divergence-free outside the origin. Testing against any $\varphi\in C_c^\infty(\R^n)$ and customarily integrating by parts on a punctured ball covering $\operatorname{supp}\varphi$ yields
\[
    \int_{\R^n}\langle A_{\mathtt r}(x,D\mathcal G_{p,q}),D\varphi\rangle\dx=\lim_{\varepsilon\searrow0}\frac{1}{n\omega_n\varepsilon^{n-1}}\int_{\partial B_\varepsilon}\varphi\,d\mathcal H^{n-1}=\varphi(0),
\]
showing that $-\divergence A_{\mathtt r}(x,D\mathcal G_{p,q})=\delta_0$ in $\mathcal D'(\R^n)$. The flux magnitude matches that of the standard $p$-Laplacian fundamental solution $\mathcal G_p$: regardless of how $p$, $q$, and $\alpha$ are related, \eqref{flux.fundamental} implies
\begin{equation}\label{fundamental.rel.gradient}
h_{\texttt r}\bigl(|x|,|D\mathcal G_{p,q}(x)|\bigr)=|D\mathcal G_p(x)|^{p-1}=\frac{1}{n\omega_n|x|^{n-1}}.
\end{equation}
For $\varepsilon>0$ and $b>0$, we therefore have
\[
\int_{B_\varepsilon}\bigl[h_{\texttt r}\bigl(|x|,|D\mathcal G_{p,q}(x)|\bigr)\bigr]^b\dx\approx\int_0^\varepsilon r^{n-1-b(n-1)}\,dr,
\]
which converges if and only if
\begin{equation}\label{fundamental.b.threshold}
b<\frac{n}{n-1};
\end{equation}
both the upper bound in Lemma~\ref{local.rH} for $p\leq n$ and the constraint on $b$ in Paragraph~\ref{exponent.b} are dictated by this exponent. Integrating \eqref{fundamental.rel.gradient} over $B_R$ shows that the radial solution achieves the exact flux scaling of Proposition~\ref{decay.measure.data} and Corollary~\ref{flux.decay}: indeed
\[
\int_{B_R}h(x,D\mathcal G_{p,q})\dx=R,\qquad R^{n-1}\mean{B_R}h(x,D\mathcal G_{p,q})\dx=\frac1{\omega_n}.
\]
The same model is also responsible for the presence of the linear Riesz potential in Theorem~\ref{main.thm}: when $\mu=\delta_0$ and $0<|x|<R$, the origin lies inside $B_\rho(x)$ only for $\rho>|x|$, and in this case
\[
\bI_1^{\delta_0}(x,R)=\int_{|x|}^R\frac{d\rho}{\rho^n}=\frac{|x|^{1-n}-R^{1-n}}{n-1}\approx_n|x|^{1-n}\approx_n h\bigl(x,D\mathcal G_{p,q}(x)\bigr)
\]
if we restrict ourselves to $0<|x|<R/2$. The linear potential in \eqref{potential.estimate.Du} therefore captures the exact singularity scale of the fundamental solution's flux. Finally, differentiating \eqref{flux.fundamental} gives
\begin{equation}\label{radial.flux.derivative}
\frac{d}{dr}\Bigl\{r^{n-1}\bigl[f_{p,q}(r)^{p-1}+r^\alpha f_{p,q}(r)^{q-1}\bigr]\Bigr\}=0,
\end{equation}
from which we deduce
\begin{equation}\label{f.logderivative}
\frac{rf_{p,q}'(r)}{f_{p,q}(r)}=-\frac{(n-1)+(n+\alpha-1)r^\alpha f_{p,q}(r)^{q-p}}{(p-1)+(q-1)r^\alpha f_{p,q}(r)^{q-p}}.
\end{equation}

\subsection{The pure phases, and their interpolation}\label{log.par}

We set, for $p>1$ and $n\in\N$ with $n\geq2$,
\begin{equation}\label{cpn}
c_p(n)=(n\omega_n)^{-\frac{1}{p-1}}.  
\end{equation}
If $a\equiv0$, \eqref{flux.fundamental} collapses to $r^{n-1}f(r)^{p-1}=\frac{1}{n\omega_n}$, reducing the profile to
\begin{equation}\label{fundamental.p.gradient}
f_p(r)=c_p(n)r^{-\frac{n-1}{p-1}},
\end{equation}
which is simply the familiar radial gradient of the $p$-Laplacian fundamental solution. Conversely, retaining only the weighted $q$-term yields $r^{n-1+\alpha}f(r)^{q-1}=\frac{1}{n\omega_n}$, giving
\begin{equation}\label{fundamental.q.gradient}
f_q(r) = c_q(n) r^{-\frac{n+\alpha-1}{q-1}}:
\end{equation}
the exponent of the radial gradient is the one obtained by replacing $n$ with $n+\alpha$ in the unweighted $q$-Laplacian profile. When both terms interact, the radial gradient transitions smoothly between these two pure-phase behaviors. The logarithmic derivative $rf'(r)/f(r)$ measures how much $f$ deviates from a pure power law (as discussed in Paragraph~\ref{N-function}), remaining identically equal to $m$ whenever $f(r)=cr^m$; we can call this quantity the logarithmic slope. The logarithmic slopes of the pure $p$-profile and the weighted $q$-profile are, respectively, 
\[
m_p:=-\frac{n-1}{p-1}\qquad\text{and}\qquad m_q:=-\frac{n+\alpha-1}{q-1}.
\]
Furthermore, introducing the ratio
\[
\lambda(r):=r^\alpha f_{p,q}(r)^{q-p} = \frac{r^\alpha f_{p,q}(r)^{q-1}}{f_{p,q}(r)^{p-1}}= \frac{r^\alpha f_{p,q}(r)^{q}}{f_{p,q}(r)^{p}},
\]
we can recast \eqref{f.logderivative} in the form
\[
\frac{rf_{p,q}'(r)}{f_{p,q}(r)} = \frac{p-1}{(p-1)+(q-1)\lambda(r)}\,m_p+\frac{(q-1)\lambda(r)}{(p-1)+(q-1)\lambda(r)}\,m_q.
\]
The logarithmic slope of the double phase fundamental gradient is thus a convex combination of the pure-phase slopes $m_p$ and $m_q$, with weights proportional to $p-1$ and $(q-1)\lambda(r)$. This shows how the profile interpolates between the two regimes: it stays close to $m_p$ when $\lambda(r)$ is small (where the $p$-term dominates the flux) and approaches $m_q$ when $\lambda(r)$ is large (where the $q$-term dominates). 

\subsection{Gradient asymptotics under the subcritical balance condition}\label{gradient.asymptotics}

We define the balance parameter
\[
\beta:= \alpha-\frac{n-1}{p-1}(q-p) = \alpha+m_p(q-p)
\]
and the comparison of the two blow-up rates yields the following equivalences:
\begin{equation}\label{fundamental.balance}
\beta>0\qquad\iff\qquad\frac{q-1}{p-1} < 1+\frac{\alpha}{n-1}\qquad\iff\qquad -m_q=\frac{n+\alpha-1}{q-1} > -m_p=\frac{n-1}{p-1}.
\end{equation}
Assuming henceforth that $\beta>0$, we introduce the normalized profile
\begin{equation}\label{normalized.fundamental.profile}
g(r) := \frac{f_{p,q}(r)}{r^{m_p}}.
\end{equation}
Substituting this into \eqref{flux.fundamental} produces
\begin{equation}\label{normalized.fundamental.equation}
g(r)^{p-1} + r^\beta g(r)^{q-1} =\frac1{n\omega_n},
\end{equation}
which immediately implies $0<g(r)\leq c_p(n)$. Furthermore, observing that
\[
0\leq \frac{1}{n\omega_n}-g(r)^{p-1} = r^\beta g(r)^{q-1} \leq c(n,p,q) r^\beta\qquad\implies\qquad g(r)\xrightarrow[r \searrow 0]{} {(n\omega_n)}^{-\frac{1}{p-1}},
\]
we deduce the limiting relation
\begin{equation}\label{eqn:fund-sol-asympt}
\lim_{x\to0} \frac{|D\mathcal G_{p,q}(x)|}{|D\mathcal G_p(x)|} = \lim_{x\to0} \frac{|D\mathcal G_{p,q}(x)|}{(n\omega_n)^{-\frac{1}{p-1}}|x|^{-\frac{n-1}{p-1}}} =1.
\end{equation}
Applying the implicit function theorem to $\mathcal F(\varepsilon,z):=z^{p-1}+\varepsilon z^{q-1}-1/(n\omega_n)$ in a neighborhood of $(0,c_p(n))$ allows us to sharpen this profile, yielding the expansion
\[
g(r)=c_p(n)-\frac{c_p(n)^{q-p+1}}{p-1}r^\beta+O(r^{2\beta})\qquad\text{as }r\searrow0. 
\]
Substituting this back into \eqref{normalized.fundamental.profile} gives the refined gradient asymptotics
\begin{equation}\label{fundamental.gradient.expansion}
f_{p,q}(r)=c_p(n)r^{-\frac{n-1}{p-1}}\left[1-\frac{c_p(n)^{q-p}}{p-1}r^\beta+O(r^{2\beta})\right]\qquad\text{as }r\searrow0,
\end{equation}
with $c_p(n)$ as defined in \eqref{cpn}.

In view of \eqref{eqn:fund-sol-asympt}, local Lebesgue integrability holds for every positive exponent below the critical threshold as in \eqref{int.sharp.Lebesgue}, while the endpoint only allows for the weak-$L^s$ integrability
\[
D\mathcal G_{p,q}\in L^{\frac{n(p-1)}{n-1},\infty}_{\rm loc}.
\]
Energy integrability is governed by the exact same asymptotics: since
\[
H(x,D\mathcal G_{p,q})\approx |x|^{-\frac{p(n-1)}{p-1}},\qquad\int_{B_\varepsilon}H(x,D\mathcal G_{p,q})\dx\approx\int_0^\varepsilon r^{-\frac{n-1}{p-1}}\,dr,
\]
the integral remains finite if and only if $p>n$ and when $p\leq n$, even the fundamental solution generated by a Dirac mass escapes the natural energy space - precisely motivating our reliance on the approximation framework of Definition~\ref{def.SOLA}.

\medskip

Along this fundamental profile, the relative contribution of the $q$-term to the flux vanishes at the origin:
\[
\frac{r^\alpha f_{p,q}(r)^{q-1}}{f_{p,q}(r)^{p-1}}=r^\alpha f_{p,q}(r)^{q-p}=r^\beta g(r)^{q-p}\xrightarrow[r\searrow0]{}0.
\]
Thus, provided $\beta>0$, the $p$-term dictates the behavior near the pole. By contrast, the geometric phase classification relies exclusively on the weight $a(x)=|x|^\alpha$: because
\[
\sup_{B_\rho(0)}a(\,\cdot\,)=\rho^\alpha=[a]_{0,\alpha}\rho^\alpha,
\]
every ball centered at the origin satisfies the $p$-phase condition \eqref{p.phase.sup.condition} for any $M\geq1$, irrespective of the sign of $\beta$.

\subsection{Asymptotics of the fundamental solution itself, under the subcritical balance condition}\label{asy.fund.sol}

Integrating \eqref{fundamental.gradient.expansion} recovers the classical trichotomy associated with the $p$-Laplacian fundamental solution. Setting
\[
C_p(n):=\frac{p-1}{n-p}(n\omega_n)^{-\frac1{p-1}}=\frac{p-1}{n-p} c_p(n)
\]
for $p\neq n$, we distinguish three regimes:
\begin{itemize}
\item if $1<p<n$, then, up to an additive normalization constant,
\begin{equation}\label{fundamental.function.p<n}
\mathcal G_{p,q}(x) \sim C_p(n) |x|^{\frac{p-n}{p-1}} \qquad\text{as }x\to0,
\end{equation}
in the sense that their quotient tends to one;
\item if $p=n$, then
\begin{equation}\label{fundamental.function.p=n}
\mathcal G_{n,q}(x) \sim c_n(n)\log\frac1{|x|} \qquad\text{as }x\to0;
\end{equation}
\item if $p>n$, expansion \eqref{fundamental.gradient.expansion} ensures $\int_0^1 f_{p,q}(r)\,dr<\infty$, meaning $\mathcal G_{p,q}$ possesses a finite limit at the origin, denoted by $\mathcal G_{p,q}(0)$. More precisely,
\begin{equation}\label{fundamental.function.p>n}
\mathcal G_{p,q}(0)-\mathcal G_{p,q}(x) \sim -C_p(n)|x|^{\frac{p-n}{p-1}} \qquad\text{as }x\to0.
\end{equation}
\end{itemize}
When $p>n$, relation \eqref{fundamental.function.p>n} identifies $(p-n)/(p-1)$ as the optimal Hölder exponent at the origin. A direct bridge to the potential estimates stems from \eqref{fundamental.gradient.inverse}, which provides, for $0<r<R$,
\[
F_{p,q}(r)-F_{p,q}(R)=\int_r^R h_{\texttt r}(\rho,\,\cdot\,)^{-1}\biggl(\frac1{n\omega_n\rho^{n-1}}\biggr)\,d\rho.
\]
Recalling that $\sup_{B_\rho(0)}a(\,\cdot\,)=\rho^\alpha$, this identity represents the annular counterpart of the scale-dependent potential introduced in Remark~\ref{rem:potential.u.sup}, up to a fixed geometric constant.

These structural asymptotics correspond directly to the endpoint scales for the solutions established in Theorem~\ref{existence.SOLA}: for $p<n$, we have
\[
\mathcal G_{p,q}\in L_{\rm loc}^{\frac{n(p-1)}{n-p},\infty},\qquad\mathcal G_{p,q}\notin L_{\rm loc}^{\frac{n(p-1)}{n-p}}.
\]
In the borderline case $p=n$, the logarithmic singularity implies that the measures of the superlevel sets decay exponentially on every fixed ball, whereas for $p>n$ the solution remains locally bounded.

\subsection{The critical and supercritical balances}
In this paragraph we relax \eqref{bounds} to examine the critical and supercritical balances: when $\beta=0$, equation \eqref{normalized.fundamental.equation} degenerates to $g(r)^{p-1}+g(r)^{q-1}=1/[n\omega_n]$, reducing $g$ to the unique positive root $g_*\in(0,c_p(n))$ of the algebraic condition
\[
g_*^{p-1}+g_*^{q-1}=\frac1{n\omega_n}.
\]
Consequently, the size of the radial gradient is
\[
f_{p,q}(r) = g_*r^{-\frac{n-1}{p-1}}.
\]
At this critical balance, the scaling exponent coincides with that of the $p$-Laplacian fundamental solution, while the prefactor is modulated by the interplay between the two growth terms. 

On the other hand, if $\beta<0$, the weighted $q$-phase takes over: defining $g_q(r):=r^{-m_q}f_{p,q}(r)$, equation \eqref{flux.fundamental} transforms into
\[
g_q(r)^{q-1} + r^\delta g_q(r)^{p-1} =\frac1{n\omega_n},\qquad\text{where}\quad \delta := n-1 - \frac{(p-1)(n+\alpha-1)}{q-1} = -\frac{p-1}{q-1}\beta >0.
\]
This implies
\[
g_q(r)\xrightarrow[r \searrow 0]{}(n\omega_n)^{-\frac{1}{q-1}} \qquad\implies\qquad f_{p,q}(r) \sim c_q(n)r^{-\frac{n+\alpha-1}{q-1}}.
\]

\section{Preliminaries}\label{notation}

\subsection{Notation}

We denote by
\[
B_r(x_0):=\{x \in \R^n \, : \,  |x-x_0|< r\}
\]
the open ball with center $x_0$ and radius $r>0$; when the center is irrelevant or clear from the context, we shall simply avoid writing it after specifying $B_r\equiv B_r(x_0)$. Its volume will be $\omega_nr^n$, where $\omega_n=|B_1|$ when $B_1$ is the unit ball in $\R^n$. For a measurable set $\mathcal B\subset\R^n$ with $0<|\mathcal B|<\infty$ and a map $g\in L^1(\mathcal B;\R^k)$, $k\geq1$, we shall denote by
\[
   (g)_{\mathcal B} \equiv \mean{\mathcal B}  g(x) \dx  := \frac{1}{|\mathcal B|}\int_{\mathcal B}  g(x) \dx
\]
its integral average.

\medskip

The Sobolev exponent $\esp^*$ is $n\esp/(n-\esp)$ if $\esp<n$ or a fixed number larger than $\esp$, in the case $\esp\geq n$. The H\"older conjugate of an exponent $\esp>1$ will be the number
\[
\esp'=\frac \esp {\esp-1}>1.
\]
Notice that $\esp\mapsto \esp'$ is a decreasing map, surjective onto $(1,+\infty)$, with $\esp=2$ as its fixed point. For $\esp>1$, $\chi_{(-\infty,2)}(\esp)$ will be the characteristic function of the open half-line $(-\infty,2)$:
\[
\chi_{(-\infty,2)}(\esp)=
\begin{cases}
 0\qquad&\text{if $\esp\geq2$}\\[2mm]
 1&\text{if $\esp<2$}
\end{cases};
\]
similarly for $\chi_{[2,+\infty)}(\esp)$. Terms multiplied by a vanishing characteristic function are understood to be omitted, without evaluating the remaining factors; powers with exponent zero are taken to equal one also at the origin. $\N$ is the set of the natural numbers, where we do not include the zero, and we call ``interval of natural numbers'' any interval of $\R$ intersected with $\N$.

\medskip

We shall denote by $c$ a general positive constant possibly varying from line to line; special occurrences will be denoted by $c_1,  c_*, \bar c$ or the like; all such constants will always be {\em greater than or equal to one} and therefore appear in the denominator in lower bounds. Relevant dependencies on parameters will be emphasized using parentheses, for instance, $c_1\equiv c_1(n,p,\nu)$ means that $c_1$  depends on $n,p$ and $\nu$; if no dependence is specified, this will usually mean that the constant is numerical and does not depend on any of the parameters, even if sometimes this will mean that the dependence of the constant is obvious from the context. We shall use the short notation $A\lesssim_S B$ (often when dealing with $N$-functions) to mean that $A\leq c B$ for a constant $c\geq 1$ depending only on the parameters in the set $S$ (where $A,B\geq0$). It will also be convenient to write $A\approx_S B$ if the two inequalities $A\lesssim_S B \lesssim_S A$ hold at the same time. When the constants are numerical,  we shall simply use the symbols $\lesssim,\approx$.

\medskip

Finally, and most importantly, we shall employ the compact notation
\[
\data:=\bigl\{n,p,q,\nu,L,\alpha,{[a]}_{0,\alpha}\bigr\}
\]
and the geometric data of the problem will be compactly collected with
\begin{equation}\label{data.gdata}
\gdata:=\bigl\{|\Omega|,|\mu|(\Omega),d_{\Omega_0},\mathcal M\bigr\},
\end{equation}
where $\Omega_0$ and $d_{\Omega_0}$ are as in \eqref{Omega0}, and $\mathcal M$ as in \eqref{energy.M}.

\medskip

For a ball $B\subset\Omega$ and a map $v\in L^1(B;\R^N)$, it is well known that the quantity 
\[
\mean{B}|v-(v)_B|\dx
\] 
is called the ($L^1$-)excess of $v$ on $B$. By the triangle inequality and Jensen's inequality, it satisfies
\begin{equation}\label{exc.elementary}
\mean{B}|v-(v)_B|\dx\leq2\mean{B}|v-k|\dx\qquad\text{for every }k\in\R^N.
\end{equation}

\subsection{The \texorpdfstring{$\boldsymbol{N}$}{N}-function setting}\label{N-function}

In the following we shall need some results which apply in a more general setting, the one of $N$- and Orlicz functions. 

\vs
We start by considering a function $g$ such that 
\begin{align}
g:[0,\infty) \to [0,\infty),\quad g(0)=0, \quad &g \in C^0([0,\infty)) \cap C^1((0,\infty)),\label{prop-g}\\
&\qquad \qquad 0< \frac1{c_g} \leq\mathcal O_g(t)=\frac{g'(t)t}{g(t)}=\frac{\frac{d}{dt} \log g(t)}{\frac d{dt} \log t}\leq c_g<\infty\label{prop-g-Orlicz}
\end{align}
for some $c_g\geq1$; the model we have in mind is $g(t)=t^\alpha$ for $\alpha>0$ and, as mentioned before, the quantity $\mathcal O_g$ quantifies how close $g$ is, pointwise, to a pure power;  the logarithmic slope $\mathcal O_g$ is identically equal to $\alpha$ on an interval if and only if $g(t)=ct^\alpha$ there for some $c>0$. This is why it coincides with the notion of elasticity used pervasively in economics. Under these assumptions $g$ is strictly increasing and surjective onto $[0,\infty)$, since $t\mapsto t^{-1/c_g}g(t)$ is non-decreasing and therefore $g(t)\geq g(1)t^{1/c_g}$ for every $t\geq1$. Therefore if $g$ is as above, its inverse function $g^{-1}:[0,\infty) \to [0,\infty)$ exists and it satisfies \eqref{prop-g}-\eqref{prop-g-Orlicz} with $c_{g^{-1}}=c_g$, as $g^{-1}\in C^1((0,\infty))$ due to $g'(t)>0$. Moreover, if $g_1$ and $g_2$ both are as $g$ and they both satisfy \eqref{prop-g}-\eqref{prop-g-Orlicz} with constants, respectively, $c_{g_1}$ and $c_{g_2}$, we have that their sum $g_1+g_2$, their product $g_1\cdot g_2$ and their composition $\tilde g=g_1\circ g_2$ satisfy \eqref{prop-g}-\eqref{prop-g-Orlicz} with constants, respectively, $c_{g_1+g_2}=\max\{c_{g_1},c_{g_2}\}$, $c_{g_1g_2}=c_{g_1}+c_{g_2}$ and $c_{g_1\circ g_2}=c_{g_1}c_{g_2}$ (see for instance \cite{BCMU}). If $g$ is as above and $0<\alpha<c_g^{-1}$, then $t\mapsto g(t)/t^\alpha$ extends continuously to $t=0$ with value zero and satisfies \eqref{prop-g}-\eqref{prop-g-Orlicz} with constant $\max\{c_g-\alpha,(c_g^{-1}-\alpha)^{-1}\}$.

\medskip

Often (but not always) a function $g$ as considered in the previous lines will be the derivative of a regular convex function. Let us therefore also consider a function $\varphi$ such that 
\begin{multline}\label{prop-varphi}
\varphi : [0,\infty) \to [0,\infty), \qquad \varphi(0)= \varphi'(0)= 0, \qquad \varphi \in C^1([0,\infty)) \cap C^2((0,\infty)), \\
0<g_0\leq \mathcal O_{\varphi'}(t)=\frac{\varphi''(t)t}{\varphi'(t)}\leq g_1<\infty
\end{multline}
- the basic example being $t^p/p$, $p>1$; $g_0=g_1=p-1$. Integrating by parts one immediately sees that
\begin{equation}\label{Or.primo}
1<g_0+1 \leq \frac{\varphi'(t)t}{\varphi(t)}\leq g_1+1
\end{equation}
holds for every $t>0$.  Thus under the assumptions in \eqref{prop-varphi} $\varphi'$ turns out to be non-decreasing, $t\mapsto \varphi(t)/t^{g_0+1}$ too, $t\mapsto \varphi(t)/t^{g_1+1}$ non-increasing,
\[
\varphi\text{\ convex \quad with} \quad \lim_{t\to 0^+}\frac{\varphi(t)}t=0,\quad \lim_{t\to\infty}\frac{\varphi(t)}t=+\infty;
\]
a continuous function over $[0,+\infty)$ satisfying the properties in the last display is usually called an $N$-function. Due to these additional properties, it is possible to define, for a function $\varphi$ satisfying \eqref{prop-varphi}, its convex conjugate $\widetilde\varphi$:
\begin{equation}\label{conv.conj}
\widetilde\varphi(t)=\sup_{s>0} \bigl[ts-\varphi(s)\bigr]=\max_{s\geq 0} \bigl[ts-\varphi(s)\bigr]=\int_0^t(\varphi')^{-1}(\tau)\,d\tau;
\end{equation}
notice that a  direct computation shows that, if $\varphi$ satisfies \eqref{prop-varphi}, then its convex conjugate $\widetilde\varphi$ satisfies
\[
\frac1{g_1} \leq \mathcal O_{\widetilde\varphi'}(t) = \frac{\widetilde\varphi''(t)t} {\widetilde\varphi'(t)} \leq \frac1{g_0} \qquad\text{for every }t>0
\]
and consequently, by \eqref{Or.primo},
\[
1< \frac{g_1+1}{g_1} \leq \mathcal O_{\widetilde\varphi}(t) = \frac{\widetilde\varphi'(t)t} {\widetilde\varphi(t)} \leq \frac{g_0+1}{g_0} <\infty.
\]
Further basic properties of the convex conjugate are the following ones, valid for every $t\in(0,+\infty)$:
\begin{equation}\label{conj.inverse}
t\leq \varphi^{-1}(t)\widetilde \varphi^{-1}(t)\leq 2t,\qquad \varphi(t)\leq \widetilde\varphi\Bigl(2\frac{\varphi(t)}{t}\Bigr),\qquad \widetilde\varphi\Bigl(\frac{\varphi(t)}{t}\Bigr)\leq\varphi(t);
\end{equation}
note that the convex conjugation is an involution on the class of convex lower semicontinuous functions. Finally, the Young-type inequality
\begin{equation}\label{Young-N}
 ab\leq \varphi\bigl(\varepsilon^{\frac1{1+g_\smallzero}}a\bigr) +\widetilde \varphi\bigl(\varepsilon^{-\frac1{1+g_\smallzero}}b\bigr)\leq \varepsilon \varphi(a)+c(g_\smallzero,g_1,\varepsilon)\widetilde \varphi(b) 
\end{equation}
for all $a,b\geq0$ and $\varepsilon\in(0,1)$, holds as a consequence of the monotonicity of $t\mapsto \varphi(t)/t^{g_0+1}$ and the estimate above.  As a reference for the properties above we suggest \cite[Chapter 8]{RAd}, while for the properties we are going to describe now a good reference is \cite{DE}; see also \cite[Section 3.2]{BCMU}.

\vs

As a subsequent step, one defines the auxiliary vector fields $A_\varphi, V_\varphi\colon  \R^\ell \to \R^\ell$, $\ell\geq1$, by
\begin{equation}\label{V.varphi}
A_\varphi(z):= \frac{\varphi'(|z|)}{|z|} \,z  ,\qquad V_\varphi(z):= \left( \frac{\varphi'(|z|)}{|z|} \right)^{\frac12}z  
\end{equation}
whenever $z \in \R^\ell\setminus\{0\}$; we set $A_\varphi(0)=V_\varphi(0)=0$. By \eqref{prop-varphi} they turn out to be bijections of $\R^\ell$. Thanks to \eqref{Or.primo} we have that $|V_\varphi(z)|^2$ is comparable to $\varphi(|z|)$, namely for every $z\in \R^\ell$
\begin{equation}\label{V-phi}
\bigl|V_\varphi(z)\bigr|^2\approx_{g_0,g_1} \varphi(|z|).
\end{equation}
$V_\varphi$ is very useful as it describes the monotonicity properties of the map $z\mapsto\frac{\varphi'(|z|)}{|z|}z$; indeed, for $z_1, z_2 \in \R^\ell$ it holds that
\begin{equation}\label{monotonicity}
 \bigl\langle \frac{\varphi'(|z_1|)}{|z_1|} z_1-\frac{\varphi'(|z_2|)}{|z_2|} z_2, z_1-z_2\bigr\rangle=\langle A_\varphi(z_1)-A_\varphi(z_2), z_1-z_2\rangle\approx_{g_0,g_1}\big|V_\varphi(z_1)-V_\varphi(z_2)\big|^2.
\end{equation}
The previous relation is a consequence of the following algebraic fact: if $\varphi$ satisfies \eqref{prop-varphi}, then
\begin{equation}\label{keep}
\varphi''(|z_1|+|z_2|){|z_1-z_2|}^2 \approx_{g_0,g_1}\frac{\varphi'(|z_1|+|z_2|)}{|z_1|+|z_2|}{|z_1-z_2|}^2 \approx_{g_0,g_1}{\big|V_{\varphi}(z_1)-V_{\varphi}(z_2)\big|}^2.
\end{equation}
We also have
\begin{equation}\label{keep.a}
\big|A_\varphi(z_1)-A_\varphi(z_2)|\approx_{g_0,g_1}\varphi''(|z_1|+|z_2|)|z_1-z_2| \approx_{g_0,g_1}\frac{\varphi'(|z_1|+|z_2|)}{|z_1|+|z_2|}|z_1-z_2| ;
\end{equation}
the last two inequalities hold for $|z_1|+|z_2|\neq 0$.  These equivalences follow by integrating the growth and ellipticity bounds for $\partial_z A_\varphi$ along the segment joining $z_1$ and $z_2$, using the segment estimate in \cite[Lemma 20]{DE} and \eqref{prop-varphi}.

\subsection{Frozen functions and basic properties}\label{frozen.functions}

We now specialize the previous facts to the functions
\begin{equation}\label{hH.frozen}
\begin{cases}
H_x(t):=\displaystyle{H(x,t)=\frac{t^p}{p}+a(x)\frac{t^q}{q}}\\[3mm] 
h_x(t):=H_x'(t)=t^{p-1}+a(x)t^{q-1}
\end{cases}
\end{equation}
for $t\in[0,\infty)$ and $x\in\Omega$; recalling the definitions in \eqref{H.def} and \eqref{h.def}, we shall indeed often regard the spatial variable as fixed. When the point $x\in\Omega$ is fixed, we shall write $h_{x}=H_{x}'$: for $t>0$,
\[
h_{x}'(t) = (p-1)t^{p-2}+a(x)(q-1)t^{q-2}
\]
and
\[
h_{x}''(t) = (p-1)(p-2)t^{p-3} + a(x)(q-1)(q-2)t^{q-3}.
\]
Both the maps $H_x$ and $h_x$ are continuous and strictly increasing on $[0,+\infty)$, vanish at zero and tend to $+\infty$ as $t\to+\infty$; their inverse functions $H_x^{-1},h_x^{-1}:[0,+\infty)\to[0,+\infty)$ are therefore well defined and strictly increasing. The functions $H_x$ and $h_x$ are also subadditive up to a multiplicative constant: for every $x\in\Omega$ and $z_1,z_2\in\mathbb R^n$,
\[
H_x(z_1+z_2)\leq 2^{q-1}\bigl[H_x(z_1)+H_x(z_2)\bigr],\qquad h_x(z_1+z_2)\leq 2^{\max\{q-2,0\}}\bigl[h_x(z_1)+h_x(z_2)\bigr].
\]
For brevity, we shall refer to these inequalities simply as ``subadditivity''. When $q\leq2$, notice that the inequality for $h$ holds with constant one. Secondly, we notice that 
\begin{equation}\label{below.function}
\max\big\{t^{p-1},a(x)t^{q-1}\big\}\leq h_{x}(t),\qquad \max\Bigl\{\frac{t^p}p,\frac{a(x)t^q}q\Bigr\}\leq H_{x}(t)
\end{equation}
from which it follows, when $a(x)>0$,
\begin{equation}\label{above.below.inverseh}
2^{-\frac1{p-1}}\min \biggl\{t^{\frac1{p-1}},\Big[\frac t{a(x)}\Big]^{\frac1{q-1}}\biggr\}\leq h^{-1}_{x}(t)\leq \min \biggl\{t^{\frac1{p-1}},\Big[\frac t{a(x)}\Big]^{\frac1{q-1}}\biggr\}
\end{equation}
and
\begin{equation}\label{above.below.inverseH}
2^{-\frac1p}\min \biggl\{{(pt)}^{\frac1p},\Big[\frac{qt}{a(x)}\Big]^{\frac1{q}}\biggr\}\leq H^{-1}_{x}(t)\leq \min \biggl\{{(pt)}^{\frac1p},\Big[\frac{qt}{a(x)}\Big]^{\frac1{q}}\biggr\}:
\end{equation}
indeed, calling $m_h,m_H$ the two minima, one has $t\leq h_x(m_h)\leq2t$, $t\leq H_x(m_H)\leq2t$ and 
\[
h_x\bigl(2^{-\frac1{p-1}}m_h\bigr)\leq \frac12 h_x(m_h)\leq t, \qquad H_x\bigl(2^{-\frac1p}m_H\bigr)\leq \frac12 H_x(m_H)\leq t.
\]
When $a(x)=0$, we have $h_x^{-1}(t)=t^{\frac1{p-1}}$ and $H_x^{-1}(t)=(pt)^{\frac1p}$ for every $t\geq0$; in this case, \eqref{above.below.inverseh} and \eqref{above.below.inverseH} remain valid with the second term omitted from each minimum. Now a direct computation gives
\begin{equation}\label{Orlicz_trivial}
p-1\leq\mathcal O_{h_x}(t)=\frac{h_x'(t)t}{h_x(t)}\leq q-1,\qquad p\leq\mathcal O_{H_x}(t)=\frac{H_x'(t)t}{H_x(t)}\leq q
\end{equation}
for every $t>0$; in particular,
\begin{equation}\label{H.h.comparison}
\frac1q t h_x(t) \leq H_x(t) \leq \frac1p t h_x(t) \qquad\text{for every }t\geq0.
\end{equation}
The fact that the maps $t\mapsto \frac{H_x(t)}{t^p}$ and $t\mapsto \frac{h_x(t)}{t^{p-1}}$ are non-decreasing, whereas $t\mapsto \frac{H_x(t)}{t^q}$ and $t\mapsto \frac{h_x(t)}{t^{q-1}}$ are non-increasing, gives
\begin{equation}\label{constant.h}
\min\{\lambda^p,\lambda^q\}\leq \frac{H_x(\lambda t)}{H_x(t)}\leq \max\{\lambda^p,\lambda^q\},\quad \min\big\{\lambda^\frac1p,\lambda^\frac1q\big\}\leq \frac{H_x^{-1}(\lambda t)}{H_x^{-1}(t)}\leq \max\big\{\lambda^\frac1p,\lambda^\frac1q\big\} 
\end{equation}
and 
\begin{equation}\label{constant.hinverse}
\min\{\lambda^{p-1},\lambda^{q-1}\}\leq \frac{h_x(\lambda t)}{h_x(t)}\leq \max\{\lambda^{p-1},\lambda^{q-1}\},\quad \min\{\lambda^\frac1{p-1},\lambda^\frac1{q-1}\}\leq \frac{h_x^{-1}(\lambda t)}{h_x^{-1}(t)}\leq \max\{\lambda^\frac1{p-1},\lambda^\frac1{q-1}\}
\end{equation}
for every $t,\lambda>0$. 

\subsection{Scaling properties of inverse compositions}\label{inverse.compositions}

We next collect some scaling properties of the compositions involving the inverse frozen functions: these estimates hold uniformly with respect to $x_0\in\Omega$. We set
\begin{equation}\label{g_x_0}
g_x := h_x\circ H_x^{-1}, \qquad G_x := H_x\circ h_x^{-1}. 
\end{equation}
Both functions are increasing on $[0,+\infty)$, and they are mutually inverse, as $G_x=g_x^{-1}$. We first study $g_x$: for $t,\lambda>0$,  using \eqref{H.h.comparison} we have
\[
\frac pq \lambda\frac{H_x^{-1}(t)}{H_x^{-1}(\lambda t)} \leq  \frac{g_x(\lambda t)}{g_x(t)} = \frac{h_x(H_x^{-1}(\lambda t))}{h_x(H_x^{-1}(t))} \leq \frac qp \lambda\frac{H_x^{-1}(t)}{H_x^{-1}(\lambda t)};
\]
inserting here the scaling properties of $H_x^{-1}$ in \eqref{constant.h} yields
\begin{equation}\label{g_x.lambda}
\frac1c \min\bigl\{ \lambda^{\frac1{p'}},\lambda^\frac1{q'} \bigr\} \leq \frac{g_x(\lambda t)}{g_x(t)} \leq c \max\bigl\{ \lambda^\frac1{p'},\lambda^\frac1{q'} \bigr\},
\end{equation}
where $c=c(p,q)\geq1$. As a consequence, for every $\varsigma>0$ and every $t_1,t_2\geq0$,
\begin{equation*}
\bigl[g_x(t_1+t_2)\bigr]^\varsigma \leq c \Bigl( \bigl[g_x(t_1)\bigr]^\varsigma + \bigl[g_x(t_2)\bigr]^\varsigma \Bigr),
\end{equation*}
where now $c=c(p,q,\varsigma)$. Indeed, by monotonicity, $g_x(t_1+t_2) \leq g_x(2\max\{t_1,t_2\})$ and \eqref{g_x.lambda}, applied with $\lambda=2$, gives the claim. We then consider $G_x$: the scaling properties of $h_x^{-1}$ give
\[
\min\bigl\{ \lambda^\frac1{p-1},\lambda^\frac1{q-1} \bigr\} \leq \frac{h_x^{-1}(\lambda t)}{h_x^{-1}(t)} \leq \max\bigl\{ \lambda^\frac1{p-1},\lambda^\frac1{q-1} \bigr\}
\]
and using again \eqref{H.h.comparison} we find
\[
\frac pq \lambda\frac{h_x^{-1}(\lambda t)}{h_x^{-1}(t)} \leq \frac{H_x(h_x^{-1}(\lambda t))}{H_x(h_x^{-1}(t))} \leq \frac qp \lambda\frac{h_x^{-1}(\lambda t)}{h_x^{-1}(t)};
\]
therefore,
\begin{equation}\label{G_x.lambda}
\frac1c \min\big\{ \lambda^{p'},\lambda^{q'} \big\} \leq \frac{G_x(\lambda t)}{G_x(t)} \leq c \max\big\{ \lambda^{p'},\lambda^{q'} \big\},
\end{equation}
for every $t>0$ and $\lambda>0$, where $c=c(p,q)\geq1$. These properties, together with those of the previous paragraph, allow us to move constants in and out of estimates, and in general we will do so without saying it every time; when the exponent of the quantity being pulled out will matter, however, we will state explicitly which property we are employing.

\subsection{Young inequalities and nonlinear vector fields}

Finally, the pointwise Young-type inequality
\begin{equation}\label{Young.H}
h_x(t)s \leq \varepsilon H_x(t) + c(p,q,\varepsilon)H_x(s)
\end{equation}
holds for every $s,t\geq0$ and every $\varepsilon>0$. Indeed, applying the classical Young inequality separately to the $p$- and $q$-terms, we obtain
\[
t^{p-1}s\leq\frac{\varepsilon}{p}t^p+c(p,\varepsilon)\frac{s^p}{p} \qquad\text{and}\qquad a(x)t^{q-1}s\leq\frac{\varepsilon}{q}a(x)t^q+c(q,\varepsilon)a(x)\frac{s^q}{q};
\]
summing the previous inequalities and renaming $\eps$ yields \eqref{Young.H}.

\vs

We shall denote by $a_\esp(\,\cdot\,)$ and $V_\esp(\,\cdot\,)$, $\esp>1$, the vector fields defined in \eqref{V.varphi} associated with $\varphi(t)= t^\esp/\esp$, namely 
\begin{equation}\label{V.fields}
a_\esp(z)=|z|^{\esp-2}z\qquad\text{and}\qquad V_\esp(z):= |z|^{\frac{\esp-2}2}z. 
\end{equation}
 Another case that will occur frequently is $\varphi(t)=\frac{t^p}p+a(x)\frac{t^q}q$, $x\in\Omega$, and in this case we shall use the shorthand notation $V_{a(x)}:=V_\varphi$. In the power case, \eqref{keep} reads as
\begin{equation}\label{diffV}
\big|V_\esp(z_1)-V_\esp(z_2)\big|\approx_\esp\bigl(|z_1|+|z_2|\bigr)^{\frac{\esp-2}{2}}|z_1-z_2|
\end{equation}
for every $z_1,z_2\in\R^\ell$, where the right-hand side is defined by continuity to be zero when $z_1=z_2=0$. Consequently, if $\esp\geq2$, then
\begin{equation}\label{diffV2}
|z_1-z_2|\lesssim_\esp\big|V_\esp(z_1)-V_\esp(z_2)\big|^{\frac2\esp},
\end{equation}
whereas, if $1<\esp<2$, then, using Young's inequality,
\begin{equation}\label{diffVp}
|z_1-z_2|\lesssim_\esp\big|V_\esp(z_1)-V_\esp(z_2)\big|^{\frac2\esp}+|z_2|^{\frac{2-\esp}{2}}\big|V_\esp(z_1)-V_\esp(z_2)\big|.
\end{equation}
Applying the previous inequalities with $\esp=p$ and $\esp=q$, and raising them to the corresponding powers, we obtain the composite estimate
\begin{multline}\label{diffVHx0}
H_x(|z_1-z_2|)\lesssim_{p,q}\big|V_{a(x)}(z_1)-V_{a(x)}(z_2)\big|^2+\chi_{(-\infty,2)}(p)|z_2|^{\frac{(2-p)p}{2}}\big|V_p(z_1)-V_p(z_2)\big|^p\\+\chi_{(-\infty,2)}(q)a(x)|z_2|^{\frac{(2-q)q}{2}}\big|V_q(z_1)-V_q(z_2)\big|^q.
\end{multline}
Young's inequality therefore yields, for every $\varepsilon\in(0,1]$,
\begin{equation}\label{diffVHx0.varepsilon}
H_{x}(|z_1-z_2|) \leq c(p,q,\varepsilon) {\big|V_{a(x)}(z_1)-V_{a(x)}(z_2)\big|}^2 +\varepsilon\, H_x(|z_2|).
\end{equation}
Finally, applying \eqref{keep} to $H_x$ and separately to the power functions $t^p/p$ and $t^q/q$, we obtain
\begin{equation}\label{Vx0}
\big|V_{a(x)}(z_1)-V_{a(x)}(z_2)\big|^2\approx_{p,q}\big|V_p(z_1)-V_p(z_2)\big|^2+a(x)\big|V_q(z_1)-V_q(z_2)\big|^2.
\end{equation}

\medskip

\begin{lemma}\label{lem:a-a}
Let $A$ satisfy \eqref{A.strong}. Then, for almost every $x\in\Omega$ and every $z_1,z_2\in\R^n$, there holds
\begin{equation}\label{flux.difference.general}
\bigl|A(x,z_1)-A(x,z_2)\bigr|\leq c\,h_x\bigl(|z_1-z_2|\bigr)+c\Big[|z_2|^{p-2}\chi_{[2,\infty)}(p)+a(x)|z_2|^{q-2}\chi_{[2,\infty)}(q)\Big]|z_1-z_2|;
\end{equation}
moreover, if $1<p\leq q\leq2$, then for every $z_1,z_2\in\mathbb R^n$ with $z_2\neq0$,
\begin{equation}\label{flux.difference.singular}
\bigl|A(x,z_1)-A(x,z_2)\bigr|\leq c\bigl[|z_2|^{p-2}+a(x)|z_2|^{q-2}\bigr]|z_1-z_2|.
\end{equation}
Finally, if $1<p<2<q$, then for every $z_1,z_2\in\mathbb R^n$ with $|z_1|+|z_2|>0$, we have
\begin{equation}\label{flux.difference.mixed}
|A(x,z_1)-A(x,z_2)|\leq c{\bigl(|z_1|+|z_2|\bigr)}^{p-2}|z_1-z_2|+c\,a(x)\Bigl[|z_1-z_2|^{q-2}+|z_2|^{q-2}\Bigr]|z_1-z_2|.
\end{equation}
All the constants depend only on $p,q$ and $L$.
\end{lemma}

\begin{proof}
Fix $x$ for which the structural assumptions hold. The case $z_1=z_2$ is trivial, so we assume $z_1\neq z_2$; by the first condition in \eqref{A.strong}, we have
\[
|A(x,z)|\leq Lh_x(|z|),\qquad|\partial_z A(x,z)|\leq L\bigl[|z|^{p-2}+a(x)|z|^{q-2}\bigr]
\]
for $z\neq0$; the first inequality also holds at $z=0$ since $A(x,0)=0$. 

\medskip

Now, if $|z_1-z_2|\geq |z_2|/2$, then $|z_1|+|z_2|\leq5|z_1-z_2|$; the first growth bound therefore gives
\begin{equation}\label{case1.AFestimate}
\bigl|A(x,z_1)-A(x,z_2)\bigr|\leq L\bigl[h_x(|z_1|)+h_x(|z_2|)\bigr]\leq c\,h_x(|z_1-z_2|), 
\end{equation}
which already proves \eqref{flux.difference.general} in this case. Suppose instead that $|z_1-z_2|<|z_2|/2$: for every $t\in[0,1]$, we have
\[
\frac{|z_2|}{2}\leq |z_2+t(z_1-z_2)|\leq \frac{3|z_2|}{2}.
\]
In particular, the segment joining $z_2$ to $z_1$ avoids the origin: the fundamental theorem of calculus and the bound on the derivative of $A$ then yield
\begin{equation}\label{A-A.derivative}
|A(x,z_1)-A(x,z_2)|\leq |z_1-z_2|\int_0^1\bigl|\partial_z A(x,z_2+t(z_1-z_2))\bigr|\,dt\leq c\bigl[|z_2|^{p-2}+a(x)|z_2|^{q-2}\bigr]|z_1-z_2|. 
\end{equation}
These two estimates hold in their respective cases for all $p,q>1$, regardless of their position relative to $2$; we now specialize them to prove the three assertions.

\medskip

To prove \eqref{flux.difference.general}, we note that for each $\esp\in\{p,q\}$ with $1<\esp<2$, the inequality $|z_1-z_2|<|z_2|$ implies $|z_2|^{\esp-2}|z_1-z_2|\leq |z_1-z_2|^{\esp-1}$; thus the terms with exponent smaller than two are bounded by $c\,h_x(|z_1-z_2|)$, while those with exponent at least two are retained in their present form. 

\medskip

To prove \eqref{flux.difference.singular}, we just employ \eqref{case1.AFestimate}, which holds when $|z_1-z_2|\geq |z_2|/2$: under this assumption we have $|z_1-z_2|^{\esp-1}\lesssim_\esp|z_2|^{\esp-2}|z_1-z_2|$ for $\esp\in\{p,q\},\esp\in(1,2]$ and, as the case $|z_1-z_2|<|z_2|/2$ has already been treated before, this concludes the proof.

\medskip

Finally, to prove \eqref{flux.difference.mixed}, if $|z_1-z_2|\geq|z_2|/2$, the estimate obtained in \eqref{case1.AFestimate}, together with $|z_1|+|z_2|\leq5|z_1-z_2|$ and $p<2$, gives
\[
|A(x,z_1)-A(x,z_2)|\leq c{\bigl(|z_1|+|z_2|\bigr)}^{p-2}|z_1-z_2|+c\,a(x)|z_1-z_2|^{q-1}.
\]
If instead $|z_1-z_2|<|z_2|/2$ then $|z_1|+|z_2|\leq(5/2)|z_2|$, so that
\[
|z_2|^{p-2}\lesssim_p(|z_1|+|z_2|)^{p-2},
\]
and combining this with \eqref{A-A.derivative} proves \eqref{flux.difference.mixed}.
\end{proof}

\medskip

We finally record the following scaling property, which holds for the model in \eqref{model.A} but also for general vector fields satisfying the differential double phase structure.

\begin{lemma}\label{lem:scaling.a}
Assume that $A\colon\Omega\times\R^n\to\R^n$ has the regularity specified in Paragraph~\ref{regular.A} and satisfies the first two conditions in \eqref{A.strong}. Then
\begin{equation}\label{scaling.a}
\big|A(x,\lambda z_1)-A(x,\lambda z_2)\big|\leq c\max\big\{|\lambda|^{p-1},|\lambda|^{q-1}\big\} \big|A(x,z_1)-A(x,z_2)\big|
\end{equation}
for every $z_1,z_2\in\R^n$, $\lambda\in\R$, and almost every $x\in\Omega$, where $c\equiv c(p,q,\nu,L)$.
\end{lemma}

\begin{proof}
For exponents $\ti>1$, the elementary estimate
\begin{equation}\label{segment.power}
\int_0^1\big|z_2+t(z_1-z_2)\big|^{\ti-2}\,dt\approx_\ti\bigl(|z_1|+|z_2|\bigr)^{\ti-2}
\end{equation}
holds whenever $|z_1|+|z_2|>0$. Set $z_t:=z_2+t(z_1-z_2)$; by the fundamental theorem of calculus and the upper bound in \eqref{A.strong}$_1$, we have
\[
\big|A(x,z_1)-A(x,z_2)\big|\leq\int_0^1\big|\partial_z A(x,z_t)\big||z_1-z_2|\,dt\leq L|z_1-z_2| \int_0^1\bigl[|z_t|^{p-2}+a(x)|z_t|^{q-2}\bigr]\,dt.
\]
On the other hand, the ellipticity assumption \eqref{A.strong}$_2$ gives
\begin{align*}
\big|A(x,z_1)-A(x,z_2)\big||z_1-z_2|&\geq\big\langle A(x,z_1)-A(x,z_2),z_1-z_2\big\rangle\\
&=\int_0^1\big\langle\partial_z A(x,z_t)(z_1-z_2),z_1-z_2\big\rangle\,dt\\
&\geq\nu|z_1-z_2|^2\int_0^1\bigl[|z_t|^{p-2}+a(x)|z_t|^{q-2}\bigr]\,dt.
\end{align*}
Consequently, by \eqref{segment.power},
\begin{equation}\label{difference.a.general}
\big|A(x,z_1)-A(x,z_2)\big|\approx_{p,q,\nu,L}\Big[\bigl(|z_1|+|z_2|\bigr)^{p-2}+a(x)\bigl(|z_1|+|z_2|\bigr)^{q-2}\Big]|z_1-z_2|.
\end{equation}
Applying the upper estimate in \eqref{difference.a.general} to $(\lambda z_1,\lambda z_2)$, we find
\begin{align*}
\big|A(x,\lambda z_1)-A(x,\lambda z_2)\big|&\leq c\Big[|\lambda|^{p-1}\bigl(|z_1|+|z_2|\bigr)^{p-2}+a(x)|\lambda|^{q-1}\bigl(|z_1|+|z_2|\bigr)^{q-2}\Big]|z_1-z_2|\\
&\leq c\max\big\{|\lambda|^{p-1},|\lambda|^{q-1}\big\}\Big[\bigl(|z_1|+|z_2|\bigr)^{p-2}+a(x)\bigl(|z_1|+|z_2|\bigr)^{q-2}\Big]|z_1-z_2|;
\end{align*}
the lower estimate in \eqref{difference.a.general} now yields \eqref{scaling.a}.
\end{proof}

\subsection{Regularity and harmonic approximation for autonomous equations}

Next, we collect two quantitative regularity results for autonomous equations with general Orlicz growth. Let $\varphi$ satisfy \eqref{prop-varphi} and consider weak solutions $w\in W^{1,\varphi}(B_{R/8})$ to
\begin{equation}\label{autonomous.equation}
-\divergence A_0(Dw)=0\qquad\text{in }B_{R/8},
\end{equation}
where $A_0\in C^0(\R^n;\R^n)\cap C^1(\R^n\setminus\{0\};\R^n)$, $A_0(0)=0$, and
\begin{equation}\label{Orlicz.structure}
|A_0(z)|+|z||\partial_z A_0(z)|\leq L\varphi'(|z|),\qquad \big\langle\partial_z A_0(z)\xi,\xi\big\rangle\geq\nu\frac{\varphi'(|z|)}{|z|}|\xi|^2
\end{equation}
for every $z\in\R^n\setminus\{0\}$ and $\xi\in\R^n$, with $0<\nu\leq L$: these assumptions are the autonomous counterparts of the strong structural conditions in \eqref{A.strong}. Moreover, for consistency with their subsequent applications, we formulate the  results only on the ball $B_{R/8}$. We impose them here since the forthcoming arguments involve pointwise values of the gradient and, in particular, require its local continuity and quantitative non-degeneracy estimates. Gradient regularity in this setting is already known from \cite[Lemma 5.1]{Lieb} and \cite[Theorem 4.1 \& Remark 5.6]{DKS}, and in full generality from \cite[Theorem 4.1]{Antonini} (see also \cite[Theorem~1.1]{OkUhlenbeck} for a recent development). With this mature theory in hand, we prove an excess-decay estimate for the flux $A_0(Dw)$ and infer pointwise non-degeneracy from integral lower bounds. Notably, every constant and exponent depends exclusively on $n,g_0,g_1,\nu$ and $L$, making them uniform across all frozen problems encountered below.

\begin{proposition}\label{varphi-harmonic.are.regular}
Suppose that $w\in W^{1,\varphi}(B_{R/8})$ is a weak energy solution to \eqref{autonomous.equation} under the assumptions \eqref{Orlicz.structure}. There exist $c\geq1$ and $\beta_0\in(0,1)$, depending only on $n,g_0,g_1,\nu$ and $L$, such that $A_0(Dw)\in C^{0,\beta_0}_{\rm loc}(B_{R/8})$,
\begin{equation}\label{lip.varphi.harmonic}
\sup_{B_{\rho/2}}|A_0(Dw)|\leq c\mean{B_\rho}|A_0(Dw)|\dx
\end{equation}
for every ball $B_\rho\subset B_{R/8}$, and
\begin{equation}\label{decay.varphi.harmonic}
\mean{B_\rho}\big|A_0(Dw)-(A_0(Dw))_{B_\rho}\big|\dx\leq c\Bigl(\frac{8\rho}{R}\Bigr)^{\beta_0}\mean{B_{R/8}}\big|A_0(Dw)-(A_0(Dw))_{B_{R/8}}\big|\dx
\end{equation}
for every ball $B_\rho$ concentric with $B_{R/8}$, with $0<\rho\leq R/8$.
\end{proposition}

\begin{proof}
All constants in this proof depend only on $n,g_0,g_1,\nu$ and $L$, unless otherwise specified. Let us stress that \eqref{keep} and \eqref{keep.a} continue to hold for $A_0$ satisfying \eqref{Orlicz.structure}, with constants depending on $g_0,g_1,\nu$ and $L$, not only for the model case in \eqref{V.varphi}, see \cite[Lemma 20]{DE}. In particular, $|A_0(z)|\approx\varphi'(|z|)$; moreover, the scaling properties of $\varphi'$ yield
\begin{equation}\label{autonomous.flux.modulus}
|A_0(z_1)-A_0(z_2)|\leq c\,\varphi'(M)\Bigl(\frac{|z_1-z_2|}{M}\Bigr)^{\min\{g_0,1\}}\qquad\text{if }|z_1|,|z_2|\leq M,\quad M>0.
\end{equation}
Indeed, when $g_0\leq1$, the upper bound in \eqref{keep.a} is at most $c\,\varphi'(M)M^{-g_0}(|z_1|+|z_2|)^{g_0-1}|z_1-z_2|$, and we use $|z_1-z_2|\leq|z_1|+|z_2|$; when $g_0\geq1$, the same bound is at most $c\,\varphi'(M)|z_1-z_2|/M$.

Theorem~4.1 from \cite{Antonini} then applies with the weight $t\mapsto\varphi'(t)/t$, whose Orlicz ratio lies between $g_0-1$ and $g_1-1$; note that the normalization required in the Theorem can be obtained by dividing both $\varphi$ and $A_0$ by $\varphi'(1)$. By that theorem, together with the usual covering and interpolation arguments for the supremum estimate, which allow us to lower the power of $Dw$ on the right-hand side, we obtain
\begin{equation}\label{autonomous.gradient.osc}
\sup_{B_{R/16}}|Dw|\leq c\biggl(\mean{B_{R/8}}|Dw|^{g_0}\dx\biggr)^{\frac1{g_\smallzero}}\quad\text{and}\quad \osc_{B_\rho}Dw\leq c\Bigl(\frac{\rho}{R}\Bigr)^\beta\mean{B_{R/32}}\bigl|Dw-(Dw)_{B_{R/32}}\bigr|\dx
\end{equation}
with $\beta\equiv \beta(n,g_0,g_1,\nu,L)\in(0,1)$, for $0<\rho\leq R/64$. To express the first estimate in terms of the flux, we use the elementary consequence $(s/t)^{g_0}\leq1+\varphi'(s)/\varphi'(t)$ of \eqref{prop-varphi}, valid for $s\geq0$ and $t>0$: choosing 
\[
t=(\varphi')^{-1}\biggl(\mean{B_{R/8}}|A_0(Dw)|\dx\biggr)
\]
(that we can suppose non-zero without loss of generality) and recalling $|A_0(z)|\approx\varphi'(|z|)$, we obtain
\begin{equation}\label{autonomous.gradient.from.flux}
\sup_{B_{R/16}}|Dw|\leq c\,(\varphi')^{-1}\biggl(\mean{B_{R/8}}|A_0(Dw)|\dx\biggr).
\end{equation}
Estimate \eqref{lip.varphi.harmonic} follows directly from \eqref{autonomous.gradient.from.flux}, applied on any ball $B_\rho\subset B_{R/8}$ with $8\rho$ in place of $R$, together with the scaling properties of $\varphi'$. Furthermore, for every $0<\rho\leq R/64$ we have
\begin{equation}\label{autonomous.flux.osc.size}
\osc_{B_\rho}A_0(Dw)\leq c\Bigl(\frac{\rho}{R}\Bigr)^{\beta_0}\mean{B_{R/8}}|A_0(Dw)|\dx,
\end{equation}
where $\beta_0=\beta\min\{g_0,1\}$; this bound is obtained by combining \eqref{autonomous.gradient.osc} and \eqref{autonomous.flux.modulus}, where we also used \eqref{autonomous.gradient.from.flux} to pass back to the average of $|A_0(Dw)|$. The claimed local H\"older continuity follows at once by applying \eqref{autonomous.flux.osc.size} on interior balls.

\medskip

Next, we need to trade the flux average in \eqref{autonomous.flux.osc.size} for its excess. Since strict monotonicity and coercivity make $A_0$ a bijection on $\R^n$, we can pick the unique vector $\xi\in\R^n$ satisfying
\[
A_0(\xi)={\bigl(A_0(Dw)\bigr)}_{B_{R/8}}.
\]
Let us fix $\varepsilon\in(0,1)$, depending solely on the structural parameters, which will be specified shortly. If
\[
\mean{B_{R/8}}|A_0(Dw)-A_0(\xi)|\dx\geq\varepsilon|A_0(\xi)|,
\]
then the triangle inequality gives
\[
\mean{B_{R/8}}|A_0(Dw)|\dx\leq\mean{B_{R/8}}|A_0(Dw)-A_0(\xi)|\dx+|A_0(\xi)|\leq\bigl(1+\varepsilon^{-1}\bigr)\mean{B_{R/8}}|A_0(Dw)-A_0(\xi)|\dx,
\]
and \eqref{decay.varphi.harmonic} follows directly from \eqref{autonomous.flux.osc.size} for $\rho\leq R/64$. Suppose instead that
\begin{equation}\label{autonomous.flux.small.excess}
\mean{B_{R/8}}|A_0(Dw)-A_0(\xi)|\dx<\varepsilon|A_0(\xi)|:
\end{equation}
again we can suppose $\xi\neq0$, and \eqref{autonomous.gradient.from.flux}, together with the triangle inequality, yields $\sup_{B_{R/16}}|Dw|\leq c|\xi|$. Consequently, $|\xi|\leq|Dw|+|\xi|\leq c|\xi|$ in $B_{R/16}$, so the lower bound in \eqref{keep.a} gives
\[
|Dw-\xi|\leq c\,\frac{|\xi|}{\varphi'(|\xi|)}|A_0(Dw)-A_0(\xi)|\qquad\text{in }B_{R/16}.
\]
Using \eqref{autonomous.gradient.osc} with $\rho=R/64$, the preceding inequality and \eqref{autonomous.flux.small.excess}, we infer
\begin{align*}
\sup_{B_{R/64}}|Dw-\xi|&\leq\osc_{B_{R/64}}Dw+\mean{B_{R/64}}|Dw-\xi|\dx\leq c\mean{B_{R/32}}|Dw-\xi|\dx\\
&\leq c\,\frac{|\xi|}{\varphi'(|\xi|)}\mean{B_{R/8}}|A_0(Dw)-A_0(\xi)|\dx\leq c\varepsilon|\xi|.
\end{align*}
Choosing $\varepsilon$ sufficiently small, we conclude that
\[
\frac{|\xi|}{2}\leq|Dw|\leq\frac{3|\xi|}{2}\qquad\text{in }B_{R/64}.
\]
The upper bound in \eqref{keep.a} can now be used on this ball with the fixed factor $\varphi'(|\xi|)/|\xi|$. Hence, for $0<\rho\leq R/64$, another application of \eqref{autonomous.gradient.osc} gives
\begin{align*}
\mean{B_\rho}|A_0(Dw)-(A_0(Dw))_{B_\rho}|\dx&\leq\osc_{B_\rho}A_0(Dw)\\
&\leq c\,\frac{\varphi'(|\xi|)}{|\xi|}\osc_{B_\rho}Dw\\
&\leq c\Bigl(\frac{\rho}{R}\Bigr)^\beta\frac{\varphi'(|\xi|)}{|\xi|}\mean{B_{R/32}}|Dw-\xi|\dx\\
&\leq c\Bigl(\frac{\rho}{R}\Bigr)^\beta\mean{B_{R/8}}|A_0(Dw)-A_0(\xi)|\dx.
\end{align*}
Since $\beta_0\leq\beta$, this proves \eqref{decay.varphi.harmonic} also in the second case. Finally, for $R/64<\rho\leq R/8$, the same estimate is straightforward by elementary properties of the excess, after increasing the constant by a factor depending only on $n$ and $\beta_0$.
\end{proof}

\begin{remark}\label{approx.vf}
{\rm The previous interior estimate for $Dw$ requires the autonomous vector field $A_0$ to be $C^1$-regular away from the origin and in our applications we will take $A_0=A(x_m,\,\cdot\,)$, where $x_m$ is a minimum point of $a(\,\cdot\,)$ on $\overline{B_R}$. This is the reason why in \eqref{A.strong} we assume $A(x,\,\cdot\,)$ to be $C^1$ away from the origin for {\em every} $x$.}
\end{remark}

The following lemma allows us to turn an integral lower bound on an arbitrarily small ball into a pointwise non-degeneracy estimate on a larger one, with constants independent of the number of intermediate scales, by the use of the quantitative control of the oscillation of a solution to \eqref{autonomous.equation}; see \cite[Lemma 2.3]{KMI1} for the $p$-Laplacian case.

\begin{lemma}\label{integral.to.pointwise}
Let $w\in W^{1,\varphi}(B_{R/8})$ be a weak solution to \eqref{autonomous.equation}; then there exist constants $c,c_{\sup}\geq1$ and $\beta\in(0,1)$, depending only on $n,g_0,g_1,\nu$ and $L$, such that for every $\mathcal A\geq1$, setting
\begin{equation}\label{bound.delta}
\delta_0=\frac{1}{c\,\mathcal A^{2/(\beta g_\smallzero)}}\leq\frac1{64},
\end{equation}
the following holds: if $\delta\in(0,\delta_0]$, $m\geq1$ is an integer and $\lambda>0$ satisfy
\begin{equation}\label{integral.bounds}
\frac{\lambda}{\mathcal A}\leq\mean{B_{\delta^m R}}\varphi'(|Dw|)\dx,\qquad \mean{B_{R/8}}\varphi'(|Dw|)\dx\leq\mathcal A\lambda,
\end{equation}
then
\begin{equation}\label{pointwise.bounds}
\frac{\lambda}{2^{g_1}\mathcal A}\leq\inf_{B_{\delta R}}\varphi'(|Dw|)\leq\sup_{B_{R/16}}\varphi'(|Dw|)\leq c_{\sup}\mathcal A\lambda.
\end{equation}
\end{lemma}

\begin{proof}
Combining \eqref{autonomous.gradient.from.flux} and \eqref{lip.varphi.harmonic} with the two-sided bound $|A_0(z)|\approx\varphi'(|z|)$, we deduce from the second estimate in \eqref{integral.bounds} that
\[
\sup_{B_{R/16}}|Dw|\leq c(\varphi')^{-1}(\mathcal A\lambda),\qquad\sup_{B_{R/16}}\varphi'(|Dw|)\leq c_{\sup}\mathcal A\lambda.
\]
Next, \eqref{autonomous.gradient.osc} on $B_{R/8}$ together with $\delta\leq1/64$ yields
\[
\osc_{B_{\delta R}}Dw\leq c\,\delta^\beta(\varphi')^{-1}(\mathcal A\lambda)\leq c\,\delta^\beta\mathcal A^{\frac2{g_\smallzero}}(\varphi')^{-1}(\lambda/\mathcal A)\leq\frac12(\varphi')^{-1}(\lambda/\mathcal A),
\]
using the scaling behavior of $(\varphi')^{-1}$ in the middle step and the choice of $\delta_0$ at the end. Meanwhile, the first bound in \eqref{integral.bounds} allows us to pick a point $x_1\in B_{\delta^mR}\subset B_{\delta R}$ where $|Dw(x_1)|\geq(\varphi')^{-1}(\lambda/\mathcal A)$. For any $x\in B_{\delta R}$, we thus get
\[
|Dw(x)|\geq|Dw(x_1)|-\osc_{B_{\delta R}}Dw\geq\frac12(\varphi')^{-1}(\lambda/\mathcal A).
\]
Invoking once more the scaling properties of $\varphi'$ gives
\[
\varphi'(|Dw(x)|)\geq2^{-g_1}\frac{\lambda}{\mathcal A},
\]
and notice that none of the constants involved depend on $m$.
\end{proof}

We now introduce the notion of harmonicity associated with the general autonomous vector fields considered above.

\begin{definition}
Let $\varphi$ satisfy \eqref{prop-varphi} and let $A_0$ satisfy \eqref{Orlicz.structure}. We say that $w\in W^{1,\varphi}_{\rm loc}(E)$ is locally $A_0$-harmonic in an open set $E\subset\R^n$ if
\[
\int_E\langle A_0(Dw),D\psi\rangle\dx=0\qquad\text{for every $\psi\in C_c^\infty(E)$}.
\]
When $A_0=A_\varphi$ as defined in \eqref{V.varphi}, we simply say that $w$ is $\varphi$-harmonic; in this case the preceding condition is equivalent to the local minimality of the functional $u\mapsto\int_E\varphi(|Du|)\dx$; in particular, the choice $\varphi(t)=t^p/p$ gives the usual notion of $p$-harmonicity.
\end{definition}

The following quantitative harmonic approximation result extends \cite[Lemma 1]{BCMgen} to general autonomous vector fields satisfying \eqref{Orlicz.structure}; see also \cite[Lemma 4.13]{HO23}. 

\begin{lemma}\label{p-harm}
Let $B\equiv B_R(x_0)\subset\R^n$ be a ball, let $\varphi$ satisfy \eqref{prop-varphi} with the normalization $\varphi(1)=1$, and let $A_0$ satisfy \eqref{Orlicz.structure}. Suppose that $v\in W^{1,\varphi}(2B)$ satisfies
\begin{equation}\label{eqn:v-int00}
\mean{2B}\varphi(|Dv|)\dx\leq\tilde c_1
\end{equation}
and
\begin{equation}\label{eqn:v-int002}
\mean{B}[\varphi(|Dv|)]^{1+\eta_0}\dx\leq\tilde c_2,
\end{equation}
where $\tilde c_1,\tilde c_2\geq1$ and $\eta_0>0$ are fixed constants. Moreover, suppose that for some $\sigma>0$ there holds
\begin{equation}\label{eqn:almost-p-harm}
\biggl|\mean{B}\langle A_0(Dv),D\psi\rangle\dx\biggr|\leq\sigma\|D\psi\|_{L^\infty(B)}\qquad\text{for every $\psi\in C_c^\infty(B)$}.
\end{equation}
Then there exists a unique $A_0$-harmonic map $w\in v+W^{1,\varphi}_0(B)$, and it satisfies
\begin{equation}\label{eqn:v-int1}
\mean{B}\varphi(|Dw|)\dx\leq c\,\tilde c_1,\qquad\mean{B}[\varphi(|Dw|)]^{1+\eta_0}\dx\leq c\,\tilde c_2,
\end{equation}
where $c$ depends only on $n,g_0,g_1,\nu,L$ and $\eta_0$. In addition, there exist $\iota\equiv\iota(n,g_0,g_1,\nu,L,\eta_0)\in(0,1)$ and $c\equiv c(n,g_0,g_1,\nu,L,\eta_0,\tilde c_1,\tilde c_2)\geq1$ such that
\begin{equation}\label{general.harmonic.approximation}
\mean{B}|V_\varphi(Dv)-V_\varphi(Dw)|^2\dx\leq c\,\sigma^{\iota}.
\end{equation}
\end{lemma}

Notice that once the ball $B$, the boundary datum $v$, and the autonomous field $A_0$ are fixed, the Dirichlet problem uniquely determines $w$, independently of $\sigma$ and of how the assumptions of the lemma are verified. In particular, for the frozen operators $A_i$ we will introduce in \eqref{vf.frozen.ai}, the resulting comparison maps will coincide whether the bounds come from the $p$-phase or the $(p,q)$-phase; therefore, the unified notation $w_1,w_2$ used later across the two regimes is unambiguous, as the freezing points and the Dirichlet data are fixed (see the beginning of Section~\ref{lin.chapter}).

\begin{proof}
Existence and uniqueness of the $A_0$-harmonic map $w\in v+W^{1,\varphi}_0(B)$ follow directly from standard monotone operator theory. Moreover, we have
\[
\mean{B}\varphi(|Dw|)\dx\leq c\mean{B}\varphi(|Dv|)\dx\leq c\,\tilde c_1,
\]
as seen by simply testing the weak formulation with $w-v$ and using the coercivity of $A_0$ together with Young's inequality. Finally, the second estimate in \eqref{eqn:v-int1} follows from the global Calder\'on-Zygmund bound in \cite[Lemma 4.5]{HO22}, applied to the Orlicz function $t\mapsto t^{1+\eta_0}$.

Proving \eqref{general.harmonic.approximation} relies on the Lipschitz truncation method of \cite[Lemma 1]{BCMgen} in the general framework of \cite[Lemma 4.13]{HO23}. Under \eqref{Orlicz.structure}, we have
\[
\langle A_0(z_1)-A_0(z_2),z_1-z_2\rangle\approx_{g_0,g_1,\nu,L}|V_\varphi(z_1)-V_\varphi(z_2)|^2,
\]
which plays the exact role of \eqref{monotonicity} for a general field $A_0$ (cf.~\cite[Lemma 20]{DE}). Everything else needed in the proof - namely, the growth of $\varphi$ and its conjugate, the Orlicz boundedness of the maximal operator, and \eqref{eqn:v-int1} - is insensitive to the specific form of $A_\varphi$. Running the same argument therefore yields \eqref{general.harmonic.approximation}, with $c$ and $\iota$ depending solely on the listed parameters.
\end{proof}

\subsection{Iteration and potential estimates}\label{dyadic.potentials}

We conclude the preliminaries by collecting a few standard tools concerning iteration arguments, dyadic summations and potential estimates. We begin with the following elementary iteration lemma; see \cite[Lemma 6.1]{Giusti}.
\begin{lemma}\label{Giusti.fundamental}
Let $\phi:[R_1,R_2]\to [0,+\infty)$ be a bounded function and assume that there exist $\mathcal A, \mathcal B\geq0$ and $\beta>0$ such that
\[
\phi(r)\leq \frac12\phi(s)+\frac{\mathcal A}{(s-r)^\beta}+\mathcal B\qquad  \text{for all $R_1\leq r<s\leq R_2$};
\]
then
\[
\phi(R_1)\leq c(\beta)\biggl[\frac{\mathcal A}{(R_2-R_1)^\beta}+\mathcal B\biggr].
\]
\end{lemma}

We shall also use the following elementary lemma, encoding a standard argument that allows us to pass from a decay estimate at a fixed ratio of radii to the corresponding estimate at arbitrary intermediate radii.

\begin{lemma}\label{geometric.decay.iteration}
Let $\Phi:(0,R_*]\to[0,\infty)$ be non-decreasing, and let $0<\sigma\leq\beta$, $\tau,\theta\in(0,1)$ and $A\geq0$. Suppose moreover that
\begin{equation}\label{decay.dyadic.classic}
\Phi(\tau r)\leq\theta\tau^\sigma\Phi(r)+Ar^\beta 
\end{equation}
for every $0<r\leq R_*$; then
\[
\Phi(\rho)\leq c\Bigl(\frac{\rho}{R}\Bigr)^\sigma\bigl[\Phi(R)+AR^\beta\bigr],\qquad \text{for every $0<\rho\leq R\leq R_*$,}
\]
where $c$ depends only on $\tau,\theta$ and $\sigma$. In the homogeneous case, let $\sigma\geq0$ and $0<r_*\leq R_*$. If $\Phi(\tau r)\leq\tau^\sigma\Phi(r)$ whenever $r_*\leq\tau r$ and $r\leq R_*$, then
\[
\Phi(\rho)\leq\tau^{-\sigma}\Bigl(\frac{\rho}{R}\Bigr)^\sigma\Phi(R)\qquad\text{for every }r_*\leq\rho\leq R\leq R_*.
\]
\end{lemma}
\begin{proof}
Fix $0<R\leq R_*$. Dividing the hypothesis  by $(\tau r)^\sigma$ and iterating at the radii $R,\tau R,\ldots,\tau^{j-1}R$, we obtain, for every $j\geq1$,
\[
\frac{\Phi(\tau^jR)}{(\tau^jR)^\sigma}\leq\theta^j\frac{\Phi(R)}{R^\sigma}+\frac{A}{\tau^\sigma}\sum_{k=0}^{j-1}\theta^{j-1-k}(\tau^kR)^{\beta-\sigma}\leq\theta^j\frac{\Phi(R)}{R^\sigma}+\frac{A}{\tau^\sigma(1-\theta)}R^{\beta-\sigma},
\]
since $\beta\geq\sigma$. The last bound also holds for $j=0$. Given $0<\rho\leq R$, choose $j\geq0$ such that $\tau^{j+1}R<\rho\leq\tau^jR$. By monotonicity, $\Phi(\rho)\leq\Phi(\tau^jR)$, and the conclusion follows from $\tau^jR<\rho/\tau$. In the homogeneous case, all the radii $\tau^kR$ with $0\leq k\leq j$ are at least $\rho\geq r_*$, so iteration gives $\Phi(\rho)\leq\Phi(\tau^jR)\leq\tau^{j\sigma}\Phi(R)\leq\tau^{-\sigma}(\rho/R)^\sigma\Phi(R)$.
\end{proof}

\medskip

Let now $x_0\in\Omega$, $B_R(x_0)\subset\Omega$, and fix $\delta\in(0,1/2]$. We consider the sequence of concentric balls with radii $R_j=\delta^jR$, $j\in\N_0$. For every $\sigma>0$, we have
\[
\sum_{j=0}^\infty R_j^\sigma=\frac{R^\sigma}{1-\delta^\sigma}.
\]
Using the monotonicity of $\varrho\mapsto|\mu|(B_\varrho(x_0))$, and splitting the integral over the intervals $[R_{j+1},R_j]$, we obtain (see, e.g., \cite{KMI1})
\begin{equation}\label{dyadic.potential.comparison}
\frac{1-\delta^{n-1}}{n-1}\sum_{j=1}^\infty\frac{|\mu|(B_{R_j}(x_0))}{R_j^{n-1}}\leq \bI_1^{|\mu|}(x_0,R)\leq\frac{\delta^{1-n}-1}{n-1}\sum_{j=0}^\infty\frac{|\mu|(B_{R_j}(x_0))}{R_j^{n-1}}.
\end{equation}
Notice that the first term of the sum on the right-hand side cannot, in general, be controlled by the potential truncated at $R$. However, starting from $R/2$ instead of $R$, we obtain
\begin{equation}\label{dyadic.potential.bound}
\sum_{j=0}^\infty\frac{|\mu|(B_{\delta^jR/2}(x_0))}{(\delta^jR/2)^{n-1}}\leq \frac{n-1}{1-2^{1-n}}\,\bI_1^{|\mu|}(x_0,R),
\end{equation}
and the constant is independent of $\delta\in(0,1/2]$. To prove this estimate, we may assume that $\bI_1^{|\mu|}(x_0,R)<\infty$, since otherwise the conclusion is immediate. By \eqref{dyadic.potential.comparison} and since $\delta\leq1/2$, we have
\begin{align*}
\sum_{j=0}^\infty\frac{|\mu|(B_{\delta^jR/2}(x_0))}{(\delta^jR/2)^{n-1}}&\leq\frac{n-1}{1-2^{1-n}}\left[\bI_1^{|\mu|}(x_0,R)-\bI_1^{|\mu|}(x_0,R/2)\right]+\frac{n-1}{1-\delta^{n-1}}\bI_1^{|\mu|}(x_0,R/2)\\
&\leq\frac{n-1}{1-2^{1-n}}\bI_1^{|\mu|}(x_0,R).
\end{align*}
We next recall the connection between Dini continuity and summability over geometric scales. Let $\varpi:[0,1]\to[0,\infty)$ be non-decreasing, with $\varpi(0)=0$. For every $R\in(0,1]$ and $\delta\in(0,1)$, we have
\begin{equation}\label{dyadic.Dini}
\frac1{\log(1/\delta)}\int_0^R\varpi(\varrho)\,\frac{d\varrho}{\varrho}\leq\sum_{j=0}^\infty\varpi(\delta^jR)\leq\varpi(R)+\frac1{\log(1/\delta)}\int_0^R\varpi(\varrho)\,\frac{d\varrho}{\varrho}.
\end{equation}
In particular, the Dini condition \eqref{Dini} is equivalent to the convergence of $\sum_{j=0}^\infty\omega(\delta^jR)$ for any fixed $R\in(0,1]$ and $\delta\in(0,1)$. Moreover, if $\int_0^1\varpi(\varrho)\,\frac{d\varrho}{\varrho}<\infty$, then, for every fixed $\delta\in(0,1)$,
\begin{equation}\label{dyadic.Dini.smallness}
\lim_{R\searrow0}\sum_{j=0}^\infty\varpi(\delta^jR)=0.
\end{equation}

\section{The balance condition and the two phases}

When $p\leq n$, the balance condition \eqref{subcritical.bound} reads $(n-1)(q-p)<\alpha(p-1)$, which in particular ensures that
\begin{equation}\label{rsmaller2}
r^{\alpha-\frac{(n-1)(q-p)}{p-1}}\leq1\qquad\text{for every }0<r\leq1.
\end{equation}
Moreover, since $\alpha\leq1$, this condition implies
\begin{equation}\label{bound.trivial}
\frac{q-1}{p-1}<1+\frac{\alpha}{n-1}\leq\frac{n}{n-1}\leq2.
\end{equation}
Thus the exponent $\es=(q-1)/(p-1)$ is admissible as an exponent for the $p$-part in \eqref{rough.compa.solopartep} - a choice we shall rely on to estimate the $q$-part in the $p$-phase. Furthermore, since $n'\leq p'$ for $p\leq n$, \eqref{bound.trivial} yields
\begin{equation}\label{bound.lasco}
q-1<n'(p-1)\leq p\leq n,\qquad\text{and hence}\qquad q-p<1.
\end{equation}
Conversely, when $p>n$, condition \eqref{bound.supercritical} guarantees that $q-p\leq\alpha p/n$, meaning the exponent of $r$ in \eqref{rsmaller2} satisfies
\[
\alpha-\frac{(n-1)(q-p)}{p-1}\geq\frac{\alpha(p-n)}{n(p-1)}>0.
\]
Thus \eqref{rsmaller2} continues to hold in this regime, including the endpoint $q/p=1+\alpha/n$, where the exponent reduces to $\alpha(p-n)/[n(p-1)]$.

\subsection{The two phases}

We call the unweighted term the $p$-part (or $p$-term) and the term weighted by $a(x)$ - even when potentially frozen at some point $\bar x$ -  the $q$-part/term. Their respective degrees of homogeneity in $z$ are $p$ and $q$ in the energy density \eqref{H.def}, and $p-1$ and $q-1$ in the flux magnitude \eqref{h.def} and the model vector field \eqref{model.A}; we use the same terminology for the comparison estimates below.

The terms ``$p$-phase'' and ``$(p,q)$-phase'' instead identify two distinct regimes on any given ball, which can be distinguished in several ways: via the supremum of $a(\,\cdot\,)$ during preliminary estimates (cf.~the forthcoming \eqref{pq.phase.sup.condition} and \eqref{p.phase.sup.condition}), via its infimum when freezing coefficients in the linearized comparisons (cf.~\eqref{two.different.cases}), or simply by evaluating $a(\,\cdot\,)$ at the center of the ball. Thanks to Hölder continuity, all these criteria remain mutually consistent up to a shift in the phase constant, even though a fixed threshold may occasionally assign differing phases to the same ball (see Remark~\ref{stability.p.phase}). For every ball $B_R\Subset\Omega$, H\"older continuity gives
\begin{equation}\label{oscillation.a}
\osc_{B_R}a \leq 2[a]_{0,\alpha}R^\alpha.
\end{equation}
In the supremum formulation, the two complementary conditions are going to be
\begin{equation}\label{p.phase.sup.condition}
\sup_{B_R}a(\,\cdot\,)\leq M[a]_{0,\alpha}R^\alpha\qquad\text{($p$-phase)}
\end{equation}
and
\begin{equation}\label{pq.phase.sup.condition}
\sup_{B_R}a(\,\cdot\,)>M[a]_{0,\alpha}R^\alpha\qquad\text{($(p,q)$-phase)}.
\end{equation}
For the latter formulation we will always take $M>2$, so that the coefficient is comparable in the whole closure of the ball to its value at the center of the ball, see \eqref{above.below}.

\begin{remark}\label{essent.constant}
{\rm
If \eqref{pq.phase.sup.condition} holds with $M>2$ then, employing \eqref{oscillation.a},
\[
\inf_{B_R}a(\,\cdot\,)\geq\sup_{B_R}a(\,\cdot\,)-2[a]_{0,\alpha}R^\alpha\geq\Bigl(1-\frac2M\Bigr)\sup_{B_R}a(\,\cdot\,)>0.
\]
Consequently, for every $x,y\in\overline{B_R}$,
\begin{equation}\label{above.below}
\frac{M-2}{M}a(y)\leq a(x)\leq\frac{M}{M-2}a(y).
\end{equation}
The comparison constants depend only on $M$ and are uniform for $M\geq M_0>2$. If instead
\[
\inf_{B_R}a(\,\cdot\,)>M[a]_{0,\alpha}R^\alpha
\]
for some $M\geq 2$, then directly
\begin{equation}\label{inf.phase.comparability}
0<\inf_{B_R}a(\,\cdot\,)\leq\sup_{B_R}a(\,\cdot\,)\leq\Bigl(1+\frac2M\Bigr)\inf_{B_R}a(\,\cdot\,)\leq 2\inf_{B_R}a(\,\cdot\,).
\end{equation}
}
\end{remark}
\begin{remark}\label{stability.phases}
{\rm
Let $B_\rho(y)\subset B_R$, with $R\leq 1$. If \eqref{pq.phase.sup.condition} holds with $M>2$, then analogously
\[
\inf_{B_\rho(y)}a(\,\cdot\,)\geq\inf_{B_R}a(\,\cdot\,)>(M-2)[a]_{0,\alpha}R^\alpha\geq(M-2)[a]_{0,\alpha}\rho^\alpha;
\]
consequently, passing from the supremum condition on $B_R$ to a sub-ball $B_\rho(y)$ introduces a controlled shift in the phase parameter. By contrast, the coefficient comparability established on $B_R$ persists on every sub-ball with identical comparison constants. Conversely, the infimum formulation avoids this degradation entirely: if $\inf_{B_R}a(\,\cdot\,)>M[a]_{0,\alpha}R^\alpha$, then
\[
\inf_{B_\rho(y)}a(\,\cdot\,)\geq\inf_{B_R}a(\,\cdot\,)>M[a]_{0,\alpha}R^\alpha\geq M[a]_{0,\alpha}\rho^\alpha.
\]
Passing to a smaller ball thus preserves the condition with the exact same constant $M$. The same reasoning applies directly to the modified thresholds $R^{\alpha-\x}$ used for any $0\leq\x<\alpha$.

}
\end{remark}

\begin{remark}\label{stability.p.phase}
{\rm
Let $B_{\rho_1}\subset B_{\rho_2}\Subset\Omega$ with $0<\rho_1\leq\rho_2\leq1$, let $M>0$, and let $0\leq\x<\alpha$, so that the phase condition is unmodified when $\x=0$. If
\[
\inf_{B_{\rho_1}}a(\,\cdot\,)\leq M[a]_{0,\alpha}\rho_1^{\alpha-\x},
\]
then the condition in \eqref{oscillation.a} implies
\[
\sup_{B_{\rho_2}}a(\,\cdot\,)\leq\inf_{B_{\rho_1}}a(\,\cdot\,)+2[a]_{0,\alpha}\rho_2^\alpha\leq(M+2)[a]_{0,\alpha}\rho_2^{\alpha-\x}.
\]
Thus, imposing a (modified) $p$-phase condition on the smaller ball implies an upper bound on the larger one with the phase constant independent of the ratio of radii. Regarding the infimum condition, the monotonicity of the infimum directly preserves the constant $M$ via domain enlargement. Conversely, if
\[
\inf_{B_{\rho_2}}a(\,\cdot\,)\leq M[a]_{0,\alpha}\rho_2^{\alpha-\x},
\]
it follows that
\[
\sup_{B_{\rho_1}}a(\,\cdot\,)\leq\inf_{B_{\rho_2}}a(\,\cdot\,)+2[a]_{0,\alpha}\rho_2^\alpha\leq(M+2)\Bigl(\frac{\rho_2}{\rho_1}\Bigr)^{\alpha-\x}[a]_{0,\alpha}\rho_1^{\alpha-\x};
\]
hence, restricting the domain in this direction introduces a dependence on the ratio of the radii.

Both conclusions hold, mutatis mutandis, if the assumptions are instead formulated in terms of the supremum. We will routinely bound $M+2\leq2M$ for $M\geq2$, and $M+2\leq3M$ when $M\geq1$. In the special case $\x=0$, the upper bound restriction $\rho_2\leq1$ is superfluous.
}
\end{remark}

\section{Some basic measure data theory elements}\label{basic.measure.data.comparisons}

Essentially all the results of this section and the following ones are local. We therefore consider balls $B_R\equiv B_R(x_0)\subset\Omega$ with $R\leq1$; when $p>n$, we additionally require $B_R\subset\Omega_1$, where $\Omega_1$ is the intermediate domain introduced after \eqref{Omega0}. In particular, these conditions are automatically satisfied for balls centered in $\Omega_0$ with $R\leq\min\{1,d_{\Omega_0}/8\}$. 

\subsection{The existence of SOLAs}
We first establish a global estimate for energy solutions with bounded data in Lemma~\ref{global.RH} and use it to prove Theorem~\ref{existence.SOLA} by approximation. Although the existence of SOLAs has recently been established in \cite{SongYoun} under differentiability assumptions on the vector field with respect to the gradient variable, the corresponding result under our weaker monotonicity assumptions, and with nonhomogeneous boundary data, does not seem to be explicitly available. We provide the details.

\begin{lemma}[Global integrability in the case $p\leq n$]\label{global.RH}
Let $u\in g+W^{1,H(\cdot)}_0(\Omega)$ be an energy solution to \eqref{CD.sola} under the assumptions of Theorem~\ref{existence.SOLA}, with the additional assumption $\mu\in L^\infty(\Omega)$, and let $\mathfrak D$ be as in \eqref{mathcalD}. For every exponent
\begin{equation}\label{ess.energy.SOLA}
\min\Bigl\{1,\frac1{p-1}\Bigr\}\leq s<\frac n{n-1},
\end{equation}
there holds
\begin{equation}\label{global.est.SOLA}
\biggl(\int_\Omega|Du|^{(p-1)s}\dx\biggr)^{\frac1s}\leq c\,\mathfrak D,
\end{equation}
with the constant $c$ depending on $n,p,q,\nu,L,|\Omega|$ and $s$. Moreover,
\begin{equation}\label{global.energy.persola}
\int_\Omega h(x,Du)\dx\leq c\biggl[\mathfrak D+\|a\|_{L^\infty(\Omega)}^{\frac1q}\mathfrak D^{\frac{p'}{q'}}\biggr],
\end{equation}
where $c$ depends only on $n,p,q,\nu,L$ and $|\Omega|$.
\end{lemma}

\begin{proof}
We test the weak formulation with the test functions
\[
\varphi_\pm:=\pm\frac{\kappa^{1-\xi}-(\kappa+(u-g)_\pm)^{1-\xi}}{\xi-1}\in L^\infty(\Omega)\cap W^{1,H(\cdot)}_0(\Omega),
\]
which are Lipschitz compositions of $u-g$ with scalar maps vanishing at zero, where $\xi>1$ and $\kappa>0$ are to be chosen below. The two test functions satisfy
\[
{\|\varphi_\pm\|}_{L^\infty(\Omega)}\leq\frac{\kappa^{1-\xi}}{\xi-1} \qquad\text{and}\qquad D\varphi_\pm = \frac{Du-Dg}{(\kappa+{(u-g)}_\pm)^\xi} \chi_{\{{(u-g)}_\pm>0\}},
\]
so we deduce
\[
\int_\Omega \frac{\langle A(x,Du),Du-Dg\rangle}{(\kappa+{(u-g)}_\pm)^\xi} \chi_{\{{(u-g)}_\pm>0\}}\dx = \int_\Omega \varphi_\pm\,d\mu \leq \frac{\kappa^{1-\xi}}{\xi-1}|\mu|(\Omega).
\]
Inserting into the previous display the estimate from below that follows from coercivity, growth and the Young inequality in \eqref{Young.H}
\[
\langle A(x,Du),Du-Dg\rangle\geq \nu H(x,Du)-L\,h(x,Du)|Dg| \geq \frac{\nu}{2}H(x,Du)-c\,H(x,Dg),
\]
we infer
\begin{align}\label{weighted.global.energy}
\int_\Omega \frac{H(x,Du)}{(\kappa+{(u-g)}_\pm)^\xi}\chi_{\{{(u-g)}_\pm>0\}}\dx&\leq c\frac{\kappa^{1-\xi}}{\xi-1}|\mu|(\Omega) +c\int_\Omega \frac{H(x,Dg)}{(\kappa+{(u-g)}_\pm)^\xi} \chi_{\{{(u-g)}_\pm>0\}}\dx\notag\\
&\leq c\frac{\kappa^{1-\xi}}{\xi-1}|\mu|(\Omega) +\frac{c}{\kappa^\xi}\int_\Omega H(x,Dg)\dx.
\end{align}
Unless otherwise stated, constants in the following argument depend on $n,p,q,\nu,L,s,|\Omega|$ and, until it is fixed, on $\xi$. Summing the two inequalities, for the two signs, in \eqref{weighted.global.energy} and recalling that $Du=Dg$ a.e. on $\{u=g\}$ since $u-g\in W^{1,1}(\Omega)$, we have
\[
\int_\Omega \frac{H(x,Du)}{(\kappa+|u-g|)^\xi}\,\dx\leq c\frac{\kappa^{1-\xi}}{\xi-1}|\mu|(\Omega)+\frac{c}{\kappa^\xi}\int_\Omega H(x,Dg)\,\dx.
\]
We now apply H\"older's inequality with conjugate exponents $(p'/s,p'/(p'-s))$, since for $s$ as in \eqref{ess.energy.SOLA} we have $s<n'\leq p'$, and we find
\[
\int_\Omega|Du|^{(p-1)s}\dx\leq c\int_\Omega{\bigl[H(x,Du)\bigr]}^{\frac{s}{p'}}\dx \leq c\biggl(\int_\Omega \frac{H(x,Du)}{(\kappa+|u-g|)^\xi}\dx\biggr)^\frac{s}{p'}\biggl(\int_\Omega (\kappa+|u-g|)^\frac{s\xi}{p'-s}\dx\biggr)^{1-\frac s{p'}}.
\]
Finally we choose 
\[
\kappa^{\frac{s \xi}{p'-s }}=\mean\Omega |u-g|^\frac{s \xi}{p'-s }\dx\quad\Longrightarrow\quad {\biggl(\int_\Omega (\kappa+|u-g|)^\frac{s \xi}{p'-s }\dx\biggr)}^{1-\frac s {p'}}\leq c\,{|\Omega|}^{1-\frac s{p'}}\kappa^{\xi\frac{s }{p'}};
\]
if $u=g$ almost everywhere, \eqref{global.est.SOLA} follows directly from H\"older's inequality; hence we may assume $\kappa>0$; moreover, the choice of $\xi$ below will ensure that $\kappa$ is finite. Using the preceding weighted estimate and Young's inequality with conjugate exponents $p$ and $p'$, for $\varepsilon\in(0,1)$ to be chosen, we obtain
\begin{align*}
\int_\Omega|Du|^{(p-1)s}\dx &\leq c\frac{\kappa^\frac{s}{p'}}{(\xi-1)^\frac{s}{p'}}\bigl[|\mu|(\Omega)\bigr]^\frac{s}{p'}+c{\biggl(\int_\Omega H(x,Dg)\dx\biggr)}^\frac{s}{p'}\\
&\leq \varepsilon \kappa^{(p-1)s}+ c_\varepsilon(\xi)\bigl[|\mu|(\Omega)\bigr]^s+c{\biggl(\int_\Omega H(x,Dg)\dx\biggr)}^\frac{s}{p'}.
\end{align*}
Now, if we are in the case $(p-1)s\geq1$ (this is the case, for instance, when $p\geq2$), we make the choice (note that $(p-1)s<n$ since $p\leq n$ and $s<p'$)
\begin{equation}\label{xi.choice.1}
\xi=\frac{n[p-(p-1)s ]}{n-(p-1)s} \qquad\Longleftrightarrow\qquad \frac{s \xi}{p'-s }=[(p-1)s]^*=\frac{n(p-1)s }{n-(p-1)s};
\end{equation}
this choice gives $\xi>1$ because $s<\frac n{n-1}$. Since $1\leq(p-1)s<p\leq n$ and $u-g\in W^{1,(p-1)s}_0(\Omega)$, Sobolev's inequality ensures that $\kappa$ is finite and yields
\begin{align*}
 \kappa^{(p-1)s}={\biggl(\mean\Omega |u-g|^{[(p-1)s]^*}\dx\biggr)}^{\frac{(p-1)s}{[(p-1)s]^*}}&\leq c\mean\Omega {|D(u-g)|}^{(p-1)s}\dx\\
 &\leq c\int_\Omega {|Du|}^{(p-1)s}\dx+c{\biggl(\int_\Omega H(x,Dg)\dx\biggr)}^\frac{s}{p'}.
\end{align*}
Hence,
\[
\int_\Omega|Du|^{(p-1)s}\dx\leq  c\,\varepsilon\int_\Omega {|Du|}^{(p-1)s}\dx+ c_\varepsilon(\xi)\bigl[|\mu|(\Omega)\bigr]^s+c{\biggl(\int_\Omega H(x,Dg)\dx\biggr)}^\frac{s}{p'};
\]
and if $\varepsilon$ is sufficiently small, depending only on $n,p,q,\nu,L,s$ and $|\Omega|$, we can reabsorb the $L^{(p-1)s}$-norm of $Du$, obtaining \eqref{global.est.SOLA} if $s\geq1/(p-1)$, with $c$ depending only on $n,p,q,\nu,L,|\Omega|$ and $s$. If $(p-1)s<1$, the assumption $p>2-\frac1n$ ensures that $\frac1{p-1}<\frac n{n-1}$ and we may therefore apply the estimate just proved with $s=\frac1{p-1}$ and use H\"older's inequality to obtain
\[
{\biggl(\int_\Omega |Du|^{(p-1)s}\dx\biggr)}^{\frac1s}\leq|\Omega|^{\frac{1-(p-1)s}{s}}{\biggl(\int_\Omega |Du|\dx\biggr)}^{p-1}\leq c\,\mathfrak D.
\]
To establish \eqref{global.energy.persola}, we note first that \eqref{global.est.SOLA} with $s=1$ controls the integral of $|Du|^{p-1}$. For the $q$-term, we take 
\[
m\in 
\begin{cases}
\displaystyle\Bigl(\max\{n',q-1\} ,\frac{n(p-1)}{n-p}\Bigr)\qquad &\text{if $p< n$}\\[3mm]
\bigl(\max\{n',q-1\} ,+\infty\bigr)\qquad &\text{if $p\geq n$}
\end{cases}
\]
and fix $1<\xi<m/(q-1)$; both exponents depend only on $n,p$ and $q$ and exist under our hypotheses. Sobolev's embedding, \eqref{global.est.SOLA} for $s=nm/[(n+m)(p-1)]$ and H\"older's inequality (to handle $Dg$) yield
\[
\biggl(\int_\Omega|u-g|^m\dx\biggr)^{\frac1m}\leq c\biggl(\int_\Omega|D(u-g)|^{\frac{nm}{n+m}}\dx\biggr)^{\frac{n+m}{nm}}\leq c\,\mathfrak D^{\frac1{p-1}},
\]
with $c=c(n,p,q,\nu,L,|\Omega|)$. Assuming without loss of generality that $\mathfrak D>0$, we set $\kappa:=\mathfrak D^{\frac1{p-1}}$ in the weighted estimate just proved. Since 
\[
|\mu|(\Omega)\leq\mathfrak D=\kappa^{p-1}\qquad\text{and}\qquad \int_\Omega H(x,Dg)\dx\leq\mathfrak D^{p'}=\kappa^p,
\]
we obtain
\[
\int_\Omega\frac{H(x,Du)}{(\kappa+|u-g|)^\xi}\dx\leq c\kappa^{1-\xi}|\mu|(\Omega)+\frac{c}{\kappa^\xi}\int_\Omega H(x,Dg)\dx\leq c\kappa^{p-\xi}.
\]
We can now estimate the $q$-part and complete the proof of \eqref{global.energy.persola}. Using the previous estimate alongside H\"older's inequality with exponents $q'$ and $q$, we find
\begin{align*}
\int_\Omega a(x)|Du|^{q-1}\dx&\leq\biggl(\int_\Omega\frac{a(x)|Du|^q}{(\kappa+|u-g|)^\xi}\dx\biggr)^{\frac1{q'}}\biggl(\int_\Omega a(x)(\kappa+|u-g|)^{\xi(q-1)}\dx\biggr)^{\frac1q}\\
&\leq c\|a\|_{L^\infty(\Omega)}^{\frac1q}\kappa^{\frac{p-\xi}{q'}}\kappa^{\frac{\xi}{q'}}=c\|a\|_{L^\infty(\Omega)}^{\frac1q}\mathfrak D^{\frac{p'}{q'}},
\end{align*}
where the last line relies on the $L^m(\Omega)$ bound for $u-g$ and the condition $\xi(q-1)<m$. The constant depends again on $n,p,q,\nu,L$ and $|\Omega|$.
\end{proof}

\begin{proof}[Proof of Theorem~\ref{existence.SOLA}]
The argument follows the lines of \cite{SongYoun}, but we outline it here by constructing a local SOLA through the approximation of the measure datum and verifying the requirements of Definition~\ref{def.SOLA}. 

We extend $\mu$ by zero outside $\Omega$ and choose a sequence of non-negative standard mollifiers $\varrho_j\in C_c^\infty(B_{\varepsilon_j})$ with $\varepsilon_j\searrow0$ and $\int_{\R^n}\varrho_j\dx=1$, setting $\mu_j:=\mu*\varrho_j$. Standard properties of convolution imply that $\mu_j\in C^\infty(\Omega)\cap L^\infty(\Omega)$, $\mu_j\dx\stackrel{*}{\rightharpoonup}\mu$ in $\mathcal M_{\rm loc}(\Omega)$, and
\[
\int_\Omega|\mu_j|\dx\leq|\mu|(\Omega).
\]
Moreover, for every ball $B_\rho(x_0)\Subset\Omega$ and all sufficiently large $j$,
\begin{equation}\label{SOLA.approx.variation}
\int_{B_\rho(x_0)}|\mu_j|\dx\leq|\mu|(B_{\rho+\varepsilon_j}(x_0)).
\end{equation}
Passing to the $\limsup$ and using the continuity from above of the total variation measure then guarantees the bound on the total variations in \eqref{SOLA.uniform.meas}. 

Next, we consider the energy solution to the approximating Dirichlet problem
\[
-\divergence A(x,Du_j)=\mu_j, \qquad u_j-g\in W^{1,H(\cdot)}_0(\Omega).
\]
Since $g\in W^{1,H(\cdot)}(\Omega)$ and $W^{1,H(\cdot)}_0(\Omega)$ is reflexive, and the structure conditions in \eqref{A.veryweak} ensure that the operator
\[
w\mapsto-\divergence A(x,Dg+Dw)
\]
is bounded, hemicontinuous, coercive, and strictly monotone from $W^{1,H(\cdot)}_0(\Omega)$ into its dual, the Browder-Minty theorem guarantees a unique energy solution $u_j\in g+W^{1,H(\cdot)}_0(\Omega)$ for each $\mu_j\in L^\infty(\Omega)\subset (W^{1,H(\cdot)}_0(\Omega))^*$. 

By Lemma~\ref{global.RH} and the bound $\|\mu_j\|_{L^1(\Omega)}\leq|\mu|(\Omega)$, the gradients of the approximating sequence are uniformly bounded in $L^t(\Omega)$ for every
\begin{equation}\label{usual.bound.weak}
1\leq t<\frac{n(p-1)}{n-1}.
\end{equation}
Their flux magnitudes are uniformly bounded in $L^1(\Omega)$, since
\[
\sup_j\int_\Omega h(x,Du_j)\dx\leq c\biggl[\mathfrak D+\|a\|_{L^\infty(\Omega)}^{\frac1q}\mathfrak D^{\frac{p'}{q'}}\biggr],
\]
with $c$ depending only on $n,p,q,\nu,L$ and $|\Omega|$; in particular, $\{u_j\}$ satisfies the uniform local flux bound in \eqref{SOLA.convergence}. Thanks to \eqref{2-1/n}, we can fix an exponent $t_0>1$ in the range \eqref{usual.bound.weak}. Poincaré's inequality then implies that $\{u_j-g\}$ is bounded in $W^{1,t_0}_0(\Omega)$. By weak compactness and the compact embedding into $L^{t_0}(\Omega)$, there exist a function $u\in g+W^{1,t_0}_0(\Omega)$ and a subsequence (not relabeled) such that
\[
u_j\rightharpoonup u\quad\text{weakly in }W^{1,t_0}(\Omega),\qquad u_j\longrightarrow u\quad\text{strongly in }L^{t_0}(\Omega)\text{ and almost everywhere in }\Omega.
\]
To prove almost everywhere convergence of the gradients, we follow the classical truncation and compactness argument of Boccardo and Gallou\"et. For $\varsigma>0$, we define $T_\varsigma(t):=\max\{\min\{t,\varsigma\},-\varsigma\}$ and we test the difference of the two approximating equations with the bounded function $T_\varsigma(u_j-u_k)$ - $u_j-u_k\in W^{1,H(\cdot)}_0(\Omega)$ so the test function also has zero trace; using $\|\mu_j\|_{L^1(\Omega)}+\|\mu_k\|_{L^1(\Omega)}\leq2|\mu|(\Omega)$, we obtain
\begin{equation}\label{SOLA.monotonicity.defect}
\int_{\{|u_j-u_k|<\varsigma\}} \langle A(x,Du_j)-A(x,Du_k),Du_j-Du_k\rangle \dx \leq 2\varsigma |\mu|(\Omega).
\end{equation}
We fix $\rho>0$, take $\lambda>\rho$ and define
\[
E^\rho_{j,k} := \bigl\{|Du_j-Du_k|\geq\rho\bigr\}, \qquad U^\varsigma_{j,k} := \bigl\{|u_j-u_k|<\varsigma\bigr\}, \qquad V^\lambda_{j,k} :=\bigl\{\max\{|Du_j|,|Du_k|\}\leq\lambda\bigr\},
\]
along with
\[
{\textnormal{\texttt{E}}}_{\rho,\lambda}(x) := \inf\bigl\{ \langle A(x,z_1)-A(x,z_2),z_1-z_2\rangle : |z_1|,|z_2|\leq\lambda,\ |z_1-z_2|\geq\rho \bigr\},\qquad x\in\Omega.
\]
Continuity and strict monotonicity of $A(x,\,\cdot\,)$ imply ${\textnormal{\texttt{E}}}_{\rho,\lambda}(x)>0$ for a.e. $x\in\Omega$, while \eqref{SOLA.monotonicity.defect} yields
\[
\int_{E^\rho_{j,k}\cap U^\varsigma_{j,k}\cap V^\lambda_{j,k}}\textnormal{\texttt{E}}_{\rho,\lambda}(x)\dx\leq2\varsigma|\mu|(\Omega).
\]
Since $\{u_j\}$ converges in $L^{t_0}(\Omega)$ (hence is Cauchy in measure) and $\{Du_j\}$ is bounded in $L^{t_0}(\Omega)$, we verify that $\{Du_j\}$ is Cauchy in measure via the standard decomposition
\[
|E^\rho_{j,k}| \leq |E^\rho_{j,k}\cap U^\varsigma_{j,k}\cap V^\lambda_{j,k}| + |\Omega\setminus U^\varsigma_{j,k}| + |\Omega\setminus V^\lambda_{j,k}|;
\]
given $\varepsilon>0$, we first pick $\lambda$ large enough so that $|\Omega\setminus V^\lambda_{j,k}|<\varepsilon$ for all $j,k$. Then, since ${\textnormal{\texttt{E}}}_{\rho,\lambda}>0$ a.e., we next choose $\delta>0$ with $|\{{\textnormal{\texttt{E}}}_{\rho,\lambda}\le\delta\}|<\varepsilon$, and then $\varsigma>0$ small enough so that
\[
|E^\rho_{j,k}\cap U^\varsigma_{j,k}\cap V^\lambda_{j,k}| \leq |\{{\textnormal{\texttt{E}}}_{\rho,\lambda}\le\delta\}| + \frac{c\varsigma|\mu|(\Omega)}{\delta} < 2\varepsilon
\]
for all $j,k$. Finally, we choose $N\in\N$ such that $|\Omega\setminus U^\varsigma_{j,k}| = |\{|u_j-u_k|\geq\varsigma\}|<\varepsilon$ whenever $j,k\geq N$. Because $\rho$ and $\varepsilon$ are arbitrary, $\{Du_j\}$ is Cauchy in measure on $\Omega$, hence converges in measure. As $\{Du_j\}$ is uniformly bounded in $L^{t_0}(\Omega)$ with $t_0>1$, Vitali's theorem yields strong convergence in $L^1(\Omega;\R^n)$ and the limit of the gradients must be $Du$, so that, up to a further subsequence, we have
\[
Du_j\longrightarrow Du\qquad\text{almost everywhere in }\Omega.
\]
Fixing now $t$ as dictated by \eqref{usual.bound.weak}, we choose $r\in(t,n'(p-1))$. Through a careful synthesis of the uniform $L^r(\Omega)$ bound on $\{Du_j\}$ with the pointwise almost everywhere convergence - invoking once more the classical machinery of Vitali's convergence theorem - we deduce that $Du_j$ converges strongly to $Du$ in $L^t(\Omega;\R^n)$. Observing next that $u_j-g\in W^{1,t}_0(\Omega)$, we turn to Poincar\'e's inequality to establish that $\{u_j-g\}$ forms a Cauchy sequence in $W^{1,t}_0(\Omega)$; its unique limit is thereby identified as $u-g$ in light of the aforementioned almost everywhere convergence of $u_j$. We thus arrive at the convergences
\[
u_j\longrightarrow u\quad\text{strongly in }W^{1,t}(\Omega),\qquad u-g\in W^{1,t}_0(\Omega),
\]
which hold collectively for every exponent $t$ satisfying \eqref{usual.bound.weak}. Proceeding further, we fix $m$ and $\xi$ precisely as in the proof of \eqref{global.energy.persola} to ensure the strict requirement $\xi(q-1)<m$, and we set $\kappa:=\mathfrak D^{\frac1{p-1}}$ under the standing assumption that $\mathfrak D>0$ (the degenerate case $\mathfrak D=0$ being trivial). Since $\|\mu_j\|_{L^1(\Omega)}\leq|\mu|(\Omega)$, both the weighted gradient estimate and the $L^m$ bound established earlier for $u_j-g$ remain uniformly valid with respect to the index $j$. Consider now an arbitrary measurable subset $E\subset\Omega$; applying H\"older's inequality leads to
\begin{align*}
\int_E a(x)|Du_j|^{q-1}\dx&\leq\biggl(\int_E\frac{a(x)|Du_j|^q}{(\kappa+|u_j-g|)^\xi}\dx\biggr)^{\frac1{q'}}\biggl(\int_E a(x)(\kappa+|u_j-g|)^{\xi(q-1)}\dx\biggr)^{\frac1q}\\
&\leq c\|a\|_{L^\infty(\Omega)}^{\frac1q}\mathfrak D^{\frac{p'}{q'}}|E|^{\frac1q(1-\frac{\xi(q-1)}{m})},
\end{align*}
where the constant $c$ depends exclusively upon $n,p,q,\nu,L$, and $|\Omega|$. Since the measure of the set appears with a strictly positive exponent, the sequence $\{a(\,\cdot\,)|Du_j|^{q-1}\}$ is equi-integrable on $\Omega$. The corresponding equi-integrability for $\{|Du_j|^{p-1}\}$ follows directly from \eqref{global.est.SOLA} by choosing any fixed index $s\in(1,\frac n{n-1})$. Hence, $\{h(\,\cdot\,,Du_j)\}$ is equi-integrable, and Vitali's theorem - combined with the pointwise convergence $Du_j\to Du$ - yields
\begin{equation}\label{SOLA.global.strong.flux}
h(\,\cdot\,,Du_j)\to h(\,\cdot\,,Du)\qquad\text{strongly in }L^1(\Omega).
\end{equation}
Standard growth and continuity properties of $A$ give the analogous strong convergence
\begin{equation}\label{SOLA.strong.A}
A(\,\cdot\,,Du_j)\to A(\,\cdot\,,Du)\qquad\text{strongly in }L^1(\Omega;\R^n).
\end{equation}
To complete the argument, we test the approximating equations against an arbitrary smooth test function $\varphi\in C_c^\infty(\Omega)$; the strong $L^1$ convergence of $A(\,\cdot\,,Du_j)$ coupled with the weak$^*$ convergence $\mu_j\dx\stackrel{*}{\rightharpoonup}\mu$ permits us to pass safely to the limit in the distributional sense, thereby recovering the identity
\[
\int_\Omega\langle A(x,Du),D\varphi\rangle\dx=\int_\Omega\varphi\,d\mu.
\]
This confirms that $u$ functions as a very weak solution to \eqref{equation}; when paired with the convergence properties established above, the uniform flux bounds, and relation \eqref{SOLA.approx.variation}, it emerges that $u$ successfully fulfills every requirement set forth in Definition~\ref{def.SOLA} and thus qualifies as a local SOLA. Furthermore, the strong convergence recorded in \eqref{SOLA.global.strong.flux} provides the necessary justification to pass the foundational energy inequality \eqref{global.energy.persola} directly to the limit.

\medskip

To prove the level estimates, test the equation for $u_j$ with $T_\varsigma(u_j-g)$ for $\varsigma>0$. Coercivity and Young's inequality give
\begin{align*}
\frac{\nu}{2}\int_{\{|u_j-g|<\varsigma\}}H(x,Du_j)\dx&\leq\int_\Omega\mu_jT_\varsigma(u_j-g)\dx+c\int_{\{|u_j-g|<\varsigma\}}H(x,Dg)\dx\\
&\leq\varsigma|\mu|(\Omega)+c\int_\Omega H(x,Dg)\dx.
\end{align*}
Here and below, $c$ depends only on $n,p,q,\nu,L$ and $|\Omega|$. Since $H(x,z)\geq|z|^p/p$, the chain rule and Poincar\'e's inequality show that $\{T_\varsigma(u_j-g)\}$ is bounded in $W^{1,p}_0(\Omega)$ for each fixed $\varsigma>0$. By almost everywhere convergence to $T_\varsigma(u-g)$ and weak lower semicontinuity, we get $T_\varsigma(u-g)\in W^{1,p}_0(\Omega)$ and
\[
\int_\Omega|DT_\varsigma(u-g)|^p\dx\leq c\,\varsigma|\mu|(\Omega)+c\int_\Omega H(x,Dg)\dx\leq c\,\bigl(\varsigma\mathfrak D+\mathfrak D^{p'}\bigr).
\]
If $\mathfrak D=0$, then $u=g$ almost everywhere and the claims are trivial; we may thus assume $\mathfrak D>0$. In particular,
\begin{equation}\label{trunc.D}
\int_\Omega|DT_\varsigma(u-g)|^p\dx\leq c\,\varsigma\mathfrak D\qquad\text{for every }\varsigma\geq\mathfrak D^{\frac1{p-1}}.
\end{equation}
Suppose first that $p<n$. For $\varsigma\geq\mathfrak D^{\frac1{p-1}}$, Sobolev's inequality and \eqref{trunc.D} yield
\[
\varsigma^p\bigl|\{|u-g|>\varsigma\}\bigr|^{1-\frac pn}\leq\biggl(\int_\Omega|T_\varsigma(u-g)|^{\frac{np}{n-p}}\dx\biggr)^{1-\frac pn}\leq c\,\varsigma\mathfrak D.
\]
This gives the claimed level estimate for $\varsigma\geq\mathfrak D^{\frac1{p-1}}$, while for $0<\varsigma<\mathfrak D^{\frac1{p-1}}$ we simply note that
\[
\varsigma^{\frac{n(p-1)}{n-p}}\bigl|\{|u-g|>\varsigma\}\bigr|\leq|\Omega|\mathfrak D^{\frac n{n-p}},
\]
covering the full range of $\varsigma$ when $p<n$. Suppose next that $p=n$. Combining the Trudinger-Moser inequality with \eqref{trunc.D} yields
\[
\int_\Omega\exp\biggl(\frac{|T_\varsigma(u-g)|^{\frac n{n-1}}}{c\,(\varsigma\mathfrak D)^{\frac1{n-1}}}\biggr)\dx\leq c|\Omega|\qquad\text{for every }\varsigma\geq\mathfrak D^{\frac1{n-1}}.
\]
Since $|T_\varsigma(u-g)|=\varsigma$ on $\{|u-g|>\varsigma\}$, we immediately get
\[
\bigl|\{|u-g|>\varsigma\}\bigr|\leq c|\Omega|\exp\biggl(-\frac{\varsigma}{c\,\mathfrak D^{\frac1{n-1}}}\biggr)\qquad\text{for every }\varsigma\geq\mathfrak D^{\frac1{n-1}},
\]
and the trivial bound $|\{|u-g|>\varsigma\}|\leq|\Omega|$ extends this estimate to $0<\varsigma<\mathfrak D^{\frac1{n-1}}$ up to enlarging $c$. Notice that when $g\equiv 0$ the level estimates rely solely on the $p$-coercivity, so the resulting integrability exponents do not depend on $q$.

Finally, the level estimates and our assumption on $q$ ensure that $u-g\in L^{q-1}(\Omega)$. Since $a\in L^\infty(\Omega)$, $g\in W^{1,H(\cdot)}(\Omega)$ and $h(\,\cdot\,,Du)\in L^1(\Omega)$, we conclude that $u-g,u\in W^{1,h(\cdot)}(\Omega)$, as desired.
\end{proof}

\begin{lemma}[Convergence of the approximations]\label{SOLA.flux.stability}
Assume that $2-1/n<p\leq n$, that \eqref{a.holder} and \eqref{bounds} hold, and that $A$ satisfies \eqref{A.veryweak}. Let $u$ be a local SOLA to \eqref{equation}, and let $\{u_j,\mu_j\}$ be an admissible approximating sequence; then
\begin{equation}\label{SOLA.strong.gradient}
u_j\longrightarrow u\quad\text{strongly in }W^{1,s}_{\rm loc}(\Omega)\qquad\text{for every }1\leq s<\frac{n(p-1)}{n-1}.
\end{equation}
Moreover,
\begin{equation}\label{SOLA.strong.flux}
h(\,\cdot\,,Du_j)\longrightarrow h(\,\cdot\,,Du)\qquad\text{and}\qquad A(\,\cdot\,,Du_j)\longrightarrow A(\,\cdot\,,Du)
\end{equation}
strongly, respectively, in $L^1(K)$ and $L^1(K;\R^n)$, for every compact set $K\Subset\Omega$.
\end{lemma}

\begin{proof}
Fix
\[
\max\{1,q-1\}<t<\frac{n(p-1)}{n-1}.
\]
For every ball $B_{2r}\Subset\Omega$, applying Lemma~\ref{local.rH} to $u_j$ with $\es=t/(p-1)$ yields
\[
\biggl(\mean{B_r}|Du_j|^t\dx\biggr)^{\frac{p-1}{t}}\leq c\mean{B_{2r}}h(x,Du_j)\dx+c\,\frac{|\mu_j|(B_{2r})}{r^{n-1}};
\]
by \eqref{SOLA.uniform.meas} and \eqref{SOLA.convergence}, the right-hand side is bounded independently of $j$, ensuring that $\{Du_j\}$ is locally bounded in $L^t(\Omega)$. Note also that weak convergence in $W^{1,1}_{\rm loc}(\Omega)$ implies $u_j\to u$ strongly in $L^1_{\rm loc}(\Omega)$.

\medskip

To show that $Du_j\to Du$ in measure, fix $B_{2r}\Subset\Omega$ and take $\eta\in C_c^\infty(B_{2r})$ with $0\leq\eta\leq1$ and $\eta\equiv1$ on $B_r$. For $\tau>0$, testing the difference of the equations for $u_j$ and $u_k$ with $\eta T_\tau(u_j-u_k)$ leads to
\begin{multline*}
\int_{B_{2r}\cap\{|u_j-u_k|<\tau\}}\big\langle A(x,Du_j)-A(x,Du_k),Du_j-Du_k\big\rangle\eta\dx\\
\leq\tau\bigl[|\mu_j|(B_{2r})+|\mu_k|(B_{2r})\bigr]+c\,\tau\|D\eta\|_{L^\infty}\int_{B_{2r}}\bigl[h(x,Du_j)+h(x,Du_k)\bigr]\dx\leq C_r\tau,
\end{multline*}
where $C_r$ does not depend on $j,k$ or $\tau$. The convergence-in-measure argument from Theorem~\ref{existence.SOLA} then applies on $B_r$ and shows that $Du_j\to Du$ in measure on $B_r$ for the whole sequence. Combining this with the uniform local $L^t$ bound, Vitali's theorem gives $Du_j\to Du$ strongly in $L^s_{\rm loc}(\Omega;\R^n)$ for all $1\leq s<t$; by Poincar\'e's inequality on balls and the $L^1_{\rm loc}(\Omega)$ convergence of $\{u_j\}$, we actually have $u_j\to u$ strongly in $W^{1,s}_{\rm loc}(\Omega)$. Letting $t\uparrow\frac{n(p-1)}{n-1}$ gives \eqref{SOLA.strong.gradient}. Finally, as $t>q-1$, the growth assumptions ensure uniform local $L^{\frac{t}{q-1}}$ bounds for $h(\,\cdot\,,Du_j)$ and $A(\,\cdot\,,Du_j)$, so \eqref{SOLA.strong.flux} follows again from convergence in measure and Vitali's theorem.
\end{proof}

Mean oscillation estimates remain stable under the approximations of Definition~\ref{def.SOLA}, enabling a direct transfer of the zero-order potential bounds from energy solutions to local SOLAs:

\begin{lemma}[Stability of oscillation estimates]\label{SOLA.oscillation.stability}
Let $u$ be a local SOLA in the sense of Definition~\ref{def.SOLA}, and suppose that the approximating sequence $\{u_j,\mu_j\}$ satisfies
\begin{equation}\label{SOLA.oscillation.approx}
\mean{B_{\theta R}(x_0)}|u_j-(u_j)_{B_{\theta R}(x_0)}|\dx\leq\frac12\mean{B_R(x_0)}|u_j-(u_j)_{B_R(x_0)}|\dx+c\,R h_{x_0}^{-1}\biggl( R\mean{B_R(x_0)}|\mu_j|\dx\biggr)
\end{equation}
for some $x_0\in\Omega$, $\theta\in(0,1/2)$, $R_*>0$ with $B_{R_*}(x_0)\Subset\Omega$, and $c>0$, and for every $j$ and every $0<R\leq R_*$. Then the same inequality holds with $u$ in place of $u_j$ and the normalized measure term $|\mu|(B_R(x_0))/|B_R(x_0)|$ in place of the integral average of $|\mu_j|$, for every $0<R\leq R_*$.
\end{lemma}

\begin{proof}
By \eqref{SOLA.convergence} and local compactness, $u_j\to u$ strongly in $L^1_{\rm loc}(\Omega)$, ensuring the convergence of mean oscillations on every ball compactly contained in $\Omega$. For every $0<R\leq R_*$ such that $|\mu|(\partial B_R(x_0))=0$, combining \eqref{SOLA.uniform.meas} with the continuity and monotonicity of $h_{x_0}^{-1}$ lets us pass to the limit in \eqref{SOLA.oscillation.approx}. Since exceptional radii form at most a countable set, the estimate extends to every $0<R\leq R_*$ by approximation from below, relying on the continuity of integral averages and the continuity from below of $\mu$.
\end{proof} 

\subsection{Local energy estimates}

For a fixed ball $B_R\Subset\Omega$, we introduce the short notation
\begin{equation}\label{dtwo}
d_1:=\mean{B_R}h(x,Du)\dx+\frac{|\mu|(B_R)}{R^{n-1}}.
\end{equation}
We next prove the following a priori estimate, a local version of Lemma~\ref{global.RH}, also covering the case $p>n$.

\begin{lemma}[Measure data reverse H\"older inequality]\label{local.rH}
Let $u\in W^{1,H(\cdot)}_{\rm loc}(\Omega)$ be a local energy solution to \eqref{equation}, under \eqref{a.holder} to \eqref{bounds} and \eqref{A.veryweak}. Then, for every ball $B_R\Subset\Omega$ and every exponent $\es$ satisfying
\begin{equation}\label{es}
\begin{cases}
\displaystyle 0<\es<\frac{n}{n-1},&\text{if }p\leq n,\\[3mm]
\displaystyle 0<\es\leq\frac{p}{p-1},&\text{if }p>n,
\end{cases}
\end{equation}
there holds
\begin{equation}\label{local.energy.SOLA}
\biggl(\mean{B_{R/2}}|Du|^{(p-1)\es}\dx\biggr)^{\frac{1}{\es}}\leq c_{rH}\mean{B_R}h(x,Du)\dx+c_{rH}\frac{|\mu|(B_R)}{R^{n-1}},
\end{equation}
where $c_{rH}$ depends only on $n,p,q,\nu,L$, and additionally on $\es$ only if $p\leq n$.
\end{lemma}

\begin{proof}
We prove a more general component-wise estimate, which will also be used later in Lemma~\ref{local.rH.pq}. Let $\esp\in\{p,q\}$, set $\bar a_p=1$, and let $\bar a_q\geq0$ be any constant such that $\bar a_q\leq a(x)$ for almost every $x\in B_R$; here and below, the subscripts $p$ and $q$ distinguish the unweighted and weighted terms, respectively, even when $p=q$. We claim that, whenever
\[
\begin{cases}
\displaystyle 0<\sigma<\frac{n}{n-1},&\qquad \text{if }\esp\leq n,\\[3mm]
\displaystyle 0<\sigma\leq\frac{\esp}{\esp-1},&\qquad \text{if }\esp>n,
\end{cases}
\]
there holds
\begin{equation}\label{component.local.rH}
\biggl(\mean{B_{R/2}}{\bigl[\bar a_\esp|Du|^{\esp-1}\bigr]}^\sigma\dx\biggr)^{\frac{1}{\sigma}}\leq c\biggl[\mean{B_R}h(x,Du)\dx+\frac{|\mu|(B_R)}{R^{n-1}}\biggr],
\end{equation}
where $c$ depends only on $n,p,q,\nu,L$ and, if $\esp\leq n$, also on $\sigma$; in particular, $c$ is independent of $\bar a_\esp$. The conclusion \eqref{local.energy.SOLA} follows by taking $\esp=p$, $\bar a_p=1$ and $\sigma=\es$. Here and below, $d_1$ is defined in \eqref{dtwo}. As \eqref{component.local.rH} is immediate if $\bar a_\esp=0$ when $\esp=q$, we may assume $\bar a_\esp>0$.

\medskip
\noindent\emph{Case $\esp\leq n$.} By H\"older's inequality, it is sufficient to consider exponents $\sigma\geq1/(\esp-1)$; note that this range is nonempty by \eqref{2-1/n}. We fix two exponents and a constant $\kappa$ that we can assume without loss of generality strictly larger than zero, since otherwise the conclusion is immediate:
\[
\zeta:=(\esp-1)\sigma,\qquad \xi:=\frac{n(\esp-\zeta)}{n-\zeta}>1,\qquad \kappa:=\biggl(\mean{B_{R/2}}{|u-(u)_{B_{R/2}}|}^{\zeta^*}\dx\biggr)^{\frac{1}{\zeta^*}}.
\]
We observe that $\zeta^*=n\zeta/(n-\zeta)$,
\[
1\leq\zeta<\esp\leq n\qquad\text{and}\qquad \frac{\xi\zeta}{\esp-\zeta}=\zeta^*.
\]
Now we define
\[
\psi(s):=\int_0^s(\kappa+|\tau|)^{-\xi}\,d\tau,\quad\text{so that}\quad |\psi(s)|\leq\frac{\kappa^{1-\xi}}{\xi-1},
\]
and choose a standard cut-off $\eta\in C_c^\infty(B_R)$ such that $0\leq\eta\leq1$, $\eta\equiv1$ on $B_{R/2}$ and $|D\eta|\leq c/R$. Testing \eqref{equation} with $\psi(u-(u)_{B_{R/2}})\eta^p$, and recalling that $\psi'(t)=(\kappa+|t|)^{-\xi}$, we obtain
\begin{align*}
\nu\int_{B_R}\frac{H(x,Du)\eta^p}{(\kappa+|u-(u)_{B_{R/2}}|)^\xi}\dx&\leq\int_{B_R}\frac{\langle A(x,Du),Du\rangle\eta^p}{(\kappa+|u-(u)_{B_{R/2}}|)^\xi}\dx\\
&=\int_{B_R}\psi(u-(u)_{B_{R/2}})\eta^p\,d\mu-p\int_{B_R}\psi(u-(u)_{B_{R/2}})\langle A(x,Du),D\eta\rangle\eta^{p-1}\dx\\
&\leq c(p,L)\frac{\kappa^{1-\xi}}{\xi-1}\biggl[|\mu|(B_R)+\|D\eta\|_{L^\infty(B_R)}\int_{B_R}h(x,Du)\dx\biggr].
\end{align*}
Using $|D\eta|\leq c/R$ and $\bar a_\esp|Du|^\esp\leq qH(x,Du)$, we conclude that
\begin{equation}\label{component.local.rH.weighted}
\mean{B_{R/2}}\frac{\bar a_\esp|Du|^\esp}{(\kappa+|u-(u)_{B_{R/2}}|)^\xi}\dx\leq c\frac{\kappa^{1-\xi}}{R}d_1.
\end{equation}
H\"older's inequality with conjugate exponents $\esp/\zeta$ and $\esp/(\esp-\zeta)$, together with the choice of $\kappa$ and \eqref{component.local.rH.weighted}, gives
\begin{align*}
\mean{B_{R/2}}{\bigl[\bar a_\esp|Du|^{\esp-1}\bigr]}^\sigma\dx&\leq\bar a_\esp^{\frac{\sigma}{\esp}}\biggl(\mean{B_{R/2}}\frac{\bar a_\esp|Du|^\esp}{(\kappa+|u-(u)_{B_{R/2}}|)^\xi}\dx\biggr)^{\frac{\zeta}{\esp}}\biggl(\mean{B_{R/2}}(\kappa+|u-(u)_{B_{R/2}}|)^{\zeta^*}\dx\biggr)^{\frac{\esp-\zeta}{\esp}}\\
&\leq c\,\bar a_\esp^{\frac{\sigma}{\esp}}{\Bigl(\frac{\kappa}{R}d_1\Bigr)}^{\frac{\zeta}{\esp}};
\end{align*}
on the other hand, Sobolev-Poincar\'e's inequality yields
\[
\frac{\kappa}{R}\leq c\biggl(\mean{B_{R/2}}|Du|^\zeta\dx\biggr)^{\frac{1}{\zeta}}=c\,\bar a_\esp^{-\frac{1}{\esp-1}}{\biggl(\mean{B_{R/2}}{\bigl[\bar a_\esp|Du|^{\esp-1}\bigr]}^\sigma\dx\biggr)}^{\frac{1}{\zeta}}.
\]
Substituting this into the previous estimate and recalling that $\zeta=(\esp-1)\sigma$, we obtain
\[
\mean{B_{R/2}}{\bigl[\bar a_\esp|Du|^{\esp-1}\bigr]}^\sigma\dx\leq c\,\bar a_\esp^{\frac{\sigma}{\esp}-\frac{\zeta}{\esp(\esp-1)}}{\biggl(\mean{B_{R/2}}{\bigl[\bar a_\esp|Du|^{\esp-1}\bigr]}^\sigma\dx\biggr)}^{\frac{1}{\esp}}d_1^{\frac{\zeta}{\esp}}=c{\biggl(\mean{B_{R/2}}{\bigl[\bar a_\esp|Du|^{\esp-1}\bigr]}^\sigma\dx\biggr)}^{\frac{1}{\esp}}d_1^{\frac{\zeta}{\esp}};
\]
dividing by the common factor when it is different from zero, we conclude with \eqref{component.local.rH} in this case.

\medskip
\noindent\emph{Case $\esp>n$.} It suffices to prove \eqref{component.local.rH} for $\sigma=\esp'$, since smaller exponents follow from H\"older's inequality. As $\bar a_\esp>0$, we have $u\in W^{1,\esp}(B_R)$ and may use Morrey-Poincar\'e's inequality. Fix $R/2\leq r<s\leq R$ and choose $\eta\in C_c^\infty(B_s)$ such that $0\leq\eta\leq1$, $\eta\equiv1$ on $B_r$ and $|D\eta|\leq c/(s-r)$. Testing \eqref{equation} with $(u-(u)_{B_s})\eta^p$, coercivity and growth give
\begin{align*}
\nu\int_{B_s}H(x,Du)\eta^p\dx&\leq\int_{B_s}\langle A(x,Du),Du\rangle\eta^p\dx\\
&=\int_{B_s}(u-(u)_{B_s})\eta^p\,d\mu-p\int_{B_s}(u-(u)_{B_s})\langle A(x,Du),D\eta\rangle\eta^{p-1}\dx\\
&\leq\|u-(u)_{B_s}\|_{L^\infty(B_s)}\biggl[|\mu|(B_s)+pL\|D\eta\|_{L^\infty(B_s)}\int_{B_s}h(x,Du)\dx\biggr],
\end{align*}
whence
\begin{align*}
\bar a_\esp\int_{B_r}|Du|^\esp\dx&\leq c\|u-(u)_{B_s}\|_{L^\infty(B_s)}\biggl[|\mu|(B_R)+\frac{1}{s-r}\int_{B_R}h(x,Du)\dx\biggr]\\
&\leq cR^{1-\frac{n}{\esp}}\biggl(\int_{B_s}|Du|^\esp\dx\biggr)^{\frac{1}{\esp}}\biggl[|\mu|(B_R)+\frac{1}{s-r}\int_{B_R}h(x,Du)\dx\biggr].
\end{align*}
Defining
\[
Z(t):=\bar a_\esp^{\esp'}\int_{B_t}|Du|^\esp\dx,
\]
and recalling the definition of $d_1$, the previous estimate can be rewritten as
\[
Z(r)\leq cR^{\frac{n}{\esp'}}\frac{R}{s-r}Z(s)^{\frac{1}{\esp}}d_1\leq\frac12Z(s)+cR^n\Bigl(\frac{R}{s-r}\Bigr)^{\esp'}d_1^{\esp'};
\]
the last inequality follows from Young's inequality with conjugate exponents $(\esp,\esp')$. At this point we can apply Lemma~\ref{Giusti.fundamental} on $[R/2,R]$ to conclude that $Z(R/2)\leq cR^nd_1^{\esp'}$, which proves \eqref{component.local.rH} at $\sigma=\esp'$. This completes the proof.
\end{proof}

\begin{remark}\label{rem:local.rH.SOLA}
{\rm
When $p\leq n$, the componentwise estimate \eqref{component.local.rH} holds for local SOLAs under the structural assumptions of Lemma~\ref{local.rH}, for both $\esp=p$ and $\esp=q$; in particular, the choice $\esp=p$ gives \eqref{local.energy.SOLA}. Indeed, we apply \eqref{component.local.rH} to an admissible approximating sequence and Fatou's lemma allows us to pass to the limit on the left-hand side for every admissible $\sigma$ (up to a subsequence), since by Lemma~\ref{SOLA.flux.stability} the gradients converge almost everywhere. On the right-hand side we can instead use \eqref{SOLA.strong.flux} and \eqref{SOLA.uniform.meas}, first for radii such that $|\mu|(\partial B_R)=0$ and then by approximation from below, keeping $\bar a_\esp$ fixed.

\medskip

Still assuming $p\leq n$, we also record a uniform local integrability bound. We let $b$ be the exponent fixed in \eqref{bbi}, and recall the definitions of $\Omega_0,\Omega_1$ and $\mathcal M$ in \eqref{Omega0}-\eqref{energy.M}. We also let $x_0\in\Omega_0$ and $0<R\leq\min\{1,d_{\Omega_0}/8\}$; we apply \eqref{local.energy.SOLA} on $B_{d_{\Omega_0}/4}(x_0)\subset\Omega_1$ with $\es=b$ to find 
\begin{equation}\label{global.hi.energy}
\int_{B_R(x_0)}|Du|^{b(p-1)}\dx\leq c\,d_{\Omega_0}^n\biggl[\frac{\mathcal M}{d_{\Omega_0}^n}+\frac{|\mu|(\Omega)}{d_{\Omega_0}^{n-1}}\biggr]^b\leq c(\data,\gdata),
\end{equation}
the dependence on $b$ being included in that on $\data$, by its fixed choice in \eqref{bbi}.
}
\end{remark}

In the $p$-phase, combining the preceding reverse H\"older inequality with the local bound \eqref{energy.M} controls the average of $h(x,Du)$ by that of $|Du|^{p-1}$ and the measure density of $\mu$ (a bound needed later in Proposition~\ref{decay.measure.data}).

\begin{lemma}[A reverse flux estimate in the $p$-phase]\label{reverse.flux.pphase}
Under the assumptions of Theorem~\ref{main.thm.2}, consider a ball $B_R=B_R(x_0)$ with $x_0\in\Omega_0$ and $0<R\leq\min\{1,d_{\Omega_0}/8\}$. If
\begin{equation}\label{p.phase.reverse.measure}
\sup_{B_R}a(\,\cdot\,)\leq M[a]_{0,\alpha}R^\alpha
\end{equation}
for some $M\geq2$, then
\begin{equation}\label{reverse.flux.pphase.est}
\mean{B_{R/2}}h(x,Du)\dx\leq c\mean{B_R}|Du|^{p-1}\dx+c\frac{|\mu|(B_R)}{R^{n-1}},
\end{equation}
where $c$ depends only on $\data,\gdata$ and $M$.
\end{lemma}

\begin{proof}
All constants in this proof depend only on $\data,\gdata$ and $M$, unless otherwise specified. We assume $q>p$, since for $q=p$ the conclusion follows directly from the fact that we are in the $p$-phase and $R\leq1$.

\medskip

\noindent{\em The case $p\leq n$.} We fix
\[
b:=\frac12\biggl[\max\biggl\{1,\frac{1}{p-1},\frac{n(q-p)}{\alpha(p-1)}\biggr\}+n'\biggr],
\]
consistently with the forthcoming \eqref{bbi}: by \eqref{2-1/n} and \eqref{bounds}, we have
\[
1<b<n',\qquad b(p-1)>1,\qquad \alpha-\frac{n(q-p)}{b(p-1)}>0.
\]
We set
\[
\es:=\frac{b}{b-\frac{q-p}{p-1}},\qquad \widehat{\es}:=\frac{\es+n'}{2};
\]
the choice of $b$ and the fact $\alpha\leq1$ imply
\[
1<\es<\widehat{\es}<n',\qquad \frac{q-p}{b(p-1)}+\frac1{\es}=1.
\]
For every $r\in[R/2,R]$, H\"older's inequality and \eqref{p.phase.reverse.measure} yield
\[
\mean{B_r}a(x)|Du|^{q-1}\dx\leq M[a]_{0,\alpha}R^\alpha\biggl(\mean{B_r}|Du|^{b(p-1)}\dx\biggr)^{\frac{q-p}{b(p-1)}}\biggl(\mean{B_r}|Du|^{\es(p-1)}\dx\biggr)^{\frac1{\es}};
\]
since $r\geq R/2$ and $R\leq1$, from \eqref{global.hi.energy} we deduce
\[
R^\alpha\biggl(\mean{B_r}|Du|^{b(p-1)}\dx\biggr)^{\frac{q-p}{b(p-1)}}\leq cR^{\alpha-\frac{n(q-p)}{b(p-1)}}\biggl(\int_{B_R}|Du|^{b(p-1)}\dx\biggr)^{\frac{q-p}{b(p-1)}}\leq c(\data,\gdata).
\]
Using H\"older's inequality also on the $p$-term gives
\begin{equation}\label{flux.to.lower.moment}
\mean{B_r}h(x,Du)\dx\leq c\biggl(\mean{B_r}|Du|^{\es(p-1)}\dx\biggr)^{\frac1{\es}},\qquad R/2\leq r\leq R.
\end{equation}
Next, fix $\theta\in(0,1)$ via
\[
\frac1{\es}=1-\theta+\frac{\theta}{\widehat{\es}}.
\]
For $R/2\leq r<s\leq R$, interpolation together with the inclusion $B_r\subset B_R$ yields
\[
\biggl(\mean{B_r}|Du|^{\es(p-1)}\dx\biggr)^{\frac1{\es}}\leq c\biggl(\mean{B_R}|Du|^{p-1}\dx\biggr)^{1-\theta}\biggl(\mean{B_r}|Du|^{\widehat{\es}(p-1)}\dx\biggr)^{\frac{\theta}{\widehat{\es}}}.
\]
Finally, covering $B_r$ with finitely many balls of radius $(s-r)/4$ having bounded overlap, whose doubles are contained in $B_s$, and applying \eqref{local.energy.SOLA} with exponent $\widehat{\es}$ on each doubled ball, in the form established for local SOLAs in Remark~\ref{rem:local.rH.SOLA}, we obtain
\[
\biggl(\mean{B_r}|Du|^{\widehat{\es}(p-1)}\dx\biggr)^{\frac1{\widehat{\es}}}\leq c\Bigl(\frac{R}{s-r}\Bigr)^n\biggl[\mean{B_s}h(x,Du)\dx+\frac{|\mu|(B_R)}{R^{n-1}}\biggr].
\]
Combining these estimates with \eqref{flux.to.lower.moment} via Young's inequality leads to
\[
\mean{B_r}h(x,Du)\dx\leq\frac12\mean{B_s}h(x,Du)\dx+c\Bigl(\frac{R}{s-r}\Bigr)^{\frac{n\theta}{1-\theta}}\mean{B_R}|Du|^{p-1}\dx+c\frac{|\mu|(B_R)}{R^{n-1}},
\]
and an application of Lemma~\ref{Giusti.fundamental} on $[R/2,R]$ then yields \eqref{reverse.flux.pphase.est}.
\medskip

\noindent{\em The case $p>n$.} By \eqref{bounds}, we have
\[
q-p\leq\frac{\alpha p}{n}\leq p,\qquad \alpha-\frac{n(q-p)}p\geq0.
\]
For every $s\in[R/2,R]$, H\"older's inequality, \eqref{p.phase.reverse.measure} and \eqref{energy.M} therefore give
\begin{align*}
\mean{B_s}a(x)|Du|^{q-p}\dx&\leq M[a]_{0,\alpha}R^\alpha\biggl(\mean{B_s}|Du|^p\dx\biggr)^{\frac{q-p}{p}}\\
&\leq cM[a]_{0,\alpha}R^{\alpha-\frac{n(q-p)}p}\biggl(\int_{\Omega_1}|Du|^p\dx\biggr)^{\frac{q-p}{p}}\leq cM[a]_{0,\alpha}R^{\alpha-\frac{n(q-p)}p}\mathcal M^{\frac{q-p}{p}}\leq c;
\end{align*}
consequently, H\"older's inequality with conjugate exponents $p$ and $p'$ yields
\[
\mean{B_s}a(x)|Du|^{q-1}\dx\leq\biggl(\mean{B_s}a(x)|Du|^{q-p}\dx\biggr)^{\frac1p}\biggl(\mean{B_s}a(x)|Du|^q\dx\biggr)^{\frac1{p'}}\leq c\biggl(\mean{B_s}H(x,Du)\dx\biggr)^{\frac1{p'}}.
\]
Using H\"older's inequality on the $p$-term as well, we get
\begin{equation}\label{flux.from.energy.pphase}
\mean{B_s}h(x,Du)\dx\leq c\biggl(\mean{B_s}H(x,Du)\dx\biggr)^{\frac1{p'}},\qquad R/2\leq s\leq R.
\end{equation}
Fix $R/2\leq r<s\leq R$ and pick a cut-off $\eta\in C_c^\infty(B_s)$ such that $0\leq\eta\leq1$, $\eta\equiv1$ on $B_r$, and $|D\eta|\leq c/(s-r)$. Following the proof of Lemma~\ref{local.rH}, we test \eqref{equation} with $(u-(u)_{B_s})\eta^p$, which belongs to the energy space, has compact support in $B_s$, and admits a continuous representative by Morrey's embedding. Retaining the full energy on the left-hand side, coercivity and growth properties yield
\begin{align*}
\int_{B_r}H(x,Du)\dx&\leq c\|u-(u)_{B_s}\|_{L^\infty(B_s)}\biggl[\frac1{s-r}\int_{B_s}h(x,Du)\dx+|\mu|(B_R)\biggr]\\
&\leq c\,\osc_{B_s}u\biggl[\frac1{s-r}\int_{B_s}h(x,Du)\dx+|\mu|(B_R)\biggr].
\end{align*}
We then set
\[
m:=\frac{p+\max\{n,p-1\}}2
\]
and define $\theta\in(0,1)$ through the relation
\[
\frac1m=\frac{1-\theta}{p-1}+\frac{\theta}{p}.
\]
Since $\max\{n,p-1\}<m<p$, combining Morrey's inequality, interpolation, and inclusion gives
\begin{align*}
\osc_{B_s}u\leq cs\biggl(\mean{B_s}|Du|^m\dx\biggr)^{\frac1m}&\leq cR\biggl(\mean{B_s}|Du|^{p-1}\dx\biggr)^{\frac{1-\theta}{p-1}}\biggl(\mean{B_s}|Du|^p\dx\biggr)^{\frac{\theta}{p}}\\
&\leq cR\biggl(\mean{B_R}|Du|^{p-1}\dx\biggr)^{\frac{1-\theta}{p-1}}\biggl(\mean{B_s}H(x,Du)\dx\biggr)^{\frac{\theta}{p}}.
\end{align*}
Inserting this into the Caccioppoli inequality above, dividing by $|B_r|$, and invoking \eqref{flux.from.energy.pphase}, we obtain
\begin{multline*}
\mean{B_r}H(x,Du)\dx\\
\leq c\,\frac{R}{s-r}\biggl(\mean{B_R}|Du|^{p-1}\dx\biggr)^{\frac{1-\theta}{p-1}}\biggl(\mean{B_s}H(x,Du)\dx\biggr)^{\frac{\theta}{p}}\biggl[\biggl(\mean{B_s}H(x,Du)\dx\biggr)^{\frac1{p'}}+\frac{|\mu|(B_R)}{R^{n-1}}\biggr].
\end{multline*}
Since
\[
\frac{\theta}{p}+\frac1{p'}=1-\frac{1-\theta}{p}<1,
\]
Young's inequality implies
\[
\mean{B_r}H(x,Du)\dx\leq\frac12\mean{B_s}H(x,Du)\dx+c\,\Bigl(\frac{R}{s-r}\Bigr)^{\gamma}\biggl[\biggl(\mean{B_R}|Du|^{p-1}\dx\biggr)^{p'}+\biggl(\frac{|\mu|(B_R)}{R^{n-1}}\biggr)^{p'}\biggr],
\]
where $\gamma = \frac{p}{1-\theta}$. Applying Lemma~\ref{Giusti.fundamental} on $[R/2,R]$, we conclude that
\[
\mean{B_{R/2}}H(x,Du)\dx\leq c\biggl[\biggl(\mean{B_R}|Du|^{p-1}\dx\biggr)^{p'}+\biggl(\frac{|\mu|(B_R)}{R^{n-1}}\biggr)^{p'}\biggr];
\]
combining this estimate with \eqref{flux.from.energy.pphase}, evaluated at $s=R/2$, proves \eqref{reverse.flux.pphase.est} also when $p>n$.
\end{proof}

In the $(p,q)$-phase, we have $a(\,\cdot\,)\approx_M a(x_0)$. This allows us to apply \eqref{component.local.rH} to the $q$-growth component and combine it with Lemma~\ref{local.rH}, obtaining a reverse H\"older estimate for the full frozen flux magnitude $h_{x_0}(Du)$. The admissible range of exponents is now determined by $q$.

\begin{lemma}[Measure data reverse H\"older inequality in the $(p,q)$-phase]\label{local.rH.pq}
Assume \eqref{a.holder} to \eqref{bounds} and \eqref{A.veryweak}, and let $u$ be a local solution to \eqref{equation} in the sense of Definition~\ref{def.solution}. Let moreover $B_R\Subset\Omega$ be a ball on which \eqref{pq.phase.sup.condition} holds for some $M>2$. Then, for every exponent $\ti$ satisfying
\begin{equation}\label{ti}
\begin{cases}
\displaystyle 0<\ti<\frac{n}{n-1},&\text{if }q\leq n,\\[2mm]
\displaystyle 0<\ti\leq\frac{q}{q-1},&\text{if }q>n,
\end{cases}
\end{equation}
there holds
\begin{equation}\label{local.energy.SOLA.pq}
\biggl(\mean{B_{R/2}}[h_{x_0}(Du)]^\ti\dx\biggr)^{\frac{1}{\ti}}\leq c_{rH}\mean{B_R}h(x,Du)\dx+c_{rH}\frac{|\mu|(B_R)}{R^{n-1}},
\end{equation}
where $c_{rH}$ depends only on $n,p,q,\nu,L,M$, and additionally on $\ti$ only if $q\leq n$.

\medskip

The same conclusion holds if the phase assumption is replaced by
\[
0<\inf_{B_R}a(\,\cdot\,)\leq\sup_{B_R}a(\,\cdot\,)\leq2\inf_{B_R}a(\,\cdot\,),
\]
with $c_{rH}$ depending only on $n,p,q,\nu,L$ and additionally on $\ti$ only if $q\leq n$.
\end{lemma}

\begin{proof}
Set $\bar a_q:=\inf_{B_R}a(\,\cdot\,)$. Under \eqref{pq.phase.sup.condition}, estimate \eqref{above.below} gives
\[
0<\bar a_q\leq a(x_0)\leq\frac{M}{M-2}\bar a_q
\]
and under the alternative assumption, the same bound holds with $M/(M-2)$ replaced by $2$. Applying \eqref{component.local.rH} with $\esp=q$ and $\sigma=\ti$, we therefore obtain
\begin{equation}\label{local.energy.SOLA.partq}
\mean{B_{R/2}}{\bigl[a(x_0)|Du|^{q-1}\bigr]}^\ti\dx\leq c\,d_1^\ti,
\end{equation}
where $d_1$ is as in \eqref{dtwo}. Since the range in \eqref{ti} is contained in that in \eqref{es}, applying \eqref{component.local.rH} with $\esp=p$ and $\sigma=\ti$ gives the corresponding bound for the $p$-component. When $p\leq n$, both applications of \eqref{component.local.rH} to the local SOLA $u$ are justified by Remark~\ref{rem:local.rH.SOLA}. Combining the two bounds yields \eqref{local.energy.SOLA.pq}.

Finally, when $q>n$, both component estimates can be applied at the fixed exponent $\ti=q'$, which is admissible also for the $p$-component, and  the smaller exponents then follow from H\"older's inequality, with constants independent of $\ti$.
\end{proof}

\section{Homogeneous problems and estimates below the energy level}\label{cor:nice-rev.sec}

We consider the homogeneous equation
\begin{equation}\label{homogeneous.equation}
 -\divergence  A(x,Dv)=0 \qquad\text{in $B_R(x_0)\subset\R^n$} 
\end{equation}
with $B_R=B_R(x_0)$ and $R\leq1$; a solution $v\in W^{1,H(\cdot)}(B_R)$ of \eqref{homogeneous.equation} is understood to satisfy the weak formulation \eqref{en.solution} of Definition~\ref{def.local.energy.solution} with $u=v$ and $\mu=0$. All the results of this section hold for irregular vector fields $A$ satisfying \eqref{A.weak}.

\subsection{The fundamental exponent}\label{exponent.b}

For $p\leq n$, we fix an exponent $b$ in the following range
\begin{equation}\label{def.b}
\mathfrak m:=\max\Big\{1,\frac{1}{p-1},\frac{(q-p)n}{\alpha(p-1)}\Big\}\leq  b<\frac{n}{n-1}.
\end{equation}

\medskip

We notice that 
\[
\frac{1}{p-1} < \frac{n}{n-1}\qquad \text{if and only if} \qquad p > 2 -\frac1n 
\]
so that the lower bound in \eqref{2-1/n} is necessary for the range in \eqref{def.b} to be nonempty. Moreover, for the bound from below for $b$, we have
\[
\frac{(q-p)n}{\alpha(p-1)}<\frac n{n-1}\qquad \text{if and only if} \qquad \frac{q-1}{p-1} < 1 +\frac\alpha n \cdot \frac n{n-1}.
\]
We are going to determine a fixed value of $b$, depending only on $n,p,q$ and $\alpha$: the midpoint of the two bounds in \eqref{def.b}, namely
\begin{equation}\label{bbi}
 b=\frac12\Bigl(\mathfrak m+\frac n{n-1}\Bigr);
\end{equation}
none of the constants below will therefore depend ultimately on $b$. Note that we already employed similar exponents in the proof of Lemma~\ref{reverse.flux.pphase}. With this just fixed choice of $b$ satisfying \eqref{def.b}, we set
\begin{equation}\label{esse}
\mathfrak s:= \begin{cases} b(p-1) &\text{if }p\leq n,\\ 
p &\text{if }p>n, \end{cases} 
\end{equation}
and
\begin{equation*}
\benergy(f,B_R(x_0)):={\|f\|}_{L^{\mathfrak s}(B_R(x_0))}
=\begin{cases}
{\|f\|}_{L^{b(p-1)}(B_R(x_0))}\qquad &\text{if $p\leq n$}\\[2mm]
{\|f\|}_{L^p(B_R(x_0))} &\text{if $p> n$}
\end{cases},
\end{equation*}
with $b>1$ defined in \eqref{bbi}; notice that $1<b(p-1)<p$. We shall occasionally use  $\benergy(f,\Omega)$ for $\Omega$ not necessarily a ball, with the same definition. We also stress that the bound from below on $b$ implies that
\begin{equation}\label{rsmaller1}
r^{\alpha-\frac{n(q-p)}{\mathfrak s}}\leq 1\qquad\qquad\text{if $r\leq 1$}, 
\end{equation}
a fact that we stress here for future reference, since it will be used countless times in the paper.

\medskip

We both collect some results present in \cite{BCMgen, CM1, DM} and prove new ones. Particular attention will be devoted to the dependence of both constants and exponents on the energy of $Dv$. In order to have a more transparent presentation, we shall assume the measure-data energy bound
\begin{equation}\label{bound.energy.mdata}
\benergy(Dv,B_R(x_0))\leq E\qquad \Bigg(\text{iff}\ \ 
\begin{cases}
{\|Dv\|}_{L^{b(p-1)}(B_R(x_0))}\leq  E\quad &\text{if $p\leq n$}\\[2mm]
{\|Dv\|}_{L^p(B_R(x_0))}\leq  E &\text{if $p> n$}
\end{cases}\Bigg)
\end{equation}
for some positive constant $E>0$.

\subsection{Basic regularity results tailored to our needs}

We start by stating a Sobolev-Poincar\'e-type result we are going to use several times: 

\begin{proposition}\label{prop:intrinsic.sobolev.poincare}
Let $a\in C^{0,\alpha}(B_R(x_0))$, with $R\leq1$, be non-negative, and assume that $p$ and $q$ satisfy \eqref{2-1/n} and \eqref{bounds}. Moreover let $\benergy$ be defined as above, with $b$ fixed as in \eqref{def.b}-\eqref{bbi} when $p\leq n$. Then there exist a constant $c=c(n,p,q,\alpha)$ and two exponents $d_1=d_1(n,p,q,\alpha)>1$, $d_2=d_2(n,p,q,\alpha)\in(0,1)$, such that every $w\in W^{1,H(\cdot)}(B_R(x_0))$ satisfies
\begin{multline}\label{Sob.Poinc.H}
{\biggl(\mean{B_R(x_0)}{\biggl[H{\Bigl(\,\cdot\,,\frac{w-(w)_{B_R(x_0)}}R\Bigr)}\biggr]}^{d_1}\dx\biggr)}^{\frac1{d_\smallone}}\\
\leq c{\Bigl(1+[a]_{0,\alpha}\bigl[\benergy(Dw,B_R(x_0))\bigr]^{q-p}\Bigr)}{\biggl(\mean{B_R(x_0)}{\bigl[H(\,\cdot\,,Dw)\bigr]}^{d_2}\dx\biggr)}^{\frac1{d_\smalltwo}}.
\end{multline}
The same estimate holds for every $w\in W^{1,H(\cdot)}_0(B_R(x_0))$, with $w$ replacing $w-(w)_{B_R(x_0)}$ on the left-hand side.
\end{proposition}

\begin{proof}
We follow the proof of \cite[Theorem 1.6]{CM1}, keeping the exponent $\mathfrak s$ defined in \eqref{esse} fixed throughout the argument. By \eqref{bounds} and, when $p\leq n$, the choice of $b$ in \eqref{bbi}, we have
\begin{equation}\label{condition.H.Sobolev}
\alpha-\frac{n(q-p)}{\mathfrak s}\geq 0.
\end{equation}
We extend $a(\,\cdot\,)$ to $\mathbb R^n$, preserving non-negativity and its H\"older seminorm, and keep the same notation. Let $B_\rho(\bar x)$ be a ball with $0<\rho\leq2$, and let $f\in L^{\mathfrak s}(B_\rho(\bar x))$ satisfy $H(\,\cdot\,,f)\in L^1(B_\rho(\bar x))$. We claim that, for a constant $c=c(n,p,q,\alpha)$,
\begin{equation}\label{H1.first}
H\bigl(\bar x,(f)_{B_\rho(\bar x)}\bigr) \leq c \left(1+{[a]}_{0,\alpha} \|f\|_{L^{\mathfrak s}(B_\rho(\bar x))}^{q-p} \right) \mean{B_\rho(\bar x)} H(x,f(x))\dx.
\end{equation}
We first observe that H\"older's inequality gives
\[
\big|(f)_{B_\rho(\bar x)}\big|\leq{\biggl(\mean{B_\rho(\bar x)}|f|^{\mathfrak s}\dx\biggr)}^{\frac1{\mathfrak s}}\leq c(n,\mathfrak s)\rho^{-\frac n{\mathfrak s}}\|f\|_{L^{\mathfrak s}(B_\rho(\bar x))}.
\]
Freezing $a(\,\cdot\,)$ at its infimum over $B_\rho(\bar x)$ and using Jensen's inequality, followed by the previous estimate, we obtain
\begin{align*}
H\bigl(\bar x,(f)_{B_\rho(\bar x)}\bigr)&\leq\mean{B_\rho(\bar x)}H(x,f(x))\dx+c[a]_{0,\alpha}\rho^\alpha\big|(f)_{B_\rho(\bar x)}\big|^q\\
&\leq c\Bigl(1+[a]_{0,\alpha}\rho^{\alpha-\frac{n(q-p)}{\mathfrak s}}\|f\|_{L^{\mathfrak s}(B_\rho(\bar x))}^{q-p}\Bigr)\mean{B_\rho(\bar x)}H(x,f(x))\dx:
\end{align*}
as $0<\rho\leq2$, the condition in \eqref{condition.H.Sobolev} gives $\rho^{\alpha-\frac{n(q-p)}{\mathfrak s}}\leq2$, and \eqref{H1.first} follows. Let now $\Omega\subset\mathbb R^n$ be, in general, a bounded open set. For $f\in L^{\mathfrak s}(\Omega)$, we denote by $M_\Omega f$ the restricted maximal function
\[
M_\Omega f(\bar x):=\sup_{\substack{B_\rho(\bar x)\subset\Omega\\0<\rho\leq2}}\mean{B_\rho(\bar x)}|f(x)|\dx,\qquad \bar x\in\Omega.
\]
For every $t\geq1$ and every $f\in L^{\mathfrak s}(\Omega)$ such that $H(\,\cdot\,,f)\in L^t(\Omega)$, we claim that, for a constant $c=c(n,p,q,\alpha,t)$,
\begin{equation}\label{maxi.H}
\int_\Omega \bigl[H(x,M_\Omega f(x))\bigr]^t\dx \leq c \left( 1+ {[a]}_{0,\alpha}^t \|f\|_{L^{\mathfrak s}(\Omega)}^{(q-p)t} \right) \int_\Omega [H(x,f(x))]^t\dx.
\end{equation}
For $t>1$ this indeed follows by taking the supremum in \eqref{H1.first}, integrating over $\Omega$ and using the Hardy-Littlewood inequality, while for $t=1$ we introduce the auxiliary exponent $1<{\textnormal{\texttt{E}}}<p$, and define $\hat p:=\frac p{\textnormal{\texttt{E}}}$, $\hat q:=\frac q{\textnormal{\texttt{E}}}$, $\hat\alpha:=\frac\alpha{\textnormal{\texttt{E}}}$, $\hat a(x):=a(x)^{1/{\textnormal{\texttt{E}}}}$, together with
\[
\widehat H(x,t) := |t|^{\widehat p} + \widehat a(x)|t|^{\widehat q};
\]
we have $[\hat a]_{0,\hat\alpha} \leq {[a]}_{0,\alpha}^{1/{\textnormal{\texttt{E}}}}$ and $[H(x,t)]^{1/{\textnormal{\texttt{E}}}} \approx_{p,q,\textnormal{\texttt{E}}} \widehat H(x,t)$. The exponent $\mathfrak s$ remains fixed, and
\[
\hat\alpha-\frac{n(\hat q-\hat p)}{\mathfrak s}=\frac{1}{\textnormal{\texttt{E}}}\Bigl(\alpha-\frac{n(q-p)}{\mathfrak s}\Bigr)\geq0.
\]
Since $\hat p>1$, we can repeat the proof of \eqref{H1.first} for $\widehat H$ and $|f|$ while preserving the same $L^{\mathfrak s}$ norm. Taking the supremum over admissible balls gives
\begin{equation}\label{point.maximal.Hhat}
\widehat H(x,M_\Omega f(x)) \leq c \left( 1+ {[a]}_{0,\alpha}^\frac1{{\textnormal{\texttt{E}}}} \|f\|_{L^{\mathfrak s}(\Omega)}^{\frac{q-p}{{\textnormal{\texttt{E}}}}} \right) M_\Omega\bigl(\widehat H(\,\cdot\,,f)\bigr)(x).
\end{equation}
Raising \eqref{point.maximal.Hhat} to the power ${\textnormal{\texttt{E}}}t$, integrating over $\Omega$, and using the strong $L^{{\textnormal{\texttt{E}}}t}$-boundedness of the maximal operator (recall that ${\textnormal{\texttt{E}}}t>1$), we obtain
\[
\int_\Omega \bigl[ \widehat H(x,M_\Omega f(x)) \bigr]^{{\textnormal{\texttt{E}}} t}\dx \leq c \left( 1+ {[a]}_{0,\alpha}^{t} \|f\|_{L^{\mathfrak s}(\Omega)}^{(q-p)t} \right) \int_\Omega [\widehat H(x,f(x))]^{{\textnormal{\texttt{E}}} t}\dx.
\]
Estimate \eqref{maxi.H} then follows from $[H(x,t)]^{1/{\textnormal{\texttt{E}}}} \approx_{p,q,\textnormal{\texttt{E}}} \widehat H(x,t)$.

We now prove a preliminary Sobolev-Poincar\'e estimate. We denote by $\widetilde D:=Dw\,\chi_{B_R(x_0)}$ the null extension of $Dw$, and write $\tilde w:=(w-(w)_{B_R(x_0)})/R$ and $M:=M_{B_{3R}(x_0)}$. By the standard potential representation and a decomposition into dyadic annuli, see \cite[Proof of Theorem 1.6]{CM1}, for almost every $\bar x\in B_R(x_0)$ and every $\varepsilon\in(0,1]$, we have
\begin{equation}\label{w.hedberg}
|\tilde w(\bar x)|\leq c(n)\varepsilon M(\widetilde D)(\bar x)+\frac{c(n)}{\varepsilon^{n-1}}\mean{B_{2R}(\bar x)}|\widetilde D(x)|\dx.
\end{equation}
Using the growth properties of $H$, we obtain
\[
H(\bar x,\tilde w(\bar x)) \leq c\,\varepsilon^p H\bigl(\bar x,M(\tilde D)(\bar x)\bigr) + c\,\varepsilon^{-q(n-1)} H\biggl(\bar x, \mean{B_{2R}(\bar x)}|\tilde D(x)|\dx\biggr);
\]
we again choose
\[
\varepsilon := \left[ \frac{1}{ H\Bigl(\bar x,M(\tilde D)(\bar x)\Bigr)}\ H\biggl(\bar x, \mean{B_{2R}(\bar x)}|\tilde D(x)|\dx \biggr)\right]^{ \frac1{p+q(n-1)} } \in(0,1],\quad\text{by}\quad \mean{B_{2R}(\bar x)}|\tilde D(x)|\dx \leq M(\tilde D)(\bar x),
\]
with the usual convention when the denominator vanishes. This yields
\begin{equation}\label{pointwise.preliminary.H}
H(\bar x,\tilde w(\bar x)) \leq c\, \Big[ H\bigl(\bar x,M(\tilde D)(\bar x)\bigr) \Big]^{\frac{q(n-1)}{p+q(n-1)}}  \biggl[ H\biggl(\bar x, \mean{B_{2R}(\bar x)}|\tilde D(x)|\dx \biggr) \biggr]^{\frac p{p+q(n-1)}}.
\end{equation}
We now set
\[
d_H := 1+\frac p{q(n-1)} = \frac{p+q(n-1)}{q(n-1)} >1 \quad\text{and}\quad \Lambda_H := 1+ {[a]}_{0,\alpha} {\bigl[ \benergy(Dw,B_R(x_0)) \bigr]}^{q-p}.
\]
Applying \eqref{H1.first} to $f=|\widetilde D|$ on $B_{2R}(\bar x)$, we obtain
\[
H\biggl( \bar x, \mean{B_{2R}(\bar x)}|\tilde D(x)|\dx\biggr) \leq c\,\Lambda_H \mean{B_{2R}(\bar x)} H(x,\tilde D(x))\dx \leq c\,\Lambda_H \mean{B_R(x_0)} H(x,Dw(x))\dx; 
\]
we also used the fact that $\tilde D=0$ outside $B_R(x_0)$ and $\benergy(\tilde D,B_{2R}(\bar x)) = \benergy(Dw,B_R(x_0))$. Raising \eqref{pointwise.preliminary.H} to the power $d_H$, inserting the last estimate, averaging over $B_R(x_0)$, and applying \eqref{maxi.H} with $t=1$ on $B_{3R}(x_0)$, we obtain
\begin{equation}\label{preliminary.Sobolev.H}
{\biggl(\mean{B_R(x_0)} [H(x,\tilde w(x))]^{d_\smallH}\dx\biggr)}^\frac1{d_\smallH} \leq c\,\Lambda_H  \mean{B_R(x_0)} H(x,Dw(x))\dx . 
\end{equation}
We finally improve \eqref{preliminary.Sobolev.H} as follows. We choose $1<{\textnormal{\texttt{E}}}_H<\min\{p,d_H\}$ and define $\bar p:=\frac p{{\textnormal{\texttt{E}}}_H}$, $\bar q:=\frac q{{\textnormal{\texttt{E}}}_H}$, $\bar\alpha:=\frac{\alpha}{{\textnormal{\texttt{E}}}_H}$, $\bar a(x):=a(x)^{1/{\textnormal{\texttt{E}}}_H}$, as well as $\bar H(x,t) := |t|^{\bar p} + \bar a(x)|t|^{\bar q}$; then, $[\bar a]_{0,\bar\alpha} \leq {[a]}_{0,\alpha}^{1/{\textnormal{\texttt{E}}}_H}$ and $[H(x,t)]^{1/{\textnormal{\texttt{E}}}_H} \approx \bar H(x,t)$. Keeping $\mathfrak s$ fixed, we have
\[
\bar\alpha-\frac{n(\bar q-\bar p)}{\mathfrak s}=\frac{1}{{\textnormal{\texttt{E}}}_H}\Bigl(\alpha-\frac{n(q-p)}{\mathfrak s}\Bigr)\geq0.
\]
We can therefore repeat the argument leading to \eqref{preliminary.Sobolev.H} with $H$ replaced by $\bar H$, retaining the same $L^{\mathfrak s}$ norm, and we find
\[
{\biggl( \mean{B_R(x_0)} [\bar H(x,\tilde w)]^{d_\smallH}\dx \biggr)}^{\frac1{d_\smallH}} \leq c \Bigl( 1+ {[a]}_{0,\alpha}^\frac1{{\textnormal{\texttt{E}}}_\smallH} {\benergy(Dw,B_R(x_0))}^\frac{q-p}{{\textnormal{\texttt{E}}}_\smallH} \Bigr) \mean{B_R(x_0)} \bar H(x,Dw)\dx.
\]
We choose $d_1:={d_H}/{{\textnormal{\texttt{E}}}_H}>1$ and $d_2:=1/{{\textnormal{\texttt{E}}}_H}\in(0,1)$. Raising the previous estimate to the power ${\textnormal{\texttt{E}}}_H$, employing the equivalence between $[H(x,t)]^{1/{\textnormal{\texttt{E}}}_\smallH}$ and $\bar H(x,t)$, and observing that
\[
{\Bigl(1+ {[a]}_{0,\alpha}^{\frac1{{\textnormal{\texttt{E}}}_\smallH}} \benergy(Dw,B_R(x_0))^{\frac{q-p}{{\textnormal{\texttt{E}}}_\smallH}}\Bigr)}^{{\textnormal{\texttt{E}}}_\smallH} \leq c\,\Lambda_H,
\]
we conclude with \eqref{Sob.Poinc.H}. When $w\in W^{1,H(\cdot)}_0(B_R(x_0))$, the same argument applies with $\tilde w(x):=w(x)/R$, since the potential estimate \eqref{w.hedberg} holds also for this choice of $\tilde w$; see also \cite[Remark 2]{CM1}.
\end{proof}

The first ``new'' result we propose as a consequence of Sobolev's inequality is an intrinsic higher integrability result for solutions to \eqref{homogeneous.equation}, which differs from that in \cite{CM1} in the dependence of the constants.

\begin{theorem}\label{hi.correct.exponent}
Let $v\in W^{1,H(\cdot)}(B_R)$ be a solution to the equation \eqref{homogeneous.equation} and let $E$ satisfy \eqref{bound.energy.mdata}. Then there exists a higher integrability exponent $\eta_0\equiv \eta_0(\data,E)>0$ such that $H(x,Dv)\in L^{1+\eta_0}_{\rm loc}(B_R)$; moreover for any ball $ B_\rho(\bar x) \subset B_R$ (not necessarily concentric with $B_R$) and for any $\theta \in[1/2,3/4]$ the local estimate
\begin{equation}\label{hi.local}
\biggl(\mean{B_{\theta\rho}(\bar x)} \bigl[H(x,Dv)\bigr]^{1+\eta_\smallzero}\dx\biggr)^{\frac{1}{1+\eta_\smallzero}}\leq c_h\mean{B_\rho(\bar x)}H(x,Dv)\dx
\end{equation}
holds for a constant $c_h$ depending on $\data$ and $E$. 
\end{theorem}

\begin{proof}
As before, we need only to replace the dependence of the constant and of the exponent with a more appropriate one. This is however immediate once we have at hand the Sobolev inequality of Proposition~\ref{prop:intrinsic.sobolev.poincare}. Indeed, the classic way to prove \eqref{hi.local} is to combine the Caccioppoli inequality with a Sobolev-type inequality: for every ball $B_{\rho}(\bar x)\subset B_R$, merging \eqref{Sob.Poinc.H} and \eqref{Caccio2} (with a suitable choice of radii) yields
\[
\mean{B_{\rho/2}(\bar x)}H(x,Dv)\dx \leq c{\biggl( \mean{B_\rho(\bar x)} [H(x,Dv)]^{d_2}\dx \biggr)}^\frac1{d_\smalltwo},
\]
where $d_2\in(0,1)$ and $c=c(\data,E)$. At this point \eqref{hi.local} follows from the standard self-improving character of reverse-H\"older inequalities (see for instance \cite[Theorem 6.6]{Giusti}).
\end{proof}

\begin{remark}
{\rm When $p>n$, we decrease $\eta_0$, if necessary, so that \eqref{hi.local} also holds with $2\eta_0$ in place of $\eta_0$ and
\[
\frac{n\eta_0}{1+\eta_0}<\frac{p-n}{p-1}.
\]
This choice depends only on $\data$ and $E$ and preserves the stated dependencies of the constants. 
}
\end{remark}

\begin{remark}
{\rm The well-known self-improving character of reverse H\"older's inequalities implies that the higher integrability estimate \eqref{hi.local} can be improved in the following way: for any $\sigma>0$ we have
\begin{equation}\label{rev.self.improved}
\mean{B_{\theta\rho}} H(x,Dv)\dx\leq c{\biggl(\mean{B_\rho}\bigl[H(x,Dv)\bigr]^\sigma\dx\biggr)}^{\frac1\sigma},  
\end{equation}
uniformly for $\theta\in [1/2,3/4]$, with the constant $c$ depending now on $\data,E$ and $\sigma$ too.  }
\end{remark}

Next, we state the Caccioppoli inequality for local energy solutions to \eqref{homogeneous.equation}. Its proof uses only the coercivity and growth bounds in \eqref{A.weak}; see also \cite[Proposition 1]{BCMgen} for the variational counterpart.

\begin{lemma}[Caccioppoli inequality]\label{Lemma.Caccio}
Let $v\in W^{1,H(\cdot)}(B_R)$ be a weak solution to \eqref{homogeneous.equation} under the assumptions \eqref{A.weak}. Then there exists a constant $c\geq1$, depending only on $n,p,q$ and $\ratio$, such that
\begin{equation}\label{Caccio}
\int_{B_\varrho}H\bigl(x,D(v-k)_\pm\bigr)\eta^q\dx\leq c\int_{B_\varrho}H\bigl(x,(v-k)_\pm|D\eta|\bigr)\dx
\end{equation}
for every ball $B_\varrho\subseteq B_R$ (not necessarily concentric), every $k\in\mathbb R$ and every cut-off $\eta\in C_c^\infty(B_\varrho)$ with $0\leq\eta\leq1$. In particular,
\begin{equation}\label{Caccio2}
\int_{B_\varrho}H(x,Dv)\eta^q\dx\leq c\int_{B_\varrho}H\bigl(x,(v-(v)_{B_\varrho})|D\eta|\bigr)\dx.
\end{equation}
\end{lemma}

\begin{proof}
We test \eqref{homogeneous.equation} with the test functions $\pm(v-k)_\pm\eta^q$; by coercivity and growth we have
\[
\int_{B_\varrho}H(x,D(v-k)_\pm)\eta^q\dx\leq\frac{qL}{\nu}\int_{B_\varrho}h(x,D(v-k)_\pm)(v-k)_\pm|D\eta|\eta^{q-1}\dx.
\]
Now we apply Young's inequality separately to the $p$- and $q$-parts and use $\eta^{q-p}\leq1$ to get
\[
\int_{B_\varrho}H(x,D(v-k)_\pm)\eta^q\dx\leq\frac12\int_{B_\varrho}H(x,D(v-k)_\pm)\eta^q\dx+c\int_{B_\varrho}H\bigl(x,(v-k)_\pm|D\eta|\bigr)\dx
\]
with $c=c(p,q,\ratio)$, so that \eqref{Caccio} follows by absorption. Adding the inequalities for the two signs and taking $k=(v)_{B_\varrho}$ gives \eqref{Caccio2}.
\end{proof}

As a direct consequence of Lemma~\ref{Lemma.Caccio}, we have the following version of the Caccioppoli inequality.
\begin{corollary}[Product Caccioppoli inequality]\label{Caccio.product}
Let $v\in W^{1,H(\cdot)}(B_R)$ be a solution to \eqref{homogeneous.equation}. For $B_\varrho\subset B_R$, $k\in\mathbb R$ and $\eta\in C_c^\infty(B_\varrho)$ as in Lemma~\ref{Lemma.Caccio}, we have
\begin{equation}\label{Caccio.product.estimate}
\int_{B_\varrho} H\bigl(x,D[(v-k)_\pm\eta^m]\bigr)\dx \leq c \int_{B_\varrho} H\bigl(x,(v-k)_\pm|D\eta|\bigr)\dx,
\end{equation}
for $m:=\left\lceil q/p\right\rceil$ and $c\equiv c(n,p,q,\ratio)$.
\end{corollary}

\begin{proof}
Since $|D\eta^m|=m\eta^{m-1}|D\eta| \leq m|D\eta|$, using the product rule and the growth properties of $H$, we obtain
\begin{align*}
H\bigl(x,D[(v-k)_\pm\eta^m]\bigr) &\leq c(p,q)\Big[ \eta^{mp}|D(v-k)_\pm|^p + (v-k)_\pm^p|D\eta|^p\\
&\hspace{4cm}+ a(x)\eta^{mq}|D(v-k)_\pm|^q + a(x)(v-k)_\pm^q|D\eta|^q \Big] \\
&\leq c(p,q)\Big[ \eta^q H\bigl(x,D(v-k)_\pm\bigr) + H\bigl(x,(v-k)_\pm|D\eta|\bigr)\Big]
\end{align*}
as $(\eta^m)^p=\eta^{mp}\leq\eta^q$ and $(\eta^m)^q=\eta^{mq}\leq\eta^q$; integrating over $B_\varrho$ and applying \eqref{Caccio} with cutoff $\eta$, we get
\begin{align*}
\int_{B_\varrho}H\bigl(x,D[(v-k)_\pm\eta^m]\bigr)\dx &\leq c \int_{B_\varrho}\eta^qH\bigl(x,D(v-k)_\pm\bigr)\dx + c \int_{B_\varrho}H\bigl(x,(v-k)_\pm|D\eta|\bigr)\dx\\
&\leq c \int_{B_\varrho}H\bigl(x,(v-k)_\pm|D\eta|\bigr)\dx. 
\end{align*}
\end{proof}

The next result we need is a $\sup$-$BMO$-estimate.
\begin{theorem}\label{sup.correct.exponent}
Let $v\in W^{1,H(\cdot)}(B_R)$ be a solution to \eqref{homogeneous.equation} and suppose that \eqref{bound.energy.mdata} holds. Let $B_\rho\equiv B_\rho(\bar x)\subset B_R$ be a ball, not necessarily concentric. Then the local estimate
\begin{equation}\label{sup.local}
\sup_{B_{\rho/2}}\biggl[\frac{v-(v)_{B_{\rho}}}\rho\biggr]_\pm \leq  c\mean{B_\rho}|Dv|\dx
\end{equation}
holds for a constant $c$ depending on $\data$ and $E$.
\end{theorem}

\begin{proof}
We follow the standard De Giorgi iteration scheme used to prove local boundedness; however, in view of the subsequent applications to measure data problems, we need an estimate with lower integrability on the right-hand side and we must pay particular attention to the dependence of the constants on the energy of $\tilde v$. Although local boundedness for solutions to \eqref{homogeneous.equation} was established in \cite{CM1}, we must modify the argument at several points. Thus, we start by rescaling
\begin{equation}\label{rescaled}
\tilde v(x)=\frac{v(\bar x+\rho x)-(v)_{B_{\rho }}}{\rho},\qquad x\in B_1\equiv B_1(0),
\end{equation}
which is a solution in $B_1$ to an equation structurally similar to \eqref{homogeneous.equation}, where the switching coefficient is
\begin{equation}\label{a.rescal}
\tilde a(x)=a(\bar x+\rho x)\qquad \text{with}\qquad  {[\tilde a]}_{C^{0,\alpha}(B_1)}\leq \rho ^\alpha{[a]}_{0,\alpha}; 
\end{equation}
we set $\tilde H(x,z)=|z|^p+\tilde a(x)|z|^q$. 

\medskip

For reasons that will become clearer later, we first focus on the range $p\leq n$: in this case, we need to use the fact that $\tilde v$ solves an equation, and to do so we fix two radii $0<r<s\leq 1$, two intermediate ones $r\leq \varrho_1<\varrho_2\leq s$ and their midpoint $\tilde \varrho=(\varrho_1+\varrho_2)/2$. Using a standard cut-off function $\zeta\in C_c^\infty(B_{\tilde\varrho})$, $0\leq\zeta\leq1$, $\zeta\equiv1$ on $B_{\varrho_1}$, $\|D\zeta\|_{L^\infty} \leq {4}/(\varrho_2-\varrho_1)$, we set $m:=\lceil q/p\rceil$, $\eta:=\zeta^m$. We notice that
\[
\eta\in C_c^\infty(B_{\tilde\varrho}), \qquad 0\leq\eta\leq1, \qquad \eta\equiv1\quad\text{on }B_{\varrho_1}, \qquad \|D\eta\|_{L^\infty} \leq \frac{c(p,q)}{\varrho_2-\varrho_1}.
\]
Using the Sobolev inequality \eqref{Sob.Poinc.H} of Proposition~\ref{prop:intrinsic.sobolev.poincare}, with $d_1=d_1(n,p,q,\alpha)>1$, we infer, for every level $k\geq0$, that
\begin{equation}\label{parz.DeG}
{\biggl( \int_{B_{\varrho_1}} \bigl[ \tilde H(x,(\tilde v\mp k)_\pm) \bigr]^{d_1}\dx \biggr)}^{\frac1{d_\smallone}} \leq c \int_{B_{\tilde\varrho}} \tilde H\big( x,D[(\tilde v\mp k)_\pm\eta] \big)\dx.
\end{equation}
Let us be more precise on the form the constant in the display above can be taken: in view of \eqref{a.rescal}, it can be taken of the form 
\[
c(n,p,q,\alpha) \Bigl(1+[\tilde a]_{C^{0,\alpha}(B_1)}\bigl[\benergy(D[(\tilde v\mp k)_\pm\eta],B_1)\bigr]^{q-p}\Bigr)\leq c(\data)\Bigl(1+\rho^\alpha\bigl\|D[(\tilde v\mp k)_\pm\eta]\bigr\|_{L^{b(p-1)}(B_1)}^{q-p}\Bigr).
\]
Now, observing that $(\tilde v\mp k)_\pm\leq \tilde v_\pm\leq |\tilde v|$, we have, employing Poincar\'e's inequality
\begin{align*}
{\big\|D[(\tilde v\mp k)_\pm\eta]\big\|}_{L^{b(p-1)}(B_1)}^{b(p-1)}&=\int_{B_1}{\big| D[(\tilde v\mp k)_\pm\eta]\big|}^{b(p-1)}\dx\\
&\leq \frac{c(n,q)}{(\varrho_2-\varrho_1)^{b(p-1)}}\biggl[\mean{B_\rho(\bar x)}{|Dv|}^{b(p-1)}\dx+\mean{B_\rho(\bar x)}{\bigg|\frac{v-(v)_{B_\rho}}{\rho}\bigg|}^{b(p-1)}\dx\biggr]\\
&\leq \frac{c(\data)E^{b(p-1)}\rho^{-n}}{(\varrho_2-\varrho_1)^{b(p-1)}}.
\end{align*}
Summarizing, the computations above show us that the constant in \eqref{parz.DeG} has the form
\[
\frac{c(\data,E)}{(\varrho_2-\varrho_1)^{q-p}}\Bigl(1+\rho^{\alpha-n\frac{q-p}{b(p-1)}}\Bigr)\leq \frac{c(\data,E)}{(\varrho_2-\varrho_1)^{q-p}}
\]
as, due to \eqref{bbi}, the exponent of $\rho\leq 1$ is positive.  Applying the product Caccioppoli inequality \eqref{Caccio.product.estimate} in $B_{\tilde\varrho}$ yields
\[
\int_{B_{\tilde\varrho}} \tilde H\big(x,D({(\tilde v\mp k)}_\pm\eta)\big)\dx \leq c(\data) \int_{B_{\tilde\varrho}} \tilde H\big(x,{(\tilde v\mp k)}_\pm|D\zeta|\big)\dx \leq \frac{c(\data)}{(\varrho_2-\varrho_1)^q} \int_{B_{\varrho_2}} \tilde H(x,{(\tilde v\mp k)}_\pm)\dx;
\]
combining this estimate with \eqref{parz.DeG} and with the preceding bound for its constant, we obtain
\begin{equation}\label{almost.final.sup}
\biggl(\int_{B_{\varrho_\smallone}} \bigl[\tilde H\bigl(x,(\tilde v\mp k)_\pm\bigr)\bigr]^{d_1}\dx\biggr)^{\frac1{d_\smallone}}\leq \frac{c(\data,E)}{{(\varrho_2-\varrho_1)}^{2q-p}}\int_{B_{\varrho_\smalltwo}} \tilde H\bigl(x,(\tilde v\mp k)_\pm\bigr)\dx.
\end{equation}
The rest of the proof is now quite standard, but we must take some care: we fix, for $\texttt d>0$ to be chosen, 
\[
k_j=2\texttt d\Bigl(1-\frac1{2^{j+1}}\Bigr)\nearrow 2\texttt d,\qquad R_j:=r+\frac{s-r}{2^j}\searrow r,\qquad\qquad j\in\N_0,
\]
we set $B_j=B_{R_j}$ and we choose in \eqref{almost.final.sup} above, for $j\in\N_0$ fixed, $k=k_j$, $\varrho_1=R_{j+1}$, $\varrho_2=R_j$ and we distinguish as usual two cases.

\medskip 

\subsubsection*{The $p$-case}

Here we suppose that 
\begin{equation}\label{Linfty.p.phase}
\sup_{B_1}\tilde a(\,\cdot\,)\leq 3[\tilde a]_{0,\alpha} 
\end{equation}
holds and we prove the following (inhomogeneous) estimate:
\begin{equation}\label{d.explicit.pre}
{\|\tilde v_\pm\|}_{L^\infty(B_r)}\leq \frac{c(\data,E)}{(s-r)^{{\textnormal{\texttt{E}}}_\smallone}}\biggl(\mean{B_s}\tilde H(x,\tilde v_\pm)\dx\biggr)^{\frac1p},
\end{equation}
for a positive exponent ${\textnormal{\texttt{E}}}_1=(2q-p)d_1/[p(d_1-1)]$ if $p\leq n$, ${\textnormal{\texttt{E}}}_1=q/p$ if $p>n$. In the case $p\leq n$, we define
\[
\Psi^\pm_j=\frac1{\texttt d^p}\int_{B_j}\tilde H\bigl(x,{(\tilde v\mp k_j)}_\pm\bigr)\dx
\]
and we estimate from below as follows:
\begin{align*}
\int_{B_{j+1}} \bigl[\tilde H\bigl(x,{(\tilde v\mp k_j)}_\pm\bigr)\bigr]^{d_1}\dx &\geq \int_{B_{j+1}\cap\{{(\tilde v\mp k_{j+1})}_\pm> 0\}}\tilde H\bigl(x,{(\tilde v\mp k_j)}_\pm\bigr)\bigl[\tilde H\bigl(x,k_{j+1}- k_j\bigr)\bigr]^{d_1-1}\dx\\
&\geq (k_{j+1}- k_j)^{p(d_1-1)} \int_{B_{j+1}}\tilde H\bigl(x,{(\tilde v\mp k_{j+1})}_\pm\bigr)\dx\\
&= 2^{-(j+1)p(d_1-1)} \texttt d^{pd_1}\Psi^\pm_{j+1},
\end{align*}
whence
\[
2^{-(j+1)p(d_1-1)} \texttt d^{pd_1}\Psi^\pm_{j+1}\leq c(\data,E)\biggl[\frac{2^{(j+1)(2q-p)}}{(s-r)^{2q-p}}\texttt d^p\Psi^\pm_{j}\biggr]^{d_1}\qquad\text{for all $j\in\N_0$,}
\]
that is
\[
 \Psi^\pm_{j+1}\leq c(\data,E)\frac{2^{j {\textnormal{\texttt{E}}}_p}}{(s-r)^{(2q-p)d_1}}{\bigl[\Psi^\pm_j\bigr]}^{d_1} 
\]
with $d_1>1$ and ${\textnormal{\texttt{E}}}_p=(2q-p)d_1+p(d_1-1)>0$; this implies, using a standard iteration procedure (cf. \cite[Lemma 7.1]{Giusti}),
\begin{equation}\label{est.sup.implicit}
\frac{1}{\texttt d^p}\int_{B_r}\tilde H(x,{(\tilde v\mp 2 \texttt d)}_\pm\bigr)\dx \leq \Psi^\pm_j\to 0\quad\text{that is}\quad \sup_{B_r} \tilde v_\pm\leq 2 \texttt d 
\end{equation}
provided that
 \[
 \Psi^\pm_0=\frac1{\texttt d^p}\int_{B_s}\tilde H\bigl(x,{(\tilde v\mp \texttt d)}_\pm\bigr)\dx\leq \varepsilon_0(\data,E)(s-r)^{{\textnormal{\texttt{E}}}_1p},
 \]
with $\varepsilon_0(\data,E)$ sufficiently small,  and this can be guaranteed by choosing
\begin{equation}\label{d.explicit}
\texttt d=\frac{c(\data,E)}{(s-r)^{\textnormal{\texttt{E}}_1}}\biggl(\mean{B_s}\tilde H(x,\tilde v_\pm)\dx\biggr)^{\frac1p}
\end{equation}
for $c(\data,E)$ sufficiently large. Now that we know that $\tilde v$ is bounded in every intermediate ball between $B_{1/2}$ and $B_1$, we can interpolate the powers; we consider radii $1/2\leq r<s\leq 3/4$ and we observe that, as $\tilde v_\pm\leq |\tilde v|$,
\begin{align*}
\biggl(\mean{B_s}\tilde v_\pm\dx\biggr)^\frac1p&\leq c(n,p,q) \biggl(\mean{B_s}\tilde v_\pm\dx\biggr)^{\frac{p-q+1}p} \biggl(\mean{B_1}|\tilde v|\dx\biggr)^{\frac{q-p}p} \\
&= c\,\biggl(\mean{B_s}\tilde v_\pm\dx\biggr)^{\frac{p-q+1}p} \biggl(\mean{B_\rho(\bar x)}\Big|\frac{v-(v)_{B_\rho(\bar x)}}\rho\Big|\dx\biggr)^{\frac{q-p}p} \\
&\leq c \biggl(\mean{B_s}\tilde v_\pm\dx\biggr)^{\frac{p-q+1}p} \biggl(\mean{B_\rho(\bar x)}|Dv|\dx\biggr)^{\frac{q-p}p} \\
&\leq c \biggl(\mean{B_s}\tilde v_\pm\dx\biggr)^{\frac{p-q+1}p} \biggl(\mean{B_\rho(\bar x)}{|Dv|}^{b(p-1)}\dx\biggr)^{\frac{q-p}{pb(p-1)}} \\
&\leq c(n,p,q,E)\rho^{-n{\frac{q-p}{pb(p-1)}} } \biggl(\mean{B_s}\tilde v_\pm\dx\biggr)^{\frac{p-q+1}p},
\end{align*}
employing \eqref{bound.energy.mdata}. This gives, using twice Young's inequality with conjugate exponents $(p,p')$ and $(p/(q-1),p/(p-q+1))$ (note that $p/(q-1)>1$, see \eqref{bound.lasco})
\begin{align}\label{reabsorb.Linfty.pphase}
{\|\tilde v_\pm\|}_{L^\infty(B_r)}&\leq  \frac{c}{(s-r)^{{\textnormal{\texttt{E}}}_1}}\biggl(\mean{B_s}\Big[\tilde v_\pm^p+\tilde a(x)\tilde v_\pm^q\Big]\dx\biggr)^\frac1p\notag\\
&\leq \frac{c}{(s-r)^{{\textnormal{\texttt{E}}}_1}}\|\tilde v_\pm\|_{L^\infty(B_s)}^{\frac1{p'}}\biggl(\mean{B_s}\tilde v_\pm\dx\biggr)^\frac1p+ \frac{c\,\rho^\frac\alpha p}{(s-r)^{{\textnormal{\texttt{E}}}_1}}\|\tilde v_\pm\|_{L^\infty(B_s)}^{\frac{q-1}p}\biggl(\mean{B_s}\tilde v_\pm\dx\biggr)^\frac1p\notag\\
&\leq \frac14\|\tilde v_\pm\|_{L^\infty(B_s)}+\frac{c}{(s-r)^{{\textnormal{\texttt{E}}}_1p}}\mean{B_s}|\tilde v|\dx+\frac{c}{(s-r)^{{\textnormal{\texttt{E}}}_1}}\rho^{\frac1p[\alpha-n{\frac{q-p}{\mathfrak s}}] }\|\tilde v_\pm\|_{L^\infty(B_s)}^{\frac{q-1}p} \biggl(\mean{B_s}|\tilde v|\dx\biggr)^{\frac{p-q+1}p}\notag\\
&\leq \frac12\|\tilde v_\pm\|_{L^\infty(B_s)}+\frac{c}{(s-r)^{\frac{{\textnormal{\texttt{E}}}_\smallone p}{p-q+1}}}\mean{B_1}|\tilde v|\dx
\end{align}
with $c\equiv c(\data,E)$; we uniformly bounded the power of the radius using \eqref{rsmaller1} and we also used that $s\geq1/2$. As all the quantities in play are now finite, the standard iteration Lemma~\ref{Giusti.fundamental} applied to the function $\phi(t)=\|\tilde v_\pm\|_{L^\infty(B_t)}$ yields
\[
{\|\tilde v_\pm\|}_{L^\infty(B_{1/2})} \leq c(\data,E) \mean{B_1}|\tilde v|\dx\leq c(\data,E) \mean{B_1}|D\tilde v|\dx
\]
since $(\tilde v)_{B_1}=0$ and we can use the $W^{1,1}$-Poincaré inequality. Recalling the definition of $\tilde v$, this proves \eqref{sup.local} when \eqref{Linfty.p.phase} is in force and $p\leq n$. When $p>n$, estimate \eqref{d.explicit.pre} follows from the Morrey-Sobolev embedding combined with the product Caccioppoli inequality \eqref{Caccio.product.estimate}. Indeed, set $\tilde\varrho:=(r+s)/2$ and let $\zeta\in C_c^\infty(B_{\tilde\varrho})$, $0\leq\zeta\leq1$, satisfy $\zeta\equiv1$ on $B_r$, $\|D\zeta\|_{L^\infty}\leq c/(s-r)$; set moreover $m:=\left\lceil  q/p\right\rceil$, $\eta:=\zeta^m$. Then, by the Morrey-Sobolev embedding,
\begin{equation}\label{pphase.pleqn}
{\|\tilde v_\pm\|}_{L^\infty(B_r)} \leq {\|\tilde v_\pm\eta\|}_{L^\infty(B_{\tilde\varrho})} \leq c(n,p) \biggl( \int_{B_{\tilde\varrho}} |D(\tilde v_\pm\eta)|^p\dx \biggr)^{\frac1p} \leq c(n,p) \biggl( \int_{B_{\tilde\varrho}} \tilde H\bigl(x,D(\tilde v_\pm\eta)\bigr)\dx \biggr)^{\frac1p}.
\end{equation}
We estimate the term on the right-hand side by applying \eqref{Caccio.product.estimate} to the rescaled solution $\tilde v$, with $k=0$: since $s-r\leq1$ and $\|D\zeta\|_{L^\infty}\leq c/(s-r)$, we conclude that
\begin{equation}\label{local.p>n}
{\|\tilde v_\pm\|}_{L^\infty(B_r)} \leq c(\data) \biggl( \int_{B_{\tilde\varrho}} \tilde H\bigl(x,\tilde v_\pm|D\zeta|\bigr)\dx \biggr)^{\frac1p} \leq \frac{c(\data)}{(s-r)^{q/p}} \biggl( \int_{B_s} \tilde H(x,\tilde v_\pm)\dx \biggr)^{\frac1p}.
\end{equation}
In this way we obtain the localized version of \eqref{d.explicit.pre}, with a constant here independent of $E$. Using \eqref{Linfty.p.phase} we have
\begin{equation}\label{bra.bound}
\mean{B_s}\widetilde H(x,\widetilde v_\pm)\dx \leq {\|\tilde v_\pm\|}_{L^\infty(B_s)}^{p-1}\mean{B_s}\tilde v_\pm\dx + c(\data)\rho^\alpha {\|\tilde v_\pm\|}_{L^\infty(B_s)}^{q-1}\mean{B_s}\tilde v_\pm\dx . 
\end{equation}
We use again Sobolev-Morrey embedding to estimate part of the $L^\infty$ norm on the right-hand side:
\begin{align*}
\rho^\alpha {\|\tilde v_\pm\|}_{L^\infty(B_s)}^{q-p}\leq \rho^\alpha {\|\tilde v\|}_{L^\infty(B_1)}^{q-p} &\leq c(\data) \rho^\alpha\biggl(\mean{B_1} {|D\tilde v|}^p\dx\biggr)^\frac{q-p}p\\
&\leq c(\data) \rho^\alpha\biggl(\mean{B_\rho(\bar x)} {|Dv|}^p\dx\biggr)^\frac{q-p}p\\
&\leq c(\data)\rho^{\alpha-n\frac{q-p}p}{\|Dv\|}_{L^p(B_R)}^{q-p}\leq c(\data,E),
\end{align*}
using again \eqref{rsmaller1} where here $\mathfrak s=p$ as $p>n$. Consequently, inserting this estimate in \eqref{bra.bound} and the latter, in turn, into \eqref{d.explicit.pre}, using again $s\geq1/2$ to enlarge the integration ball, we obtain
\[
{\|\tilde v_\pm\|}_{L^\infty(B_r)}  \leq \frac{c(\data,E)}{(s-r)^{q/p}} {\|\tilde v_\pm\|}_{L^\infty(B_s)}^\frac1{p'}\biggl( \mean{B_1}\tilde v_\pm\dx  \biggr)^{\frac1p}\leq \frac12{\|\tilde v_\pm\|}_{L^\infty(B_s)} + \frac{c(\data,E)}{(s-r)^q}\mean{B_1}|\tilde v|\dx,
\]
using Young's inequality with conjugate exponents $p'$ and $p$. Now we can reabsorb the $L^\infty$ norm on the right-hand side and use Poincar\'e's inequality exactly as in the case $p\leq n$, and this proves \eqref{sup.local} when \eqref{Linfty.p.phase} holds.

\medskip

\subsubsection*{The $(p,q)$-case}

We now consider the complementary condition to \eqref{Linfty.p.phase}, namely, we assume that
\begin{equation}\label{Linfty.pq.phase}
\sup_{B_1}\tilde a(\,\cdot\,)> 3[\tilde a]_{0,\alpha}.
\end{equation}
Then $\tilde a(x)\approx \tilde a(0)$ for every $x\in B_1$, as in \eqref{above.below}. As before, we first consider the case $p\leq n$; here the iteration must be performed in a slightly different form. We set $\tilde H_0(\xi):=\tilde H(0,\xi)$, define
\[
\Psi^\pm_j=\frac1{\tilde H_0(\texttt d)}\int_{B_j}\tilde H\bigl(x,{(\tilde v\mp k_j)}_\pm\bigr)\dx,
\]
and proceed as follows:
\begin{align*}
\int_{B_{j+1}} \bigl[\tilde H\bigl(x,{(\tilde v\mp k_j)}_\pm\bigr)\bigr]^{d_1}\dx &\geq \frac1c\int_{B_{j+1}\cap\{(\tilde v\mp k_{j+1})_\pm> 0\}} \tilde H\bigl(x,{(\tilde v\mp k_j)}_\pm\bigr) \bigl[\tilde H\bigl(0,k_{j+1}- k_j\bigr)\bigr]^{d_1-1}\dx\\
& \geq \frac1c\bigl[\tilde H_0(2^{-(j+1)}\texttt d)\bigr]^{d_1-1} \int_{B_{j+1}}\tilde H\bigl(x, {(\tilde v\mp k_{j+1})}_\pm\bigr)\dx \\
&\geq \frac1c\,2^{-(j+1)q(d_1-1)} \bigl[\tilde H_0(\texttt d)\bigr]^{d_1}\Psi^\pm_{j+1}.
\end{align*}
Hence
\[
2^{-(j+1)q(d_1-1)}\bigl[\tilde H_0(\texttt d)\bigr]^{d_1}\Psi^\pm_{j+1} \leq c(\data,E)\biggl[\frac{2^{j(2q-p)}}{(s-r)^{2q-p}}\tilde H_0(\texttt d)\Psi^\pm_j\biggr]^{d_1}
\]
for every $j\in\N_0$, that is,
\[
\Psi^\pm_{j+1} \leq c(\data,E)\frac{2^{j{\textnormal{\texttt{E}}}_q}}{(s-r)^{(2q-p)d_1}} \bigl[\Psi^\pm_j\bigr]^{d_1},
\]
with ${\textnormal{\texttt{E}}}_q=(2q-p)d_1+q(d_1-1)>0$. Since $d_1>1$, this again implies that
\[
\sup_{B_r} \tilde v_\pm\leq 2\texttt d,
\]
provided that
\[
\frac1{\tilde H_0(\texttt d)} \int_{B_s}\tilde H\bigl(x, (\tilde v\mp \texttt d)_\pm\bigr)\dx \leq \varepsilon_0(\data,E)(s-r)^{\frac{(2q-p)d_\smallone}{d_\smallone-1}},
\]
with $\varepsilon_0(\data,E)$ sufficiently small. This can be ensured by choosing
\[
\texttt d=\tilde H_0^{-1}\biggl(\frac{c(\data,E)}{(s-r)^{(2q-p)d'_\smallone}}\mean{B_s}\tilde H(x,\tilde v_\pm)\dx \biggr),
\]
with $c(\data,E)$ large. Therefore,
\[
\sup_{B_r}\tilde H(0,\tilde v_\pm) \leq \frac{c(\data,E)}{(s-r)^{(2q-p)d'_\smallone}}\mean{B_s}\tilde H(x,\tilde v_\pm)\dx \leq \frac{c(\data,E)}{(s-r)^{(2q-p)d'_\smallone}}\mean{B_s}\tilde H(0,\tilde v_\pm)\dx,
\]
and in particular $\tilde v$ is bounded in any $B_r$ for $r\in (0,1)$. Now we restrict ourselves to radii $1/2\leq r<s\leq 3/4$ and we suppose $\|\tilde v_\pm\|_{L^\infty(B_s)}\neq 0$, otherwise what we are going to prove is trivial, and employing Young's inequality for Young functions - we here denote by $\hat H_0$ the convex conjugate of $\tilde H_0$ - we have, using the monotonicity of $t\mapsto\widetilde H_0(t)/t$, where we set $\tilde H_0(t)/t=0$ if $t=0$,
\begin{align}\label{reabsorb.Linfty.pqphase}
\tilde H_0\bigl(\|\tilde v_\pm\|_{L^\infty(B_r)}\bigr)&\leq \frac c{(s-r)^{(2q-p)d'_\smallone}}\mean{B_s}\frac{\tilde H_0(\tilde v_\pm)}{\tilde v_\pm}\tilde v_\pm\dx\notag\\
&\leq \frac c{(s-r)^{(2q-p)d'_\smallone}}\,\frac{\tilde H_0({\|\tilde v_\pm\|}_{L^\infty(B_s)})}{{\|\tilde v_\pm\|}_{L^\infty(B_s)}}\mean{B_s}\tilde v_\pm\dx\notag\\
&\leq \frac12\hat H_0\biggl(\frac{\tilde H_0({\|\tilde v_\pm\|}_{L^\infty(B_s)})}{{\|\tilde v_\pm\|}_{L^\infty(B_s)}}\biggr)+c\,\tilde H_0\biggl(\frac1{(s-r)^{(2q-p)d'_\smallone}}\mean{B_s}\tilde v_\pm\dx\biggr)\notag\\
&\leq \frac12\tilde H_0\bigl({\|\tilde v_\pm\|}_{L^\infty(B_s)}\bigr)+\frac c{(s-r)^{(2q-p)d'_\smallone q}}\,\tilde H_0\biggl(\mean{B_1}\tilde v_\pm\dx\biggr),
\end{align}
by the standard property of convex conjugates, see \eqref{Young-N} and the last inequality in \eqref{conj.inverse}; again this yields, using the scaling properties of $\tilde H_0$ and Poincar\'e's inequality in $W^{1,1}(B_1)$,
\[
\tilde H_0\bigl(\|\tilde v_\pm\|_{L^\infty(B_{1/2})}\bigr)\leq c\,\tilde H_0\biggl(\mean{B_1}\tilde v_\pm\dx\biggr)\leq \tilde H_0\biggl(c\mean{B_1}|\tilde v|\dx\biggr)\leq \tilde H_0\biggl(c\mean{B_1}|D\tilde v|\dx\biggr)
\]
and the claim follows when \eqref{Linfty.pq.phase} is in force and $p\leq n$, since the map $\xi\mapsto \tilde H_0(\xi)$ is strictly increasing. For the range $p>n$, we observe that the map $t\mapsto \tilde H_0(t^{1/p})=t+\tilde  a(0)t^{q/p}$ is convex. Hence, by Jensen's inequality, the Sobolev-Morrey embedding and the scaling properties of $\tilde H_0$, using the same notation introduced before \eqref{pphase.pleqn},
\[
\tilde  H_0\bigl({\|\tilde  v_\pm\|}_{L^\infty(B_r)} \bigr) \leq \tilde  H_0\bigl({\|\tilde  v_\pm\eta\|}_{L^\infty(B_{\tilde \varrho})} \bigr) \leq \tilde  H_0\biggl( c(n,p) \biggl( \mean{B_{\tilde \varrho}}{\big|D(\tilde v_\pm\eta)\big|}^p\dx \biggr)^{\frac1p} \biggr) \leq c\mean{B_{\tilde \varrho}} \tilde  H_0\bigl(D(\tilde v_\pm\eta)\bigr)\dx,
\]
with $c\equiv c(n,p,q)$. By the comparability of $\tilde  H(x,\,\cdot\,)$ and $\tilde  H_0$ in the $(p,q)$-phase, and applying the product Caccioppoli inequality \eqref{Caccio.product.estimate}, we find
\[
\mean{B_{\tilde \varrho}} \tilde  H_0\bigl(D(\eta\tilde  v_\pm)\bigr)\dx \leq c \mean{B_{\tilde \varrho}} \tilde  H\bigl(x,D(\eta\tilde  v_\pm)\bigr)\dx \leq c \mean{B_{\tilde \varrho}} \tilde  H\bigl(x,\tilde  v_\pm|D\zeta|\bigr)\dx \leq \frac{c}{(s-r)^q} \mean{B_s}\tilde  H_0(\tilde  v_\pm)\dx,
\]
$c\equiv c(\data)$; consequently,
\begin{equation}\label{local.sup.pq.supercritical}
\tilde  H_0\bigl( {\|\tilde  v_\pm\|}_{L^\infty(B_r)} \bigr) \leq \frac{c(\data)}{(s-r)^q} \mean{B_s}\tilde  H_0(\tilde  v_\pm)\dx.
\end{equation}
We now perform the absorption: we have, very similarly to \eqref{reabsorb.Linfty.pqphase},
\begin{align*}
\tilde  H_0\bigl( {\|\tilde  v_\pm\|}_{L^\infty(B_r)} \bigr)&\leq\frac{c(\data)}{(s-r)^q} \frac{ \tilde  H_0\bigl( {\|\tilde  v_\pm\|}_{L^\infty(B_s)} \bigr) }{ {\|\tilde  v_\pm\|}_{L^\infty(B_s)} } \mean{B_s}\tilde  v_\pm\dx\\
&\leq \frac12 \tilde  H_0\bigl( {\|\tilde  v_\pm\|}_{L^\infty(B_s)} \bigr) + \frac{c(\data)}{(s-r)^{q^2}} \tilde  H_0\biggl( \mean{B_1}|\tilde  v|\dx \biggr).
\end{align*}
Iteration Lemma~\ref{Giusti.fundamental} therefore gives
\[
\tilde  H_0\bigl( {\|\tilde  v_\pm\|}_{L^\infty(B_{1/2})} \bigr) \leq c(\data) \tilde  H_0\biggl( \mean{B_1}|\tilde  v|\dx \biggr)\leq  \tilde  H_0\biggl( c(\data)\mean{B_1}|\tilde  v|\dx \biggr)
\]
and we conclude since $\tilde  H_0$ is strictly increasing on $[0,\infty)$, and using again the $W^{1,1}$-Poincaré inequality.
\end{proof}

\begin{corollary}\label{corol.supBMO}
Under the same assumptions as in Theorem~\ref{sup.correct.exponent}, there holds
\begin{equation}\label{sup.local.mean}
\frac{\osc_{B_{\rho/2}(\bar x)}v}{\rho/2}\leq4\sup_{B_{\rho/2}(\bar x)}\biggl|\frac{v-(v)_{B_\rho(\bar x)}}{\rho}\biggr|\leq c(\data,E)\mean{B_\rho(\bar x)}|Dv|\dx.
\end{equation}
\end{corollary}

\begin{proof}
The first inequality follows from the elementary estimate $\osc_Bv\leq2\sup_B|v-k|$, applied with $k=(v)_{B_\rho(\bar x)}$. The second inequality follows by applying \eqref{sup.local} to the positive and negative parts.
\end{proof}

The proof of Theorem~\ref{sup.correct.exponent}, before the final application of Poincar\'e's inequality, provides an $L^1$-$L^\infty$ estimate for the oscillation of homogeneous solutions. Combining this estimate with \eqref{Caccio2}, we obtain the following reverse Poincar\'e inequality that will be useful in the proof of Theorem~\ref{main.thm.2}.

\begin{corollary}[$L^1$-Caccioppoli inequality]\label{Caccio.L1}
Let $v\in W^{1,H(\cdot)}(B_R)$ be a solution to \eqref{homogeneous.equation}, with $R\leq1$, and suppose that \eqref{bound.energy.mdata} holds. Then, for every pair of concentric balls $B_r(\bar x)\subset B_s(\bar x)\subset B_R$, with $0<r<s$, and every $k\in\R$, there holds
\begin{equation}\label{Caccio.L1.estimate}
\int_{B_r(\bar x)}|Dv|\dx\leq\frac{c}{s-r}\int_{B_s(\bar x)}|v-k|\dx,
\end{equation}
where $c$ depends only on $\data$ and $E$.
\end{corollary}

\begin{proof}
Let $B_\rho(y)\subset B_R$ be an arbitrary ball. From the proof of Theorem~\ref{sup.correct.exponent} (before the final application of Poincar\'e's inequality), we already know that
\begin{equation}\label{sup.L1.preliminary}
\sup_{B_{\rho/2}(y)}\bigl|v-(v)_{B_\rho(y)}\bigr|\leq c\mean{B_\rho(y)}|v-(v)_{B_\rho(y)}|\dx.
\end{equation}
We use \eqref{Caccio2} with a cut-off function supported in $B_{\rho/2}(y)$ such that $\eta\equiv 1$ on $B_{\rho/4}(y)$ and $|D\eta|\leq c/\rho$; combining the resulting bound with \eqref{sup.L1.preliminary} gives
\begin{equation}\label{Caccio.L1.energy}
\mean{B_{\rho/4}(y)}H(x,Dv)\dx\leq c\biggl[{\biggl(\frac1\rho\mean{B_\rho(y)}\bigl|v-(v)_{B_\rho(y)}\bigr|\dx\biggr)}^p+\sup_{B_\rho(y)}a(\,\cdot\,){\biggl(\frac1\rho\mean{B_\rho(y)}\bigl|v-(v)_{B_\rho(y)}\bigr|\dx\biggr)}^q\biggr].
\end{equation}
We now separate the two phases in $B_\rho(y)$. If the $p$-phase condition \eqref{p.phase.sup.local} holds with $M=3$, then Poincar\'e's and H\"older's inequalities together with \eqref{bound.energy.mdata} imply
\[
\frac1\rho\mean{B_\rho(y)}\bigl|v-(v)_{B_\rho(y)}\bigr|\dx\leq c\mean{B_\rho(y)}|Dv|\dx\leq c\rho^{-\frac n{\mathfrak s}}E,
\]
where $\mathfrak s$ is as in \eqref{esse}. By \eqref{condition.H.Sobolev}, it follows that
\[
\sup_{B_\rho(y)}a(\,\cdot\,){\biggl(\frac1\rho\mean{B_\rho(y)}\bigl|v-(v)_{B_\rho(y)}\bigr|\dx\biggr)}^{q-p}\leq c[a]_{0,\alpha}E^{q-p}\rho^{\alpha-\frac{n(q-p)}{\mathfrak s}}\leq c(\data,E);
\]
plugging this bound into \eqref{Caccio.L1.energy} and applying H\"older's inequality leads to
\begin{equation}\label{Caccio.L1.oneball}
\mean{B_{\rho/4}(y)}|Dv|\dx\leq \frac c\rho\mean{B_\rho(y)}|v-(v)_{B_\rho(y)}|\dx.
\end{equation}
In the $(p,q)$-phase, \eqref{above.below} gives $\sup_{B_\rho(y)}a(\,\cdot\,)\leq c\inf_{B_\rho(y)}a(\,\cdot\,)$; setting $H_\rho(t):=t^p+\inf_{B_\rho(y)}a(\,\cdot\,) t^q$, Jensen's inequality and \eqref{Caccio.L1.energy} yield
\[
H_\rho\biggl(\mean{B_{\rho/4}(y)}|Dv|\dx\biggr)\leq\mean{B_{\rho/4}(y)}H_\rho(|Dv|)\dx\leq cH_\rho\biggl(\frac1\rho\mean{B_\rho(y)}\bigl|v-(v)_{B_\rho(y)}\bigr|\dx\biggr).
\]
Inequality \eqref{Caccio.L1.oneball} follows in this case from the scaling properties of $H_\rho(\,\cdot\,)$. In both phases, the elementary estimate \eqref{exc.elementary} for the excess leads to
\begin{equation}\label{caccioL1_before.cover}
\int_{B_{\rho/4}(y)}|Dv|\dx\leq\frac c\rho\int_{B_\rho(y)}|v-k|\dx
\end{equation}
for every $k\in\R$, with $c=c(\data,E)$.

\medskip

To complete the proof, take $B_r(\bar x)\subset B_s(\bar x)\subset B_R$ and set $\rho:=(s-r)/4$. We cover $B_r(\bar x)$ with balls $B_{\rho/4}(y_j)$, $y_j\in B_r(\bar x)$, in such a way that the dilated balls $B_\rho(y_j)$ have finite overlap bounded only in terms of $n$. Since $B_\rho(y_j)\subset B_s(\bar x)$, summing \eqref{caccioL1_before.cover} over $j$ gives \eqref{Caccio.L1.estimate}.
\end{proof}

We now record some useful estimates adapted to the different phases.

\subsection{Estimates in the \texorpdfstring{$\boldsymbol{p}$}{p}-phase}

Let $B_\rho\equiv B_\rho(\bar x)\subset B_R$ be a ball, not necessarily concentric with $B_R$. We say that the $p$-phase condition holds in $B_\rho$ if
\begin{equation}\label{p.phase.sup.local}
\sup_{x\in B_\rho}a(x)\leq M{[a]}_{0,\alpha}\rho^\alpha
\end{equation}
for some constant $M\geq0$. The following result provides a quantitative reverse H\"older inequality in the $p$-phase, with the dependence of the constant required in the study of double phase problems with measure data.

\begin{proposition}[A reverse H\"older inequality in the $p$-phase]\label{cor:nice.reverse}
Let $v\in W^{1,H(\cdot)}(B_R)$ be a solution to \eqref{homogeneous.equation} and suppose that the energy bound \eqref{bound.energy.mdata} holds;  let $B_\rho\subset B_R $ be a ball with radius $\rho \leq 1$ and such that the inequality in \eqref{p.phase.sup.local} is satisfied for some $M \geq 0$. Then, for every $\theta \in [1/2,3/4]$, we have
\begin{equation}\label{nice.reverse.2}
\mean{B_{\theta\rho}} H(x,Dv)\dx \leq c(\data,E,M)\biggl(\mean{B_\rho} |Dv|^{b(p-1)}\dx \biggr)^{\frac{p}{b(p-1)}}
\end{equation}
if $p\leq n$, for $b$ as in \eqref{def.b}-\eqref{bbi} and
\begin{equation}\label{nice.reverse.3}
\mean{B_{\theta\rho}} H(x,Dv)\dx \leq c(\data,E,M)\biggl(\mean{B_\rho} |Dv|^{p-1}\dx \biggr)^{p'}
\end{equation}
if $p>n$.
\end{proposition}

\begin{proof}
We prove \eqref{nice.reverse.2} for $\theta=3/4$ only; the estimates for $\theta\in [1/2,3/4)$ follow simply by enlarging the integration ball. 

\subsubsection*{Case $p\leq n$} For $b$ as in the statement, we choose in \eqref{rev.self.improved}
\[
\sigma=\frac{b(p-1)}q\leq \frac{b(p-1)}p
\]
and this gives for $c$ depending on $\data$ and $E$
\begin{equation*}
\mean{B_{3\rho/4}} H(x,Dv)\dx\leq c{\biggl(\mean{B_\rho}|Dv|^{p\sigma}\dx\biggr)}^{\frac 1\sigma}+c{\biggl(\mean{B_\rho}\bigl[a(x)|Dv|^q\bigr]^\frac{b(p-1)}{q}\dx\biggr)}^{\frac q{b(p-1)}}.
\end{equation*}
The first term on the right-hand side is estimated simply by applying H\"older's inequality:
\[
{\biggl(\mean{B_\rho}|Dv|^{\frac pq b(p-1)}\dx\biggr)}^{\frac qp\frac p{b(p-1)}}\leq {\biggl(\mean{B_\rho}|Dv|^{b(p-1)}\dx\biggr)}^{\frac p{b(p-1)}}.
\]
For the second, we use the fact that we are in the $p$-phase, that is \eqref{p.phase.sup.local}, \eqref{rsmaller1} and \eqref{bound.energy.mdata}
\begin{align*}
\biggl(\mean{B_\rho}\bigl[a(x)|Dv|^q\bigr]^\frac{b(p-1)}{q}\dx\biggr)^{\frac q{b(p-1)}}&\leq  c\,M {[a]}_{0,\alpha} \rho^\alpha\biggl(\mean{B_\rho}|Dv|^{b(p-1)}\dx\biggr)^{\frac q{b(p-1)}}\notag\\
&\leq  c\,M {[a]}_{0,\alpha} \rho^{\alpha-\frac{n(q-p)}{b(p-1)}}{\|Dv\|}_{L^{b(p-1)}(B_R(x_0))}^{q-p}\biggl(\mean{B_\rho}|Dv|^{b(p-1)}\dx\biggr)^{\frac p{b(p-1)}}\\
&\leq  c(\data,E)M \biggl(\mean{B_\rho}|Dv|^{b(p-1)}\dx\biggr)^{\frac p{b(p-1)}}.
\end{align*}
\subsubsection*{Case $p> n$} Here in \eqref{rev.self.improved} we choose $\sigma=(p-1)/q$: we obtain
\begin{align}\label{partial.reverse.p>n}
\mean{B_{3\rho/4}} H(x,Dv)\dx&\leq c\biggl(\mean{B_\rho}|Dv|^{\frac{p(p-1)}q}\dx\biggr)^{\frac q{p-1}}+c\biggl(\mean{B_\rho}\bigl[a(x)|Dv|^q\bigr]^\frac{p-1}{q}\dx\biggr)^{\frac q{p-1}}\notag\\
&\leq c\biggl(\mean{B_\rho}|Dv|^{p-1}\dx\biggr)^{p'}+c\,\rho^\alpha\biggl(\mean{B_\rho}{|Dv|}^{p-1}\dx\biggr)^{\frac q{p-1}}
\end{align}
for $c\equiv c(\data,E, M)$, since $(p-1)p/q\leq p-1$. Now we split and use H\"older's inequality to obtain
\begin{align*}
\biggl(\mean{B_\rho}{|Dv|}^{p-1}\dx\biggr)^{\frac q{p-1}}&\leq \biggl(\mean{B_\rho}{|Dv|}^p\dx\biggr)^{\frac{q-p}p}\biggl(\mean{B_\rho}{|Dv|}^{p-1}\dx\biggr)^{p'}\\
&\leq c(n,p,q)\rho^{-n\frac{q-p}p}E^{q-p}\biggl(\mean{B_\rho}{|Dv|}^{p-1}\dx\biggr)^{p'}, 
\end{align*}
and the proof is completed by inserting the last inequality into \eqref{partial.reverse.p>n}, thanks to \eqref{rsmaller1}.
\end{proof}

The previous result estimates the natural energy in terms of a lower integrability level of the gradient. We now state a converse estimate, which will be used repeatedly to recover the corresponding lower-order quantities for the energy.

\begin{theorem}[A non-homogeneous H\"older inequality in the $p$-phase]
Let $v\in W^{1,H(\cdot)}(B_R)$ be a solution to \eqref{homogeneous.equation}, $b$ as in \eqref{def.b}-\eqref{bbi} when $p\leq n$, and suppose that the assumptions of Proposition~\ref{cor:nice.reverse} are satisfied. Then we have 
\begin{equation}\label{directnice.reverse.2}
\mean{B_\rho}|Dv|^{b(p-1)}\dx +{\biggl(\mean{B_\rho}a(x)|Dv|^{q-1}\dx\biggr)}^b\leq  c(\data,E,M)\biggl(\mean{B_\rho}H(x,Dv)\dx\biggr)^{\frac{b}{p'}} 
\end{equation}
if $p\leq n$ and
\begin{equation}\label{directnice.reverse.3}
\mean{B_\rho} h(x,Dv)\dx \leq c(\data,E,M)\biggl(\mean{B_\rho} H(x,Dv)\dx \biggr)^{\frac1{p'}}
\end{equation}
if $p>n$.
\end{theorem}

\begin{proof}
We prove the two cases separately.
\subsubsection*{Case $p\leq n$}
Since $b(p-1)<p$, H\"older's inequality and the bound $|Dv|^p\leq p\,H(x,Dv)$ give
\begin{equation}\label{holder.b.easy}
\mean{B_\rho}|Dv|^{b(p-1)}\dx\leq c(n,p,q,\alpha)\biggl(\mean{B_\rho}H(x,Dv)\dx\biggr)^{\frac b{p'}}:
\end{equation}
for the second term, we use H\"older's inequality, \eqref{p.phase.sup.local}, and the bound $q-p\leq b(p-1)$ following from \eqref{def.b}. Recalling \eqref{bound.energy.mdata}, we obtain
\begin{align}\label{recover.est}
\mean{B_\rho}a(x)|Dv|^{q-1}\dx&\leq\biggl(\mean{B_\rho}a(x)|Dv|^{q-p}\dx\biggr)^{\frac1p}\biggl(\mean{B_\rho}a(x)|Dv|^q\dx\biggr)^{\frac1{p'}}\notag\\
&\leq c\rho^{\frac\alpha p}\biggl(\mean{B_\rho}|Dv|^{b(p-1)}\dx\biggr)^{\frac{q-p}{pb(p-1)}}\biggl(\mean{B_\rho}H(x,Dv)\dx\biggr)^{\frac1{p'}}\notag\\
&\leq c\rho^{\frac1p[\alpha-\frac{n(q-p)}{b(p-1)}]}E^{\frac{q-p}p}\biggl(\mean{B_\rho}H(x,Dv)\dx\biggr)^{\frac1{p'}}\notag\\
&\leq c\biggl(\mean{B_\rho}H(x,Dv)\dx\biggr)^{\frac1{p'}},
\end{align}
where the last inequality follows from \eqref{rsmaller1}, and $c$ depends only on $\data,E$ and $M$. Raising this estimate to the power $b$ and combining it with \eqref{holder.b.easy} proves \eqref{directnice.reverse.2}.

\subsubsection*{Case $p>n$} 
The proof in this case is very similar to that in the case $p\leq n$, with some elements replacing others. We use H\"older's inequality as in \eqref{holder.b.easy} for $b=1$ for the $p$-part and we estimate, by $a(x)|Dv|^{q-1}=[a(x)]^{1/p}|Dv|^{q/p-1}[a(x)|Dv|^q]^{1/p'}$, $q-p\leq p$ by \eqref{bound.supercritical}, \eqref{rsmaller1} and \eqref{bound.energy.mdata}
\begin{align}\label{q.part.useful.later}
\mean{B_{\rho}}a(x)|Dv|^{q-1}\dx&\leq c\biggl(\mean{B_{\rho}}a(x)|Dv|^{q-p}\dx\biggr)^{\frac1p}\biggl(\mean{B_{\rho}}a(x)|Dv|^q\dx\biggr)^{\frac1{p'}}\notag\\
&\leq c \,\rho^{\frac \alpha p} \biggl(\mean{B_{\rho}}|Dv|^p\dx\biggr)^{\frac{q-p}{p^\smalltwo}}\biggl(\mean{B_{\rho}}H(x,Dv)\dx\biggr)^{\frac1{p'}}\notag\\
&\leq c \,\rho^{\frac 1 p[\alpha-n \frac{q-p}p]} \biggl(\int_{B_{\rho}}|Dv|^p\dx\biggr)^{\frac{q-p}{p^\smalltwo}}\biggl(\mean{B_{\rho}}H(x,Dv)\dx\biggr)^{\frac1{p'}}\notag\\
&\leq c E^{\frac{q-p}p}\biggl(\mean{B_{\rho}}H(x,Dv)\dx\biggr)^{\frac{1}{p'}},
\end{align}
for the $q$-part, to get \eqref{directnice.reverse.3}.
\end{proof}

As a corollary of \eqref{nice.reverse.2}, we finally get a Caccioppoli inequality in the $p$-phase with the correct dependence of the constant; compare with \cite[Lemma 10.1]{CM1}.
\begin{lemma}[Almost-standard Caccioppoli estimate - Version I]\label{AsCaccioppoli}
Let $v\in W^{1,H(\cdot)}(B_R)$ be a solution to \eqref{homogeneous.equation}, suppose that the energy bound \eqref{bound.energy.mdata} holds and let $B_\rho\equiv B_\rho(\bar x)\subset B_R $ be a ball with radius $\rho \leq 1$ and such that the inequality in \eqref{p.phase.sup.local} holds for some $M \geq 0$. Then
\[
\mean{B_{\rho/4}(\bar x)}H(x,Dv)\dx\leq c\mean{B_{\rho/2}(\bar x)}\bigg|\frac{v-(v)_{B_{\rho/2}(\bar x)}}{\rho}\bigg|^p\dx
\]
with $c$ depending on $\data,E$ and $M$.
\end{lemma}

\begin{proof}

We start from the standard Caccioppoli estimate in \eqref{Caccio2}: taking a standard cut-off function $\eta\in C_c^\infty(B_{\rho/2}(\bar x))$, $0\leq\eta\leq1$, $\eta\equiv1$ on $B_{\rho/4}(\bar x)$ and $\|D\eta\|_{L^\infty} \leq c/\rho$. This leads to 
\[
\mean{B_{\rho/4}(\bar x)}H(x,Dv)\dx\leq c\mean{B_{\rho/2}(\bar x)}H\left( x,\frac{v-(v)_{B_{\rho/2}(\bar x)}}{\rho} \right)\dx
\]
 and we use \eqref{p.phase.sup.local} to pointwise estimate the right-hand side: for $x\in B_{\rho/2}(\bar x)$,
\begin{equation}\label{H.parz.est}
H\biggl(x,\frac{v-{(v)}_{B_{\rho/2}(\bar x)}}{\rho}\biggr)\leq \bigg|\frac{v-{(v)}_{B_{\rho/2}(\bar x)}}{\rho}\bigg|^p\Big(1+ c(p,q) M{[a]}_{0,\alpha}\rho^\alpha\sup_{B_{\rho/2}(\bar x)}{\bigg|\frac{v-{(v)}_{B_{\rho/2}(\bar x)}}{\rho}\bigg|}^{q-p} \Big).
\end{equation}
Using the $\sup$-estimate in \eqref{sup.local.mean} and H\"older's inequality gives
\[
\sup_{B_{\rho/2}(\bar x)}{\bigg|\frac{v-{(v)}_{B_{\rho/2}(\bar x)}}{\rho}\bigg|}^{q-p} \leq c(\data,E){\biggl(\mean{B_\rho(\bar x)}{|Dv|}^{\mathfrak s}\dx\biggr)}^\frac{q-p}{\mathfrak s}.
\]
Inserting this estimate into \eqref{H.parz.est} and using \eqref{bound.energy.mdata}, we obtain
\[
\rho^\alpha\sup_{B_{\rho/2}(\bar x)}{\bigg|\frac{v-{(v)}_{B_{\rho/2}(\bar x)}}{\rho}\bigg|}^{q-p}\leq c(\data,E)\rho^\alpha {\biggl(\mean{B_\rho(\bar x)}{|Dv|}^{\mathfrak s}\dx\biggr)}^\frac{q-p}{\mathfrak s} \leq c(\data,E)\rho^{\alpha -n\frac{q-p}{\mathfrak s}}E^{q-p},
\]
and the proof is concluded thanks to \eqref{rsmaller1}.

\end{proof}

The utility of the next Lemma will become clearer later on.

\begin{lemma}[Almost-standard Caccioppoli estimate - Version II]\label{AsCaccioppoli2}
Under the same assumptions as in Lemma~\ref{AsCaccioppoli}, there holds
\begin{equation}\label{cac.k.correct.dependence}
\mean{B_{\rho/4}(\bar x)}H\bigl(x,D{(v-k)}_\pm\bigr)\dx\leq c\mean{B_{\rho/2}(\bar x)}\frac{{(v-k)}^p_\pm}{\rho^p}\dx 
\end{equation}
for every level 
\[
k\in \Big[\sup_{B_{\rho/2}(\bar x)}v-\frac12\osc_{B_{\rho/2}(\bar x)}v,\sup_{B_{\rho/2}(\bar x)}v \Big]\quad\text{if ``$+$''},\qquad k\in \Big[\inf_{B_{\rho/2}(\bar x)}v,\inf_{B_{\rho/2}(\bar x)}v+\frac12\osc_{B_{\rho/2}(\bar x)}v\Big]\quad\text{if ``$-$''},
\]
with $c$ depending on $\data,E$ and $M$. More generally, for every $0<r<s\leq\rho/2$, we have
\[
\int_{B_r(\bar x)}H\bigl(x,D(v-k)_\pm\bigr)\dx\leq c\,\frac{\rho^{q-p}}{(s-r)^q}\int_{B_s(\bar x)}|(v-k)_\pm|^p\dx,
\]
with the same dependence of the constant.
\end{lemma}

\begin{proof}
If $\osc_{B_{\rho/2}(\bar x)}v=0$, the conclusion is immediate; otherwise, for the prescribed levels we have
\[
0\leq(v-k)_\pm\leq\frac12\osc_{B_{\rho/2}(\bar x)}v.
\]
By \eqref{sup.local.mean}, H\"older's inequality, \eqref{p.phase.sup.local}, \eqref{bound.energy.mdata} and \eqref{rsmaller1}, we obtain
\[
\sup_{B_{\rho/2}(\bar x)}a(\,\cdot\,)\biggl[\frac{\osc_{B_{\rho/2}(\bar x)}v}{\rho}\biggr]^{q-p}\leq c(\data,E,M).
\]
Fix $0<r<s\leq\rho/2$. Applying \eqref{Caccio} with a cut-off supported in $B_s(\bar x)$, equal to one on $B_r(\bar x)$ and satisfying $|D\eta|\leq c/(s-r)$, gives
\begin{align*}
\int_{B_r(\bar x)}H\bigl(x,D(v-k)_\pm\bigr)\dx&\leq c\int_{B_s(\bar x)}\biggl[\frac{|(v-k)_\pm|^p}{(s-r)^p}+a(x)\frac{|(v-k)_\pm|^q}{(s-r)^q}\biggr]\dx\\
&\leq c\frac{\rho^{q-p}}{(s-r)^q}\int_{B_s(\bar x)}|(v-k)_\pm|^p\dx,
\end{align*}
where we also used $s-r\leq\rho$. Taking $r=\rho/4$ and $s=\rho/2$, and dividing by the volume of $B_{\rho/4}$, proves \eqref{cac.k.correct.dependence}.
\end{proof}

\subsection{Estimates in the \texorpdfstring{$\boldsymbol{(p,q)}$}{(p,q)}-phase}

We need the following reverse H\"older inequality in the $(p,q)$-phase: we assume that the complementary $(p,q)$-phase condition holds, namely
\begin{equation}\label{pq.phase.sup.local}
\sup_{x \in B_\rho}\, a(\,\cdot\,) > M{[a]}_{0, \alpha}\rho^\alpha 
\end{equation}
for $M>2$ to be chosen; notice that if \eqref{pq.phase.sup.local} holds, then we have \eqref{above.below} at our disposal. 
\begin{proposition}\label{cor:nice-rev.pq}
Let $v\in W^{1,H(\cdot)}(B_R)$ be a solution to \eqref{homogeneous.equation}; let $E$ be as in \eqref{bound.energy.mdata} and $B_\rho \equiv B_\rho(\bar x)\subset B_R$. Let moreover the inequality in \eqref{pq.phase.sup.local} be satisfied for some number $M>2$.  Then there exists a constant $c$, depending on $\data,E$ and $M$, such that
\begin{equation}\label{nice.reverse.2.pq}
\mean{B_{\rho/2}} H(x,Dv)\dx \leq c\bigl(H_{\bar x}\circ h_{\bar x}^{-1}\bigr)\biggl(\mean{B_\rho} h(x,Dv)\dx \biggr).
\end{equation}
The phase assumption may alternatively be replaced by the existence of a constant $a_*>0$ such that
\[
a_*\leq a(x)\leq2a_*\qquad\text{for every }x\in B_\rho.
\]
Setting
\[
H_*(t):=\frac{t^p}{p}+a_*\frac{t^q}{q},\qquad h_*(t):=t^{p-1}+a_*t^{q-1},
\]
the same estimate holds with $H_*\circ h_*^{-1}$ in place of $H_{\bar x}\circ h_{\bar x}^{-1}$, with a constant depending only on $\data$ and $E$.
\end{proposition}

\begin{proof}
We start from the reverse H\"older inequality \eqref{rev.self.improved} that, together with \eqref{above.below} that implies $a(x)\approx_M a(\bar x)$, gives
\[
\mean{B_{\rho/2}} H(x,Dv)\dx \leq c(\data,E,M,\sigma){\biggl(\mean{B_\rho}\bigl[H(\bar x,Dv)\bigr]^\sigma\dx\biggr)}^{\frac1\sigma}.
\]
Fixing $\sigma$, depending on $p$ and $q$, sufficiently small so that the map $t\mapsto [(H_{\bar x}\circ h_{\bar x}^{-1})(t)]^\sigma$ is concave, see \eqref{Psi.small.concave} in the Appendix, concludes the proof, recalling Jensen's inequality: indeed
\[
\mean{B_\rho}\bigl[H(\bar x,Dv)\bigr]^\sigma\dx=\mean{B_\rho}\bigl[H_{\bar x}(Dv)\bigr]^\sigma\dx\leq  \biggl[\bigl(H_{\bar x}\circ h_{\bar x}^{-1}\bigr){\biggl(\mean{B_\rho} h_{\bar x}(Dv)\dx \biggr)}\biggr]^\sigma
\]
and using again \eqref{above.below} and the properties in \eqref{constant.h}-\eqref{constant.hinverse}
\[
{\biggl(\mean{B_\rho}\bigl[H(\bar x,Dv)\bigr]^\sigma\dx\biggr)}^{\frac1\sigma}\leq \bigl(H_{\bar x}\circ h_{\bar x}^{-1}\bigr){\biggl(\mean{B_\rho} h_{\bar x}(Dv)\dx \biggr)}\leq c(p,q,M)\bigl(H_{\bar x}\circ h_{\bar x}^{-1}\bigr)\biggl(\mean{B_\rho} h(x,Dv)\dx \biggr).
\]
The alternative formulation follows by the same argument, using $H_*\leq H(x,\,\cdot\,)\leq2H_*$ and $h_*\leq h(x,\,\cdot\,)\leq2h_*$ throughout $B_\rho$.
\end{proof}

\subsection{Rough comparison estimates}

We now establish comparison estimates between a solution to \eqref{equation} and the corresponding homogeneous comparison map. We first work with a local energy solution $u\in W^{1,H(\cdot)}_{\rm loc}(\Omega)$, under \eqref{a.holder} to \eqref{bounds} and \eqref{A.weak}; when $p\leq n$, we additionally assume that $\mu\in L^\infty_{\rm loc}(\Omega)$.

\medskip

Fix $B_R\equiv B_R(x_0)$ with $x_0\in\Omega_0$ and $0<R\leq\min\{1,d_{\Omega_0}/8\}$, so that $B_R\subset\Omega_1$. We define $v$ as the unique weak solution to
\begin{equation}\label{def.v}
\begin{cases}
-\divergence A(x,Dv)=0&\text{in $B_{R/2}$},\\[2mm]
v\in u+W^{1,H(\cdot)}_0(B_{R/2}).
\end{cases}
\end{equation}
For the moment, we fix an energy bound $E\geq1$ such that
\[
\benergy(Dv,B_{R/2})\leq E,
\]
that is, \eqref{bound.energy.mdata} holds with $R/2$ in place of $R$. Accordingly, the estimates for homogeneous solutions established in Section~\ref{cor:nice-rev.sec} will be applied to $v$ with outer radius $R/2$. In Paragraph~\ref{choice.E}, we shall show that $E$ can be chosen depending only on $\data$ and $\gdata$, uniformly over the balls under consideration.

For a local SOLA, this construction is applied to an admissible approximating sequence $\{u_j,\mu_j\}$, yielding comparison maps $v_j\in u_j+W^{1,H(\cdot)}_0(B_{R/2})$. The comparison estimates below are then applied to the pairs $(u_j,v_j)$; the uniform choice of $E$ along the approximations will also be justified in Paragraph~\ref{choice.E}.

\medskip

Thanks to the fact that for $p\leq n$ the inequality in \eqref{bounds} is strict, we can choose two exponents $\overline \x\equiv \overline \x(\data),\underline \x\equiv \underline \x(\data,E)$ such that
\[
\overline \x=\frac12{\biggl(\alpha-\frac{(n-1)(q-p)}{p-1}\biggr)},\qquad \underline \x=\frac12{\biggl(\alpha-(q-p)\frac{n}{p(1+\eta_0)}\biggr)},
\]
where $\eta_0$ is the positive exponent from Theorem~\ref{hi.correct.exponent}. We thus choose a generic non-negative exponent such that
\begin{equation}\label{iccs}
 \x\in
\begin{cases}
[0,\overline \x]\qquad &\text{if $p\leq n$}\\[2mm]
[0,\underline \x] &\text{if $p>n$}
\end{cases};
\end{equation}
notice that thanks to the choices of $\overline\x, \underline\x$ and the bound \eqref{bounds}, we have 
\begin{equation}\label{bound.with.x}
0<s_0:=\begin{cases}
\displaystyle{ \alpha-\overline\x-(q-p)\frac{n-1}{p-1}}\qquad &\text{if $p\leq n$}\\[3mm]  
\displaystyle{ \alpha-\underline \x-(q-p)\frac{n}{p(1+\eta_0)}}\qquad &\text{if $p>n$}
\end{cases}.
\end{equation}
With this ``choice'' of $\x$, and for a constant $M\geq1$ to be fixed, we set
\[
a_i(B_R):=\inf_{{B_R}}a(\,\cdot\,)=\min_{\overline{B_R}}a(\,\cdot\,);
\]
we say that we are in the modified (if $\x>0$) $p$-phase if
\begin{equation}\label{p.phase}
a_i(B_R)\leq M[a]_{0,\alpha}R^{\alpha-\x},
\end{equation}
whereas we say that we are in the modified $(p,q)$-phase if the complementary condition
\begin{equation}\label{pq.phase}
a_i(B_R)>M[a]_{0,\alpha}R^{\alpha-\x}
\end{equation}
holds. By Remark~\ref{stability.p.phase}, applied with equal radii, the $p$-phase condition gives
\begin{equation}\label{sup-inf_withx}
\sup_{B_R}a(\,\cdot\,)\leq(M+2)[a]_{0,\alpha}R^{\alpha-\x}\leq3M[a]_{0,\alpha}R^{\alpha-\x},
\end{equation}
where we used $M\geq1$. In the $(p,q)$-phase, \eqref{oscillation.a} instead gives, for every $x\in\overline{B_R}$,
\begin{equation}\label{comparability.a}
a_i(B_R)\leq a(x)\leq a_i(B_R)+2[a]_{0,\alpha}R^\alpha\leq\Bigl(1+\frac{2R^\x}{M}\Bigr)a_i(B_R);
\end{equation}
in particular,
\begin{equation}\label{relative.oscillation.a}
|a(x)-a_i(B_R)|\leq\frac{2R^\x}{M}a_i(B_R).
\end{equation}
For $M\geq2$, since $R\leq1$ and $\x\geq0$, \eqref{comparability.a} yields
\[
a_i(B_R)\leq a(x)\leq2a_i(B_R)\qquad\text{for every }x\in\overline{B_R},
\]
uniformly in $M$. The sharper relative oscillation estimate \eqref{relative.oscillation.a} will be used in the freezing argument.

\begin{remark}
{\rm The comparison estimates of this Section remain valid with $\x=0$. Moreover, the uniform choice of $E$ in Paragraph~\ref{choice.E} follows from Remark~\ref{rem:singular.p.component} and the estimates for $u$, without using any phase assumption or the higher integrability exponent $\eta_0$. Once $E$ has been fixed in terms of $\data$ and $\gdata$, Theorem~\ref{hi.correct.exponent} provides an exponent $\eta_0$ depending only on these data, and the same is consequently true of $\underline\x$. Thus there is no circular dependence in the choice of these parameters. We retain the full range of $\x$ in \eqref{iccs} because positive values will be needed in subsequent applications, in particular in Section~\ref{lin.chapter}, where the factor $R^\x/M$ in the freezing estimates must tend to zero as $R\searrow0$, with $M$ fixed.}
\end{remark}

The first comparison estimate we need is basic:
\begin{lemma}\label{comparison.uv.weighted}
Assume \eqref{a.holder} to \eqref{bounds} and \eqref{A.weak}. Let $u$ be a local energy solution to \eqref{equation} in the sense of Definition~\ref{def.local.energy.solution}, with the convention stated after Definition~\ref{def.solution}. Let $v$ be the homogeneous comparison map defined in \eqref{def.v} on a ball $B_R=B_R(x_0)$, with $x_0\in\Omega_0$ and $0<R\leq\min\{1,d_{\Omega_0}/8\}$. Then there exists a constant $c$, depending only on $n,p,q$ and $\nu$, such that for every $\xi>1$ and $\kappa>0$ we have
\begin{equation}\label{compa.V.gen}
\mean{B_{R/2}}\frac{|V_p(Du)-V_p(Dv)|^2+a(x)|V_q(Du)-V_q(Dv)|^2}{(\kappa+|u-v|)^\xi}\dx\leq c\frac{\kappa^{1-\xi}}{\xi-1}\frac{|\mu|(B_{R/2})}{(R/2)^n}.
\end{equation}
Moreover, let $\bar a_p=1$ and let $\bar a_q\geq0$ be any constant satisfying $\bar a_q\leq a(x)$ for almost every $x\in B_{R/2}$; for $\esp\in\{p,q\}$ we have
\begin{equation}\label{compa.V.component.L1}
\bar a_\esp\biggl(\mean{B_{R/2}}|V_\esp(Du)-V_\esp(Dv)|^{\frac2\esp}\dx\biggr)^\esp\leq c\,\frac{|\mu|(B_{R/2})}{(R/2)^{n-1}}\mean{B_{R/2}}|Du-Dv|\dx,
\end{equation}
where $c$ depends only on $n,p,q$ and $\nu$. Finally, for each $\esp\in\{p,q\}$ and every exponent $\ti$ satisfying
\[
\begin{cases}
\displaystyle \frac1{\esp-1}\leq\ti<n'&\text{if $\esp<2$},\\[2mm]
0<\ti<n'&\text{if $2\leq\esp\leq n$},\\[2mm]
0<\ti\leq\esp'&\text{if $\esp>n$},
\end{cases}
\]
we also have
\begin{multline}\label{compa.q.component.ell}
\biggl(\mean{B_{R/2}}\bigl[\bar a_\esp|Du-Dv|^{\esp-1}\bigr]^\ti\dx\biggr)^{\frac1\ti}\leq c\,\frac{|\mu|(B_{R/2})}{(R/2)^{n-1}}\\
+c\,\chi_{(-\infty,2)}(\esp)\biggl[\frac{|\mu|(B_{R/2})}{(R/2)^{n-1}}\biggr]^{\esp-1}\biggl(\mean{B_{R/2}}\bigl[\bar a_\esp|Du|^{\esp-1}\bigr]^\ti\dx\biggr)^{\frac{2-\esp}{\ti}}.
\end{multline}
The second term on the right-hand side is omitted when $\esp\geq2$. Here $c$ depends only on $n,p,q,\nu$ and $\ti$ if $\esp\leq n$, and only on $n$ and $\esp$ and $\nu$ if $\esp>n$; in both cases it is independent of $\bar a_\esp$. This estimate requires only the lower bound $\bar a_q\leq a(x)$ already assumed above, without any $(p,q)$-phase condition.
\end{lemma}

\begin{proof}
Estimate \eqref{compa.V.gen} follows from the standard weighted testing argument for the $p$-Laplacian (see \cite[Lemma 2.4]{KMI1} or \cite[Lemma 1]{KMGuide}), employing the monotonicity assumption in \eqref{A.weak} and the equivalence \eqref{diffV}. For the proof of \eqref{compa.V.component.L1},  for $\esp\in\{p,q\}$, we choose
\[
\kappa=\frac R2\mean{B_{R/2}}|Du-Dv|\dx>0,\qquad \xi=\frac{n(\esp-1)}{n-1}>1,
\]
where the last inequality follows from \eqref{2-1/n}. We can suppose $\kappa\neq 0$ as otherwise the conclusion is immediate. Since $u-v\in W^{1,1}_0(B_{R/2})$, Sobolev's inequality gives
\[
\biggl(\mean{B_{R/2}}{\bigl(\kappa+|u-v|\bigr)}^{\frac n{n-1}}\dx\biggr)^{\frac{n-1}{n}}\leq c(n)\kappa.
\]
Recalling that $\xi/(\esp-1)=n/(n-1)$, H\"older's inequality with conjugate exponents $\esp,\esp'$, followed by \eqref{compa.V.gen}, yields
\begin{align*}
\bar a_\esp\biggl(\mean{B_{R/2}}|V_\esp(Du)-V_\esp(Dv)|^{\frac2\esp}\dx\biggr)^\esp&\leq\mean{B_{R/2}}\frac{\bar a_\esp|V_\esp(Du)-V_\esp(Dv)|^2}{(\kappa+|u-v|)^\xi}\dx\biggl(\mean{B_{R/2}}(\kappa+|u-v|)^{\frac n{n-1}}\dx\biggr)^{\esp-1}\\
&\leq c\,\frac{\kappa^{1-\xi}}{\xi-1}\frac{|\mu|(B_{R/2})}{(R/2)^n}\kappa^\xi\leq c\,\frac{|\mu|(B_{R/2})}{(R/2)^{n-1}}\mean{B_{R/2}}|Du-Dv|\dx.
\end{align*}
Here we used $\bar a_p=1$ and $\bar a_q\leq a(x)$ to retain the corresponding component in \eqref{compa.V.gen}. Since $\xi$ is fixed in terms of $n$ and $\esp$, the constant depends only on $n,p,q$ and $\nu$. To prove \eqref{compa.q.component.ell}, fix $\esp\in\{p,q\}$. We may assume $\bar a_\esp>0$, since otherwise the conclusion is immediate; in particular, $u-v\in W^{1,\esp}_0(B_{R/2})$. First suppose that $\esp\leq n$. It suffices to consider $\ti\geq1/(\esp-1)$: when $\esp\geq2$, the remaining exponents will follow by H\"older's inequality, while for $\esp<2$ this lower bound is already required in the statement. Since $\ti<n'\leq\esp'$, we have $1\leq\ti(\esp-1)<\esp\leq n$. We choose
\[
\kappa=\frac R2\biggl(\mean{B_{R/2}}|Du-Dv|^{\ti(\esp-1)}\dx\biggr)^{\frac1{\ti(\esp-1)}},\qquad \xi=\frac{n[\esp-\ti(\esp-1)]}{n-\ti(\esp-1)},
\]
and we may suppose $\kappa>0$, as otherwise the conclusion is immediate. Our choice gives
\[
\xi-1=\frac{(\esp-1)[n-(n-1)\ti]}{n-\ti(\esp-1)}>0,\qquad \frac{\xi\ti}{\esp'-\ti}=[\ti(\esp-1)]^*,
\]
where the strict inequality follows from $\ti<n'$. By Sobolev's inequality and the choice of $\kappa,\xi$, we have
\[
\biggl(\mean{B_{R/2}}(\kappa+|u-v|)^{\frac{\xi\ti}{\esp'-\ti}}\dx\biggr)^{\frac{\esp'-\ti}{\xi\ti}}\leq c\kappa
\]
and H\"older's inequality with conjugate exponents $\esp'/\ti$ and $\esp'/(\esp'-\ti)$, followed by \eqref{compa.V.gen}, yields
\begin{align*}
\bar a_\esp^\ti\mean{B_{R/2}}|V_\esp(Du)-V_\esp(Dv)|^{\frac{2\ti}{\esp'}}\dx&\leq\bar a_\esp^{\frac\ti\esp}\biggl(\mean{B_{R/2}}\frac{\bar a_\esp|V_\esp(Du)-V_\esp(Dv)|^2}{(\kappa+|u-v|)^\xi}\dx\biggr)^{\frac\ti{\esp'}}\\
&\pushright{\times\,\biggl(\mean{B_{R/2}}(\kappa+|u-v|)^{\frac{\xi\ti}{\esp'-\ti}}\dx\biggr)^{\frac{\esp'-\ti}{\esp'}}}\\
&\leq c\,\bar a_\esp^{\frac\ti\esp}\kappa^{\frac{(1-\xi)\ti}{\esp'}}\biggl[\frac{|\mu|(B_{R/2})}{(R/2)^n}\biggr]^{\frac\ti{\esp'}}\kappa^{\frac{\xi\ti}{\esp'}}=c\,\bar a_\esp^{\frac\ti\esp}\kappa^{\frac\ti{\esp'}}\biggl[\frac{|\mu|(B_{R/2})}{(R/2)^n}\biggr]^{\frac\ti{\esp'}}.
\end{align*}
Here we used $\bar a_p=1$ and $\bar a_q\leq a(x)$ to retain the corresponding component in \eqref{compa.V.gen}; the constant $c$ depends only on $n,p,q,\nu$ and $\ti$ and is independent of $\bar a_\esp$. Substituting the definition of $\kappa$ into the preceding estimate, we obtain
\[
\bar a_\esp^\ti\mean{B_{R/2}}|V_\esp(Du)-V_\esp(Dv)|^{\frac{2\ti}{\esp'}}\dx\leq c\biggl(\mean{B_{R/2}}\bigl[\bar a_\esp|Du-Dv|^{\esp-1}\bigr]^\ti\dx\biggr)^{\frac1\esp}\biggl[\frac{|\mu|(B_{R/2})}{(R/2)^{n-1}}\biggr]^{\frac\ti{\esp'}}.
\]
If $\esp\geq2$, \eqref{diffV2} and Young's inequality with conjugate exponents $\esp,\esp'$ yield
\begin{align*}
\mean{B_{R/2}}\bigl[\bar a_\esp|Du-Dv|^{\esp-1}\bigr]^\ti\dx&\leq c\,\bar a_\esp^\ti\mean{B_{R/2}}|V_\esp(Du)-V_\esp(Dv)|^{\frac{2\ti}{\esp'}}\dx\\
&\leq c\biggl(\mean{B_{R/2}}\bigl[\bar a_\esp|Du-Dv|^{\esp-1}\bigr]^\ti\dx\biggr)^{\frac1\esp}\biggl[\frac{|\mu|(B_{R/2})}{(R/2)^{n-1}}\biggr]^{\frac\ti{\esp'}}\\
&\leq\frac12\mean{B_{R/2}}\bigl[\bar a_\esp|Du-Dv|^{\esp-1}\bigr]^\ti\dx+c\biggl[\frac{|\mu|(B_{R/2})}{(R/2)^{n-1}}\biggr]^\ti;
\end{align*}
absorbing the first term on the right-hand side and taking the power $1/\ti$ proves \eqref{compa.q.component.ell} when $2\leq\esp\leq n$ and $\ti\geq1/(\esp-1)$, while the remaining exponents $0<\ti<1/(\esp-1)$ follow by H\"older's inequality. Suppose now that $\esp<2$: applying \eqref{diffVp} with $z_1=Dv$ and $z_2=Du$, followed by H\"older's inequality with conjugate exponents $2/\esp$ and $2/(2-\esp)$, gives
\begin{align*}
\mean{B_{R/2}}\bigl[\bar a_\esp|Du-Dv|^{\esp-1}\bigr]^\ti\dx&\leq c\,\bar a_\esp^\ti\mean{B_{R/2}}|V_\esp(Du)-V_\esp(Dv)|^{\frac{2\ti}{\esp'}}\dx\\
&\quad+c\biggl(\bar a_\esp^\ti\mean{B_{R/2}}|V_\esp(Du)-V_\esp(Dv)|^{\frac{2\ti}{\esp'}}\dx\biggr)^{\frac\esp2}\biggl(\mean{B_{R/2}}\bigl[\bar a_\esp|Du|^{\esp-1}\bigr]^\ti\dx\biggr)^{\frac{2-\esp}2}.
\end{align*}
Using the estimate for the $V_\esp$-difference obtained before distinguishing the cases, and appealing to Young's inequality first with conjugate exponents $\esp,\esp'$ and then with exponents $2,2$, we infer
\begin{align*}
&\mean{B_{R/2}}\bigl[\bar a_\esp|Du-Dv|^{\esp-1}\bigr]^\ti\dx\\
&\quad\leq\frac14\mean{B_{R/2}}\bigl[\bar a_\esp|Du-Dv|^{\esp-1}\bigr]^\ti\dx+c\biggl[\frac{|\mu|(B_{R/2})}{(R/2)^{n-1}}\biggr]^\ti\\
&\quad\qquad\qquad+c\biggl[\frac{|\mu|(B_{R/2})}{(R/2)^{n-1}}\biggr]^{\frac{\ti(\esp-1)}2}\biggl(\mean{B_{R/2}}\bigl[\bar a_\esp|Du-Dv|^{\esp-1}\bigr]^\ti\dx\biggr)^{\frac12}\biggl(\mean{B_{R/2}}\bigl[\bar a_\esp|Du|^{\esp-1}\bigr]^\ti\dx\biggr)^{\frac{2-\esp}2}\\
&\quad\leq\frac12\mean{B_{R/2}}\bigl[\bar a_\esp|Du-Dv|^{\esp-1}\bigr]^\ti\dx+c\biggl[\frac{|\mu|(B_{R/2})}{(R/2)^{n-1}}\biggr]^\ti+c\biggl[\frac{|\mu|(B_{R/2})}{(R/2)^{n-1}}\biggr]^{\ti(\esp-1)}\biggl(\mean{B_{R/2}}\bigl[\bar a_\esp|Du|^{\esp-1}\bigr]^\ti\dx\biggr)^{2-\esp};
\end{align*}
reabsorbing the first term on the right-hand side proves \eqref{compa.q.component.ell} also when $\esp<2$. Suppose finally that $\esp>n$. Since $u-v\in W^{1,\esp}_0(B_{R/2})$, Morrey's inequality gives
\[
\|u-v\|_{L^\infty(B_{R/2})}\leq c(n,\esp)\frac R2\biggl(\mean{B_{R/2}}|Du-Dv|^\esp\dx\biggr)^{\frac1\esp}.
\]
Testing the difference of \eqref{equation} and \eqref{def.v} with $u-v$, using the monotonicity assumption in \eqref{A.weak}, $\bar a_p=1$, $\bar a_q\leq a(x)$ and $\esp>n\geq2$, we obtain
\[
\bar a_\esp\mean{B_{R/2}}|Du-Dv|^\esp\dx\leq c\,\frac{|\mu|(B_{R/2})}{(R/2)^n}\|u-v\|_{L^\infty(B_{R/2})}\leq c\,\frac{|\mu|(B_{R/2})}{(R/2)^{n-1}}\biggl(\mean{B_{R/2}}|Du-Dv|^\esp\dx\biggr)^{\frac1\esp},
\]
where $c$ depends only on $n,\esp$ and $\nu$ and is independent of $\bar a_\esp$. Consequently, trivially also when $Du=Dv$ almost everywhere,
\[
\biggl(\mean{B_{R/2}}\bigl[\bar a_\esp|Du-Dv|^{\esp-1}\bigr]^{\esp'}\dx\biggr)^{\frac1{\esp'}}=\bar a_\esp\biggl(\mean{B_{R/2}}|Du-Dv|^\esp\dx\biggr)^{\frac1{\esp'}}\leq c\,\frac{|\mu|(B_{R/2})}{(R/2)^{n-1}}.
\]
H\"older's inequality yields \eqref{compa.q.component.ell} for every $0<\ti\leq\esp'$, with the same constant, completing the proof.
\end{proof}

As a first consequence, we prove the following comparison estimate involving only the $p$-part.

\begin{lemma}\label{lem:rough.compa.solopartep}
Let $u$ and $v$ be as in Lemma~\ref{comparison.uv.weighted}. Then, for $\es$ satisfying \eqref{es},
\begin{equation}\label{rough.compa.solopartep}
{\biggl(\mean{B_{R/2}}|Du-Dv|^{(p-1)\es} \dx\biggr)}^\frac1\es\leq c\,\frac{|\mu|(B_{R})}{R^{n-1}}+c\,\chi_{(-\infty,2)}(p)\biggl[\frac{|\mu|(B_R)}{R^{n-1}}\biggr]^{p-1}{\biggl(\mean{B_R}h(x,Du)\dx\biggr)}^{2-p}
\end{equation}
for a constant $c$ depending on $\data$ and $\es$ if $p\leq n$ and only on $n,p$ and $\nu$ if $p>n$.
\end{lemma}

\begin{proof}
If $p\geq2$, the conclusion follows directly from \eqref{compa.q.component.ell} with $\esp=p$, $\bar a_p=1$ and $\ti=\es$, using $|\mu|(B_{R/2})\leq|\mu|(B_R)$ and absorbing the fixed ratio of the radii into the constant. Suppose that $p<2$ and first consider $\es\geq1/(p-1)$. Applying \eqref{compa.q.component.ell} with the same choices, followed by \eqref{local.energy.SOLA}, gives
\begin{align*}
\biggl(\mean{B_{R/2}}|Du-Dv|^{(p-1)\es}\dx\biggr)^{\frac1\es}&\leq c\,\frac{|\mu|(B_R)}{R^{n-1}}+c\biggl[\frac{|\mu|(B_R)}{R^{n-1}}\biggr]^{p-1}\biggl(\mean{B_{R/2}}|Du|^{(p-1)\es}\dx\biggr)^{\frac{2-p}{\es}}\\
&\leq c\,\frac{|\mu|(B_R)}{R^{n-1}}+c\biggl[\frac{|\mu|(B_R)}{R^{n-1}}\biggr]^{p-1}\biggl[\mean{B_R}h(x,Du)\dx+\frac{|\mu|(B_R)}{R^{n-1}}\biggr]^{2-p}\\
&\leq c\,\frac{|\mu|(B_R)}{R^{n-1}}+c\biggl[\frac{|\mu|(B_R)}{R^{n-1}}\biggr]^{p-1}\biggl(\mean{B_R}h(x,Du)\dx\biggr)^{2-p}.
\end{align*}
This proves \eqref{rough.compa.solopartep} for $\es\geq1/(p-1)$. Since $1/(p-1)<n'$ by \eqref{2-1/n}, the remaining exponents follow by applying this estimate with $\es=1/(p-1)$ and using H\"older's inequality.
\end{proof}

\begin{remark}\label{rem:singular.p.component}
{\rm A quick inspection of the proof shows that it is possible to avoid the quantity $h_x(\,\cdot\,)$ on the right-hand side in the case $p<2$; we therefore also have
\[
\mean{B_{R/2}}|Du-Dv|^{(p-1)\es} \dx\leq c\biggl[\frac{|\mu|(B_R)}{R^{n-1}}\biggr]^\es+c\,\chi_{(-\infty,2)}(p)\biggl[\frac{|\mu|(B_R)}{R^{n-1}}\biggr]^{(p-1)\es}{\biggl(\mean{B_R}|Du|^{(p-1)\es}\dx\biggr)}^{2-p}
\]
for $\es$ as in \eqref{es}, with the additional restriction $\es\geq1/(p-1)$ when $p<2$, and a constant depending only on $\data$ and $\es$. This estimate will be useful later, when fixing the constant $E$.
}
\end{remark}

Next we suppose we are in the $(p,q)$-phase, that is we suppose here \eqref{pq.phase} holds for some $M\geq2$ and some $\x\geq0$ as in \eqref{iccs}, which implies, in particular, that \eqref{comparability.a} holds in $B_R$. 

\begin{lemma}\label{lem:rough.compa.pq}
Let $u$ and $v$ be as in Lemma~\ref{comparison.uv.weighted}; suppose that \eqref{pq.phase} holds on $B_R$ for some $M\geq2$ and some $\x$ as in \eqref{iccs}. Then for every $\ti$ satisfying \eqref{ti} we have
\begin{equation}\label{rough.compa.pqphase}
{\biggl(\mean{B_{R/2}}\bigl[h_{x_\smallzero}(Du-Dv)\bigr]^\ti \dx\biggr)}^\frac1\ti\leq c\,\frac{|\mu|(B_R)}{R^{n-1}}+c\,d_{p,q}
\end{equation}
with 
\begin{equation}\label{d_pq}
d_{p,q}=\sum_{\esp\in\{p,q\}}\chi_{(-\infty,2)}(\esp)\biggl[\frac{|\mu|(B_R)}{R^{n-1}}\biggr]^{\esp-1}\biggl(\mean{B_R}h(x,Du)\dx\biggr)^{2-\esp}
\end{equation}
and for a constant $c$ depending only on $\data$ and, if $q\leq n$, also on $\ti$. In particular, $c$ is independent of $M\geq2$.
\end{lemma}

\begin{proof}
By \eqref{comparability.a}, the assumption $M\geq2$ gives $a_i(B_R)\leq a(x)\leq2a_i(B_R)$ throughout $B_R$. In particular,
\[
\frac12 h_{x_0}(z)\leq h(x,z)\leq2h_{x_0}(z)\qquad\text{for every }x\in B_R,\ z\in\mathbb R^n:
\]
thus all coefficient comparisons are uniform in $M$, and Lemma~\ref{local.rH.pq} applies in its alternative formulation. Since Lemma~\ref{lem:rough.compa.solopartep} controls the $p$-component, it remains to estimate the $q$-component. We show that 
\begin{multline}\label{rough.compa.qcomponent}
\mean{B_{R/2}}\bigl[a(x_0)|Du-Dv|^{q-1}\bigr]^\ti\dx\leq c\biggl[\frac{|\mu|(B_R)}{R^{n-1}}\biggr]^\ti\\
+c\,\chi_{(-\infty,2)}(q)\biggl[\frac{|\mu|(B_R)}{R^{n-1}}\biggr]^{(q-1)\ti}\biggl(\mean{B_R}h(x_0,Du)\dx\biggr)^{(2-q)\ti},
\end{multline}
for $\ti$ as in \eqref{ti} and  $c$ depending only on $\data$ and, if $q\leq n$, also on $\ti$. This estimate is also of independent interest.

\subsubsection*{Range $q\geq2$}
If $q\leq n$, let $0<\ti<n'$; if $q>n$, just take $\ti=q'$ (H\"older's inequality indeed would give \eqref{rough.compa.pqphase} for every $0<\ti\leq q'$): for these exponents, \eqref{rough.compa.qcomponent} follows by considering \eqref{compa.q.component.ell} with $\esp=q$ and $\bar a_q=a_i(B_R)$, using \eqref{comparability.a} and the fact $|\mu|(B_{R/2})\leq|\mu|(B_R)$. For the estimate of the $p$-part, we use Lemma~\ref{lem:rough.compa.solopartep} that applies with $\es=\ti$: it is immediate to check that if $q\leq n$ $\ti$ satisfies the range condition in \eqref{es}, while for $q>n$ we need to use $q'<n'$ and $q'\leq p'$. \eqref{rough.compa.pqphase} simply follows by combining the two component estimates and appealing to \eqref{comparability.a}, and the dependencies of the constant stated are a consequence of the dependencies of the constants in \eqref{rough.compa.solopartep} and \eqref{rough.compa.qcomponent}.

\subsubsection*{Range $q<2$}
First suppose that $\ti\geq1/(q-1)$. Applying \eqref{compa.q.component.ell} with $\esp=q$ and $\bar a_q=a_i(B_R)$, using \eqref{comparability.a} and taking the power $\ti$, we obtain
\[
\mean{B_{R/2}}\bigl[a(x_0)|Du-Dv|^{q-1}\bigr]^\ti\dx\leq c\biggl[\frac{|\mu|(B_{R/2})}{(R/2)^{n-1}}\biggr]^\ti+c\biggl[\frac{|\mu|(B_{R/2})}{(R/2)^{n-1}}\biggr]^{\ti(q-1)}\biggl(\mean{B_{R/2}}\bigl[h_{x_0}(Du)\bigr]^\ti\dx\biggr)^{2-q},
\]
where $c$ depends only on $\data$ and $\ti$ and it is independent of $M$. By \eqref{local.energy.SOLA.pq} and \eqref{comparability.a}, we have
\begin{align*}
\biggl[\frac{|\mu|(B_{R/2})}{(R/2)^{n-1}}\biggr]^{\ti(q-1)}\biggl(\mean{B_{R/2}}\bigl[h_{x_0}(Du)\bigr]^\ti&\dx\biggr)^{2-q}\\
&\leq c\biggl[\frac{|\mu|(B_R)}{R^{n-1}}\biggr]^{\ti(q-1)}\biggl[\mean{B_R}h(x_0,Du)\dx+\frac{|\mu|(B_R)}{R^{n-1}}\biggr]^{\ti(2-q)}\\
&\leq c\biggl[\frac{|\mu|(B_R)}{R^{n-1}}\biggr]^\ti+c\biggl[\frac{|\mu|(B_R)}{R^{n-1}}\biggr]^{\ti(q-1)}\biggl(\mean{B_R}h(x_0,Du)\dx\biggr)^{\ti(2-q)}.
\end{align*}
Together with $|\mu|(B_{R/2})\leq|\mu|(B_R)$, this proves \eqref{rough.compa.qcomponent} for $\ti\geq1/(q-1)$. Since $1/(q-1)<n'$ by \eqref{2-1/n}, the remaining exponents follow by applying this estimate with $\ti=1/(q-1)$ and using H\"older's inequality. Combining \eqref{rough.compa.qcomponent} with Lemma~\ref{lem:rough.compa.solopartep} and \eqref{comparability.a} proves \eqref{rough.compa.pqphase} also when $q<2$.

\end{proof}

\begin{corollary}\label{comparison.pq.frozen}
Under the assumptions of Lemma~\ref{lem:rough.compa.pq}, we have
\begin{equation*}
\mean{B_{R/2}}|Du-Dv| \dx\leq c\,h_{x_\smallzero}^{-1}\biggl(\frac{|\mu|(B_R)}{R^{n-1}}+d_{p,q}\biggr)
\end{equation*}
for a constant depending on $\data$, with $d_{p,q}$ defined in \eqref{d_pq}.
\end{corollary}

\begin{proof}
We choose $\ti:=\max\{1,1/(p-1)\}$: this exponent is admissible in \eqref{ti}. Indeed, this is immediate when $p\geq2$. If $p<2$ and $q\leq n$, then $\ti=1/(p-1)<n'$ by \eqref{2-1/n}, while if $p<2$ and $q>n$ \eqref{bounds} together with $\alpha\leq1$ and \eqref{2-1/n} imply
\[
(2-p)(q-1)<\frac{p-1}{n-1}\leq p-1,
\]
which implies $q-1<q(p-1)$ whence $\ti<q'$. Since $\ti(p-1)\geq1$, moreover, Lemma~\ref{convexity.h.power} ensures that $t\mapsto[h_{x_\smallzero}(t)]^\ti$ is convex. Jensen's inequality and \eqref{rough.compa.pqphase} therefore give
\[
h_{x_\smallzero}\biggl(\mean{B_{R/2}}|Du-Dv|\dx\biggr)\leq\biggl(\mean{B_{R/2}}[h_{x_\smallzero}(Du-Dv)]^\ti\dx\biggr)^{\frac1\ti}\leq c\biggl[\frac{|\mu|(B_R)}{R^{n-1}}+d_{p,q}\biggr]
\]
and applying $h_{x_\smallzero}^{-1}$ and using the scaling inequality in \eqref{constant.hinverse} proves the conclusion. Since $\ti$ depends only on $p$, the constant depends only on $\data$.
\end{proof}

\begin{remark}\label{rem:comparison.epsilon}
{\rm Notice that in all cases, regardless of whether $p$ or $q$ are smaller or larger than two, thanks to Young's inequality (only if $p$ and/or $q$ are less than two) and \eqref{comparability.a} in the $(p,q)$ phase, we have two comparison estimates in the following form:
\[
\mean{B_{R/2}}|Du-Dv|^{(p-1)\es} \dx\leq c_\varepsilon\biggl[\frac{|\mu|(B_R)}{R^{n-1}}\biggr]^\es+\varepsilon^\es{\biggl(\mean{B_R}h(x,Du)\dx\biggr)}^\es
\]
for $\es$ as in \eqref{es} and 
\[
\mean{B_{R/2}}h_{x_\smallzero}(Du-Dv) \dx\leq c_\varepsilon\frac{|\mu|(B_R)}{R^{n-1}}+\varepsilon \mean{B_R}h(x,Du)\dx
\]
if \eqref{pq.phase} holds with $M\geq2$. In both estimates, $\varepsilon\in(0,1]$ and $c_\varepsilon$ depends only on $\data,\varepsilon$, and additionally on $\es$ in the first estimate when $p\leq n$.
}
\end{remark}

\subsection{The choice of \texorpdfstring{$E$}{E}}\label{choice.E}

We now show that the constant $E$ in \eqref{bound.energy.mdata} can be chosen uniformly over the balls under consideration. In the case $p\leq n$, we have for $c\equiv c(\data)$
\begin{align*}
\int_{B_{R/2}}|Dv|^{(p-1)b}\dx&\leq c\int_{B_{R/2}}|Du|^{(p-1)b}\dx+c\int_{B_{R/2}}|Du-Dv|^{(p-1)b}\dx \\
&\leq c\int_{B_R}|Du|^{(p-1)b}\dx+c\int_{B_{R/2}}|Du-Dv|^{(p-1)b}\dx.
\end{align*}
Here and below, $x_0\in\Omega_0$ and $0<R\leq\min\{1,d_{\Omega_0}/8\}$. The first term is controlled by \eqref{global.hi.energy}, see Remark~\ref{rem:local.rH.SOLA}; the second by the form of the comparison estimate in Remark~\ref{rem:singular.p.component}:
\[
\int_{B_{R/2}}|Du-Dv|^{(p-1)b}\dx\leq c\,R^n\biggl[\frac{|\mu|(B_R)}{R^{n-1}}\biggr]^b+c\,\chi_{(-\infty,2)}(p)R^n\biggl[\frac{|\mu|(B_R)}{R^{n-1}}\biggr]^{(p-1)b}{\biggl(\mean{B_R}|Du|^{(p-1)b}\dx\biggr)}^{2-p};
\]
we have, since $R\leq 1$,
\[
R^n\biggl[\frac{|\mu|(B_R)}{R^{n-1}}\biggr]^b\leq R^{n-b(n-1)} \bigl[|\mu|(\Omega)\bigr]^b\leq c(\data,\gdata)
\]
and, when $p<2$,
\begin{align*}
R^n\biggl[\frac{|\mu|(B_R)}{R^{n-1}}\biggr]^{(p-1)b}&{\biggl(\mean{B_R}|Du|^{(p-1)b}\dx\biggr)}^{2-p}\\
&\leq c\,R^{n-(n-1)(p-1)b-n(2-p)}\bigl[|\mu|(\Omega)\bigr]^{(p-1)b}{\biggl(\int_{B_R}|Du|^{(p-1)b}\dx\biggr)}^{2-p}\\
&\leq c(\data,\gdata)\,R^{(p-1)[n-(n-1)b]}\leq c(\data,\gdata),
\end{align*}
using again \eqref{global.hi.energy}. If $p>n$, on the other hand and similarly, recalling the expression for the global energy of $Du$ in \eqref{energy.M}
\[
\int_{B_{R/2}}|Dv|^p\dx\leq c\,\mathcal M+c\int_{B_{R/2}}|Du-Dv|^p\dx.
\]
The comparison term is estimated by Remark~\ref{rem:singular.p.component} with $\es=p'$:
\[
\int_{B_{R/2}}|Du-Dv|^p\dx\leq c\,R^n\biggl[\frac{|\mu|(B_R)}{R^{n-1}}\biggr]^{p'}\leq c\,R^{\frac{p-n}{p-1}} \bigl[|\mu|(\Omega)\bigr]^{p'} \leq c(\data,\gdata),
\]
as $R\leq 1$ and $p-n>0$. Therefore, we can choose $E\equiv E(\data,\gdata)\geq1$ independently of $x_0\in\Omega_0$ and $0<R\leq\min\{1,d_{\Omega_0}/8\}$. When $p\leq n$, this also applies uniformly to the homogeneous comparison maps associated via \eqref{def.v} to an admissible approximating sequence $\{u_j,\mu_j\}$ for a local SOLA. Indeed, for each fixed $x_0\in\Omega_0$, Lemma~\ref{SOLA.flux.stability} and \eqref{SOLA.uniform.meas} give
\[
\int_{\Omega_1}h(x,Du_j)\dx\leq1+\mathcal M,\qquad |\mu_j|(B_{d_{\Omega_0}/4}(x_0))\leq1+|\mu|(\Omega)
\]
 after discarding finitely many terms. If we apply Lemma~\ref{local.rH} with $\es=b$ on $B_{d_{\Omega_0}/4}(x_0)$ and repeat the preceding argument, we find
\[
\benergy(Dv_j,B_{R/2}(x_0))\leq E
\]
for every remaining index $j$ and every $0<R\leq\min\{1,d_{\Omega_0}/8\}$. Note that $E$ may be enlarged, still depending only on $\data$ and $\gdata$; while the number of initial terms we discard may depend on $x_0$, both $E$ and the range of radii do not. As a result, $\eta_0$ and $\underline\x$ are determined solely in terms of $\data$ and $\gdata$, and the same applies to the constants in Theorem~\ref{hi.correct.exponent} and any constant depending only on $\data$ and $E$ is thus uniform along the sequence.

\section{Morrey decay estimates}

In this section we establish Morrey decay estimates for the energy of solutions to \eqref{homogeneous.equation} starting from local H\"older continuity under the structural assumptions \eqref{A.weak}. An intrinsic counterpart for minimizers of double phase functionals appears in \cite{BCMgen}; our goal here is to trace explicitly how the energy bound fixed in Paragraph~\ref{choice.E} governs the constants involved. Throughout the section, $\alpha_0\in(0,1)$ denotes a H\"older exponent for solutions to \eqref{homogeneous.equation}, whose qualitative existence is guaranteed by \cite{CM1}. Once the energy bound is prescribed, both $\alpha_0$ and the multiplicative constants in the decay estimates can be selected solely in terms of $\data$ and $\gdata$, while keeping track of the exact decay exponents in each regime. Notice, however, that the two phase hypotheses are not complementary: the $p$-phase condition is enforced on the inner ball, whereas the $(p,q)$-phase condition is placed on the outer one (see \eqref{p.phase.energy.decay} and \eqref{pq.phase.energy.decay}, respectively). This deliberate asymmetry will prove essential in the exit-time arguments that follow.

\medskip

\subsection{Morrey decay at the energy level from H\"older continuity} 

Throughout this section, every homogeneous solution $v$ to \eqref{homogeneous.equation} in $B_R(x_0)$ (with $R\leq1$) is required to satisfy \eqref{bound.energy.mdata} in terms of the constant $E\equiv E(\data,\gdata)\geq1$ fixed in Paragraph~\ref{choice.E}; whenever $v$ is defined exclusively on $B_{R/2}(x_0)$, the energy bound is naturally understood to hold on that reduced domain. While this condition follows as an immediate consequence of Paragraph~\ref{choice.E} for the comparison maps introduced in \eqref{def.v}, it is taken as a standing hypothesis for general homogeneous solutions.

\medskip 

We begin by proving an oscillation decay estimate with constants depending only on $\data$ and $E$. Fix a ball $B_\rho\equiv B_\rho(x_1)\Subset B_R(x_0)$ and assume first that
\[
\sup_{B_\rho}a(\,\cdot\,)\leq4[a]_{0,\alpha}\rho^\alpha.
\]
Under this condition, the arbitrary-radius version of Lemma~\ref{AsCaccioppoli2}, applied on $B_{\rho/2}$, permits us to carry out the De Giorgi iteration of \cite[Lemmas 10.2-10.3]{CM1}, whence
\[
\osc_{B_{\rho/8}}v\leq\vartheta\osc_{B_{\rho/2}}v\leq\vartheta\osc_{B_\rho}v
\]
for some $\vartheta\in(0,1)$ depending only on $\data$ and $E$; when $p\geq n$, the same iteration goes through by using the Sobolev embedding into $L^{2p}$ in place of $L^{p^*}$. Suppose now that
\[
\sup_{B_\rho}a(\,\cdot\,)>4[a]_{0,\alpha}\rho^\alpha;
\]
then $a(x)$ and $a(x_1)$ are comparable on $B_\rho$ by \eqref{above.below}. Using $v-w$ as a test function in the homogeneous equation (for $v-w$ having compact support in $B_\rho$) and applying Young's inequality, we see that $v$ is a quasi-minimizer of the functional with integrand $H(x,\,\cdot\,)$, with quasi-minimality constant depending only on $p,q,\nu$ and $L$; because of the comparability of $a(\,\cdot\,)$, the same property holds for the frozen integrand $H(x_1,\,\cdot\,)$. We can therefore apply \cite[Theorem 11.3]{CM1}, whose constants do not depend on the frozen coefficient or on the specific quasi-minimizer, and choose $\tau_0\in(0,\min\{1/8,\vartheta\})$, depending only on $\data$ and $E$, so that in both cases
\[
\osc_{B_{\tau_0\rho}}v\leq\vartheta\osc_{B_\rho}v.
\]
Since $0<\tau_0<\vartheta<1$, we have $\alpha_0:=\log\vartheta/\log\tau_\smallzero\in(0,1)$ and $\tau_0^{\alpha_0}=\vartheta$. Applying Lemma~\ref{geometric.decay.iteration} (in the homogeneous case, to the non-decreasing function $r\mapsto\osc_{B_r(x_1)}v$ with $\tau=\tau_0$ and $\sigma=\alpha_0$) then yields
\begin{equation}\label{oscillation.v}
\osc_{B_{\rho_\smallone}}v\leq c\Bigl(\frac{\rho_1}{\rho_2}\Bigr)^{\alpha_0}\osc_{B_{\rho_\smalltwo}}v,\qquad 0<\rho_1\leq\rho_2\leq\rho_3,
\end{equation}
for concentric balls
\begin{equation}\label{double.ball.inclusion}
B_{\rho_\smallone}\equiv B_{\rho_\smallone}(x_1)\subset B_{\rho_\smalltwo}\equiv B_{\rho_\smalltwo}(x_1)\subset B_{\rho_\smallthree}\equiv B_{\rho_\smallthree}(x_1)\Subset B_R(x_0),
\end{equation}
where $c$ and $\alpha_0$ depend only on $\data$ and $E$, and hence only on $\data$ and $\gdata$ under our standing choice of $E$.

\begin{lemma}[Morrey estimate at the natural integrability in the $p$-phase]\label{decay.energy.p}
Let $v\in W^{1,H(\cdot)}(B_R(x_0))$ be a solution to \eqref{homogeneous.equation}, suppose that \eqref{bound.energy.mdata} holds and let $B_{\rho_\smallone}, B_{\rho_\smalltwo}$ be two balls with center $x_1$, as in \eqref{double.ball.inclusion}. Suppose that, for some $M\geq2$,
\begin{equation}\label{p.phase.energy.decay}
\sup_{B_{\rho_\smallone}}a(\,\cdot\,) \leq M{[a]}_{0,\alpha}\rho_1^\alpha.
\end{equation}
Then there exists a constant $c\geq1$, depending on $\data,\gdata$ and $M$, such that
\begin{equation}\label{decay.energy.p.est}
\mean{B_{\rho_\smallone}(x_1)}H(x,Dv)\dx\leq c\Bigl(\frac{\rho_1}{\rho_2}\Bigr)^{-(1-\alpha_\smallzero)p}\mean{B_{\rho_\smalltwo}(x_1)}H(x,Dv)\dx .
\end{equation}
\end{lemma}

\begin{proof}
We prove \eqref{decay.energy.p.est} only for $\rho_1\leq\rho_2/4$; otherwise, the estimate simply follows by enlarging the integration ball, since $\rho_1/\rho_2>1/4$. By Remark~\ref{stability.p.phase} and \eqref{p.phase.energy.decay}, we have
\[
\sup_{B_{4\rho_\smallone}}a(\,\cdot\,)\leq(M+2)[a]_{0,\alpha}(4\rho_1)^\alpha.
\]
We can therefore apply Lemma~\ref{AsCaccioppoli} with $\rho=4\rho_1$, since $B_{4\rho_1}\subseteq B_{\rho_2}$. Recalling our choice of $E$, we obtain
\[
\mean{B_{\rho_\smallone}}H(x,Dv)\dx \leq c(\data,\gdata,M)\mean{B_{2\rho_\smallone}}\bigg|\frac{v-(v)_{B_{2\rho_\smallone}}}{2\rho_1}\bigg|^p\dx \leq \frac{c}{\rho_1^p}\bigl[\osc_{B_{2\rho_\smallone}}v\bigr]^p.
\]
Using the oscillation decay estimate between $B_{2\rho_1}$ and $B_{\rho_2/2}$ and then the $\sup$-BMO estimate \eqref{sup.local.mean} now yields
\begin{align*}
\mean{B_{\rho_\smallone}}H(x,Dv)\dx&\leq c \Bigl(\frac{\rho_1}{\rho_2}\Bigr)^{-(1-\alpha_\smallzero) p} \Big[\frac{\osc_{B_{\rho_\smalltwo/2}}v}{\rho_2/2}\Big]^p\\
&\leq c \Bigl(\frac{\rho_1}{\rho_2}\Bigr)^{-(1-\alpha_\smallzero)p} {\biggl(\mean{B_{\rho_\smalltwo}}|Dv|\dx\biggr)}^p\\
&\leq c \Bigl(\frac{\rho_1}{\rho_2}\Bigr)^{-(1-\alpha_\smallzero)p} \mean{B_{\rho_\smalltwo}}H(x,Dv)\dx;
\end{align*}
in the last line we used Jensen's inequality and the lower bound $H(x,z)\geq |z|^p/p$.
\end{proof}

\begin{lemma}[Morrey estimate at the natural integrability in the $(p,q)$-phase]\label{decay.energy.pq}
Let $v\in W^{1,H(\cdot)}(B_R(x_0))$ be a solution to \eqref{homogeneous.equation} under the structural assumptions \eqref{2-1/n} to \eqref{bounds} and \eqref{A.weak}, and let $B_{\rho_\smallone},B_{\rho_\smalltwo}$ be as in \eqref{double.ball.inclusion}. Suppose that, for some $M>2$,
\begin{equation}\label{pq.phase.energy.decay}
\sup_{B_{\rho_\smallthree}}a(\,\cdot\,)>M{[a]}_{0,\alpha}\rho_3^\alpha;
\end{equation}
then there exists a constant $c\geq1$, depending on $\data,\gdata$ and $M$, such that
\begin{equation}\label{decay.energy.pq.est}
\mean{B_{\rho_\smallone}(x_1)}H(x,Dv)\dx\leq c\Bigl(\frac{\rho_1}{\rho_2}\Bigr)^{-(1-\alpha_\smallzero)q}\mean{B_{\rho_\smalltwo}(x_1)}H(x,Dv)\dx\qquad \text{for every $\rho_1\leq \rho_2$.}
\end{equation}

\end{lemma}

\begin{proof}
We first consider the case $\rho_1\leq\rho_2/4$: by \eqref{pq.phase.energy.decay} and Remark~\ref{essent.constant}, we have $a(x)\approx_M a(x_1)$ and $H(x,z)\approx_M H(x_1,z)$ for every $x\in B_{\rho_2}$ and $z\in\mathbb R^n$ and by the Caccioppoli inequality \eqref{Caccio2} and the monotonicity of $t\mapsto H(x_1,t)$ we can estimate
\[
\mean{B_{\rho_1}}H(x,Dv)\dx\leq c\mean{B_{2\rho_1}}H\biggl(x_1,\frac{v-(v)_{B_{2\rho_1}}}{\rho_1}\biggr)\dx\leq c\,H\biggl(x_1,\frac{\osc_{B_{2\rho_1}}v}{\rho_1}\biggr).
\]
 Using \eqref{oscillation.v} between $B_{2\rho_1}$ and $B_{\rho_2/2}$, followed by Corollary~\ref{corol.supBMO} and the standard Poincar\'e inequality, we obtain
\[
\frac{\osc_{B_{2\rho_1}}v}{\rho_1}\leq c\Bigl(\frac{\rho_1}{\rho_2}\Bigr)^{-(1-\alpha_0)}\frac{\osc_{B_{\rho_2/2}}v}{\rho_2/2}\leq c\Bigl(\frac{\rho_1}{\rho_2}\Bigr)^{-(1-\alpha_0)}\mean{B_{\rho_2}}|Dv|\dx;
\]
combining these estimates with the scaling property \eqref{constant.h}, Jensen's inequality, and the comparability of the frozen and non-frozen integrands gives
\begin{align*}
\mean{B_{\rho_1}}H(x,Dv)\dx&\leq c\,H\biggl(x_1,\Bigl(\frac{\rho_1}{\rho_2}\Bigr)^{-(1-\alpha_0)}\mean{B_{\rho_2}}|Dv|\dx\biggr)\\
&\leq c\Bigl(\frac{\rho_1}{\rho_2}\Bigr)^{-(1-\alpha_0)q}H\biggl(x_1,\mean{B_{\rho_2}}|Dv|\dx\biggr)\\
&\leq c\Bigl(\frac{\rho_1}{\rho_2}\Bigr)^{-(1-\alpha_0)q}\mean{B_{\rho_2}}H(x_1,Dv)\dx\\
&\leq c\Bigl(\frac{\rho_1}{\rho_2}\Bigr)^{-(1-\alpha_0)q}\mean{B_{\rho_2}}H(x,Dv)\dx
\end{align*}
with a constant $c$ depending only on $\data,\gdata$ and $M$. For $\rho_2/4<\rho_1\leq\rho_2$, estimate \eqref{decay.energy.pq.est} follows directly by enlarging the integration ball and absorbing the volume ratio into the constant.

\end{proof}

Merging the last two results gives the following Morrey estimate for the gradient of solutions of equations with merely measurable structure; we explicitly write this result also for future reference.

\begin{proposition}\label{decay.energy.level}
Let $v\in W^{1,H(\cdot)}(B_R(x_0))$ be a solution to \eqref{homogeneous.equation}. Then there exist an exponent $\alpha_0\in(0,1)$ (which is actually the H\"older exponent of solutions to \eqref{homogeneous.equation}) and a constant $c\geq1$, both depending on $\data$ and $\gdata$, such that the decay estimate
\begin{equation}\label{decay.v.H}
\mean{B_{\rho_\smallone}(x_1)}H(x,Dv)\dx\leq c\,\Bigl(\frac{\rho_1}{\rho_2}\Bigr)^{-(1-\alpha_\smallzero)q}\mean{B_{\rho_\smalltwo}(x_1)}H(x,Dv)\dx
\end{equation}
holds for every pair of concentric balls $B_{\rho_\smallone}(x_1)\subset B_{\rho_\smalltwo}(x_1)\subset B_R(x_0)$.
\end{proposition}

\begin{proof}
In this proof we omit the centers of the balls, which are all equal to $x_1$ unless otherwise stated. We apply Lemmas~\ref{decay.energy.p} and~\ref{decay.energy.pq} with $M=3$.

\medskip

For $\rho_2$ as in the statement, we set $\varrho_j=4^{-j}\rho_2$,  $B_j=B_{\varrho_j}(x_1)$, $j\in\N_0$, and we define
\[
\mathcal I:=\big\{j\in\N_0: \sup_{B_j}\, a(\,\cdot\,) > 3{[a]}_{0, \alpha}\varrho_j^\alpha\big\}\subset\N_0,\qquad i_e=\min \mathcal I,
\]
where we agree that $\min\emptyset=+\infty$. Again, we will now prove \eqref{decay.v.H} only for $\rho_1\in(0,\rho_2/4]$.

\subsubsection*{Case I} If $i_e=0$ (that is, the $(p,q)$-phase condition holds on $B_{\rho_2}(x_1)$) then we are essentially done, since
\[
\sup_{B_0} a(\,\cdot\,)=\sup_{B_{\rho_\smalltwo}(x_1)} a(\,\cdot\,) > 3{[a]}_{0, \alpha}\rho_2^\alpha\qquad \Longrightarrow \qquad \eqref{pq.phase.energy.decay}\text{ holds in $B_{\rho_\smalltwo}(x_1)$ for $M=3$}
\]
and therefore for every radius $\rho_1\leq \rho_2/4$ (actually, for every $\rho_1\leq \rho_2$) we can directly use \eqref{decay.energy.pq.est}, which is \eqref{decay.v.H}.
\subsubsection*{Case II} 
We consider here the case $i_e=+\infty$, which means that all the balls are in the $p$-phase. Clearly, for $\rho_1\in(0,\rho_2/4]$ there is an index $\ell\in\N$ such that $\varrho_{\ell+1}<\rho_1\leq \varrho_\ell$. Moreover, the fact that $\mathcal I=\emptyset$ means that it holds
\[
\sup_{B_j}a(\,\cdot\,) \leq 3{[a]}_{0, \alpha}\varrho_j^\alpha\quad\text{for every $j\in\N_0$}\qquad \Longrightarrow\qquad \sup_{B_{\varrho_\ell}}\, a(\,\cdot\,) \leq  3{[a]}_{0, \alpha}\varrho_\ell^\alpha,
\]
then we can apply \eqref{decay.energy.p.est} with $\rho_1=\varrho_\ell$ to get 
\begin{equation}\label{decay.intermediate}
\mean{B_{\varrho_\ell}}H(x,Dv)\dx \leq c\Bigl(\frac{\varrho_\ell}{\rho_2}\Bigr)^{-(1-\alpha_\smallzero) p}\mean{B_{\rho_\smalltwo}}H(x,Dv)\dx.
\end{equation}
Combining this estimate with
\[
\mean{B_{\rho_\smallone}}H(x,Dv)\dx\leq 4^n\mean{B_{\varrho_\ell}}H(x,Dv)\dx,\qquad  \Bigl(\frac{\varrho_\ell}{\rho_2}\Bigr)^{-(1-\alpha_\smallzero)p}\leq \Bigl(\frac{\rho_1}{\rho_2}\Bigr)^{-(1-\alpha_\smallzero)q}
\]
yields what we wanted.

\subsubsection*{Case III} 
If $i_e\in\N$, then for every $j\in\{0,1,\dots,i_e-1\}$ we have
\begin{equation}\label{p.here}
\sup_{B_{\varrho_j}}\, a(\,\cdot\,)\leq  3{[a]}_{0, \alpha}\varrho_j^\alpha
\end{equation}
while 
\begin{equation}\label{pq.here}
\sup_{B_{\varrho_{i_e}}}\, a(\,\cdot\,)>  3{[a]}_{0, \alpha}\varrho_{i_e}^\alpha\qquad\Longrightarrow \qquad \frac12 a(x_1)\leq a(\,\cdot\,)\leq \frac32a(x_1)\quad\text{in $B_{\varrho_{i_e}}$};
\end{equation}
thus, the balls $B_0,\dots,B_{i_e-1}$ are in the $p$-phase, whereas $B_{i_e}$ is in the $(p,q)$-phase.

\medskip 

Again we define $\bar \ell\in\N$ as the natural number such that $\varrho_{\bar \ell+1}< \rho_1 \leq \varrho_{\bar \ell}$, and the proof now follows by mixing appropriately {\em Cases I \& II}: if $\bar \ell\in\{1,\dots,i_e-1\}$ then we mimic {\em Case II} to infer
\begin{align*}
\mean{B_{\rho_\smallone}}H(x,Dv)\dx&\leq c\mean{B_{\varrho_{\bar \ell}}}H(x,Dv)\dx\\
 &\leq c\Bigl(\frac{\rho_1}{\rho_2}\Bigr)^{-(1-\alpha_\smallzero)p}\mean{B_{\rho_\smalltwo}}H(x,Dv)\dx\leq c\Bigl(\frac{\rho_1}{\rho_2}\Bigr)^{-(1-\alpha_\smallzero)q}\mean{B_{\rho_\smalltwo}}H(x,Dv)\dx 
\end{align*}
using \eqref{decay.energy.p.est} as before. On the other hand, if $\bar\ell\geq i_e$, for the transition between $\rho_1$ and $\varrho_{i_e}$ we have by \eqref{decay.energy.pq.est}, applied with $\rho_2=\varrho_{i_e}$
\[
\mean{B_{\rho_\smallone}}H(x,Dv)\dx \leq c\Bigl(\frac{\rho_1}{\varrho_{i_e}}\Bigr)^{-(1-\alpha_\smallzero)q} \mean{B_{\varrho_{i_e}}}H(x,Dv)\dx.
\]
Since for the transition between $\varrho_{i_e}$ and $\rho_2$ we know that in all the balls $B_0=B_{\rho_2},B_1,\dots, B_{i_e-1}$, in particular in $B_{i_e-1}$, we are in the $p$-phase, we can repeat the argument leading to \eqref{decay.intermediate} and infer
\[
\mean{B_{\varrho_{i_e}}}H(x,Dv)\dx \leq c\Bigl(\frac{\varrho_{i_e-1}}{\rho_2}\Bigr)^{-(1-\alpha_\smallzero)p}\mean{B_{\rho_\smalltwo}}H(x,Dv)\dx\leq c\Bigl(\frac{\varrho_{i_e}}{\rho_2}\Bigr)^{-(1-\alpha_\smallzero)q}\mean{B_{\rho_\smalltwo}}H(x,Dv)\dx.
\]
The conclusion follows by merging the last two estimates. 

If the outer ball is not compactly contained in $B_R(x_0)$, the conclusion follows by approximating its radius from below, keeping the inner radius fixed, and passing to the limit. The case of equal radii is immediate.
\end{proof}

\subsection{Preliminary harmonic-type approximation in the super-dimensional regimes}\label{par.harmonic}

The energy decay estimates \eqref{decay.energy.p.est} and \eqref{decay.energy.pq.est} are consequences of the weak assumptions \eqref{A.weak}. To iterate the resulting estimate for the flux magnitude, we need the Morrey exponent in the relevant phase to exceed one. The exponents
\[
n-(1-\alpha_0)(p-1)\qquad\text{and}\qquad n-(1-\alpha_0)(q-1)
\]
exceed one whenever $p\leq n$ and $q\leq n$, respectively, since $\alpha_0>0$. When $p>n$ in the $p$-phase, the first exponent exceeds one if and only if $\alpha_0>(p-n)/(p-1)$, and similarly for the second one. The H\"older estimate \eqref{oscillation.v} provides a positive exponent $\alpha_0$ depending only on $\data$ and $E$, but does not guarantee these stronger lower bounds.

\medskip

For this reason, throughout this part we assume that $A(\,\cdot\,,\,\cdot\,)$ satisfies the stronger structural conditions in \eqref{A.strong}, noting that it remains unclear whether the decay estimate \eqref{decay.h.Du} persists for rough vector fields under \eqref{A.weak} in the super-dimensional regimes, as we still lack a sufficiently strong quantitative regularity theory for gradients of solutions to the corresponding homogeneous problems with merely measurable coefficients. Before completing the decay theory for the measure-data problem, we isolate the harmonic approximation results needed in the super-dimensional regimes, focusing on homogeneous maps and their approximation by solutions of autonomous equations obtained by freezing $A$ at a minimum point of $a$ when $p>n$ in the $p$-phase or $q>n$ in the $(p,q)$-phase. These comparison estimates are kept in energy-normalized form because they serve here solely to recover the missing homogeneous decay for $H(x,Dv)$; once \eqref{decay.h.Du} is available in the full range \eqref{2-1/n}, we shall return to them and further estimate their right-hand sides in order to deduce the complete comparisons in terms of $u$ and $\mu$.

\medskip

We begin with a preliminary energy decay estimate for $u$ in the case $p>n$, which will ensure uniform control of the coefficients after energy normalization.

\begin{lemma}\label{energy.decay.supercritical}
Let $n<p\leq q$ and let $u$ be a local energy solution to \eqref{equation} under the assumptions \eqref{a.holder}, \eqref{bound.supercritical} and \eqref{A.weak}. Fix $\Omega_0\Subset\Omega$ as in \eqref{Omega0} and take $\eta_0$ as in the remark following Theorem~\ref{hi.correct.exponent}. There exists a constant $c$, depending only on $\data$ and $\gdata$, such that
\[
\int_{B_R(x_0)}H(x,Du)\dx\leq cR^{\frac{n\eta_\smallzero}{1+\eta_\smallzero}}
\]
for every $x_0\in\Omega_0$ and $0<R\leq\min\{1,d_{\Omega_0}/8\}$.
\end{lemma}
\begin{proof}
Let $v$ be the comparison map defined in \eqref{def.v}: since $q\geq p>n$ with $n\geq2$, it is possible to test the difference of the equations for $u$ and $v$ with the test function $u-v\in W^{1,H(\cdot)}_0(B_{R/2})$, so that we have the chain of inequalities
\begin{align*}
\mean{B_{R/2}}\bigl[|Du-Dv|^p+a(x)|Du-Dv|^q\bigr]\dx&\leq c\,\frac{|\mu|(B_R)}{R^n}\sup_{B_{R/2}}|u-v|\\
&\leq c\,\frac{|\mu|(B_R)}{R^{n-1}}\biggl(\mean{B_{R/2}}|Du-Dv|^p\dx\biggr)^{\frac1p}\\
&\leq\frac12\mean{B_{R/2}}|Du-Dv|^p\dx+c\biggl[\frac{|\mu|(B_R)}{R^{n-1}}\biggr]^{p'}
\end{align*}
from our assumptions \eqref{A.weak}, Morrey's and Young's inequalities. Hence
\begin{equation}\label{rough.compa.supercritical.energy}
\mean{B_{R/2}}\bigl[|Du-Dv|^p+a(x)|Du-Dv|^q\bigr]\dx\leq c\biggl[\frac{|\mu|(B_R)}{R^{n-1}}\biggr]^{p'},
\end{equation}
where $c$ depends only on $n,p,q$ and $\nu$ and multiplying by $|B_{R/2}|$ yields
\[
\int_{B_{R/2}}H(x,Du-Dv)\dx\leq cR^{\frac{p-n}{p-1}}\bigl[|\mu|(B_R)\bigr]^{p'}.
\]
Testing \eqref{def.v} with $v-u$, and appealing to the coercivity and growth assumptions in \eqref{A.weak} together with Young's inequality, we obtain
\begin{align*}
\nu\int_{B_{R/2}}H(x,Dv)\dx&\leq\int_{B_{R/2}}\langle A(x,Dv),Dv\rangle\dx=\int_{B_{R/2}}\langle A(x,Dv),Du\rangle\dx\\
&\leq L\int_{B_{R/2}}\bigl[|Dv|^{p-1}+a(x)|Dv|^{q-1}\bigr]|Du|\dx\\
&\leq\frac{\nu}{2}\int_{B_{R/2}}H(x,Dv)\dx+c\int_{B_{R/2}}H(x,Du)\dx,
\end{align*}
whence
\[
\int_{B_{R/2}}H(x,Dv)\dx\leq c\int_{B_{R/2}}H(x,Du)\dx,
\]
where $c$ depends only on $p,q,\nu$ and $L$. For every $\tau\in(0,1/4]$, H\"older's inequality and \eqref{hi.local}, applied with $2\eta_0$ in place of $\eta_0$, give
\begin{align*}
\int_{B_{\tau R}}H(x,Dv)\dx&\leq |B_{\tau R}|^{\frac{2\eta_0}{1+2\eta_0}}\biggl(\int_{B_{R/4}}\bigl[H(x,Dv)\bigr]^{1+2\eta_0}\dx\biggr)^{\frac1{1+2\eta_0}}\\
&\leq c\tau^{\frac{2n\eta_0}{1+2\eta_0}}\int_{B_{R/2}}H(x,Dv)\dx\leq c\tau^{\frac{2n\eta_0}{1+2\eta_0}}\int_{B_{R/2}}H(x,Du)\dx,
\end{align*}
where the last inequality follows from the preceding energy comparison and $c$ depends only on $\data$ and $\gdata$. Combining the preceding estimate with \eqref{rough.compa.supercritical.energy}, and employing $H(x,z_1+z_2)\leq 2^{q-1}[H(x,z_1)+H(x,z_2)]$, we obtain
\begin{align*}
\int_{B_{\tau R}}H(x,Du)\dx&\leq c\int_{B_{\tau R}}H(x,Dv)\dx+c\int_{B_{R/2}}H(x,Du-Dv)\dx\\
&\leq c\tau^{\frac{2n\eta_\smallzero}{1+2\eta_\smallzero}}\int_{B_R}H(x,Du)\dx+cR^{\frac{p-n}{p-1}}\bigl[|\mu|(\Omega)\bigr]^{p'},
\end{align*}
where $c$ depends only on $\data$ and $\gdata$, independently of $x_0,R$ and $\tau$. This fact lets us choose $\tau\equiv \tau(\data,\gdata)\in(0,1/4]$ so small that
\[
c\tau^{\frac{2n\eta_\smallzero}{1+2\eta_\smallzero}}\leq\frac12\tau^{\frac{n\eta_\smallzero}{1+\eta_\smallzero}},
\]
where $c$ is the constant in the preceding estimate; with this choice of $\tau$, we can apply Lemma~\ref{geometric.decay.iteration} to the non-decreasing function $\Phi(r)=\int_{B_r(x_0)}H(x,Du)\dx$, with $\theta=1/2$, $\sigma=\frac{n\eta_\smallzero}{1+\eta_\smallzero}$, $\beta=\frac{p-n}{p-1}$ and $A=c[|\mu|(\Omega)]^{p'}$. Since $\sigma<\beta$ by our choice of $\eta_0$, we obtain
\[
\int_{B_\rho}H(x,Du)\dx\leq c\Bigl(\frac{\rho}{R}\Bigr)^{\frac{n\eta_\smallzero}{1+\eta_\smallzero}}\biggl[\int_{B_R}H(x,Du)\dx+R^{\frac{p-n}{p-1}}\bigl[|\mu|(\Omega)\bigr]^{p'}\biggr]
\]
for every $0<\rho\leq R\leq\min\{1,d_{\Omega_0}/8\}$, with $c$ depending only on $\data$ and $\gdata$. Finally, taking $R=\min\{1,d_{\Omega_0}/8\}$, we have $B_R(x_0)\subset\Omega_1$ and $\int_{B_R}H(x,Du)\dx\leq\mathcal M$ by \eqref{energy.M}. This proves the asserted decay estimate, uniformly for $x_0\in\Omega_0$.
\end{proof}

As before, we distinguish according to whether the minimum of $a(\,\cdot\,)$ over $B_R$ is small or large at the relevant scale. We first consider $p>n$ and assume that \eqref{p.phase} holds with $\x=0$. We take $x_m\in\overline{B_R}$ such that
\begin{equation}\label{frozen.minimum}
a(x_m)=a_i(B_R):=\inf_{B_R}a(\,\cdot\,)
\end{equation}
and define the frozen autonomous field by setting for all $z\in \R^n\setminus\{0\}$
\begin{equation}\label{vf.frozen.ai}
A_i(z)=A(x_m,z)\qquad \text{and $A_i(0)=0$}. 
\end{equation}
The field $A_i$ is autonomous, with ellipticity and growth governed by
\begin{equation}\label{inf.smaller.H}
H_i(z):=\frac{|z|^p}p+a_i(B_R)\frac{|z|^q}q,\qquad h_i(z):=|z|^{p-1}+a_i(B_R)|z|^{q-1}.
\end{equation}
We shall apply Lemma~\ref{p-harm} to the frozen field $A_i$, with energy function $H_i$, after a suitable normalization. To this end, we take $\varphi\in C^\infty_c(B_{R/8})$ and use the weak formulation of \eqref{def.v} to write
\begin{align}\label{quasi.approx.p.harmonic}
\biggl|\mean{B_{R/8}}\langle A_i(Dv),D\varphi\rangle\dx\biggr|&=\biggl|\mean{B_{R/8}}\langle A_i(Dv)-A(x,Dv),D\varphi\rangle\dx\biggr| \notag\\
&\leq\mean{B_{R/8}}|A(x_m,Dv)-A(x,Dv)|\,|D\varphi|\dx\notag \\
&\leq{\|D\varphi\|}_{L^\infty(B_{R/8})}\mean{B_{R/8}}|A(x_m,Dv)-A(x,Dv)|\dx\notag\\
&\leq c(p,q,L){\|D\varphi\|}_{L^\infty(B_{R/8})}\biggl[\omega(R)\mean{B_{R/8}}h(x,Dv)\dx+2\mean{B_{R/8}}a(x)|Dv|^{q-1}\dx\biggr]:
\end{align}
in the last inequality we used the third condition in \eqref{A.strong}, the minimality of $a(x_m)$ and the fact $\omega(|x-x_m|)\leq(9/8)\omega(R)$ for every $x\in B_{R/8}$, which follows by concavity, since $x_m\in\overline{B_R}$. We estimate the second average as in the first line of \eqref{q.part.useful.later}, obtaining
\begin{align}\label{oscillation.q.estimate}
\mean{B_{R/8}} a(x)|Dv|^{q-1}\dx&\leq\biggl(\mean{B_{R/8}}a(x)|Dv|^{q-p}\dx\biggr)^{\frac1p}\biggl(\mean{B_{R/8}}a(x)|Dv|^q\dx\biggr)^{\frac1{p'}}\notag\\
&\leq c\,R^{\frac{\alpha-\x}{p}}\biggl(\mean{B_{R/8}}|Dv|^{p(1+\eta_0)}\dx\biggr)^{\frac{q-p}{p^\smalltwo(1+\eta_0)}}\biggl(\mean{B_{R/8}}H(x,Dv)\dx\biggr)^{\frac1{p'}}.
\end{align}
In the second inequality we used the phase assumption \eqref{p.phase}, H\"older's inequality with $q-p\leq p<p(1+\eta_0)$, and $a(x)|z|^q\leq qH(x,z)$. Here $\eta_0$ is the exponent in \eqref{hi.local}, which depends only on $\data$ and $\gdata$ by the choice of $E$ in Paragraph~\ref{choice.E}. For $p>n$ and $0<R\leq\min\{1,d_{\Omega_0}/8\}$, the higher integrability estimate \eqref{hi.local}, the energy comparison between $u$ and $v$ established above, and Lemma~\ref{energy.decay.supercritical} yield
\begin{align}\label{bound.hi.Dv.harmonic}
\biggl(\mean{B_{R/4}}\bigl[H(x,Dv)\bigr]^{1+\eta_\smallzero}\dx\biggr)^{\frac{q-p}{p^\smalltwo(1+\eta_\smallzero)}}&\leq c\biggl(\mean{B_{R/2}}H(x,Dv)\dx\biggr)^{\frac{q-p}{p^\smalltwo}}\notag\\
&\leq c\biggl(\frac1{R^n}\int_{B_{R/2}}H(x,Du)\dx\biggr)^{\frac{q-p}{p^\smalltwo}}\leq c\,R^{-\frac{n(q-p)}{p^\smalltwo(1+\eta_\smallzero)}},
\end{align}
where $c$ depends only on $\data$ and $\gdata$. Now we can easily estimate the first averaged integral, as by H\"older's inequality and the preceding estimate we have
\begin{equation}\label{decay.without.decay}
\mean{B_{R/8}}h(x,Dv)\dx\leq \biggl(\mean{B_{R/8}}{|Dv|}^p\dx\biggr)^\frac{p-1}p+\mean{B_{R/8}} a(x)|Dv|^{q-1}\dx\leq c\bigl(1+R^{\frac{s_\smallzero}{p}}\bigr)E_1^{p-1}\leq c\,E_1^{p-1} 
\end{equation}
for $c\equiv c(\data,\gdata,M)$, if we remember the definition of $s_0>0$ in \eqref{bound.with.x}, the fact $R\leq 1$ and we set
\begin{equation}\label{def.E1}
E_1= \biggl(\mean{B_{R/4}}H(x,Dv)\dx\biggr)^{\frac1p}.
\end{equation}
Inserting the last considerations into \eqref{quasi.approx.p.harmonic} we get 
\begin{equation}\label{quasi.approx.p.harmonic.final}
\biggl|\mean{B_{R/8}}\langle A_i(Dv),D\varphi\rangle\dx\biggr|\leq c(\data,\gdata,M){\|D\varphi\|}_{L^\infty(B_{R/8})}\Bigl[\omega(R)+R^{\frac{s_\smallzero}p}\Bigr]E_1^{p-1}.
\end{equation}
We now renormalize the previous estimate; after noticing that we can suppose without loss of generality $E_1>0$ (otherwise the estimate \eqref{funV.DvDw} below would be trivial, as $E_1=0$ would imply $v$ constant, and in this case $w$ would be the same constant), we define
\begin{equation}\label{rescaled.v.Ai}
\tilde v:=\frac{v}{E_1},\qquad\tilde a_i:=a_i(B_R)E_1^{q-p},\qquad\tilde A_i(z):=\frac{A_i(E_1z)}{h_i(E_1)},\qquad z\in\mathbb R^n,
\end{equation}
recalling the definition of $h_i(\,\cdot\,)$ in \eqref{inf.smaller.H} and, accordingly, we introduce the rescaled functions
\begin{equation}\label{rescaled.hi}
\tilde h_i(z):=\frac{h_i(E_1z)}{h_i(E_1)}=\frac{|z|^{p-1}+\tilde a_i |z|^{q-1}}{1+\tilde a_i} 
\end{equation}
as
\[
h_i(E_1)=E_1^{p-1}+a_i(B_R)E_1^{q-1}=E_1^{p-1}(1+a_i(B_R)E_1^{q-p})=E_1^{p-1}(1+\tilde a_i).
\]
By the $p$-phase condition \eqref{p.phase} and the estimates leading to \eqref{quasi.approx.p.harmonic.final}, we have
\begin{equation}\label{rescaled.ai.bound}
\tilde a_i\leq c(\data,\gdata,M)R^{s_{\smallzero}}\leq c(\data,\gdata,M), 
\end{equation}
as $R\leq1$. Indeed, recalling \eqref{def.E1} and the energy estimate used in \eqref{bound.hi.Dv.harmonic}, we have
\[
E_1^p=\mean{B_{R/4}}H(x,Dv)\dx\leq\frac{c}{R^n}\int_{B_{R/2}}H(x,Du)\dx\leq cR^{-\frac{n}{1+\eta_\smallzero}}.
\]
Consequently, by \eqref{p.phase}, \eqref{iccs} and \eqref{bound.with.x},
\[
\tilde a_i=a_i(B_R)E_1^{q-p}\leq M[a]_{0,\alpha}R^{\alpha-\x}E_1^{q-p}\leq cR^{\alpha-\x-\frac{n(q-p)}{p(1+\eta_\smallzero)}}\leq cR^{s_\smallzero},
\]
where $c$ depends only on $\data,\gdata$ and $M$. In particular, we have
\begin{equation}\label{h_i=p-1}
h_i(E_1)=E_1^{p-1}(1+\tilde a_i)\approx_{\data,\gdata,M}E_1^{p-1}. 
\end{equation}
Finally, we set
\begin{equation}\label{rescaled.Hi}
\tilde H_i(z):=\frac{1}{\kappa_i}\int_0^{|z|}\tilde h_i(\tau)\,d\tau\approx_{p,q}\frac{|z|^p+\tilde a_i|z|^q}{1+\tilde a_i}\qquad\text{with}\qquad\kappa_i=\int_0^1\frac{h_i(E_1\tau)}{h_i(E_1)}\,d\tau=\frac{H_i(E_1)}{E_1h_i(E_1)};
\end{equation}
thus $\tilde H_i(1)=1$ and $1/q\leq\kappa_i\leq1/p$. For $z\in\R^n\setminus\{0\}$ the structural assumptions on $A_i$ yield the following bounds for every $\xi\in\R^n$:
\begin{equation}\label{rescaling.pphase}
\begin{cases}
\displaystyle\big\langle\partial_z\tilde A_i(z)\xi,\xi\big\rangle\geq\nu\frac{\tilde h_i(|z|)}{|z|}|\xi|^2\\[2mm]
\displaystyle|\partial_z\tilde A_i(z)|\leq L\frac{\tilde h_i(|z|)}{|z|}
\end{cases};
\end{equation}
thus $\widetilde A_i$ satisfies the structural assumptions of Lemma~\ref{p-harm} relative to $\widetilde H_i$, with constants depending only on the original structural data. The remaining assumptions required in Lemma~\ref{p-harm} are also preserved by the same normalization: \eqref{quasi.approx.p.harmonic.final} together with \eqref{h_i=p-1} gives
\begin{align*}
\biggl|\mean{B_{R/8}}\big\langle\tilde A_i(D\tilde v),D\varphi\big\rangle\dx\biggr|&=\frac1{h_i(E_1)}\biggl|\mean{B_{R/8}}\langle A_i(Dv),D\varphi\rangle\dx\biggr|\\
&\leq c(\data,\gdata,M){\|D\varphi\|}_{L^\infty(B_{R/8})}\Bigl[\omega(R)+R^\frac{s_\smallzero}p\Bigr], 
\end{align*}
where $s_0>0$ is defined in \eqref{bound.with.x}; the preceding display gives \eqref{eqn:almost-p-harm}. To verify \eqref{eqn:v-int00} with $B=B_{R/8}$, we use the identities $\tilde H_i(D\tilde v)=H_i(Dv)/H_i(E_1)$ and $E_1^p=\mean{B_{R/4}}H(x,Dv)\dx$, together with $H_i(E_1)\geq E_1^p/p$, to obtain
\begin{equation}\label{tildev.bounded.1}
\mean{B_{R/4}}\tilde H_i(D\tilde v)\dx=\frac1{H_i(E_1)}\mean{B_{R/4}}H_i(Dv)\dx\leq\frac p{E_1^p}\mean{B_{R/4}}H(x,Dv)\dx=p.
\end{equation}
Finally, \eqref{eqn:v-int002} follows from the higher integrability estimate \eqref{hi.local}:
\begin{align}\label{tildev.bounded.2}
\mean{B_{R/8}}[\tilde H_i(D\tilde v)]^{1+\eta_0}\dx&\leq\frac{p^{1+\eta_0}}{E_1^{p(1+\eta_0)}}\mean{B_{R/8}}[H(x,Dv)]^{1+\eta_0}\dx\notag\\
&\leq\frac c{E_1^{p(1+\eta_0)}}\biggl(\mean{B_{R/4}}H(x,Dv)\dx\biggr)^{1+\eta_0}\leq c(\data,\gdata).
\end{align}
Thus Lemma~\ref{p-harm} guarantees the existence of a function $\tilde w\in \tilde v+W^{1,\tilde H_i(\cdot)}_0(B_{R/8})$, $\tilde A_i$-harmonic in $B_{R/8}$, such that 
\begin{equation}\label{ts:sob-dist-h-u}
\mean{B_{R/8}}{\big|V_{\tilde H_i}(D\tilde v)-V_{\tilde H_i}(D\tilde w)\big|}^2\dx \leq c{\Bigl[\omega(R)+R^\frac{s_\smallzero}p\Bigr]}^{\iota},
\end{equation}
where $V_{\tilde H_i}$ is defined by \eqref{V.varphi} with $\varphi=\tilde H_i$. Now we notice that by the definition of the normalization constant $\kappa_i$ in \eqref{rescaled.Hi} we have, for $t\geq 0$,
\[
\widetilde H_i(t)=\frac{1}{\kappa_i}\int_0^t\frac{h_i(E_1\tau)}{h_i(E_1)}\,d\tau=\frac{H_i(E_1t)}{H_i(E_1)}\qquad\implies\qquad\tilde H_i'(t)=\frac{E_1h_i(E_1t)}{H_i(E_1)};
\]
both maps vanish at the origin. For $z\in\mathbb R^n\setminus\{0\}$, \eqref{V.varphi} and the preceding formula for $\tilde H_i'$ give
\begin{equation}\label{V.Hi.scaling}
 V_{\tilde H_i}(z)=\sqrt{\frac{E_1h_i(E_1|z|)}{H_i(E_1)|z|}}\,z=\frac{1}{\sqrt{H_i(E_1)}}\,V_{H_i}(E_1z)
\end{equation}
so that, setting $w=E_1\tilde w\in v+W^{1,H_i(\cdot)}_0(B_{R/8})$, we obtain
\begin{equation}\label{V.Hi.scaling.difference}
{\big|V_{\tilde H_i}(D\tilde v)-V_{\tilde H_i}(D\tilde w)\big|}^2=\frac{1}{H_i(E_1)}{\big| V_{H_i}(Dv)-V_{H_i}(Dw)\big|}^2\approx_{\data,\gdata,M}\frac{1}{E_1^p}{\big| V_{H_i}(Dv)-V_{H_i}(Dw)\big|}^2
\end{equation}
using also \eqref{h_i=p-1} and \eqref{rescaled.Hi}; therefore \eqref{ts:sob-dist-h-u} translates into
\begin{align}\label{funV.DvDw}
\mean{B_{R/8}}{\big|V_{a(x_m)}(Dv)-V_{a(x_m)}(Dw)\big|}^2\dx &\leq c {\Bigl[\omega(R)+R^\frac{s_\smallzero}p\Bigr]}^{\iota}E_1^p\notag\\
&\leq c{\Bigl[\omega(R)+R^\frac{s_\smallzero}p\Bigr]}^{\iota}\mean{B_{R/4}}H(x,Dv)\dx,
\end{align}
with $c$ depending on $\data,\gdata$ and $M$, after recalling the definition of $E_1$ in \eqref{def.E1}. Moreover, testing the equation for $w$ with $w-v$ and using the structural assumptions on $A_i$ together with Young's inequality yields
\begin{equation}\label{energy.pphase}
\mean{B_{R/8}}H_i(Dw)\dx\leq c\mean{B_{R/8}}H_i(Dv)\dx\leq c\mean{B_{R/4}}H(x,Dv)\dx,
\end{equation}
where $c$ depends only on $n,p,q,\nu$ and $L$.

\medskip

\begin{remark}\label{frozen.p.preliminary}
{\rm Here we set $p>n$ and $\x=0$; from \eqref{oscillation.q.estimate}, the higher integrability \eqref{hi.local}, and Lemma~\ref{energy.decay.supercritical} we deduce \eqref{quasi.approx.p.harmonic.final}, and thus \eqref{funV.DvDw}. Notice that the same normalization is adopted in Proposition~\ref{decay.superdimensional}, checking the needed energy bounds on the corresponding scales. Later, using Proposition~\ref{decay.measure.data}, we will return to \eqref{quasi.approx.p.harmonic} to get sharper comparison estimates.}
\end{remark}

\medskip

Assume now that \eqref{pq.phase} holds with $M\geq2$. Since $R\leq1$ and $\x\geq0$, from \eqref{comparability.a} we have
\[
a_i(B_R)\leq a(x)\leq2a_i(B_R)\qquad\text{for every }x\in\overline{B_R},
\]
so that freezing $a(\,\cdot\,)$ at a minimum point gives
\begin{equation}\label{uniform.comparability.a}
H_i(z)\leq H(x,z)\leq2H_i(z),\qquad h_i(z)\leq h(x,z)\leq2h_i(z)
\end{equation}
for every $x\in\overline{B_R}$ and $z\in\mathbb R^n$. Note that these inequalities hold both when $\x=0$, as in the present preliminary step, and for the positive values of $\x$ considered later on.

\medskip

We are going to prove that there is an $A_i$-harmonic map $w\in v+W^{1,H_i(\cdot)}_0(B_{R/8})$ that quantitatively approximates $v$, with an error that can be made arbitrarily small by choosing $M$ large and $R$ small (this latter smallness requirement is solely due to the presence of the modulus $\omega(\,\cdot\,)$ and is absent in the model case $A_{\texttt d \texttt p}(\,\cdot\,,\,\cdot\,)$ in \eqref{model.A}, for instance). To this end, we again take a test function $\varphi\in C^\infty_c(B_{R/8})$ and proceed similarly to what was done in \eqref{quasi.approx.p.harmonic}: by our assumption \eqref{A.strong}$_3$ and \eqref{comparability.a} we have
\begin{align}\label{quasi.pq.harm}
\biggl|\mean{B_{R/8}}\langle A_i(Dv),D\varphi\rangle\dx\biggr|&\leq\mean{B_{R/8}}|A(x_m,Dv)-A(x,Dv)|\,|D\varphi|\dx\notag \\
&\leq c\,{\|D\varphi\|}_{L^\infty(B_{R/8})}\biggl[\omega(R)\mean{B_{R/8}}h(x_m,Dv)\dx+\mean{B_{R/8}}|a(x)-a(x_m)||Dv|^{q-1}\dx\biggr]\notag\\
&\leq c\,{\|D\varphi\|}_{L^\infty(B_{R/8})}\biggl[\omega(R)\mean{B_{R/8}}h(x_m,Dv)\dx+\frac{R^\x}M\mean{B_{R/8}}a(x_m)|Dv|^{q-1}\dx\biggr]\notag\\
&\leq c\,{\|D\varphi\|}_{L^\infty(B_{R/8})}\Bigl[\omega(R)+\frac{R^\x}M\Bigr]\bigl(h_{x_m}\circ H_{x_m}^{-1}\bigr)\biggl(\mean{B_{R/8}}H_{x_m}(Dv)\dx\biggr)
\end{align}
with $c$ that depends only on $\data$, uniformly for $M\geq2$, by \eqref{uniform.comparability.a} and the concavity of $\omega$; we also used the averaged estimate \eqref{averaged.Psi.inverse} with $\el=1$ and \eqref{constant.hinverse}. In this case we set 
\begin{equation}\label{def.E2}
E_2= H_{x_m}^{-1}\biggl(\mean{B_{R/4}}H_{x_m}(Dv)\dx\biggr) =H_i^{-1}\biggl(\mean{B_{R/4}}H_i(Dv)\dx\biggr) 
\end{equation}
by \eqref{frozen.minimum} (again supposing without loss of generality $E_2>0$) and we define the rescaled quantities $\tilde v, \tilde h_i, \tilde H_i, \tilde A_i$ as in \eqref{rescaled.v.Ai}, \eqref{rescaled.hi} and \eqref{rescaled.Hi} but replacing $E_1$ with $E_2$. Since $\tilde H_i(D\tilde v)=H_i(Dv)/H_i(E_2)$, the definition of $E_2$ gives
\[
\mean{B_{R/4}}\tilde H_i(D\tilde v)\dx=\frac1{H_i(E_2)}\mean{B_{R/4}}H_i(Dv)\dx=1.
\]
This verifies \eqref{eqn:v-int00} with $B=B_{R/8}$. Moreover, Theorem~\ref{hi.correct.exponent}, with the choice of Paragraph~\ref{choice.E}, and \eqref{uniform.comparability.a} yield
\begin{align}\label{p-frozen.comparision}
\mean{B_{R/8}} {\bigl[\tilde H_i(D\tilde v)\bigr]}^{1+\eta_0}\dx&\leq  \frac{1}{[H_i(E_2)]^{1+\eta_0}}\mean{B_{R/8}} {[H(x,Dv)]}^{1+\eta_0}\dx\notag\\
&\leq c \frac{1}{[H_i(E_2)]^{1+\eta_0}}\biggl(\mean{B_{R/4}} H(x,Dv)\dx\biggr)^{1+\eta_0}\notag\\
&\leq c \frac{1}{[H_i(E_2)]^{1+\eta_0}}\biggl(\mean{B_{R/4}} H_i(Dv)\dx\biggr)^{1+\eta_0}\leq c(\data,\gdata).
\end{align}
Moreover, \eqref{quasi.pq.harm} can be rewritten as
\begin{multline*}
\biggl|\mean{B_{R/8}}\langle A_i(Dv),D\varphi\rangle\dx\biggr|\leq c(\data)\Bigl[\omega(R)+\frac{R^\x}M\Bigr]{\|D\varphi\|}_{L^\infty(B_{R/8})}h_i(E_2)\\
\Longrightarrow\qquad \biggl|\mean{B_{R/8}}\langle \tilde A_i(D\tilde v),D\varphi\rangle\dx\biggr|\leq c(\data)\Bigl[\omega(R)+\frac{R^\x}M\Bigr]{\|D\varphi\|}_{L^\infty(B_{R/8})},
\end{multline*}
due to the scaling property \eqref{g_x.lambda} and thus we can apply Lemma~\ref{p-harm} with $B=B_{R/8}$, $\varphi=\tilde H_i$ and $A_0=\tilde A_i$ satisfying \eqref{rescaling.pphase} and $\sigma=c\,[\omega(R)+R^\x/M]$ obtaining therefore the existence of a $\tilde A_i$-harmonic map $\tilde w\in \tilde v+W^{1,\tilde H_i(\cdot)}_0(B_{R/8})$ such that
\begin{equation}\label{pq.harm.rescaled}
\mean{B_{R/8}}\tilde H_i(D\tilde w)\dx\leq c(\data), \qquad \mean{B_{R/8}}{\big|V_{\tilde H_i}(D\tilde v)-V_{\tilde H_i}(D\tilde w)\big|}^2\dx \leq c{\Bigl[\omega(R)+\frac{R^\x}M\Bigr]}^{\iota}
\end{equation}
with $c\equiv c(\data,\gdata)$; if we set $w=\tilde w E_2\in v+W^{1,H_i(\cdot)}_0(B_{R/8})$, using the exact identity \eqref{V.Hi.scaling} with $E_2$ in place of $E_1$ to scale back, we deduce
\begin{equation}\label{pq-frozen.energy}
 \mean{B_{R/8}}H_i(Dw)\dx\leq c(\data) \mean{B_{R/4}}H_i(Dv)\dx
\end{equation}
and
\begin{equation}\label{pq-frozen.comparision}
\mean{B_{R/8}}{\big|V_{H_i}(Dv)-V_{H_i}(Dw)\big|}^2\dx \leq c{\Bigl[\omega(R)+\frac{R^\x}M\Bigr]}^{\iota}\mean{B_{R/4}}H_{x_m}(Dv)\dx 
\end{equation}
with $c,\iota$ depending on $\data$ and $\gdata$, but not on $M$.

\subsection{Morrey decay at the energy level - the superdimensional case} 

Before tackling the decay estimates for solutions to \eqref{equation} under \eqref{A.strong}, we establish their counterparts for homogeneous comparison maps in the super-dimensional regimes - that is, when $p>n$ in the $p$-phase or $q>n$ in the $(p,q)$-phase. Under analogous conditions, the methods behind \cite[Theorem~2]{BCMgen} yield, in our notation, bounds of the type
\[
\mean{B_{\rho_\smallone}(x_1)}H(x,Dv)\dx \leq c\Bigl(\frac{\rho_1}{\rho_2}\Bigr)^{-\varsigma} \mean{B_{\rho_\smalltwo}(x_1)}H(x,Dv)\dx
\]
for every $\varsigma\in(0,n]$. We combine the comparison estimates of Paragraph~\ref{par.harmonic} with the interior Lipschitz bound for the frozen comparison map $w$ to obtain a one-step energy estimate for $v$, which we then iterate. This approach makes the dependence on $\benergy(Dv,B_R(x_0))$ explicit. Since $E$ has already been fixed in terms of $\data$ and $\gdata$, it introduces no additional dependence. 

\medskip 

\begin{proposition}\label{decay.superdimensional}
Let $v\in W^{1,H(\cdot)}(B_R)$ solve \eqref{homogeneous.equation} under \eqref{a.holder} to \eqref{bounds} and \eqref{A.strong}, and assume that the energy bound is chosen as in Paragraph~\ref{choice.E}. Let $B_{\rho_\smallone},B_{\rho_\smalltwo}$ be two balls centered at $x_1$ as in \eqref{double.ball.inclusion}. For every $\varsigma\in(0,n]$, there exist constants $M_\varsigma\geq2$ and $c$, depending on $\data,\gdata,\varsigma$, and a radius $R_\varsigma\ll1$, depending on $\data,\gdata,\varsigma$ and $\omega(\,\cdot\,)$, such that the following holds. If $R\leq R_\varsigma$ and either $p>n$ with
\begin{equation}\label{superdimensional.pphase}
\inf_{B_{\rho_\smallone}}\, a(\,\cdot\,) \leq M_\varsigma{[a]}_{0, \alpha} \rho_1^\alpha,
\end{equation}
or $q>n$ and 
\begin{equation}\label{superdimensional.pqphase}
\inf_{B_{\rho_\smallthree}}\, a(\,\cdot\,) > M_\varsigma{[a]}_{0, \alpha} \rho_3^\alpha
\end{equation}
for some ball $B_{\rho_\smallthree}$ such that $B_{\rho_2}(x_1)\subset B_{\rho_3}(x_1)\subset B_R(x_0)$, then
\begin{equation}\label{decay.energy.superdimensional}
\mean{B_{\rho_\smallone}(x_1)}H(x,Dv)\dx\leq c \Bigl(\frac{\rho_1}{\rho_2}\Bigr)^{-\varsigma}\mean{B_{\rho_\smalltwo}(x_1)}H(x,Dv)\dx 
\end{equation}
for a constant depending on $\data,\gdata$ and $\varsigma$.
\end{proposition}

\begin{proof}
Fix $\varsigma\in(0,n]$. We first consider the case $q>n$, assuming that \eqref{pq.phase} holds on $B_{\rho_\smallthree}$ with $\x=0$ and a phase constant $M\geq2$ to be chosen. By Remark~\ref{stability.phases}, the infimum formulation of the $(p,q)$-phase condition is preserved under restriction with the same constant; hence, for every $0<\varrho\leq\rho_2$,
\begin{equation}\label{pq.phase.generic.scale}
a_i(B_\varrho)>M[a]_{0,\alpha}\varrho^\alpha.
\end{equation}
For each such $\varrho$, we choose $x_{m,\varrho}\in\overline{B_\varrho(x_1)}$ such that $a(x_{m,\varrho})=a_i(B_\varrho)$ and set
\begin{equation}\label{frozen.iterated}
a_{i,\varrho}:=a(x_{m,\varrho}),\qquad A_i(z):=A(x_{m,\varrho},z),\qquad H_i(z):=H(x_{m,\varrho},z).
\end{equation}
By \eqref{pq-frozen.comparision}, there exists an $A_i$-harmonic map $w\in v+W^{1,H_i}_0(B_{\varrho/8}(x_1))$ such that
\begin{align}\label{pq.comparison.generic.scale}
\mean{B_{\varrho/8}}\big|V_{a_{i,\varrho}}(Dv)-V_{a_{i,\varrho}}(Dw)\big|^2\dx &\leq  c{\Bigl[\omega(\varrho)+\frac1M\Bigr]}^{\iota}\mean{B_{\varrho/2}}H_i(Dv)\dx\notag\\
& =: c\,{\mathcal S}_{p,q}(\varrho,M)\mean{B_{\varrho/2}}H_i(Dv)\dx,
\end{align}
where ${\mathcal S}_{p,q}(\varrho,M):=[\omega(\varrho)+1/M]^\iota$ and $c,\iota$ depend only on $\data$ and $\gdata$. We recall that the approximation result is being applied here with $\x=0$. Using the triangle inequality and \eqref{V-phi}, applied with the choice $\varphi=H_i$, we get
\begin{equation}\label{elementary.Hi.V}
H_i(z_1) \leq 2H_i(z_2)+2\big|V_{a_{i,\varrho}}(z_1)-V_{a_{i,\varrho}}(z_2)\big|^2
\end{equation}
for all $z_1,z_2\in\mathbb R^n$. Consequently, for every $0<r\leq\varrho/16$, appealing to the uniform comparability \eqref{uniform.comparability.a} (which in turn follows from \eqref{pq.phase.generic.scale}), the previous estimate \eqref{elementary.Hi.V} and \eqref{pq.comparison.generic.scale}, we find
\begin{align}\label{preliminary.decay.pq.super}
\mean{B_r}H(x,Dv)\dx\leq c\mean{B_r}H_i(Dv)\dx&\leq c\mean{B_r}H_i(Dw)\dx+c\mean{B_r}\big|V_{a_{i,\varrho}}(Dv)-V_{a_{i,\varrho}}(Dw)\big|^2\dx \notag\\
&\leq c\sup_{B_{\varrho/16}}H_i(Dw)+c\left(\frac{\varrho}{r}\right)^n{\mathcal S}_{p,q}(\varrho,M)\mean{B_{\varrho/2}}H_i(Dv)\dx .
\end{align}
Since $w$ is $A_i$-harmonic, the Lipschitz estimate \eqref{lip.varphi.harmonic} and the energy comparison \eqref{pq-frozen.energy} between $w$ and $v$ yield
\[
\sup_{B_{\varrho/16}}H_i(Dw)\leq c\mean{B_{\varrho/8}}H_i(Dw)\dx\leq c\mean{B_{\varrho/4}}H_i(Dv)\dx\leq c\mean{B_{\varrho/2}}H(x,Dv)\dx .
\]
Inserting this estimate into \eqref{preliminary.decay.pq.super} and employing once more $H_i\leq H$ and enlarging the domain on the right-hand side gives
\begin{equation}\label{one.step.pq.super}
\mean{B_r}H(x,Dv)\dx\leq \bar c\Bigl[1+\left(\frac{\varrho}{r}\right)^n{\mathcal S}_{p,q}(\varrho,M)\Bigr]\mean{B_\varrho}H(x,Dv)\dx
\end{equation}
for all $0<r\leq\varrho/16$ and with a constant $\bar c \equiv \bar c(\data,\gdata)$. We now fix $\theta\in(0,1/16]$, depending only on $\data,\gdata$ and $\varsigma$, so small that
\begin{equation}\label{choice.theta.pq.super}
\bar c\leq \frac12\theta^{-\varsigma};
\end{equation}
then we first choose $M_\varsigma\geq4$ sufficiently large and then $R_1\in(0,1]$ sufficiently small so that
\begin{equation}\label{choice.M.R.pq.super}
\bar c\,\theta^{-n}{\mathcal S}_{p,q}(R_1,M_\varsigma)=\bar c\,\theta^{-n}{\Bigl[\omega(R_1)+\frac{1}{M_\varsigma}\Bigr]}^{\iota}\leq \frac12\theta^{-\varsigma}.
\end{equation}
Consequently, whenever $\varrho\leq R_1$, for the choice of $M_\varsigma$ made above, estimate \eqref{one.step.pq.super}, applied with $r=\theta\varrho$, yields
\begin{equation}\label{one.step.pq.super.tau}
\mean{B_{\theta \varrho}}H(x,Dv)\dx\leq\theta^{-\varsigma}\mean{B_\varrho}H(x,Dv)\dx \qquad\text{for every $0<\varrho\leq\rho_2$.}
\end{equation}
Multiplying \eqref{one.step.pq.super.tau} by $|B_{\theta\varrho}|$, we obtain $\Phi(\theta\varrho)\leq\theta^{n-\varsigma}\Phi(\varrho)$ for $\Phi(r)$ defined by the integral over $B_r(x_1)$ of $H(x,Dv)$. If $\rho_2\leq R_1$, the homogeneous case of Lemma~\ref{geometric.decay.iteration}, applied between $\rho_1$ and $\rho_2$ with $\sigma=n-\varsigma$, gives, after division by $|B_{\rho_1}|$,
\begin{equation}\label{partial.decay.1}
\mean{B_{\rho_1}}H(x,Dv)\dx\leq\theta^{\varsigma-n}\Bigl(\frac{\rho_1}{\rho_2}\Bigr)^{-\varsigma}\mean{B_{\rho_2}}H(x,Dv)\dx.
\end{equation}
We now turn to the regime $p>n$ under \eqref{superdimensional.pphase}, setting $M=M_\varsigma$. Throughout the remainder of the argument, the constants $c$ depend on $\data,\gdata$ and $M$. For any radius $t\in[\rho_1,\rho_2]$, Remark~\ref{stability.p.phase} guarantees that on the ball $B_t$ the problem is in the $p$-phase, with the supremum of $a(\,\cdot\,)$ being at most $(M+2)[a]_{0,\alpha}t^\alpha$, so that we are allowed  to invoke \eqref{nice.reverse.3} on $B_{\rho_2}$ with phase constant $M+2$, yielding
\[
\mean{B_{\rho_2/2}}H(x,Dv)\dx\leq c\biggl(\mean{B_{\rho_2}}|Dv|^{p-1}\dx\biggr)^{p'}\leq c\,\rho_2^{-n},
\]
where the last inequality reflects the energy bound $\|Dv\|_{L^p(B_R)}\leq E$. An application of the higher integrability estimate \eqref{hi.local} on $B_{\rho_2/2}$ therefore secures the reference bound at the fixed scale $\rho_2$:
\[
\biggl(\int_{B_{\rho_2/4}}[H(x,Dv)]^{1+\eta_0}\dx\biggr)^{\frac1{1+\eta_\smallzero}}\leq c\,\rho_2^{-\frac{n\eta_\smallzero}{1+\eta_\smallzero}}.
\]
Let now $\varrho\in[\rho_1,\rho_2/4]$ be an intermediate radius, and choose a point $x_{m,\varrho}\in\overline{B_\varrho(x_1)}$ realizing the minimum $a(x_{m,\varrho})=a_i(B_\varrho)$. Once $a_{i,\varrho},A_i$, and $H_i$ are defined as in \eqref{frozen.iterated}, we set
\[
E_\varrho:=\biggl(\mean{B_{\varrho/4}}H(x,Dv)\dx\biggr)^{\frac1p}.
\]
We may assume $E_\varrho>0$ without loss of generality - otherwise $v$ is constant on $B_{\varrho/4}$ and the comparison below holds trivially with $w=v$. Combining Hölder's inequality with the integral bound at scale $\rho_2/4$ immediately gives
\[
E_\varrho^p+\biggl(\mean{B_{\varrho/8}}[H(x,Dv)]^{1+\eta_0}\dx\biggr)^{\frac1{1+\eta_\smallzero}}\leq c\varrho^{-\frac n{1+\eta_\smallzero}}\rho_2^{-\frac{n\eta_\smallzero}{1+\eta_\smallzero}}
\]
and at the level of the switching coefficient, this translates into the decay estimate
\[
\sup_{B_\varrho}a(\,\cdot\,)E_\varrho^{q-p}\leq c\varrho^{\alpha-\frac{n(q-p)}{p(1+\eta_\smallzero)}}\rho_2^{-\frac{n\eta_\smallzero(q-p)}{p(1+\eta_\smallzero)}}=c\Bigl(\frac{\varrho}{\rho_2}\Bigr)^\gamma\rho_2^{\alpha-\frac{n(q-p)}p}\leq c\Bigl(\frac{\varrho}{\rho_2}\Bigr)^\gamma.
\]
Here $\gamma:=\alpha-n(q-p)/[p(1+\eta_0)]=2s_0>0$, with $s_0$ as in \eqref{bound.with.x}. The last inequality follows from \eqref{bound.supercritical}, since $\alpha-n(q-p)/p\geq0$ and $\rho_2\leq1$. By the same token, repeating the Hölder argument of \eqref{oscillation.q.estimate} alongside our higher integrability bound leads to
\[
\mean{B_{\varrho/8}}a(x)|Dv|^{q-1}\dx\leq c\Bigl(\frac{\varrho}{\rho_2}\Bigr)^{\frac\gamma p}E_\varrho^{p-1},\qquad\mean{B_{\varrho/8}}h(x,Dv)\dx\leq cE_\varrho^{p-1}.
\]
Inserting these estimates into the computation for \eqref{quasi.approx.p.harmonic} - now carried out at scale $\varrho$ - leads to
\[
\biggl|\mean{B_{\varrho/8}}\langle A_i(Dv),D\varphi\rangle\dx\biggr|\leq c\Bigl[\omega(\varrho)+\Bigl(\frac{\varrho}{\rho_2}\Bigr)^{\frac{s_\smallzero}p}\Bigr]E_\varrho^{p-1}\|D\varphi\|_{L^\infty(B_{\varrho/8})}
\]
for every $\varphi\in C_c^\infty(B_{\varrho/8})$, since $\gamma\geq s_0$ and $\varrho/\rho_2\leq1$. To conclude, we rescale as in \eqref{rescaled.v.Ai}-\eqref{rescaled.Hi} with $E_\varrho$ in place of $E_1$:  we still have the equivalences 
\[
h_i(E_\varrho)\approx_{\data,\gdata,M}E_\varrho^{p-1}\qquad\text{and}\qquad H_i(E_\varrho)\approx_{\data,\gdata,M}E_\varrho^p
\]
because of the coefficient bound that forces $a_{i,\varrho}E_\varrho^{q-p}\leq c$, and the definition of $E_\varrho$ and an application of \eqref{hi.local} between $B_{\varrho/8}$ and $B_{\varrho/4}$ further ensure that
\[
\mean{B_{\varrho/4}}\tilde H_i(D\tilde v)\dx\leq p,\qquad\mean{B_{\varrho/8}}[\tilde H_i(D\tilde v)]^{1+\eta_0}\dx\leq c(\data,\gdata).
\]
Therefore we can again apply Lemma~\ref{p-harm} on $B=B_{\varrho/8}$ with $A_0=\tilde A_i$ and this  yields an $A_i$-harmonic map $w\in v+W^{1,H_i}_0(B_{\varrho/8}(x_1))$ satisfying
\[
\mean{B_{\varrho/8}}\big|V_{a_{i,\varrho}}(Dv)-V_{a_{i,\varrho}}(Dw)\big|^2\dx\leq c\,{\mathcal S}_p(\varrho)\mean{B_{\varrho/4}}H(x,Dv)\dx
\]
after having rescaled back to the original variables via \eqref{V.Hi.scaling.difference}, where we have set
\[
{\mathcal S}_p(\varrho):=\Bigl[\omega(\varrho)+\Bigl(\frac{\varrho}{\rho_2}\Bigr)^{\frac{s_\smallzero}p}\Bigr]^\iota.
\]
Keep in mind that $\rho_2$ remains fixed here; the constants $c=c(\data,\gdata,M)$ and $\iota=\iota(\data,\gdata)\in(0,1)$ are independent of both $\varrho$ and $\rho_2$. Let $\rho_1\leq r\leq\varrho/32$: by \eqref{superdimensional.pphase} and \eqref{oscillation.a}, we have
\[
\sup_{B_{2r}}a(\,\cdot\,)\leq\inf_{B_{\rho_1}}a(\,\cdot\,)+2[a]_{0,\alpha}(2r)^\alpha\leq(M+2)[a]_{0,\alpha}(2r)^\alpha.
\]
This allows us to apply \eqref{nice.reverse.3} on $B_{2r}$ with $\theta=1/2$ and phase constant $M+2$, together with H\"older's and then the triangle inequality, to estimate
\[
\mean{B_r}H(x,Dv)\dx\leq c\mean{B_{2r}}|Dv|^p\dx\leq c\sup_{B_{\varrho/16}}|Dw|^p+c\mean{B_{2r}}|Dv-Dw|^p\dx.
\]
Now we estimate the first term using the Lipschitz estimate \eqref{lip.varphi.harmonic}, while for the second one we can use \eqref{diffV2} with $\esp=p$ and \eqref{Vx0} with $x=x_{m,\varrho}$, since $p>n\geq2$:
\begin{align*}
\mean{B_r}H(x,Dv)\dx&\leq c\mean{B_{\varrho/8}}H_i(Dw)\dx+c\mean{B_{2r}}\big|V_{a_{i,\varrho}}(Dv)-V_{a_{i,\varrho}}(Dw)\big|^2\dx\\
&\leq c\left[1+\left(\frac{\varrho}{r}\right)^n{\mathcal S}_p(\varrho)\right]\mean{B_{\varrho/4}}H(x,Dv)\dx.
\end{align*}
In order to deduce the inequality in the last line, we kept in mind that $B_{2r}\subset B_{\varrho/16}$, the bound for the frozen energy in \eqref{energy.pphase} and the comparison estimate obtained above. We therefore have obtained that
\begin{equation}\label{one.step.p.super}
\mean{B_r}H(x,Dv)\dx\leq\tilde c\Bigl[1+\left(\frac{\varrho}{r}\right)^n{\mathcal S}_p(\varrho)\Bigr]\mean{B_\varrho}H(x,Dv)\dx
\end{equation}
for every $r$ with $\rho_1\leq r\leq\varrho/32$, where $\tilde c$ depends only on $\data,\gdata$ and $M$. Since $M=M_\varsigma$ has already been fixed in terms of $\data,\gdata$ and $\varsigma$, we choose $\vartheta\in(0,1/32]$, depending only on these quantities, so small that $\tilde c\leq\vartheta^{-\varsigma}/2$. We next choose $\kappa\in(0,1/4]$, depending only on $\data,\gdata$ and $\varsigma$, and then $R_2\in(0,1]$, additionally depending on $\omega(\,\cdot\,)$, such that
\begin{equation}\label{choice.M.R.p.super}
\tilde c\,\vartheta^{-n}\kappa^{\frac{s_\smallzero}p\iota}\leq\frac14\vartheta^{-\varsigma},\qquad\tilde c\,\vartheta^{-n}[\omega(R_2)]^\iota\leq\frac14\vartheta^{-\varsigma}.
\end{equation}
Since $0<\iota<1$, these choices ensure that
\[
\tilde c\,\vartheta^{-n}{\mathcal S}_p(\varrho)\leq\tilde c\,\vartheta^{-n}\Bigl([\omega(R_2)]^\iota+\kappa^{\frac{s_\smallzero}p\iota}\Bigr)\leq\frac12\vartheta^{-\varsigma}
\]
whenever $\varrho\leq\min\{\kappa\rho_2,R_2\}$. Therefore, choosing $r=\vartheta\varrho$ in \eqref{one.step.p.super} gives
\begin{equation}\label{one.step.p.super.tau}
\mean{B_{\vartheta\varrho}}H(x,Dv)\dx\leq\vartheta^{-\varsigma}\mean{B_\varrho}H(x,Dv)\dx
\end{equation}
whenever $\rho_1\leq\vartheta\varrho$ and $\varrho\leq\min\{\kappa\rho_2,R_2\}$. Assume $\rho_2\leq R_2$. If $\rho_1\leq\kappa\rho_2$, multiplying \eqref{one.step.p.super.tau} by $|B_{\vartheta\varrho}|$ gives $\Phi(\vartheta\varrho)\leq\vartheta^{n-\varsigma}\Phi(\varrho)$ whenever $\rho_1\leq\vartheta\varrho$ and $\varrho\leq\kappa\rho_2$. The homogeneous case of Lemma~\ref{geometric.decay.iteration}, applied with $r_*=\rho_1$ and $R=\kappa\rho_2$, gives
\[
\Phi(\rho_1)\leq\vartheta^{\varsigma-n}\Bigl(\frac{\rho_1}{\kappa\rho_2}\Bigr)^{n-\varsigma}\Phi(\kappa\rho_2).
\]
Using $\Phi(\kappa\rho_2)\leq\Phi(\rho_2)$ and dividing by $|B_{\rho_1}|$, we obtain
\begin{equation}\label{partial.decay.2}
\mean{B_{\rho_1}}H(x,Dv)\dx\leq(\vartheta\kappa)^{\varsigma-n}\Bigl(\frac{\rho_1}{\rho_2}\Bigr)^{-\varsigma}\mean{B_{\rho_2}}H(x,Dv)\dx.
\end{equation}
If $\rho_1>\kappa\rho_2$, the same estimate follows directly from $\Phi(\rho_1)\leq\Phi(\rho_2)$ and $(\rho_2/\rho_1)^{n-\varsigma}\leq\kappa^{\varsigma-n}$. Finally, we conclude by setting $R_\varsigma=\min\{R_1,R_2\}$ and taking the constant in \eqref{decay.energy.superdimensional} as the maximum between the constants in \eqref{partial.decay.1} and \eqref{partial.decay.2}: this completes the proof, since the conditions \eqref{choice.M.R.pq.super} and \eqref{choice.M.R.p.super} are both satisfied by $R_\varsigma$.
\end{proof}

\subsection{Morrey decay below the energy level}

Our next goal is the flux decay estimate of Proposition~\ref{decay.measure.data}, which we prove for local energy solutions and extend to local SOLAs by approximation in Section~\ref{proofs}. To this end, we need to lower the integrability exponent in the estimates obtained above, passing from the natural energy quantities to those arising in the SOLA comparison arguments; we therefore first explicitly record in Lemmas~\ref{decay.1} and~\ref{decay.2} the resulting two-scale estimates, corresponding to the two possible phase configurations. In both lemmas, the balls $B_{\rho_\smallone}, B_{\rho_\smalltwo}, B_{\rho_\smallthree}$ are concentric with center $x_1$ and satisfy \eqref{double.ball.inclusion} and $\alpha_0$ denotes the decay exponent provided by Proposition~\ref{decay.energy.level}. 

\begin{lemma}[Morrey estimate below the natural exponent in the $p$-phase]\label{decay.1}
Let $v\in W^{1,H(\cdot)}(B_R(x_0))$ be a solution to \eqref{homogeneous.equation} under the assumptions \eqref{2-1/n}, \eqref{a.holder} to \eqref{bounds} and the structural assumptions \eqref{A.weak} if $p\leq n$, \eqref{A.strong} if $p>n$, let $B_{\rho_\smallone}, B_{\rho_\smalltwo}$ be as before and suppose moreover that $\rho_1\leq \rho_2/2$. Suppose that there exists a constant $M\geq 2$ such that 
\begin{equation}\label{p.phase.sup.local.help.1}
\inf_{B_{\rho_\smallone}}\, a(\,\cdot\,) \leq M{[a]}_{0, \alpha} \rho_1^\alpha:
\end{equation}
\begin{itemize}
 \item if $p\leq n$, then there exists a constant $c_1\geq1$, depending on $\data,\gdata$ and $M$, such that the following two-scale decay estimate holds true:
\begin{multline}\label{decay.v.h}
\mean{B_{\rho_\smallone}(x_1)}|Dv|^{b(p-1)}\dx+{\biggl(\mean{B_{\rho_\smallone}(x_1)}a(x)|Dv|^{q-1}\dx\biggr)}^b\\
\leq c_1\Bigl(\frac{\rho_1}{\rho_2}\Bigr)^{-b(1-\alpha_\smallzero)(p-1)}\mean{B_{\rho_\smalltwo}(x_1)}|Dv|^{b(p-1)}\dx,
\end{multline}
where $b\in(1,n')$ is defined in \eqref{def.b}-\eqref{bbi};
\item if $p>n$, for every $\varsigma\in(0,n]$, let $M_\varsigma$ and $R_\varsigma$ be the constants provided by Proposition~\ref{decay.superdimensional}. If
\[
M\leq M_\varsigma\qquad\text{and}\qquad R\leq R_\varsigma,
\]
then there exists a constant $c_{1,\varsigma}\equiv c_{1,\varsigma}(\data,\gdata,M,\varsigma)\geq1$ such that
\begin{equation}\label{decay.v.h.p>n}
\mean{B_{\rho_\smallone/4}(x_1)}h(x,Dv)\dx \leq c_{1,\varsigma}\Bigl(\frac{\rho_1}{\rho_2}\Bigr)^{-\frac{\varsigma}{p'}}\mean{B_{\rho_\smalltwo}(x_1)}|Dv|^{p-1}\dx.
\end{equation}
\end{itemize}
\end{lemma}

\begin{proof}
We note that by \eqref{p.phase.sup.local.help.1} and H\"older continuity, see Remark~\ref{stability.p.phase}, for every $\rho_1\leq s\leq\rho_2$ we have $\sup_{B_s}a(\,\cdot\,)\leq \inf_{B_{\rho_1}}a(\,\cdot\,)+2[a]_{0,\alpha}s^\alpha\leq 2M[a]_{0,\alpha}s^\alpha$; all estimates requiring the $p$-phase condition in terms of the supremum apply with $2M$ in place of $M$.

\subsubsection*{Case $p\leq n$} We start from the inequality \eqref{directnice.reverse.2}, which we combine with  \eqref{decay.energy.p.est} in order to estimate its right-hand side, as $\rho_1\leq \rho_2/2$. We also use the reverse inequality \eqref{nice.reverse.2} to conclude:
\begin{align}\label{second.second}
\mean{B_{\rho_\smallone}}|Dv|^{b(p-1)}\dx+\biggl(\mean{B_{\rho_\smallone}}a(x)|Dv|^{q-1}\dx\biggr)^b&\leq c\biggl(\mean{B_{\rho_\smallone}}H(x,Dv)\dx\biggr)^{\frac{b}{p'}}\notag\\
&\leq c\Bigl(\frac{\rho_1}{\rho_2}\Bigr)^{-b(1-\alpha_\smallzero)(p-1)}\biggl(\mean{B_{\rho_\smalltwo/2}}H(x,Dv)\dx\biggr)^{\frac{b}{p'}}\notag\\
&\leq c\Bigl(\frac{\rho_1}{\rho_2}\Bigr)^{-b(1-\alpha_\smallzero)(p-1)}\mean{B_{\rho_\smalltwo}}|Dv|^{b(p-1)}\dx,
\end{align}
recalling that now $E$ has been fixed in Paragraph~\ref{choice.E} as a constant depending on $\data$ and $\gdata$.

\subsubsection*{Case $p>n$} 
We apply \eqref{directnice.reverse.3} on $B_{\rho_\smallone}$, where the supremum condition holds with constant $2M$ and then \eqref{decay.energy.superdimensional} between the radii $\rho_1$ and $\rho_2/2$, thanks to the infimum bound \eqref{p.phase.sup.local.help.1}, to get
\begin{align*}
\mean{B_{\rho_\smallone/4}}h(x,Dv)\dx\leq c\biggl(\mean{B_{\rho_\smallone}}H(x,Dv)\dx\biggr)^{\frac1{p'}}&\leq c\Bigl(\frac{\rho_1}{\rho_2}\Bigr)^{-\frac{\varsigma}{p'}}\biggl(\mean{B_{\rho_\smalltwo/2}}H(x,Dv)\dx\biggr)^{\frac1{p'}}\\
&\leq c\Bigl(\frac{\rho_1}{\rho_2}\Bigr)^{-\frac{\varsigma}{p'}}\mean{B_{\rho_\smalltwo}}|Dv|^{p-1}\dx,
\end{align*}
where the last inequality follows from \eqref{nice.reverse.3} on $B_{\rho_\smalltwo}$ with $\theta=1/2$, and $c$ depends only on $\data,\gdata,M$ and $\varsigma$.
\end{proof}

\begin{lemma}[Morrey estimate below the natural exponent in the $(p,q)$ phase]\label{decay.2}
Let $v\in W^{1,H(\cdot)}(B_R(x_0))$ be a solution to \eqref{homogeneous.equation} under the assumptions \eqref{2-1/n}, \eqref{a.holder} to \eqref{bounds} and \eqref{A.weak} if $q\leq n$, \eqref{A.strong} if $q>n$. Moreover, suppose that there exists a constant $M>2$ such that 
\begin{equation}\label{p.phase.sup.local.help.2}
\inf_{B_{\rho_\smallthree}}\, a(\,\cdot\,) > M{[a]}_{0, \alpha} \rho_3^\alpha.
\end{equation}
Then
\begin{itemize}
\item if $q\leq n$, there exists a constant $c_2\geq1$, depending only on $\data,\gdata$ and $M$, such that the decay estimate
\begin{equation}\label{decay.v.h.pqleqn}
\mean{B_{\rho_\smallone}(x_1)}h(x,Dv)\dx\leq c_2\Bigl(\frac{\rho_1}{\rho_2}\Bigr)^{-(1-\alpha_0)(q-1)}\mean{B_{\rho_\smalltwo}(x_1)}h(x,Dv)\dx
\end{equation}
holds for every pair of concentric balls as in \eqref{double.ball.inclusion}; 
\item if $q>n$, for every $\varsigma\in(0,n]$, let $M_\varsigma$ and $R_\varsigma$ be the constants provided by Proposition~\ref{decay.superdimensional}. If
\[
M\geq M_\varsigma \qquad\text{and}\qquad R\leq R_\varsigma,
\]
then there exists a constant $c_{2,\varsigma}\equiv c_{2,\varsigma}(\data,\gdata,M,\varsigma)\geq1$ such that
\begin{equation}\label{decay.v.h.pq>n}
\mean{B_{\rho_\smallone}(x_1)}h(x,Dv)\dx\leq c_{2,\varsigma}\Bigl(\frac{\rho_1}{\rho_2}\Bigr)^{-\frac{\varsigma}{q'}}\mean{B_{\rho_\smalltwo}(x_1)}h(x,Dv)\dx.
\end{equation}
\end{itemize}
\end{lemma}

\begin{proof}
Once fixed $B_{\rho_\smalltwo}(x_1)$, we prove \eqref{decay.v.h.pqleqn} and \eqref{decay.v.h.pq>n} only for $\rho_1\leq \rho_2/8$, as otherwise the estimate is trivial.

Similarly to what was done in Paragraph~\ref{par.harmonic} to deduce \eqref{uniform.comparability.a}, we can prove that the $(p,q)$-phase condition holds on every concentric sub-ball of $B_{\rho_3}$, with the same parameter $M$ and, in particular,
\[
\frac12 H(x_1,z)\leq H(x,z)\leq2H(x_1,z),\qquad\frac12 h(x_1,z)\leq h(x,z)\leq2h(x_1,z),
\]
 for every $x\in B_{\rho_3}$ and $z\in\mathbb R^n$. We set
\[
{\textnormal{\texttt{E}}}:=
\begin{cases}
q(1-\alpha_0), & \text{if }q\leq n\\[2mm]
\varsigma, & \text{if }q>n
\end{cases}.
\]
By \eqref{decay.energy.pq.est} in the first case and \eqref{decay.energy.superdimensional} in the second, we have
\[
\mean{B_{\rho_1}}H(x,Dv)\dx\leq c\,\Bigl(\frac{\rho_1}{\rho_2}\Bigr)^{-{\textnormal{\texttt{E}}}}\mean{B_{\rho_2/2}}H(x,Dv)\dx,
\]
with a constant depending on $\data,\gdata$ and $M$ in the first case and also on $\varsigma$ in the second. Thanks to \eqref{p.phase.sup.local.help.2} and Remark~\ref{essent.constant}, we have $a(\,\cdot\,)\approx_M a(x_1)$ in $B_{\rho_3}$. Using this comparability together with the averaged estimate \eqref{averaged.Psi.inverse} with $\el=1$, we infer
\begin{align}\label{first.first}
\mean{B_{\rho_\smallone}}h(x_1,Dv)\dx&\leq c\bigl(h_{x_\smallone}\circ H_{x_\smallone}^{-1}\bigr)\biggl(c\Bigl(\frac{\rho_1}{\rho_2}\Bigr)^{-{\textnormal{\texttt{E}}}}\mean{B_{\rho_2/2}}H(x,Dv)\dx\biggr)\notag\\
&\leq c\Bigl(\frac{\rho_1}{\rho_2}\Bigr)^{-\frac{{\textnormal{\texttt{E}}}}{q'}}\bigl(h_{x_\smallone}\circ H_{x_\smallone}^{-1}\bigr)\biggl(\mean{B_{\rho_2/2}}H(x,Dv)\dx\biggr).
\end{align}
In the last line we used the scaling property \eqref{g_x.lambda}, which allows us to extract the factor $c(\rho_1/\rho_2)^{-{\textnormal{\texttt{E}}}}\geq1$ from the argument of $h_{x_\smallone}\circ H_{x_\smallone}^{-1}$, raising it to the power $1/q'$. To estimate the right-hand side, we choose $\sigma\equiv\sigma(p,q)>0$ sufficiently small so that $t\mapsto[(H_{x_1}\circ h_{x_1}^{-1})(t)]^\sigma$ is concave, as guaranteed by \eqref{Psi.small.concave}, without any restriction on $q-p$. Applying \eqref{rev.self.improved} with this choice, we have
\[
\mean{B_{\rho_\smalltwo/2}} H(x,Dv)\dx\leq c(\data,\gdata,M)\biggl(\mean{B_{\rho_\smalltwo}}\bigl[H(x_1,Dv)\bigr]^\sigma\dx\biggr)^{\frac1\sigma},  
\]
using again that $H(x_1,Dv)\approx_M H(x,Dv)$ in $B_{\rho_\smalltwo}$. We ultimately use concavity  to get
\[
\biggl(\mean{B_{\rho_\smalltwo}}\bigl[H(x_1,Dv)\bigr]^\sigma\dx\biggr)^{\frac1\sigma}\leq \bigl(H_{x_\smallone}\circ h_{x_\smallone}^{-1}\bigr)\biggl(\mean{B_{\rho_\smalltwo}}h(x_1,Dv)\dx\biggr).
\]
Combining the last estimate with \eqref{first.first} gives the conclusion as the multiplicative constants can be absorbed by \eqref{g_x.lambda}:
\[
\bigl(h_{x_\smallone}\circ H_{x_\smallone}^{-1}\bigr) \biggl(c\bigl(H_{x_\smallone}\circ h_{x_\smallone}^{-1}\bigr)\biggl(\mean{B_{\rho_\smalltwo}}h(x_1,Dv)\dx\biggr)\biggr)\leq c(p,q,M)\mean{B_{\rho_\smalltwo}}h(x,Dv)\dx.
\]
\end{proof}

\subsection{Morrey decay of the flux magnitude}

The proposition below establishes a decay estimate for the averaged flux magnitude, framing the measure density $|\mu|(B_R)/R^{n-1}$ as its natural error term. Serving as a cornerstone for the potential estimates, this bound supplies the scale-dependent threshold $\lambda\lesssim R^{-(n-1)}$ governing the flux levels in the upcoming iteration, provided the energy and total variation of the datum remain controlled on a fixed reference ball (see \eqref{R.lambda}). Furthermore, the estimate allows us to sharpen our earlier comparisons between homogeneous solutions and their frozen autonomous counterparts by reformulating the approximation error through quantities suitable for the iteration; compare \eqref{funV.DvDw} and \eqref{pq-frozen.comparision} with the refined bounds \eqref{DvDw.energy.final} and \eqref{compa.vw.V.pq}, respectively. 

Finding its roots in the scaling of the $p$-Laplacian fundamental solution (see \eqref{fundamental.rel.gradient} and the earlier discussion), the exponent $n-1$ remains entirely independent of $p$, even within the super-dimensional regime $p>n$. Such behavior hints that the estimate should persist under the weaker conditions of \eqref{A.weak}, given that \eqref{decay.h.Du} extends to $p$-Laplacian-type equations with merely measurable coefficients (\cite[(8.7)]{MingioneCZ2007} and \cite[Proposition 9.3]{MinMil}). Our proof nevertheless requires \eqref{A.strong} in the superdimensional regimes, namely when $p>n$ in the $p$-phase or $q>n$ in the $(p,q)$-phase: in these cases, the homogeneous estimates available under \eqref{A.weak} do not yield the decay needed for our argument, so we use the approximation results established above.

\begin{proposition}\label{decay.measure.data}
Let $u\in W^{1,H(\cdot)}_{\rm loc}(\Omega)$ be a local energy solution to \eqref{equation} under the structural assumptions \eqref{a.holder} to \eqref{bounds} and \eqref{A.strong}; when $p\leq n$, assume additionally that $\mu\in L^\infty_{\rm loc}(\Omega)$. Fix $\Omega_0\Subset\Omega$ as in \eqref{Omega0}, and let $\Omega_1$ and $\mathcal M$ be defined by \eqref{Omega1} and \eqref{energy.M}. There exists a radius $R_*\in(0,\min\{1,d_{\Omega_0}/8\}]$, depending only on $\data,\gdata$ and $\omega(\,\cdot\,)$, such that for every $x_0\in\Omega_0$ and $0<R\leq R_*$, the following estimate holds for every $0<\rho\leq R$:
\begin{equation}\label{decay.h.Du}
\rho^{n-1}\mean{B_\rho(x_\smallzero)}h(x,Du)\dx\leq c(\data,\gdata)\biggl[R^{n-1}\mean{B_R(x_\smallzero)}h(x,Du)\dx +|\mu|\bigl(B_R(x_0)\bigr) \biggr].
\end{equation}
\end{proposition}

\begin{proof}
We prove \eqref{decay.h.Du} only for $\rho\leq R/2$, since the remaining case follows by enlarging the integration ball. We modify the iteration argument introduced in \cite{MingioneCZ2007, MinMil}. For every $r\in(0,R]$, let $v_r\in u+W^{1,H(\cdot)}_0(B_{r/2}(x_0))$ be the homogeneous comparison map defined by \eqref{def.v}, with $r$ in place of $R$. By Paragraph~\ref{choice.E}, these maps satisfy \eqref{bound.energy.mdata} on their respective domain balls with the same constant $E\equiv E(\data,\gdata)$, independently of $r$.

\medskip

For $\tau\in(0,1/64)$ to be chosen we consider the radii, the balls and the comparison maps
\[
R_k=\tau^k R/2,\qquad  \lambda B_k=B_{\lambda R_k}(x_0),\qquad v_k=v_{R_k}\qquad k\in\N_0, \lambda>0.
\]
We set, for $M_\varsigma$ provided by Proposition~\ref{decay.superdimensional},
\[
\begin{cases}
\varsigma_q:=q'(1-\alpha_0),\quad M_q:=M_{\varsigma_q},\qquad&\text {if $q>n$}\\[1mm]
\varsigma_p:=p'(1-\alpha_0),\quad M_p:=M_{\varsigma_p}&\text {if $p>n$}
\end{cases},
\]
and
\[
M:=
\begin{cases}
3, & q\leq n,\\[1mm]
M_q, & p\leq n<q\\[1mm]
\max\{M_p,M_q\}, & p>n
\end{cases}.
\]
We set
\[
R_*:=
\begin{cases}
\min\{1,d_{\Omega_0}/8\}, & q\leq n\\[1mm]
\min\{1,d_{\Omega_0}/8,R_{\varsigma_q}\}, & p\leq n<q\\[1mm]
\min\{1,d_{\Omega_0}/8,R_{\varsigma_p},R_{\varsigma_q}\}, & p>n
\end{cases},
\]
where the radii are provided by Proposition~\ref{decay.superdimensional}. We then set 
\[
\mathcal I:=\left\{k\in\mathbb N_0:\inf_{B_k}a(\,\cdot\,)>M[a]_{0,\alpha}R_k^\alpha\right\}
\]
and we define the exit time $i_e=\min\mathcal I$ (agreeing that $i_e=+\infty$ if $\mathcal I=\emptyset$). We also set, in order to shorten the notation, 
\[
{\texttt d}_p:=
\begin{cases}
n-(1-\alpha_0)(p-1), & p\leq n\\[2mm]
n-(1-\alpha_0), & p>n
\end{cases},
\qquad{\texttt d}_{p,q}:=
\begin{cases}
n-(1-\alpha_0)(q-1), & q\leq n\\[2mm]
n-(1-\alpha_0), & q>n
\end{cases};
\]
these exponents satisfy 
\begin{equation}\label{bound.below.alphazero}
\texttt d:=\min\big\{{\texttt d}_p,{\texttt d}_{p,q}\big\}\geq1+\alpha_0>1,  
\end{equation}
since $\alpha_0>0$. Finally, we introduce a parameter $\varepsilon\in(0,1)$ to be chosen later; we are going to use the comparison estimates of Lemmas~\ref{lem:rough.compa.solopartep} and~\ref{lem:rough.compa.pq} in the form of Remark~\ref{rem:comparison.epsilon}.

\medskip

Now we analyze separately various cases: 
\subsubsection*{$\boldsymbol{\bullet}$ {\bf After the exit:} $\boldsymbol{0\leq i_e<\infty}$, {\bf indices} $\boldsymbol{k\geq i_e}$} 
We can apply Lemma~\ref{decay.2} choosing $\rho_3=3R_k/8$: taking $\varsigma=\varsigma_q$ when $q>n$, employing \eqref{decay.v.h.pqleqn} or \eqref{decay.v.h.pq>n} (depending on whether we are in the case $q\leq n$ or $q>n$) with $\rho_1=R_{k+1}$ and $\rho_2=R_k/4$, we have for every $k\geq i_e$
\[
\mean{B_{k+1}}h(x,Dv_k)\dx\leq c(4\tau)^{{\texttt d}_{p,q}-n}\mean{B_k/4}h(x,Dv_k)\dx\leq c\,\tau^{{\texttt d}_{p,q}-n}\mean{B_k/2}h(x,Dv_k)\dx,
\]
for a constant $c$ depending on $\data$ and $\gdata$. Hence, appealing twice to the comparison estimate in Lemma~\ref{lem:rough.compa.pq}-Remark~\ref{rem:comparison.epsilon}, we have
\begin{align}\label{after.exit.decay}
\mean{B_{k+1}}h(x,Du)\dx&\leq c\mean{B_{k+1}}h_{x_\smallzero}(Dv_k)\dx+\frac{c(n,p,q)}{\tau^n}\mean{B_k/2}h_{x_\smallzero}(Du-Dv_k)\dx\notag\\
&\leq c\,\tau^{{\texttt d}_{p,q}-n}\mean{B_k/2}h_{x_\smallzero}(Dv_k)\dx +c_{\tau,\varepsilon}\frac{|\mu|(B_k)}{R_k^{n-1}}+\frac\varepsilon{\tau^n}\mean{B_k}h(x,Du)\dx\notag\\
&\leq \Big[c_1\tau^{{\texttt d}_{p,q}-n}+\frac\varepsilon{\tau^n}\Big]\mean{B_k}h(x,Du)\dx+c_{\tau,\varepsilon}\frac{|\mu|(B_R)}{R_k^{n-1}}
\end{align}
with $c_1\equiv c_1(\data,\gdata)$ and $c_{\tau,\varepsilon}\equiv c_{\tau,\varepsilon}(\data,\gdata,\tau,\varepsilon)$.

\subsubsection*{$\boldsymbol{\bullet}$ {\bf Before the exit:} $\boldsymbol{1\leq i_e\leq \infty}$, {\bf indices} $\boldsymbol{0\leq k< i_e-1}$} (When $i_e=+\infty$, this range includes every $k\in\mathbb N_0$); for every such $k$, both $B_k$ and $B_{k+1}$ are in the $p$-phase with constant $M$: indeed, both indices are strictly smaller than $i_e$, so this follows from the definition of the exit time. Moreover, Remark~\ref{stability.p.phase} gives
\[
\sup_{2B_{k+1}}a(\,\cdot\,)\leq \inf_{B_{k+1}}a(\,\cdot\,)+3[a]_{0,\alpha}R_{k+1}^\alpha\leq (M+3)[a]_{0,\alpha}R_{k+1}^\alpha.
\]
We start from \eqref{reverse.flux.pphase.est}, which is therefore available on $2B_{k+1}$ with phase parameter $M+3$, and use the triangle inequality and enlarge the integration ball to infer
\begin{align}\label{primary.pphase.decay}
\mean{B_{k+1}}h(x,Du)\dx&\leq c\mean{2B_{k+1}}|Du|^{p-1}\dx+c\frac{|\mu|(2B_{k+1})}{{(2R_{k+1})}^{n-1}}\notag\\
&\leq c\mean{2B_{k+1}}|Dv_k|^{p-1}\dx+c\,(4\tau)^{-n}\mean{B_k/2}|Du-Dv_k|^{p-1}\dx+c_\tau\frac{|\mu|(B_k)}{R_k^{n-1}}, 
\end{align}
with $c\equiv c(\data,\gdata)$ and $c_\tau\equiv c_\tau(\data,\gdata,\tau)$. Now in the sub-case $p\leq n$, we estimate the first term on the right-hand side by applying \eqref{decay.v.h}  from Lemma~\ref{decay.1} with $\rho_1=2R_{k+1}, \rho_2=R_k/4$ (we have $\rho_1\leq \rho_2/2$ as $\tau\leq 1/64$) and this, also using H\"older's inequality and then enlarging the averaging ball on the right-hand side, gives
\begin{align*}
\mean{2B_{k+1}}|Dv_k|^{p-1}\dx\leq {\biggl(\mean{2B_{k+1}}|Dv_k|^{(p-1)b}\dx\biggr)}^\frac1b&\leq c\,(8\tau)^{\texttt d_p-n}\biggl(\mean{B_k/4}|Dv_k|^{b(p-1)}\dx\biggr)^{\frac1b}\\
&\leq c\,\tau^{\texttt d_p-n}\biggl(\mean{B_k/2}|Dv_k|^{b(p-1)}\dx\biggr)^{\frac1b},
\end{align*}
again for a constant depending on $\data$ and $\gdata$. Merging with \eqref{primary.pphase.decay} and using the triangle inequality again, we find
\begin{align*}
\mean{B_{k+1}}h(x,Du)\dx&\leq c\,\tau^{{\texttt d}_p-n}{\biggl(\mean{B_k/2}|Dv_k|^{(p-1)b}\dx\biggr)}^\frac1b+c(4\tau)^{-n}\mean{B_k/2}|Du-Dv_k|^{p-1}\dx+c_\tau\frac{|\mu|(B_k)}{R_k^{n-1}}\\
&\leq c\,\tau^{{\texttt d}_p-n}{\biggl(\mean{B_k/2}|Du|^{(p-1)b}\dx\biggr)}^\frac1b\\
&\hspace{4,5cm}+c\,(4\tau)^{-n}\biggl(\mean{B_k/2}|Du-Dv_k|^{(p-1)b}\dx\biggr)^\frac1b+c_\tau\frac{|\mu|(B_k)}{R_k^{n-1}}.
\end{align*}
Finally, using the reverse H\"older inequality of Lemma~\ref{local.rH} and Lemma~\ref{lem:rough.compa.solopartep}-Remark~\ref{rem:comparison.epsilon} for an appropriate choice of a small multiple of $\varepsilon>0$ (in both cases the choice $\es=b$ is allowed by $b<n'$, see \eqref{es}, the assumption for both the results), we get
\begin{align}\label{dp.pn}
\mean{B_{k+1}}h(x,Du)\dx &\leq c\bigl[c_{rH}\tau^{{\texttt d}_p-n}+\varepsilon\bigr]\mean{B_k}h(x,Du)\dx+c_{\tau,\varepsilon}\frac{|\mu|(B_k)}{R_k^{n-1}}\notag\\
&\leq \bigl[c\,\tau^{{\texttt d}_p-n}+\varepsilon\bigr]\mean{B_k}h(x,Du)\dx+c_{\tau,\varepsilon}\frac{|\mu|(B_k)}{R_k^{n-1}} 
\end{align}
after relabeling $\varepsilon$, with $c$ depending on $\data$ and $\gdata$ and $c_{\tau,\varepsilon}$ depending on $\data,\gdata,\tau$ and $\varepsilon$. In the subcase $p>n$, we set
\[
\varsigma_*:=
\begin{cases}
\varsigma_p, & \text{if $M_p\geq M_q$}\\
\varsigma_q, & \text{if $M_q>M_p$}
\end{cases};
\]
with this choice we have $M=M_{\varsigma_*}$ and $\varsigma_*/p'\leq1-\alpha_0$. We apply \eqref{decay.v.h.p>n} with $\varsigma=\varsigma_*$, $\rho_1=8R_{k+1}$, and $\rho_2=R_k/4$. Since $\tau\leq1/64$, these radii satisfy $\rho_1\leq\rho_2/2$; moreover, the $p$-phase condition on $B_{\rho_1}$ follows from that on $B_{k+1}$. Remembering that $\texttt d_p-n=-(1-\alpha_0)$ and enlarging the averaging ball on the right-hand side, we get
\[
\mean{2B_{k+1}}h(x,Dv_k)\dx\leq c(32\tau)^{\texttt d_p-n}\mean{B_k/4}|Dv_k|^{p-1}\dx\leq c\,\tau^{\texttt d_p-n}\mean{B_k/2}|Dv_k|^{p-1}\dx.
\]
Starting from \eqref{primary.pphase.decay}, we can now repeat the argument of the preceding case with $b=1$, using \eqref{local.energy.SOLA} and the comparison estimate from Remark~\ref{rem:comparison.epsilon} with $\es=1$. Equation \eqref{dp.pn} thus holds also when $p>n$.

\subsubsection*{$\boldsymbol{\bullet\ 1\leq i_e< \infty}$, {\bf index} $\boldsymbol{k= i_e-1}$}
Here we simply estimate by monotonicity 
\begin{equation}\label{change.phase.iteration}
\int_{B_{i_e}}h(x,Du)\dx \leq \int_{B_{i_e-1}}h(x,Du)\dx.
\end{equation}

\medskip

Now, recalling the definition of $\texttt d$ in \eqref{bound.below.alphazero}, combining \eqref{after.exit.decay}, \eqref{dp.pn} and \eqref{change.phase.iteration} yields for every $k\in\mathbb N_0$
\[
\int_{B_{k+1}}h(x,Du)\dx \leq \Bigl[\bigl[c_{\rm ite}+\delta_{k,i_e-1}\tau^{-\texttt d}\bigr]\tau^{\texttt d}+\varepsilon\Bigr]\int_{B_k}h(x,Du)\dx +c_{\tau,\varepsilon}|\mu|(B_R)R_k,
\]
with $c_{\rm ite}$ depending on $\data$ and $\gdata$; $\delta_{k,i_e-1}$ is the Kronecker delta (here $\delta_{k,i_e-1}=1$ if $1\leq i_e<\infty$ and $k=i_e-1$, and $\delta_{k,i_e-1}=0$ otherwise). Now first we choose $\tau\equiv \tau(\data,\gdata)$ small enough so that
\[
c_{\rm ite}\tau^{\frac{{\texttt d}-1}2}\leq \frac12,\qquad \text{then}\qquad \varepsilon=\frac{\tau^{\texttt d}}2,
\]
so that, setting 
\[
\bar {\texttt d}=\frac{\texttt d+1}2>1 \qquad \text{and}\qquad \varphi(r)=\int_{B_r(x_0)}h(x,Du)\dx,\qquad 0<r\leq R,
\]
we have $\tau^{\texttt d}\leq \tau^{\bar {\texttt d}}$, $c_{\rm ite}\tau^{\texttt d}\leq \tau^{\bar {\texttt d}}/2$ and thus
\[
\varphi(R_{k+1}) \leq \bigl[\tau^{\bar {\texttt d}}+\delta_{k,i_e-1}\bigr]\varphi(R_k)+c\,{|\mu|(B_R)}R_k,
\]
with $c$ now depending only on $\data$ and $\gdata$. Since $\delta_{k,i_e-1}$ can be nonzero for at most one index, any product of $m$ consecutive coefficients is bounded by $(1+\tau^{-\bar{\texttt d}})\tau^{m\bar{\texttt d}}$. As $\tau$ and $\bar{\texttt d}$ depend only on $\data$ and $\gdata$, we can reabsorb $1+\tau^{-\bar{\texttt d}}$ into the constant and obtain by iteration, for every $k\geq1$,
\begin{align*}
\varphi(R_k) &\leq  c\,\tau^{k \bar {\texttt d}}\varphi(R) +c\,{|\mu|(B_R)}\sum_{j=1}^k\tau^{(j-1)\bar{\texttt d}+(k-j)}R \\
&\leq   c\,\tau^{k \bar {\texttt d}}\varphi(R) +c\,{|\mu|(B_R)}\tau^{k-1}\sum_{j=0}^\infty\tau^{j(\bar{\texttt d}-1)}R \leq  c\,\tau^{k \bar {\texttt d}}\varphi(R) +c\,{|\mu|(B_R)}\tau^kR
\end{align*}
since the series is convergent being $\bar{\texttt d}-1>0$ and $\tau$ and $\bar{\texttt d}$ constants depending on $\data$ and $\gdata$. For $0<\rho\leq R/2$, choose $k\geq0$ such that $R_{k+1}<\rho\leq R_k$. As in the proof of Lemma~\ref{geometric.decay.iteration}, monotonicity gives $\varphi(\rho)/\rho\leq\varphi(R_k)/(\tau R_k)$. Since $R_k=\tau^kR/2$ and $\tau^{k\bar{\texttt d}}\leq\tau^k$, the preceding estimate yields $\varphi(\rho)/\rho\leq c[\varphi(R)/R+|\mu|(B_R)]$, which proves \eqref{decay.h.Du}.
\end{proof}

As an immediate consequence, we obtain a uniform bound for the flux density on every fixed compactly contained subset of $\Omega$.

\begin{corollary}\label{flux.decay}
Let $u$ and $\Omega_0$ be as in Proposition~\ref{decay.measure.data}. There exist a reference radius $\rho_{\rm ref}>0$, depending only on $\data,\gdata$ and $\omega(\,\cdot\,)$, and a constant $c_{\rm ref}=c_{\rm ref}(\data,\gdata,\rho_{\rm ref})$ such that
\begin{equation}\label{local.decay.ultimate}
\rho^{n-1}\mean{B_\rho(x_0)}h(x,Du)\dx+|\mu|(B_\rho(x_0))\leq c_{\rm ref}
\end{equation}
for every $x_0\in\Omega_0$ and $0<\rho\leq\rho_{\rm ref}$.
\end{corollary}

\begin{proof}
Let $R_*$ be the radius from Proposition~\ref{decay.measure.data}; we fix
\[
\rho_{\rm ref}:=\min\biggl\{R_*,\frac{d_{\Omega_0}}8\biggr\}.
\]
By the definition of $\Omega_1$, we have $B_{\rho_{\rm ref}}(x_0)\subset\Omega_1$ for every $x_0\in\Omega_0$. Therefore, \eqref{decay.h.Du}, applied with $R=\rho_{\rm ref}$, gives for every $0<\rho\leq\rho_{\rm ref}$
\begin{align*}
\rho^{n-1}\mean{B_\rho(x_0)}h(x,Du)\dx+|\mu|(B_\rho(x_0))&\leq c\biggl[\rho_{\rm ref}^{n-1}\mean{B_{\rho_{\rm ref}}(x_0)}h(x,Du)\dx+|\mu|(B_{\rho_{\rm ref}}(x_0))\biggr]\\
&\leq c\biggl[\frac1{\rho_{\rm ref}}\int_{\Omega_1}h(x,Du)\dx+|\mu|(\Omega_1)\biggr]\leq c\biggl[\frac{\mathcal M}{\rho_{\rm ref}}+|\mu|(\Omega)\biggr]
\end{align*}
by the definition of $\mathcal M$ in \eqref{energy.M}. If we denote by $c_{\rm ref}$ the quantity in the last line, then it depends only on $\data,\gdata$ and $\rho_{\rm ref}$. In particular, $\rho_{\rm ref}$ and $c_{\rm ref}$ are fixed independently of the smaller radii considered in the subsequent arguments.
\end{proof}

\subsection{Other comparison results}\label{Other.comparison}

Let $C_{\rm ref}:=\max\{1,c_{\rm ref}\}$. We fix a radius $R_{\rm cmp}\leq\min\{1,\rho_{\rm ref}\}$, depending only on $\data,\gdata$ and $\omega(\,\cdot\,)$, such that $\omega(R_{\rm cmp})\leq1$. When $p\leq n$, we additionally require
\[
C_{\rm ref}^{\frac{q-p}{p-1}}R_{\rm cmp}^{s_0(1-\frac1p)}\leq1,
\]
where $s_0>0$ is defined in \eqref{bound.with.x}. This choice is independent of the phase constant $M$ and will be retained in the following comparison estimates.

\medskip

For the remainder of this section, let $u$ and $\Omega_0$ be as in Proposition~\ref{decay.measure.data}. For each ball $B_R=B_R(x_0)$ with $x_0\in\Omega_0$ and $0<R\leq R_{\rm cmp}$, let $v\in u+W^{1,H(\cdot)}_0(B_{R/2})$ denote the homogeneous comparison map from \eqref{def.v}. We first use \eqref{decay.h.Du} to estimate $h(x,Du-Dv)$ in the $p$-phase.

\begin{lemma}\label{lem:rough.compa}
If \eqref{p.phase} holds on $B_R$ for some $M\geq1$ and $\x$ as in \eqref{iccs}, then
\begin{equation}\label{rough.compa.pphase}
\mean{B_{R/2}}h(x,Du-Dv)\dx\leq c\frac{|\mu|(B_R)}{R^{n-1}}+c\,d_{p,q},
\end{equation}
where $d_{p,q}$ is defined in \eqref{d_pq} and $c$ depends only on $\data,\gdata$ and $M$.
\end{lemma}

\begin{proof}
We distinguish two cases:

\medskip

\subsubsection{Case $p>n$} By H\"older's inequality and \eqref{rough.compa.supercritical.energy},
\begin{equation}\label{rough.compa.supercritical.p}
\mean{B_{R/2}}{|Du-Dv|}^{p-1}\dx\leq c\,\frac{|\mu|(B_R)}{R^{n-1}}.
\end{equation}
To estimate the $q$-part, we write $a(x)=a_i(B_R)+[a(x)-a_i(B_R)]$ and first consider the constant part, to estimate which we apply \eqref{compa.q.component.ell} with $\bar a_q=a_i(B_R)$ and $\ti=1$:
\begin{equation}\label{rough.compa.supercritical.q.constant}
a_i(B_R)\mean{B_{R/2}}|Du-Dv|^{q-1}\dx\leq c\frac{|\mu|(B_R)}{R^{n-1}},
\end{equation}
with $c$ depending only on $n,p,q$ and $\nu$. We next estimate the oscillating part of the coefficient: recalling that $0\leq a(x)-a_i(B_R)\leq c[a]_{0,\alpha}R^\alpha$ and $a(x)-a_i(B_R)\leq a(x)$ in $B_{R/2}$, H\"older's inequality and \eqref{rough.compa.supercritical.energy} yield
\begin{align*}
\mean{B_{R/2}}&\bigl[a(x)-a_i(B_R)\bigr]|Du-Dv|^{q-1}\dx\\
&\leq\biggl(\mean{B_{R/2}}\bigl[a(x)-a_i(B_R)\bigr]|Du-Dv|^{q-p}\dx\biggr)^{\frac1p}\biggl(\mean{B_{R/2}}\bigl[a(x)-a_i(B_R)\bigr]|Du-Dv|^q\dx\biggr)^\frac1{p'}\\
&\leq cR^{\frac\alpha p}\biggl(\mean{B_{R/2}}|Du-Dv|^p\dx\biggr)^{\frac{q-p}{p^\smalltwo}}\biggl(\mean{B_{R/2}}a(x)|Du-Dv|^q\dx\biggr)^\frac1{p'}\\
&\leq cR^{\frac\alpha p}\biggl[\frac{|\mu|(B_R)}{R^{n-1}}\biggr]^{1+\frac{q-p}{p(p-1)}}\leq c|\mu|(\Omega)^{\frac{q-p}{p(p-1)}}R^{\frac1p[\alpha-\frac{(n-1)(q-p)}{p-1}]}\frac{|\mu|(B_R)}{R^{n-1}}\leq c\frac{|\mu|(B_R)}{R^{n-1}},
\end{align*}
where we used $q-p\leq p$ and, by \eqref{bound.supercritical},
\[
\alpha-\frac{(n-1)(q-p)}{p-1}\geq\frac{\alpha(p-n)}{n(p-1)}>0.
\]
Combining this estimate with \eqref{rough.compa.supercritical.p} and \eqref{rough.compa.supercritical.q.constant}, we conclude with \eqref{rough.compa.pphase}.

\subsubsection{Case $p\leq n$} Set $\es=(q-1)/(p-1)$. Since $1\leq\es<n'$ by \eqref{bound.trivial}, we may apply \eqref{rough.compa.solopartep} with this exponent to obtain
\[
\biggl(\mean{B_{R/2}}|Du-Dv|^{q-1}\dx\biggr)^{1/\es}\leq c\biggl[\frac{|\mu|(B_R)}{R^{n-1}}+d_{p,q}\biggr]\leq c\,d_1,
\]
where the last inequality follows from the definition of $d_1$ in \eqref{dtwo}, that of $d_{p,q}$ in \eqref{d_pq} and Young's inequality. If $d_1=0$, this already gives $Du=Dv$ almost everywhere in $B_{R/2}$. Otherwise, by \eqref{local.decay.ultimate}, \eqref{bound.with.x} and the choice of $R_{\rm cmp}$, we have
\begin{equation}\label{meas.partial.global}
R^{\alpha-\x}d_1^{\es-1}\leq C_{\rm ref}^{\es-1}R^{\alpha-\x-(n-1)(\es-1)}\leq C_{\rm ref}^{\es-1}R^{s_0}\leq R^{s_0/p}\leq1.
\end{equation}
Using H\"older's inequality for the $p$-part and \eqref{sup-inf_withx} for the $q$-part, we obtain
\begin{align*}
\mean{B_{R/2}}h(x,Du-Dv)\dx
&\leq\biggl(\mean{B_{R/2}}|Du-Dv|^{q-1}\dx\biggr)^{1/\es}+cR^{\alpha-\x}\mean{B_{R/2}}|Du-Dv|^{q-1}\dx\\
&\leq c\bigl(1+R^{\alpha-\x}d_1^{\es-1}\bigr)\biggl[\frac{|\mu|(B_R)}{R^{n-1}}+d_{p,q}\biggr].
\end{align*}
Estimate \eqref{rough.compa.pphase} now follows from \eqref{meas.partial.global}, with $c$ depending only on $\data,\gdata$ and $M$.
\end{proof}

\subsection{Refined comparison estimate for the frozen map}\label{paragraph.p}

We now return to the regular reference map introduced in Paragraph~\ref{par.harmonic}, which is $A_i$-harmonic in $B_{R/8}$ and satisfies an interior gradient bound on $B_{R/16}$: the uniform decay estimate established in Corollary~\ref{flux.decay} allows us to refine the preceding comparison estimates in a form suitable for the iteration procedure developed in the next section. From now until the end of this section, we assume the stronger structural conditions \eqref{A.strong}.

\medskip

We anticipate that we may further decrease $R_{\rm cmp}$ in what follows. Recalling the choices of $x_m$, $A_i$ and $H_i$ as in \eqref{frozen.minimum}-\eqref{inf.smaller.H}, we observe that, when \eqref{p.phase} holds, \eqref{quasi.approx.p.harmonic} gives
\begin{equation}\label{approx.A_i.harm.pphase}
\biggl|\mean{B_{R/8}}\langle A_i(Dv),D\varphi\rangle\dx\biggr|\leq c\,{\|D\varphi\|}_{L^\infty(B_{R/8})}\biggl[\omega(R)\mean{B_{R/8}}h(x,Dv)\dx+2\mean{B_{R/8}}a(x)|Dv|^{q-1}\dx\biggr]
\end{equation}
for every $\varphi\in C^\infty_c(B_{R/8})$, $c\equiv c(p,q,L)$; see Remark~\ref{frozen.p.preliminary}.

\medskip

We use $d_1$ as defined in \eqref{dtwo} and recall that $E_1$ is defined in \eqref{def.E1}. If $E_1=0$, then $v$ is constant in $B_{R/4}$ and we take $w$ to be the same constant in $B_{R/8}$; all the estimates below are immediate. If $d_1=0$, then $u$ and its homogeneous comparison map $v$ are constant, so this case is also immediate. We therefore assume $E_1,d_1>0$.

We first consider the case $p\leq n$. By \eqref{local.decay.ultimate}, we have $d_1\leq C_{\rm ref}R^{-(n-1)}$. Let $b$ be as in \eqref{def.b}-\eqref{bbi}; appealing to \eqref{local.energy.SOLA}, \eqref{rough.compa.solopartep} and Young's inequality when $p<2$, we obtain
\[
\biggl(\mean{B_{R/2}}|Dv|^{b(p-1)}\dx\biggr)^{\frac1b}\leq c\,d_1,
\]
where $c$ depends only on $\data$. We now apply \eqref{rev.self.improved} on $B_{R/2}$ with $\sigma=b(p-1)/q$: since $p\leq q$, H\"older's inequality and \eqref{sup-inf_withx} give
\begin{align}\label{est.E1}
E_1^p=\mean{B_{R/4}}H(x,Dv)\dx&\leq c\biggl[\biggl(\mean{B_{R/2}}|Dv|^{b(p-1)}\dx\biggr)^{\frac p{b(p-1)}}+\sup_{B_R}a(\,\cdot\,)\,\biggl(\mean{B_{R/2}}|Dv|^{b(p-1)}\dx\biggr)^{\frac q{b(p-1)}}\biggr]\notag\\
&\leq c(\data,\gdata,M)\Bigl[1+R^{\alpha-\x}d_1^{\frac{q-p}{p-1}}\Bigr]d_1^{p'}\leq c(\data,\gdata,M)d_1^{p'},
\end{align}
where the last inequality follows from \eqref{meas.partial.global}. In particular,
\begin{equation}\label{bound.energetic}
\mean{B_{R/8}}|Dv|^p\dx\leq c(n,p)E_1^p\leq c(\data,\gdata,M)d_1^{p'}.
\end{equation}
We can therefore estimate the second integral in \eqref{approx.A_i.harm.pphase}. Since $q-1<p$ by \eqref{bound.lasco}, H\"older's inequality and \eqref{sup-inf_withx} yield
\begin{align}\label{first.part.second.integral}
\mean{B_{R/8}}a(x)|Dv|^{q-1}\dx&\leq(M+2)[a]_{0,\alpha}R^{\alpha-\x}\biggl(\mean{B_{R/8}}|Dv|^p\dx\biggr)^{\frac{q-1}p}\notag\\
&\leq cR^{\alpha-\x}E_1^{p-1}\biggl(\mean{B_{R/8}}|Dv|^p\dx\biggr)^{\frac{q-p}p}\leq cR^{\alpha-\x}d_1^{\frac{q-p}{p-1}}E_1^{p-1}\leq c\,R^{\frac{s_\smallzero}p}E_1^{p-1},
\end{align}
with $c\equiv c(\data,\gdata,M)$. The same estimates also give
\[
a_i(B_R)E_1^{q-p}\leq cR^{\alpha-\x}d_1^{\frac{q-p}{p-1}}\leq c\,R^{\frac{s_\smallzero}p}\leq c.
\]
Consequently, the normalization bounds \eqref{h_i=p-1} and \eqref{V.Hi.scaling.difference} remain valid in the present case.

If $p>n$, we instead use \eqref{oscillation.q.estimate} and \eqref{bound.hi.Dv.harmonic}: although the preliminary argument was written for $\x=0$, it applies to every admissible $\x$ by replacing $R^{\alpha/p}$ with $R^{(\alpha-\x)/p}$. By \eqref{bound.with.x}, this gives
\[
\mean{B_{R/8}}a(x)|Dv|^{q-1}\dx\leq c\,R^{\frac{s_\smallzero}p}E_1^{p-1},
\]
again with $c\equiv c(\data,\gdata,M)$. The bound $a_i(B_R)E_1^{q-p}\leq c$ follows in this case from \eqref{rescaled.ai.bound}. Therefore, setting
\begin{equation}\label{s1}
s_1:=\frac{s_0}p,
\end{equation}
we obtain in all cases that
\[
\mean{B_{R/8}}h(x,Dv)\dx\leq\biggl(\mean{B_{R/8}}|Dv|^p\dx\biggr)^{\frac1{p'}}+c\,R^{s_1}E_1^{p-1}\leq c(1+R^{s_1})E_1^{p-1}\leq cE_1^{p-1},
\]
with $c\equiv c(\data,\gdata,M)$, since $s_1>0$ and $R\leq1$. Inserting these estimates into \eqref{approx.A_i.harm.pphase} yields
\[
\biggl|\mean{B_{R/8}}\langle A_i(Dv),D\varphi\rangle\dx\biggr|\leq c(\data,\gdata,M)\Bigl[\omega(R)+R^{s_1}\Bigr]E_1^{p-1}\|D\varphi\|_{L^\infty(B_{R/8})}
\]
for every $\varphi\in C_c^\infty(B_{R/8})$. If we set $\tilde v:=v/E_1$ as in \eqref{rescaled.v.Ai}-\eqref{rescaled.Hi}, the above inequality rescales to
\[
\biggl|\mean{B_{R/8}}\langle\tilde A_i(D\tilde v),D\varphi\rangle\dx\biggr|\leq c(\data,\gdata,M)\Bigl[\omega(R)+R^{s_1}\Bigr]\|D\varphi\|_{L^\infty(B_{R/8})},
\]
using the lower bound $h_i(E_1)\geq E_1^{p-1}$. The structural conditions for $\tilde A_i$ then follow directly from \eqref{rescaling.pphase}, and the energy bounds on $B_{R/4}$ and $B_{R/8}$ follow as in \eqref{tildev.bounded.1} and \eqref{tildev.bounded.2}, respectively. Lemma~\ref{p-harm} therefore applies with $B=B_{R/8}$, $\varphi=\tilde H_i$, and $A_0=\tilde A_i$. It provides an $\tilde A_i$-harmonic map $\tilde w\in\tilde v+W^{1,\tilde H_i}_0(B_{R/8})$ satisfying
\[
\mean{B_{R/8}}\big|V_{\tilde H_i}(D\tilde v)-V_{\tilde H_i}(D\tilde w)\big|^2\dx\leq c(\data,\gdata,M)\Bigl[\omega(R)+R^{s_1}\Bigr]^\iota
\]
for an exponent $\iota\in(0,1)$ depending only on $\data,\gdata$. Scaling back via $w=E_1\tilde w$ yields an $A_i$-harmonic map in $v+W^{1,H_i}_0(B_{R/8})$ which, by uniqueness, coincides with the comparison map constructed earlier whenever both regimes overlap. Finally, since $a_i(B_R)E_1^{q-p}\leq c(\data,\gdata,M)$, it follows that $H_i(E_1)\leq c(\data,\gdata,M)E_1^p$, and \eqref{V.Hi.scaling.difference} accordingly implies
\begin{equation}\label{funV.DvDw.allcases}
\mean{B_{R/8}}\big|V_{H_i}(Dv)-V_{H_i}(Dw)\big|^2\dx\leq c\Bigl[\omega(R)+R^{s_1}\Bigr]^\iota E_1^p,
\end{equation}
with $c$ depending only on $\data,\gdata$ and $M$. We finally claim that
\begin{equation}\label{en.rh.1}
E_1^p=\mean{B_{R/4}}H(x,Dv)\dx\leq c\biggl[\mean{B_R}h(x,Du)\dx+\frac{|\mu|(B_R)}{R^{n-1}}\biggr]^{p'},
\end{equation}
where $c\equiv c(\data,\gdata,M)$. For $p\leq n$, this is precisely \eqref{est.E1}. If $p>n$, we apply \eqref{rev.self.improved} with $\sigma=1/p'$. Using $[H(x,z)]^{1/p'}\leq c(p)h(x,z)$ and subadditivity, we obtain
\[
E_1^p\leq c\biggl(\mean{B_{R/2}}h(x,Du)\dx\biggr)^{p'}+c\biggl(\mean{B_{R/2}}h(x,Du-Dv)\dx\biggr)^{p'}
\]
with $c\equiv c(\data,E)$. Estimate \eqref{en.rh.1} then follows by using \eqref{rough.compa.pphase}, recalling that here $d_{p,q}=0$. Combining \eqref{funV.DvDw.allcases} and \eqref{en.rh.1}, we conclude that
\begin{equation}\label{DvDw.energy.final}
\mean{B_{R/8}}\big|V_{H_i}(Dv)-V_{H_i}(Dw)\big|^2\dx\leq c\Bigl[\omega(R)+R^{s_1}\Bigr]^\iota\biggl[\mean{B_R}h(x,Du)\dx+\frac{|\mu|(B_R)}{R^{n-1}}\biggr]^{p'}.
\end{equation}

\medskip 

We next convert \eqref{DvDw.energy.final} into a direct comparison estimate for the gradient difference in terms of the frozen energy $H_i(\,\cdot\,)$; this passage entails a loss in the exponent of the smallness factor when $p<2$.

\begin{proposition}\label{p-harm.rough1}
Let $u$ and $v$ be as in Lemma~\ref{lem:rough.compa}, with \eqref{p.phase} holding on $B_R$ for some $M\geq1$. Then there exists an $A_i$-harmonic map $w\in v+W^{1,H_i(\cdot)}_0(B_{R/8})$, with $A_i,H_i$ as in \eqref{vf.frozen.ai}-\eqref{inf.smaller.H}, such that
\begin{align}\label{compa.Dw.Lp}
\mean{B_{R/8}}H_i(Dv-Dw)\dx&\leq c\Bigl[\omega(R)+R^{s_\smallone}\Bigr]^{\iota \min\{1,\frac p2\}}\mean{B_{R/4}} H(x,Dv)\dx\notag\\
&\leq c\Bigl[\omega(R)+R^{s_\smallone}\Bigr]^{\iota \min\{1,\frac p2\}}\biggl[\mean{B_R}h(x,Du)\dx+\frac{|\mu|(B_R)}{R^{n-1}}\biggr]^{p'}.
\end{align}
The constant $c$ depends on $\data,\gdata$ and $M$, $\iota$ is the exponent appearing in Lemma~\ref{p-harm} and $s_1$ is defined in \eqref{s1}.
\end{proposition}

\begin{proof}
We use \eqref{diffVHx0} with the choice $x=x_m$ as chosen in \eqref{frozen.minimum}; therefore, in view of \eqref{DvDw.energy.final}, the only remaining terms we need to estimate are
\[
I=\mean{B_{R/8}}|Dv|^{\frac{(2-p)p}{2}}\big|V_p(Dv)-V_p(Dw)\big|^p\dx
\]
when $p<2$ and 
\[
II=\mean{B_{R/8}}|Dv|^{\frac{(2-q)q}{2}}a(x_m)\big|V_q(Dv)-V_q(Dw)\big|^q\dx
\]
when $q<2$. The first one is treated by using H\"older's inequality with conjugate exponents $(2/p,2/(2-p))$ and \eqref{Vx0} for $x=x_m$:
\begin{align}\label{pphase.compaDvDw.partial}
I&\leq \biggl(\mean{B_{R/8}}{\big|V_p(Dv)-V_p(Dw)\big|}^2\dx\biggr)^{\frac p2}\biggl(\mean{B_{R/8}}|Dv|^p\dx\biggr)^{\frac{2-p}{2}}\notag\\
&\leq c\,\Bigl[\omega(R)+R^{s_1}\Bigr]^{\iota\frac p2}E_1^{p\cdot\frac p2+p\cdot \frac{2-p}{2}}\leq c\,\Bigl[\omega(R)+R^{s_1}\Bigr]^{\iota\frac p2}\biggl[\mean{B_R}h(x,Du)\dx+c\,\frac{|\mu|(B_R)}{R^{n-1}}\biggr]^{p'}
\end{align}
using \eqref{funV.DvDw.allcases} and \eqref{en.rh.1}. For the second one, we estimate employing again H\"older's inequality and \eqref{Vx0},
\begin{equation*}
II\leq \biggl(\mean{B_{R/8}}a(x_m){\big|V_q(Dv)-V_q(Dw)\big|}^2\dx\biggr)^{\frac q2}\biggl(\mean{B_{R/8}}a(x)|Dv|^q\dx\biggr)^{\frac{2-q}{2}}\leq c\, \Bigl[\omega(R)+R^{s_1}\Bigr]^{\iota\frac q2}E_1^{p\cdot\frac q2+p\cdot \frac{2-q}{2}};
\end{equation*}
bounding $E_1^p$ by \eqref{en.rh.1} and combining these estimates with \eqref{DvDw.energy.final} yields the conclusion.
\end{proof}

\begin{corollary}\label{corollary.compa.rough.DuDw}
For $w\in v+W^{1,H_i(\cdot)}_0(B_{R/8})$ the  $A_i$-harmonic map appearing in Proposition~\ref{p-harm.rough1}, assuming that \eqref{p.phase} holds, we have
\begin{equation}\label{compa.Dw.Lh}
\mean{B_{R/16}}h(x,Dv-Dw)\dx\leq c\Bigl[\omega(R)+R^{s_\smallone}\Bigr]^{\frac{\iota(p-1)}{\max \{p,2\}}}\biggl[\mean{B_R} h(x,Du)\dx+\frac{|\mu|(B_R)}{R^{n-1}}\biggr],
\end{equation}
the constant $c$ depending on $\data,\gdata$ and $M$. The exponents $\iota$ and $s_1$ are as in Proposition~\ref{p-harm.rough1}.
\end{corollary}

\begin{proof}
Clearly, using H\"older's inequality, the fact $t^p\leq p\,H_i(t)$ and \eqref{compa.Dw.Lp}, we find
\begin{align*}
\mean{B_{R/16}}|Dv-Dw|^{p-1}\dx &\leq c(n,p) {{\biggl(\mean{B_{R/8}}|Dv-Dw|^p\dx\biggr)}}^{\frac1{p'}} \\
&\leq c\Bigl[\omega(R)+R^{s_\smallone}\Bigr]^{\frac{\iota}{p'} \min\{1,\frac p2\}}\biggl[\mean{B_R}h(x,Du)\dx+\frac{|\mu|(B_R)}{R^{n-1}}\biggr].
\end{align*}
In the case $p\leq n$, for the $q$-part of the flux, we have by \eqref{p.phase}-\eqref{sup-inf_withx}, H\"older's inequality, which is applicable by \eqref{bound.lasco}, and \eqref{compa.Dw.Lp}
\[
\mean{B_{R/16}}a(x)|Dv-Dw|^{q-1}\dx\leq c\,R^{\alpha-\x}\Bigl[\omega(R)+R^{s_\smallone}\Bigr]^{\frac{\iota(q-1)}{\max\{2,p\}}}\biggl[\mean{B_R} h(x,Du)\dx+\frac{|\mu|(B_R)}{R^{n-1}}\biggr]^{\frac{q-1}{p-1}} 
\]
with $c\equiv c(\data,\gdata,M)$. By \eqref{meas.partial.global} and \eqref{s1}, we have $R^{\alpha-\x}d_1^{(q-p)/(p-1)}\leq R^{s_1}\leq1$. Moreover, since $q\geq p$ and $\omega(R)+R^{s_1}\leq2$, we have
\[
\bigl[\omega(R)+R^{s_1}\bigr]^{\frac{\iota(q-1)}{\max\{2,p\}}}\leq c(p,q)\bigl[\omega(R)+R^{s_1}\bigr]^{\frac{\iota(p-1)}{\max\{2,p\}}}.
\]
Combining these bounds with the estimate for the $p$-part proves \eqref{compa.Dw.Lh} when $p\leq n$. If $p>n$, we need to use the Lipschitz character of $w$, in particular  Proposition~\ref{varphi-harmonic.are.regular}, which, together with Young's inequality and the energy estimate for $w$, leads to
\[
\|Dw\|_{L^\infty(B_{R/16})}^p \leq c\mean{B_{R/8}}H_i(Dw)\dx \leq c\,E_1^p.
\]
Therefore, by \eqref{p.phase}-\eqref{sup-inf_withx}, the previous estimate and the higher integrability bound used in \eqref{bound.hi.Dv.harmonic},
\[
\mean{B_{R/16}}a(x)|Dw|^{q-1}\dx\leq cR^{\alpha-\x}E_1^{q-1}\leq cR^{\alpha-\x-\frac{n(q-p)}{p(1+\eta_\smallzero)}}E_1^{p-1}\leq cR^{s_1}E_1^{p-1},
\]
where the last inequality follows from the definition of $s_1$ in \eqref{s1}. Combining this with the estimate for $Dv$ obtained from \eqref{oscillation.q.estimate} and \eqref{bound.hi.Dv.harmonic}, we find
\begin{align*}
\mean{B_{R/16}}a(x)|Dv-Dw|^{q-1}\dx&\leq c(q)\mean{B_{R/16}}a(x)|Dv|^{q-1}\dx+c(q)\mean{B_{R/16}}a(x)|Dw|^{q-1}\dx\\
&\leq cR^{s_1}E_1^{p-1}\leq c\Bigl[\omega(R)+R^{s_1}\Bigr]^{\frac{\iota(p-1)}p}\biggl[\mean{B_R}h(x,Du)\dx+\frac{|\mu|(B_R)}{R^{n-1}}\biggr];
\end{align*}
here we used that $R\leq1$, $\iota(p-1)/p<1$, and the definition of $E_1$ in \eqref{en.rh.1}. This proves \eqref{compa.Dw.Lh} also in the case $p>n$.
\end{proof}

\medskip

We now assume that we are in the $(p,q)$-phase, that is, \eqref{pq.phase} holds for some $M\geq 2$, and we refine the estimate \eqref{pq.harm.rescaled} obtained starting from \eqref{quasi.pq.harm}. The construction of the frozen comparison map indeed extends to the case $p\leq n$: the phase comparability and Theorem~\ref{hi.correct.exponent}, with the energy bound chosen in Paragraph~\ref{choice.E}, yield the same normalized energy and higher integrability estimates. Thus, with $E_2$ defined in \eqref{def.E2}, the preceding argument provides an $A_i$-harmonic map $w\in v+W^{1,H_i(\cdot)}_0(B_{R/8})$ satisfying \eqref{pq-frozen.energy} and \eqref{pq-frozen.comparision}, with constants independent of $M$. Now we continue by estimating the term on the right-hand side of \eqref{pq-frozen.comparision}, namely $H_{x_m}(E_2)$: by \eqref{uniform.comparability.a} and the alternative formulation of Proposition~\ref{cor:nice-rev.pq}, applied with $\rho=R/2$ and $a_*=a_i(B_R)$, we have
\begin{equation}\label{estimate.E2}
H_{x_m}(E_2)\leq c\bigl(H_{x_m}\circ h_{x_m}^{-1}\bigr)\biggl(\mean{B_{R/2}} h(x,Dv)\dx \biggr); 
\end{equation}
since by subadditivity and Lemma~\ref{lem:rough.compa.pq} (in the form of Remark~\ref{rem:comparison.epsilon} with $\varepsilon=1$) there holds
\[
\mean{B_{R/2}} h(x,Dv)\dx\leq c\mean{B_R} h(x,Du)\dx+c\,\frac{|\mu|(B_R)}{R^{n-1}},
\]
we conclude with
\begin{equation}\label{compa.vw.V.pq}
 \mean{B_{R/8}}{\big|V_{H_i}(Dv)-V_{H_i}(Dw)\big|}^2\dx  \leq c {\Bigl[\omega(R)+\frac{R^\x}M\Bigr]}^{\iota}\bigl(H_{x_m}\circ h_{x_m}^{-1}\bigr)\biggl(\mean{B_R} h(x,Du)\dx+\frac{|\mu|(B_R)}{R^{n-1}} \biggr)
\end{equation}
with $c$ depending on $\data$ and $\gdata$, but not on $M$. As for Proposition~\ref{p-harm.rough1}, we have the analogue

\begin{proposition}\label{pq-harm.rough1}
If \eqref{pq.phase} holds on $B_R$ for some $M\geq2$ and $\x$ as in \eqref{iccs}, then there exists an $A_i$-harmonic map $w\in v+W^{1,H_i(\cdot)}_0(B_{R/8})$, where $A_i$ and $H_i$ are as in \eqref{vf.frozen.ai} and \eqref{inf.smaller.H}, such that
\begin{equation}\label{compa.Dw.Lp.pq}
\mean{B_{R/8}}H_i(Dv-Dw)\dx\leq c{\Bigl[\omega(R)+\frac{R^\x}M\Bigr]}^{\iota\min\{1,\frac p2\}}\bigl(H_{x_m}\circ h_{x_m}^{-1}\bigr)\biggl(\mean{B_R} h(x,Du)\dx+\frac{|\mu|(B_R)}{R^{n-1}}\biggr)
\end{equation}
with $c$ depending on $\data$ and $\gdata$, but not on $M$.
\end{proposition}

\begin{proof}
In all three cases we use \eqref{diffVHx0} and the quantity $d_1$ defined in \eqref{dtwo}. If $p\geq2$ the estimate is immediate from \eqref{compa.vw.V.pq}. If $q< 2$, using twice H\"older's inequality as in \eqref{pphase.compaDvDw.partial}, we have
\begin{multline*}
\mean{B_{R/8}}H_i(Dv-Dw)\dx\\ \leq c\mean{B_{R/8}}{\big|V_{H_i}(Dv)-V_{H_i}(Dw)\big|}^2\dx+c\biggl(\mean{B_{R/8}}{\big|V_p(Dv)-V_p(Dw)\big|}^2\dx\biggr)^\frac p2\biggl(\mean{B_{R/8}}|Dv|^p\dx\biggr)^{\frac{2-p}{2}}\\
+c\biggl(\mean{B_{R/8}}a(x_m){\big|V_q(Dv)-V_q(Dw)\big|}^2\dx\biggr)^\frac q2\biggl(\mean{B_{R/8}}a(x_m)|Dv|^q\dx\biggr)^{\frac{2-q}{2}}=:I+II_p+II_q.
\end{multline*}
It remains to estimate $II_p$ and $II_q$: we have, recalling the definition of $E_2$ in \eqref{def.E2} and its estimate \eqref{estimate.E2}
\[
II_p\leq c {\Bigl[\omega(R)+\frac{R^\x}M\Bigr]}^{\iota\frac p2}\Big[\bigl(H_{x_m}\circ h_{x_m}^{-1}\bigr)(d_1)\Big]^\frac p2\Big[H_{x_m}(E_2)\Big]^\frac{2-p}2\leq c  {\Bigl[\omega(R)+\frac{R^\x}M\Bigr]}^{\iota\frac p2} \bigl(H_{x_m}\circ h_{x_m}^{-1}\bigr)(d_1)
\]
and analogously
\[
II_q\leq c {\Bigl[\omega(R)+\frac{R^\x}M\Bigr]}^{\iota\frac q2} \bigl(H_{x_m}\circ h_{x_m}^{-1}\bigr)(d_1),
\]
both the constants depending on $\data$ and $\gdata$ only. Finally, if $p<2\leq q$ then
\[
\mean{B_{R/8}}H_{x_m}(Dv-Dw)\dx \leq c\,I+c\,II_p
\]
and again the conclusion follows.
\end{proof}

Using the averaged estimate \eqref{averaged.Psi.inverse} and the scaling property \eqref{g_x.lambda}, we infer the following corollary:

\begin{corollary}\label{corollary.compa.rough.DuDw.pq}
Let $w\in v+W^{1,H_i(\cdot)}_0(B_{R/8})$ be the $A_i$-harmonic map from Proposition~\ref{pq-harm.rough1}, $\x$ as in \eqref{iccs} and $\iota\in(0,1)$ the exponent from Lemma~\ref{p-harm}. Then, for every 
\begin{equation}\label{el}
1\leq \el\leq \frac q{q-1}  
\end{equation}
we have
\begin{equation}\label{compa.Dw.Lh.pq}
{\biggl(\mean{B_{R/8}}\bigl[h_{x_m}(Dv-Dw)\bigr]^{\el}\dx\biggr)}^\frac1\el\leq c{\Bigl[\omega(R)+\frac{R^\x}M\Bigr]}^{\frac{\iota(p-1)}{\max\{p,2\}}}\biggl[\mean{B_R} h(x,Du)\dx+\frac{|\mu|(B_R)}{R^{n-1}}\biggr],
\end{equation}
the constant $c$ depending on $\data,\gdata$ and $\mathfrak l$ but not on $M$.
\end{corollary}

\subsection{Final rough estimates}

Merging all the previous intermediate results, we obtain the following final comparison proposition: 
\begin{proposition}\label{rough.comparison.DuDw}
Let $u$, $B_R$ and $v$ be as fixed at the beginning of Paragraph~\ref{Other.comparison}. Let $x_m$ be chosen as in \eqref{frozen.minimum}, and let $A_i$ and $H_i$ be given by \eqref{vf.frozen.ai} and \eqref{inf.smaller.H}, respectively. Then there exist an $A_i$-harmonic map $w\in v+W^{1,H_i(\cdot)}_0(B_{R/8})$ and a positive exponent
\begin{equation}\label{texttt.E}
 {\textnormal{\texttt{E}}}\equiv {\textnormal{\texttt{E}}}(\data,\gdata)\in(0,1)
\end{equation}
such that the following comparison estimates hold:
\begin{itemize}
\item Let \eqref{p.phase} hold for some $M\geq1$; then, with $c\equiv c(\data,\gdata, M)$,
\begin{equation}\label{compa.Du-Dw}
\mean{B_{R/16}}h(x,Du-Dw)\dx\leq c\biggl[\Bigl[\omega(R)+R^{s_\smallone}\Bigr]^{\textnormal{\texttt{E}}}\mean{B_R} h(x,Du)\dx+\frac{|\mu|(B_R)}{R^{n-1}}+d_{p,q}\biggr]. 
\end{equation}
Moreover, let $\es$ be as in \eqref{es}; with $c\equiv c(\data,\gdata,\es,M)$ if $p\leq n$, $c\equiv c(\data,\gdata,M)$ if $p>n$, it holds
\begin{equation}\label{rough.comparison.p.component.DuDw}
\biggl(\mean{B_{R/16}}|Du-Dw|^{(p-1)\es} \dx\biggr)^\frac1\es\leq c\biggl[\Bigl[\omega(R)+R^{s_\smallone}\Bigr]^{\textnormal{\texttt{E}}}\mean{B_R} h(x,Du)\dx+\frac{|\mu|(B_R)}{R^{n-1}}+d_{p,q}\biggr].
\end{equation}
\item Let \eqref{pq.phase} hold for some $M\geq 2$; then, for $\bi$ satisfying
\begin{equation}\label{bi}
0< \bi<\min\Bigl\{\frac n{n-1}, \frac q{q-1}\Bigr\},
\end{equation}
it holds
\begin{equation}\label{compa.Du-Dw,pq}
\biggl(\mean{B_{R/16}}{\bigl[h_{x_\smallzero}(Du-Dw)\bigr]}^\bi\dx\biggr)^\frac1\bi\leq c\biggl[{\Bigl[\omega(R)+\frac{R^\x}M\Bigr]}^{\textnormal{\texttt{E}}}\mean{B_R} h(x,Du)\dx+\frac{|\mu|(B_R)}{R^{n-1}}+d_{p,q}\biggr]
\end{equation}
with $c\equiv c(\data,\gdata,\bi)$.
\end{itemize}
$d_{p,q}$ is defined in \eqref{d_pq}; $\x$ is any exponent in \eqref{iccs} and $s_1$ is defined in \eqref{s1}.
\end{proposition}

\begin{proof}
We set 
\[
{\textnormal{\texttt{E}}}:=\frac{\iota(p-1)}{\max \{p,2\}};
\]
 for the estimates available on $B_{R/8}$ we still estimate on $B_{R/16}$ for the sake of readability, absorbing the fixed volume ratios into the constants. \eqref{compa.Du-Dw} follows from sub-additivity, using Lemma~\ref{lem:rough.compa} and Corollary~\ref{corollary.compa.rough.DuDw};  \eqref{rough.comparison.p.component.DuDw} follows from Lemma~\ref{lem:rough.compa.solopartep} and Proposition~\ref{p-harm.rough1} together with H\"older's inequality. \eqref{compa.Du-Dw,pq} follows by merging Lemma~\ref{lem:rough.compa.pq} and Corollary~\ref{corollary.compa.rough.DuDw.pq}.
\end{proof}

\begin{remark}\label{rem.obvious.Young}
{\rm Young's inequality gives
\[
d_{p,q}\leq c\mean{B_R}h(x,Du)\dx+c\,\frac{|\mu|(B_R)}{R^{n-1}}
\]
for a numerical constant $c$; this estimate will be used when applying Proposition~\ref{rough.comparison.DuDw} in a rougher form.
}
\end{remark}

\section{Linearized comparison estimates}\label{lin.chapter}

In this Section we prove the inductive step that will lead to Theorem~\ref{main.thm}. We work with a local energy solution $u\in W^{1,H(\cdot)}_{\rm loc}(\Omega)$ to \eqref{equation}, under the structural assumptions \eqref{A.strong}, \eqref{a.holder} to \eqref{bounds}; when $p\leq n$, we additionally assume that $\mu\in L^\infty_{\rm loc}(\Omega)$. For local SOLAs, we apply the arguments below to admissible approximating solutions, as explained after \eqref{def.v}. We shall assume a number of inequalities that will follow from the inductive hypothesis. Recalling the reference radius $\rho_{\rm ref}$ and the constant $c_{\rm ref}$ provided by Corollary~\ref{flux.decay}, together with the comparison radius $R_{\rm cmp}$ fixed above, we consider a ball $B_R\equiv B_R(x_0)$ with $x_0\in\Omega_0$ and $0<R\leq R_0\leq R_{\rm cmp}\leq\rho_{\rm ref}$, so that $B_R(x_0)\subset B_{R_0}(x_0)\Subset\Omega_1$, with $\Omega_1$ as defined in \eqref{Omega1}. Throughout this Section, $R_0$ will be further decreased several times, while $\rho_{\rm ref}$ and $c_{\rm ref}$ remain fixed. For a constant $\delta\in(0,1/32)$ to be defined, we assume that
\begin{equation}\label{basic.bounds}
\vartheta\lambda\leq\mean{B_{\delta R}}h(x,Du)\dx,\qquad\mean{B_R}h(x,Du)\dx\leq\lambda,\qquad\frac{|\mu|(B_R)}{R^{n-1}}\leq\lambda
\end{equation}
hold for some parameter $\lambda\geq1$ and some constant $\vartheta\in (0,1/2]$. Moreover, in order to deal at the same time with the two different cases, we set
\[
R_1=R,\qquad R_2=\delta R;
\] 
indeed, we are going to define comparison maps at two scales.

\medskip

Note that in view of the first of the bounds in \eqref{basic.bounds}, Corollary~\ref{flux.decay} ensures that
\begin{equation}\label{R.lambda}
\lambda\leq\frac{c_{\rm ref}}{\vartheta\delta^{n-1}}R^{-(n-1)}.
\end{equation}
From now on, we fix
\[
\x:=
\begin{cases}
\overline\x,& p\leq n,\\
\underline\x,& p>n,
\end{cases}
\]
where the positive exponents are defined before \eqref{iccs}; in particular, $\x$ depends only on $\data$ and $\gdata$. Clearly one of the following two inequalities holds:
\begin{equation}\label{two.different.cases}
a_i(B_R)\leq2[a]_{0,\alpha}R^{\alpha-\x}\qquad\text{or}\qquad a_i(B_R)>2[a]_{0,\alpha}R^{\alpha-\x}.
\end{equation}
If {\bf the first inequality} holds, that is, if we are in the $p$-phase, we take as regular comparison maps the $A_i$-harmonic functions given by Proposition~\ref{rough.comparison.DuDw}. In particular, we set for $k=1,2$:
\begin{itemize}
 \item $v_k\in u+W^{1,H(\cdot)}_0(B_{R_k/2})$ as the solution to the homogeneous problem \eqref{def.v} over $B_{R_k/2}$; 
 \item $w_k\in v_k+W^{1,H_{i,k}(\cdot)}_0(B_{R_k/8})$ as the $A_{i,k}$-harmonic function given by Proposition~\ref{rough.comparison.DuDw} over $B_{R_k/8}$; 
 \end{itemize}
here we have set, for $k=1,2$
\[
a_i(B_{R_k})=\min_{\overline{B_{R_k}}}a(\,\cdot\,)=a(x_{m,k}),\qquad A_{i,k}(z):=A(x_{m,k},z)
\]
\begin{equation}\label{h_ik}
H_{i,k}(z):=\frac{|z|^p}{p}+a(x_{m,k})\frac{|z|^q}{q},\qquad h_{i,k}(z):=|z|^{p-1}+a(x_{m,k})|z|^{q-1},\qquad  z\in\mathbb R^n.
\end{equation}
If \textbf{the second inequality} in \eqref{two.different.cases} holds, then $B_R$ is in the $(p,q)$-phase and for each $k\in\{1,2\}$, we take
\begin{itemize}
\item $v_k\in u+W^{1,H(\cdot)}_0(B_{R_k/2})$ to be the weak solution of \eqref{def.v} in $B_{R_k/2}$,
\item $w_k\in v_k+W^{1,H_{i,k}(\cdot)}_0(B_{R_k/8})$ to be the $A_{i,k}$-harmonic map given by Proposition~\ref{rough.comparison.DuDw}.
\end{itemize}
For each fixed $k$, the two descriptions give the same pair $v_k,w_k$: indeed, the Dirichlet problem \eqref{def.v} with boundary datum $u$ on $\partial B_{R_k/2}$ uniquely determines $v_k$. We use the same minimizing point $x_{m,k}$ in both descriptions, so the frozen field $A_{i,k}$ and the trace of $v_k$ on $\partial B_{R_k/8}$ uniquely determine the same map $w_k$. For clarity, we state the definition of $v_k$ separately in each phase and the same phase determines which estimates in Proposition~\ref{rough.comparison.DuDw} apply.

We stress that the phase conditions also apply on the smaller ball $B_{\delta R}$: indeed, in the $p$-phase, Remark~\ref{stability.p.phase} gives
\begin{equation}\label{p-phase.deltaR}
a_i(B_{\delta R})\leq4[a]_{0,\alpha}R^{\alpha-\x}=M_\delta[a]_{0,\alpha}(\delta R)^{\alpha-\x},\qquad M_\delta:=4\delta^{-(\alpha-\x)}.
\end{equation}
In the $(p,q)$-phase, Remark~\ref{stability.phases} instead gives
\[
a_i(B_{\delta R})\geq a_i(B_R)>2[a]_{0,\alpha}R^{\alpha-\x}\geq2[a]_{0,\alpha}(\delta R)^{\alpha-\x},
\]
since $\alpha-\x>0$. Thus the comparison estimates at scale $\delta R$ apply with $M=M_\delta$ in the $p$-phase and with $M=2$ in the $(p,q)$-phase; in the former case, their constants may also depend on $\delta$.

\medskip

We further assume that, for some $K\geq1$ to be chosen,
\begin{equation}\label{non.degeneracy}
\frac{\lambda}{K}\leq \inf_{B_{\delta R}}h_{i,1}(Dw_1),\qquad\sup_{B_{R/16}}h_{i,1}(Dw_1)\leq K\lambda:
\end{equation}
the role of and the reason for these non-degeneracy bounds will become clear when we implement the iteration scheme; see Section~\ref{proofs}.
\medskip

For every $t\geq0$ we define the reference functions
\[
\bar h(t):=h_{i,1}(t)=t^{p-1}+a_i(B_R)t^{q-1},\qquad\qquad \bar H(t):=H_{i,1}(t)=\frac{t^p}{p}+a_i(B_R)\frac{t^q}{q}.
\]
In the $(p,q)$-phase, \eqref{comparability.a} gives $\bar h(t)\leq h_{x_0}(t)\leq2\bar h(t)$ for every $t\geq0$. Consequently,
\[
h_{x_0}^{-1}(s)\leq\bar h^{-1}(s)\qquad\text{for every }s\geq0.
\]
Thus the scalar comparison estimates stated with the coefficient frozen at the center also give the corresponding upper bounds with the coefficient frozen at its minimum.

\begin{remark}[Normalization in the $p$-phase]\label{rem:p.phase.normalization}
{\rm Suppose that $p\leq n$ and that the first inequality in \eqref{two.different.cases} holds. By \eqref{p-phase.deltaR}, \eqref{R.lambda}, \eqref{iccs}-\eqref{bound.with.x} and the choice of $R_{\rm cmp}$, recalling that $R\leq R_{\rm cmp}$ and $s_1=s_0/p$, for $k\in\{1,2\}$ we have
\[
a_i(B_{R_k})\lambda^{\frac{q-p}{p-1}}\leq c\,C_{\rm ref}^{\frac{q-p}{p-1}}R^{\alpha-\x-\frac{(n-1)(q-p)}{p-1}}=c\,C_{\rm ref}^{\frac{q-p}{p-1}}R^{s_0}\leq cR^{s_1}\leq c,
\]
where $c$ depends only on $\data,\delta$ and $\vartheta$. Evaluating $h_{i,k}$ at $\lambda^{1/(p-1)}$ and using the scaling properties of its inverse and \eqref{H.h.comparison}, we obtain
\[
h_{i,k}^{-1}(\lambda)\approx\bar h^{-1}(\lambda)\approx\lambda^{\frac1{p-1}},\qquad \bigl(H_{i,k}\circ h_{i,k}^{-1}\bigr)(\lambda)\approx\bigl(\bar H\circ\bar h^{-1}\bigr)(\lambda)\approx\lambda^{p'}.
\]
The constants in the last display depend only on $\data,\gdata,\delta$ and $\vartheta$ too.
}
\end{remark}

We reduce $R_0$ so that
\begin{equation}
\omega(R)+R^{s_1}+R^\x\leq1\qquad\text{for every $R\leq R_0$,}
\end{equation}
and introduce the error density $d_0$, retaining the notation $d_1$ from \eqref{dtwo}. Here ${\textnormal{\texttt{E}}}\in(0,1)$ is the exponent from Proposition~\ref{rough.comparison.DuDw}. By \eqref{basic.bounds} and the preceding choice of $R_0$, we have
\begin{equation}\label{dzero}
0<d_0:=\Bigl[\omega(R)+R^{\min\{s_\smallone,\x\}}\Bigr]^{\textnormal{\texttt{E}}}\mean{B_R}h(x,Du)\dx+\frac{|\mu|(B_R)}{R^{n-1}}\leq d_1\leq2\lambda,
\end{equation}
\begin{equation}\label{done}
\vartheta\delta^n\lambda\leq d_1\leq2\lambda,
\end{equation}
and
\begin{equation}\label{dpq}
d_{p,q}\leq2\chi_{(-\infty,2)}(p)\,\lambda\biggl[\frac{|\mu|(B_R)}{\lambda R^{n-1}}\biggr]^{p-1}\leq2\lambda.
\end{equation}
Although the terms in $d_{p,q}$ can be made small relative to $\lambda$, their sublinear dependence on the measure density prevents their summation over scales from being controlled by the Riesz potential alone. The purpose of the following linearization argument is therefore to eliminate these terms and obtain a final comparison estimate with linear dependence on the measure density; the ultimate goal of this section is thus the forthcoming Proposition~\ref{def.comparison.pphase}.

\medskip

We now collect several comparison estimates that will be used below. 

\begin{lemma}[Comparison estimates for the variable growth function]\label{lem:comparison.variable}
In the setting above, let $w_1$ and $w_2$ be the frozen comparison maps associated with $B_R$ and $B_{\delta R}$, respectively. After reducing $R_0$ if necessary, we have
\begin{equation}\label{Dw_1-Dw_2.variable}
\sum_{k\in\{1,2\}}\mean{B_{R_k/16}}h(x,Du-Dw_k)\dx+\mean{B_{\delta R/16}}h(x,Dw_1-Dw_2)\dx\leq c\bigl[d_0+d_{p,q}\bigr],
\end{equation}
where $d_0$ is defined in \eqref{dzero}, with ${\textnormal{\texttt{E}}}$ the common exponent from Proposition~\ref{rough.comparison.DuDw}. The constant $c$ depends only on $\data,\gdata$ and $\delta$.
\end{lemma}

\begin{proof}
Applying Proposition~\ref{rough.comparison.DuDw} at scales $R$ and $\delta R$, and using coefficient comparability in the $(p,q)$-phase, gives
\[
\sum_{k=1}^2\mean{B_{R_k/16}}h(x,Du-Dw_k)\dx \leq c\bigl[d_0+d_{p,q}\bigr].
\]
Here we enlarged the integrals and measure term on the right-hand side to $B_R$, including the resulting factors into the constant depending on $\delta$, and used $\omega(\delta R)\leq\omega(R)$. Finally, by subadditivity,
\begin{align*}
\mean{B_{\delta R/16}}h(x,Dw_1-Dw_2)\dx&\leq c\sum_{k\in\{1,2\}}\mean{B_{\delta R/16}}h(x,Du-Dw_k)\dx\\
&\leq c\,\delta^{-n}\mean{B_{R/16}}h(x,Du-Dw_1)\dx+c\mean{B_{\delta R/16}}h(x,Du-Dw_2)\dx\\
&\leq c\bigl[d_0+d_{p,q}\bigr].
\end{align*}
\end{proof}

\begin{lemma}[Comparison at the smaller scale]\label{lem:comparison.small.scale}
In the setting above, let $v_2$ and $w_2$ be the comparison maps associated with the outer ball $B_{\delta R}$ and defined on $B_{\delta R/2}$ and $B_{\delta R/8}$, respectively. Let moreover ${\textnormal{\texttt{E}}}$ be the exponent from Proposition~\ref{rough.comparison.DuDw}. Then, for every
\[
0<\bi<\min\left\{\frac n{n-1},\frac q{q-1}\right\},
\]
there hold
\begin{equation}\label{Dv_2-Dw_2.hbi}
{\biggl(\mean{B_{\delta R/8}}{\bigl[\bar h(Dv_2-Dw_2)\bigr]}^\bi\dx\biggr)}^\frac1{\bi}\leq c\Bigl[\omega(R)+R^{\min\{s_\smallone,\x\}}\Bigr]^{\textnormal{\texttt{E}}} d_1,
\end{equation}
\begin{equation}\label{Du-Dv_2}
{\biggl(\mean{B_{\delta R/2}}{\bigl[\bar h(Du-Dv_2)\bigr]}^\bi\dx\biggr)}^\frac1\bi\leq c\biggl[\frac{|\mu|(B_R)}{R^{n-1}}+d_{p,q}\biggr],
\end{equation}
and
\begin{equation}\label{Dw_1-Dw_2}
\sum_{k\in\{1,2\}}{\biggl(\mean{B_{R_k/8}}{\bigl[\bar h(Du-Dw_k)\bigr]}^\bi\dx\biggr)}^\frac1\bi+{\biggl(\mean{B_{\delta R/8}}{\bigl[\bar h(Dw_1-Dw_2)\bigr]}^\bi\dx\biggr)}^\frac1{\bi}\leq c\bigl[d_0+d_{p,q}\bigr].
\end{equation}
The constants depend on $\data,\gdata,\delta$ and  $\bi$; the constants in \eqref{Dv_2-Dw_2.hbi} and \eqref{Dw_1-Dw_2} also depend on $\vartheta$.
\end{lemma}

\begin{proof}
By Jensen's inequality, it suffices to consider exponents $\bi\geq1$. We start by proving 
\begin{equation}\label{Dv_2-Dw_2}
\mean{B_{\delta R/8}}\bar H(Dv_2-Dw_2)\dx\leq c\bigl(\bar H\circ\bar h^{-1}\bigr)\biggl(\Bigl[\omega(R)+R^{\min\{s_\smallone,\x\}}\Bigr]^{\textnormal{\texttt{E}}} d_1\biggr)
\end{equation}
for a constant depending on $\data,\gdata,\delta$ and $\vartheta$. Using the monotonicity of $\inf a(\,\cdot\,)$ with respect to set inclusion, we have
\begin{equation}\label{Dv_2-Dw_2.part1}
\mean{B_{\delta R/8}}\bar H(Dv_2-Dw_2)\dx\leq \mean{B_{\delta R/8}}H_{i,2}(Dv_2-Dw_2)\dx.
\end{equation}
Now we need to distinguish the phases: if the second inequality in \eqref{two.different.cases} holds, then we use the comparison in Proposition~\ref{pq-harm.rough1}, the monotonicity of $\omega(\,\cdot\,)$ and \eqref{G_x.lambda}
\begin{align}\label{Dv_2-Dw_2.part2}
\mean{B_{\delta R/8}}H_{i,2}(Dv_2-Dw_2)\dx &\leq c\bigl(H_{i,2}\circ h_{i,2}^{-1}\bigr)\biggl(\Bigl[\omega(\delta R)+(\delta R)^\x\Bigr]^{\textnormal{\texttt{E}}}\biggl[\mean{B_{\delta R}} h(x,Du)\dx+\frac{|\mu|(B_{\delta R})}{(\delta R)^{n-1}}\biggr]\biggr)\notag\\
&\leq c\bigl(H_{i,2}\circ h_{i,2}^{-1}\bigr)\biggl(\Bigl[\omega(R)+R^\x\Bigr]^{\textnormal{\texttt{E}}} d_1\biggr)
\end{align}
with $c\equiv c(\data,\gdata,\delta)$; we then show how to replace the reference energy at the smaller scale replacing $H_{i,2}\circ h_{i,2}^{-1}$ with  $\bar H\circ\bar h^{-1}$. Since $a_i(B_R)\leq a_i(B_{\delta R})$, we have $h_{i,2}^{-1}(t)\leq\bar h^{-1}(t)$: but recalling from \eqref{H.h.comparison} that $H_{i,2}(s)\approx sh_{i,2}(s)$ (and the same for $\bar h$ and $\bar H$), we infer
\begin{equation}\label{Dv_2-Dw_2.part3}
\bigl(H_{i,2}\circ h_{i,2}^{-1}\bigr)(t)\leq\frac{t}p\,h_{i,2}^{-1}(t)\leq\frac{t}p\,\bar h^{-1}(t)\leq\frac qp\,\bigl(\bar H\circ\bar h^{-1}\bigr)(t) 
\end{equation}
for every $t\geq0$ and the proof, in this case, is complete using again \eqref{G_x.lambda}. If the first inequality in \eqref{two.different.cases} holds, we go back to the first line of the comparison estimate \eqref{compa.Dw.Lp} with $\delta R$ replacing $R$
\begin{equation}\label{Dv_2-Dw_2.part4}
\mean{B_{\delta R/8}}H_{i,2}(Dv_2-Dw_2)\dx\leq c\Bigl[\omega(\delta R)+(\delta R)^{s_\smallone}\Bigr]^{{\textnormal{\texttt{E}}} p'}\mean{B_{\delta R/4}} H(x,Dv_2)\dx 
\end{equation}
 and we estimate more carefully the quantity
\[
E_{1,\delta}:= \biggl(\mean{B_{\delta R/4}}H(x,Dv_2)\dx\biggr)^{\frac1p},
\]
or, better, $\bar h(E_{1,\delta})=E_{1,\delta}^{p-1}\bigl(1+a_i(B_R) E_{1,\delta}^{q-p}\bigr)$. The bound \eqref{en.rh.1}, applied at scale $\delta R$, gives
\[
E_{1,\delta}^{p-1}\leq c(\data,\gdata,\delta)d_1,
\]
for $p\leq n$, using $d_1\leq2\lambda$ and Remark~\ref{rem:p.phase.normalization}, we obtain
\[
a_i(B_R)E_{1,\delta}^{q-p}\leq c\,a_i(B_R)d_1^{\frac{q-p}{p-1}}\leq c\,a_i(B_R)\lambda^{\frac{q-p}{p-1}}\leq c(\data,\gdata,\vartheta,\delta);
\]
for $p>n$ instead, the argument leading to \eqref{rescaled.ai.bound}, applied at scale $\delta R$, gives the same bound. Thus, with $c\equiv c(\data,\gdata,\delta,\vartheta)$
\[
\bar h(E_{1,\delta})\leq c\,E_{1,\delta}^{p-1}\leq c\, d_1\quad\implies\quad E_{1,\delta}\leq \bar h^{-1}(c\,d_1),
\]
using the energy estimate \eqref{en.rh.1}, applied with radius $\delta R$. We hence obtain
\[
E_{1,\delta}^p\leq p\,\bar H(E_{1,\delta})\leq p(\bar H\circ\bar h^{-1})(c\,d_1)\leq c(\bar H\circ\bar h^{-1})(d_1) 
\]
and, appealing to \eqref{G_x.lambda},
\begin{align}\label{Dv_2-Dw_2.part5}
\Bigl[\omega(\delta R)+(\delta R)^{s_\smallone}\Bigr]^{{\textnormal{\texttt{E}}} p'}E_{1,\delta}^p&\leq c\Bigl[\omega(R)+R^{\min\{s_\smallone,\x\}}\Bigr]^{{\textnormal{\texttt{E}}} p'}(\bar H\circ\bar h^{-1})(d_1)\notag\\
&\leq c\,(\bar H\circ\bar h^{-1})\biggl(\Bigl[\omega(R)+R^{\min\{s_\smallone,\x\}}\Bigr]^{\textnormal{\texttt{E}}} d_1\biggr)
\end{align}
with $c\equiv c(\data,\gdata,\delta,\vartheta)$. Merging \eqref{Dv_2-Dw_2.part1} with \eqref{Dv_2-Dw_2.part2}-\eqref{Dv_2-Dw_2.part3} or \eqref{Dv_2-Dw_2.part4}-\eqref{Dv_2-Dw_2.part5} completes the proof of \eqref{Dv_2-Dw_2}. By H\"older's inequality and the growth properties of $\bar H$ and $\bar h$, we have
\[
\biggl(\mean{B_{\delta R/8}}[\bar h(Dv_2-Dw_2)]^\bi\dx\biggr)^{\frac1\bi}\leq c\,(\bar h\circ\bar H^{-1})\biggl(\mean{B_{\delta R/8}}\bar H(Dv_2-Dw_2)\dx\biggr)\leq c\Bigl[\omega(R)+R^{\min\{s_\smallone,\x\}}\Bigr]^{\textnormal{\texttt{E}}} d_1,
\]
for every $1\leq\bi<\min\{n',q'\}$; we used \eqref{Dv_2-Dw_2} and \eqref{g_x.lambda}, and this also proves \eqref{Dv_2-Dw_2.hbi}.

\medskip

Finally, to prove \eqref{Dw_1-Dw_2} we only need to show that 
\[
\biggl(\mean{B_{R_k/8}}\bigl[\bar h(Du-Dw_k)\bigr]^\bi\dx\biggr)^{\frac1\bi}\leq c\bigl[d_0+d_{p,q}\bigr]\qquad\text{for }k=1,2
\]
for a constant as in the statement of the Lemma, since a final application of subadditivity on $B_{\delta R/8}$ would establish the estimate for $Dw_1-Dw_2$ in \eqref{Dw_1-Dw_2}: restricting the integral average of $[\bar h(Du-Dw_1)]^\bi$ from $B_{R/8}$ to $B_{\delta R/8}$ would cost at most a geometric factor $\delta^{-n}$.

To prove the previous two estimates, we notice that Lemma~\ref{lem:rough.compa.solopartep} applied at radius $R_k$ controls the $p$-part of $\bar h(Du-Dv_k)$ for $k=1,2$. As for the $q$-part, since $a_i(B_R)\leq a(x)$ in $B_{R_k}$, we can draw upon \eqref{compa.q.component.ell} with $\esp=q$ and $\bar a_q=a_i(B_R)$. When $q\geq2$, we simply set $\ti=\bi$; when $q<2$, we take $\ti=\max\{\bi,1/(q-1)\}<n'$, appeal to \eqref{component.local.rH} with the same coefficient and exponent for the $Du$-term, and step back to $\bi$ via Hölder's inequality. Combining both contributions and enlarging the integration domains on the right-hand side to $B_R$ yields \eqref{Du-Dv_2} for $k=2$, together with the corresponding bound for $Du-Dv_1$ on $B_{R/2}$, where $c$ depends only on $\data,\delta$, and $\bi$. Together with \eqref{Dv_2-Dw_2.hbi} and its counterpart at scale $R$, these bounds give - using the triangle inequality and the definition of $d_0$ in \eqref{dzero} - \eqref{Dw_1-Dw_2}.

\end{proof}

The next lemma gives comparison estimates between $Du$, $Dv_2$ and $Dw_2$; those are a linearized version of the comparison estimates in Propositions~\ref{p-harm.rough1}-\ref{pq-harm.rough1}.

\begin{lemma}[A refined comparison for $Dv_2-Dw_2$]\label{compa.Dv_2-Dw_2}
Suppose that \eqref{basic.bounds} and \eqref{non.degeneracy} are in force for some parameter $\lambda\geq1$ and some constant $K\geq1$. The following linearized comparison estimate holds, for a constant $c$ depending on $\data,\gdata,\delta$ and $\vartheta$:
\[
\mean{B_{\delta R/8}}|Dv_2-Dw_2|\dx\leq c\,\frac{\bar h^{-1}(\lambda)}\lambda\biggl[\Bigl[\omega(R)+R^{\min\{s_\smallone,\x\}}\Bigr]^\frac{{\textnormal{\texttt{E}}}}{q-1}\mean{B_R} h(x,Du)\dx+\frac{|\mu|(B_R)}{R^{n-1}}\biggr].
\]
\end{lemma}

\begin{proof}
By convexity and \eqref{Dv_2-Dw_2} we have
\[
\bar H\biggl(\mean{B_{\delta R/8}}|Dv_2-Dw_2|\dx\biggr)\leq \mean{B_{\delta R/8}}\bar H(Dv_2-Dw_2)\dx\leq c \,\bigl(\bar H\circ\bar h^{-1}\bigr)\biggl(\Bigl[\omega(R)+R^{\min\{s_1,\x\}}\Bigr]^{\textnormal{\texttt{E}}} d_1\biggr);
\]
from this, using \eqref{constant.hinverse}, we find
\[
\mean{B_{\delta R/8}}|Dv_2-Dw_2|\dx\leq c\,\bar h^{-1}\biggl(\Bigl[\omega(R)+R^{\min\{s_1,\x\}}\Bigr]^{\textnormal{\texttt{E}}} d_1\biggr)\leq c\Bigl[\omega(R)+R^{\min\{s_\smallone,\x\}}\Bigr]^\frac{{\textnormal{\texttt{E}}}}{q-1}\bar h^{-1}(d_1);
\]
the fact
\[
\bar h^{-1}(d_1)\leq c(n,p,\delta,\vartheta) \frac{\bar h^{-1}(\lambda)}\lambda d_1
\]
is immediate from the two-sided bound for $d_1$ in \eqref{done} and the monotonicity of $\bar h^{-1}$; this concludes the proof.
\end{proof}

We now use the non-degeneracy of $w_1$ to convert the preceding flux comparisons into an $L^1$ estimate for $Du-Dv_2$, with the linearization factor $\bar h^{-1}(\lambda)/\lambda$.

\begin{lemma}[A refined comparison Lemma for $Du-Dv_2$]\label{lem:refined.1}
Suppose that \eqref{basic.bounds} holds for some parameter $\lambda\geq1$ and  that \eqref{non.degeneracy} holds for some constant $K\geq1$. Then 
\begin{equation}\label{nice.compa.pphase}
\mean{B_{\delta R/8}}|Du-Dv_2|\dx\leq c\,\frac{\bar h^{-1}(\lambda)}\lambda\biggl[\Bigl[\omega(R)+R^{\min\{s_\smallone,\x\}}\Bigr]^{\textnormal{\texttt{E}}}\mean{B_R} h(x,Du)\dx+\frac{|\mu|(B_R)}{R^{n-1}}\biggr]
\end{equation}
for a constant $c$ depending on $\data,\gdata,\delta,\vartheta$ and $K$. ${\textnormal{\texttt{E}}}$ is the exponent in \eqref{texttt.E}.
\end{lemma}

\begin{proof}
We consider separately three different cases. 

\subsubsection*{The case $2< p\leq q$}

We notice that in this case $d_{p,q}=0$. We fix two exponents whose values are
\begin{equation}\label{xi.pq}
\bi_*:=\frac12\big(1+\min\{n',q'\}\big),\qquad \xi=1+2\gamma\quad\text{with}\quad \gamma=\min\Big\{\frac 1{2(n-1)},\frac1{2(q-1)},\frac{n'}{1+\frac{q-1}{n-1}},\frac{\bi_*-1}{4}\Big\},
\end{equation}
whence
\begin{equation}\label{bound.1+gamma}
1<1+\gamma<\xi<\bi_*< \min\big\{n',q'\big\}.
\end{equation}
We use \eqref{non.degeneracy} as follows: recalling the monotonicity and the sub-additivity of the map $t\mapsto [\bar h(t)/t]^{1+\gamma}$, we obtain
\begin{align}\label{est.L1.pq}
\mean{B_{\delta R/8}}|Du-Dv_2|\dx & \leq c\,\Big[\frac{\bar h^{-1}(\lambda)}\lambda\Big]^{1+\gamma}\mean{B_{\delta R/8}}\Big[\frac{\bar h(Dw_1)}{|Dw_1|}\Big]^{1+\gamma}|Du-Dv_2|\dx\notag\\
&\leq c\,\Big[\frac{\bar h^{-1}(\lambda)}\lambda\Big]^{1+\gamma}\mean{B_{\delta R/8}}\Big[\frac{\bar h(Dw_1-Dw_2)}{|Dw_1-Dw_2|}\Big]^{1+\gamma}|Du-Dv_2|\dx\notag\\
&\qquad\qquad+c\,\Big[\frac{\bar h^{-1}(\lambda)}\lambda\Big]^{1+\gamma}\mean{B_{\delta R/8}}\Big[\frac{\bar h(Dv_2-Dw_2)}{|Dv_2-Dw_2|}\Big]^{1+\gamma}|Du-Dv_2|\dx\notag\\
&\qquad\qquad\qquad\qquad+c\,\Big[\frac{\bar h^{-1}(\lambda)}\lambda\Big]^{1+\gamma}\mean{B_{\delta R/2}}\Big[\frac{\bar h(Dv_2)}{|Dv_2|}\Big]^{1+\gamma}|Du-Dv_2|\dx\notag\\
& \leq c\,\Big[\frac{\bar h^{-1}(\lambda)}\lambda\Big]^{1+\gamma}\bigl[I_p+I_q+II_p+II_q+III_p+III_q\bigr],
\end{align}
where $c\equiv c(n,p,q,K)$. With 
\[
\bar a_p:=1\qquad\text{and}\qquad \bar a_q:=a_i(B_R),
\]
we set
\begin{equation}\label{I_est}
I_\esp:=\mean{B_{\delta R/8}}{\bar a}_{\esp}^{1+\gamma}{|Dw_1-Dw_2|}^{(\esp-2)(1+\gamma)}|Du-Dv_2|\dx,\qquad \esp \in\{p,q\};
\end{equation}
similarly we define $II_p,II_q$ replacing $Dw_1-Dw_2$ in $I_p,I_q$ with $Dv_2-Dw_2$, and $III_p,III_q$, replacing $Dw_1-Dw_2$ with $Dv_2$ and taking the average over the larger ball $B_{\delta R/2}$. We have, using H\"older's inequality with conjugate exponents $(\frac{\esp-1}{\esp-2},\esp-1)$,  \eqref{below.function} and \eqref{Dw_1-Dw_2} and \eqref{Du-Dv_2} (the use of the first one is allowed by \eqref{bound.1+gamma}), and finally the estimate for $d_0$ in \eqref{dzero},
\begin{align*}
I_\esp&\leq c(n,q,\delta){\bar a}_{\esp}^{1+\gamma}\biggl(\mean{B_{\delta R/8}}{|Dw_1-Dw_2|}^{(\esp-1)(1+\gamma)}\dx\biggr)^{\frac{\esp-2}{\esp-1}}\biggl(\mean{B_{\delta R/2}}{|Du-Dv_2|}^{\esp-1} \dx\biggr)^{\frac1{\esp-1}}\notag\\
&\leq c\,{{\bar a}_{\esp}}^{\frac\gamma{\esp-1}}\biggl(\mean{B_{\delta R/8}}{\bigl[ \bar h(Dw_1-Dw_2)\bigr]}^{1+\gamma}\dx\biggr)^{\frac{\esp-2}{\esp-1}}\biggl(\mean{B_{\delta R/2 }}\bar h(Du-Dv_2) \dx\biggr)^{\frac1{\esp-1}}\notag\\
&\leq c\,{{\bar a}_{\esp}}^{\frac\gamma{\esp-1}}d_0^{(1+\gamma)\frac{\esp-2}{\esp-1}}{\biggl[\frac{|\mu|(B_R)}{R^{n-1}}\biggr]}^{\frac{1}{\esp-1}}\leq c\,{{\bar a}_{\esp}}^{\frac\gamma{\esp-1}}d_0^{1+\gamma\frac{\esp-2}{\esp-1}}\leq c\,{{\bar a}_{\esp}}^{\frac\gamma{\esp-1}}d_0\lambda^{\gamma\frac{\esp-2}{\esp-1}}
\end{align*}
with $c\equiv c(\data,\gdata,\delta,\vartheta)$. Thus
\begin{align}\label{I_pq}
 \mean{B_{\delta R/8}}\Big[\frac{\bar h(Dw_1-Dw_2)}{|Dw_1-Dw_2|}\Big]^{1+\gamma}|Du-Dv_2|\dx&\leq c(n,q)\bigl[ I_p+I_q \bigr]\notag\\
& \leq c\,d_0\Big[\frac{\lambda}{\lambda^\frac1{p-1}}+\frac{{[\bar a_q]}^{\frac1{q-1}}\lambda}{\lambda^{\frac{1}{q-1}}}\Big]^\gamma \leq c\,d_0\Big[\frac{\lambda}{\bar h^{-1}(\lambda)}\Big]^\gamma
\end{align}
thanks to \eqref{above.below.inverseh} applied at $x=x_{m,1}$, so that $h_{x_{m,1}}=\bar h$. Similarly, replacing \eqref{Dw_1-Dw_2} with \eqref{Dv_2-Dw_2.hbi} and employing the estimate 
\[
\max\biggl\{\frac{|\mu|(B_R)}{R^{n-1}},\Bigl[\omega(R)+R^{\min\{s_1,\x\}}\Bigr]^{\textnormal{\texttt{E}}}d_1\biggr\}\leq d_0\leq2\lambda
\]
as in \eqref{dzero}, we have
\begin{align}\label{II_pq}
\mean{B_{\delta R/8}}\Big[\frac{\bar h(Dv_2-Dw_2)}{|Dv_2-Dw_2|}\Big]^{1+\gamma}|Du-Dv_2|\dx&\leq c(n)\bigl[II_p+II_q\bigr]\notag\\
&\leq c\sum_{\esp\in\{p,q\}}{{\bar a}_{\esp}}^{\frac\gamma{\esp-1}}\Bigl[\omega(R)+R^{\min\{s_\smallone,\x\}}\Bigr]^{\textnormal{\texttt{E}}(1+\gamma)\frac{\esp-2}{\esp-1}}d_1^{(1+\gamma)\frac{\esp-2}{\esp-1}}d_0^{\frac{1}{\esp-1}}\notag\\
&\leq c\sum_{\esp\in\{p,q\}}{{\bar a}_{\esp}}^{\frac\gamma{\esp-1}}d_0^{1+\gamma\frac{\esp-2}{\esp-1}}\leq c\,d_0\Big[\frac{\lambda}{\bar h^{-1}(\lambda)}\Big]^\gamma.
\end{align}
The estimates for $III_p$ and $III_q$ require more effort. For $\esp \in\{p,q\}, \esp\geq 2$ we start by observing that from \eqref{diffV} we have
\[
|Dv_2|^{\esp-2}|Du-Dv_2|\leq  |Dv_2|^{\frac{\esp-2}2}\bigl[|Du|+|Dv_2|\bigr]^{\frac{\esp-2}2} |Du-Dv_2|\leq c(\esp) |Dv_2|^{\frac{\esp-2}2} \big|V_\esp(Du)-V_\esp(Dv_2)\big| ; 
\]
we estimate, for $\eta_{\esp}>0$ that will be chosen below
\begin{align*}
III_\esp&\leq c\,{\bar a}_{\esp}^{1+\gamma}\mean{B_{\delta R/2}}{|Dv_2|}^{\frac{\esp-2}2\xi}{(\eta_\esp+|u-v_2|)}^{\frac\xi2}\cdot\frac{|V_\esp(Du)-V_\esp(Dv_2)|}{{(\eta_\esp+|u-v_2|)}^{\frac\xi2}}\dx\notag\\ 
&\leq c\,\biggl( \mean{B_{\delta R/2}}{\bar a}_{\esp}^\xi{|Dv_2|}^{(\esp-2)\xi}(\eta_\esp+|u-v_2|)^\xi\dx\biggr)^{\frac12}\biggl(\mean{B_{\delta R/2}}\frac{{\bar a}_{\esp}|V_\esp(Du)-V_\esp(Dv_2)|^2}{(\eta_\esp+|u-v_2|)^\xi}\dx\biggr)^{\frac12}\notag\\
&\leq c\biggl(\mean{B_{\delta R/2}}{\bigl[{\bar a}_{\esp}{|Dv_2|}^{\esp-2}\bigr]}^\xi(\eta_\esp+|u-v_2|)^\xi\dx\biggr)^{\frac12}\eta_\esp^{\frac{1-\xi}2}\biggl[\frac{|\mu|(B_R)}{R^n}\biggr]^{\frac 12}
\end{align*}
for $c\equiv c(n,p,q,\nu,\delta)$; indeed from Lemma~\ref{comparison.uv.weighted} and $\bar a_q\leq a(\,\cdot\,)$ in $B_{\delta R/2}$, we have for $\esp \in\{p,q\}$
\[
\mean{B_{\delta R/2}}\frac{{\bar a}_{\esp}|V_\esp(Du)-V_\esp(Dv_2)|^2}{(\eta_\esp+|u-v_2|)^\xi}\dx\leq c\,\eta_\esp^{1-\xi}\frac{|\mu|(B_{\delta R/2})}{(\delta R/2)^n}\leq c\,\eta_\esp^{1-\xi}\frac{|\mu|(B_R)}{R^n},
\]
for $\xi>1$ as in \eqref{xi.pq} and the constant having the same dependencies as before. We take 
\begin{equation}\label{def.eta_t}
\eta_\esp^\xi:=\Big[\frac{\bar h^{-1}(\lambda)}\lambda\Big]^\xi\mean{B_{\delta R/2}}{\bigl[{\bar a}_{\esp}{|Dv_2|}^{\esp-2}\bigr]}^\xi|u-v_2|^\xi\dx+\varepsilon^\xi R^\xi\Big[\frac{\bar h^{-1}(\lambda)}\lambda\Big]^{(1+\gamma)\xi},
\end{equation}
with $\varepsilon>0$ that will eventually be sent to zero at the end of the proof; hence $\eta_\esp^\xi>0$. Now we want to estimate $III_\esp$ in terms of $\eta_\esp^\xi$: we first use \eqref{Dw_1-Dw_2} and \eqref{Du-Dv_2} to obtain
\begin{align*}
\sum_{\esp\in\{p,q\}}\mean{B_{\delta R/2}}\bigl[{\bar a}_{\esp}|Dw_1&-Dv_2|^{\esp-2}\bigr]^\xi\dx\\
&\leq c\sum_{\esp\in\{p,q\}}{\bar a}_{\esp}^{\frac\xi{\esp-1}}{\biggl(\mean{B_{\delta R/2}}{\bigl[\bar h(Dw_1-Dv_2)\bigr]}^\xi\dx\biggr)}^{\frac{\esp-2}{\esp-1}}\\
&\leq c\sum_{\esp\in\{p,q\}}{\bar a}_{\esp}^{\frac\xi{\esp-1}}\biggl[\mean{B_{R/8}}{\bigl[\bar h(Du-Dw_1)\bigr]}^\xi\dx+\mean{B_{\delta R/2}}{\bigl[\bar h(Du-Dv_2)\bigr]}^\xi\dx\biggr]^{\frac{\esp-2}{\esp-1}}\\
&\leq c\sum_{\esp\in\{p,q\}}{\bar a}_{\esp}^{\frac\xi{\esp-1}}\biggl[d_0+\frac{|\mu|(B_R)}{R^{n-1}}\biggr]^{\xi\frac{\esp-2}{\esp-1}}\\
&\leq c\Big[d_0^{\frac{p-2}{p-1}}+\bar a_q^{\frac1{q-1}}d_0^{\frac{q-2}{q-1}}\Big]^\xi\leq c\Big[\lambda^{\frac{p-2}{p-1}}+\bar a_q^{\frac1{q-1}}\lambda^{\frac{q-2}{q-1}}\Big]^\xi\leq c\Big[\frac{\lambda}{\bar h^{-1}(\lambda)}\Big]^\xi,
\end{align*}
using the upper bound in \eqref{dzero} and estimating as in \eqref{I_pq}, by \eqref{above.below.inverseh}. Consequently, by sub-additivity and \eqref{non.degeneracy},
\begin{align*}
\mean{B_{\delta R/2}}\bigl[{\bar a}_{\esp}|Dv_2|^{\esp-2}\bigr]^\xi\dx&\leq c(n,q)\mean{B_{\delta R/2}}\bigl[{\bar a}_{\esp}|Dw_1|^{\esp-2}\bigr]^\xi\dx+c(n,q)\mean{B_{\delta R/2}}\bigl[{\bar a}_{\esp}|Dw_1-Dv_2|^{\esp-2}\bigr]^\xi\dx\\
&\leq c\Big[\sup_{B_{\delta R/2}}\frac{\bar h(Dw_1)}{|Dw_1|}\Big]^\xi+c\mean{B_{\delta R/2}}\bigl[{\bar a}_{\esp}|Dw_1-Dv_2|^{\esp-2}\bigr]^\xi\dx\leq c\Big[\frac{\lambda}{\bar h^{-1}(\lambda)}\Big]^\xi,
\end{align*}
for a constant $c$ depending on $\data,\gdata,\delta,\vartheta$ and $K$. With the choices of $\eta_\esp$ just made, this gives
\begin{align*}
\mean{B_{\delta R/2}}\bigl[{\bar a}_{\esp}|Dv_2|^{\esp-2}\bigr]^\xi(\eta_\esp+|u-v_2|)^\xi\dx&\leq c(\gamma)\eta_\esp^\xi\mean{B_{\delta R/2}}\bigl[{\bar a}_{\esp}|Dv_2|^{\esp-2}\bigr]^\xi\dx+c(\gamma)\eta_\esp^\xi\Big[\frac{\lambda}{\bar h^{-1}(\lambda)}\Big]^\xi\\
&\leq c\,\eta_\esp^\xi\Big[\frac{\lambda}{\bar h^{-1}(\lambda)}\Big]^\xi.
\end{align*}
Inserting this into the preceding estimate for $III_\esp$, we finally obtain
\begin{equation}\label{III2.pq}
III_\esp\leq c\,\eta_\esp^{\frac12}\Big[\frac{\lambda}{\bar h^{-1}(\lambda)}\Big]^{\frac\xi2}\biggl[\frac{|\mu|(B_R)}{R^n}\biggr]^{\frac12}.
\end{equation}
We estimate now $\eta_\esp$; to do so, we use subadditivity in the following way
\begin{multline*}
\mean{B_{\delta R/2}}{\bigl[{\bar a}_{\esp}{|Dv_2|}^{\esp-2}\bigr]}^\xi|u-v_2|^\xi\dx \leq c\mean{B_{\delta R/2}}{\bigl[{\bar a}_{\esp}{|Dw_1|}^{\esp-2}\bigr]}^\xi|u-v_2|^\xi\dx\\
+c\mean{B_{\delta R/2}}{\bigl[{\bar a}_{\esp}{|Du-Dv_2|}^{\esp-2}\bigr]}^\xi|u-v_2|^\xi\dx+c\mean{B_{\delta R/2}}{\bigl[{\bar a}_{\esp}{|Du-Dw_1|}^{\esp-2}\bigr]}^\xi|u-v_2|^\xi\dx\\
=c(n,q)\bigl[IV_\esp+V_\esp+VI_\esp\bigr].
\end{multline*}
The estimate for the second term involves H\"older's and Sobolev-Poincar\'e inequalities; recalling the definition of $\bi_*$ in \eqref{xi.pq}, we set, for $\esp\in\{p,q\}$,
\[
\zeta_{\xi,\esp}:=
\begin{cases}
\displaystyle{\frac{[\bi_*(\esp-1)]^*}{[\bi_*(\esp-1)]^*-\xi}}&\text{if $\esp\leq n$},\\[3mm]
\displaystyle{\frac{n(\esp-1)}{n(\esp-1)-\xi(n-1)}}&\text{if $\esp>n$}.
\end{cases}
\]
Some remarks are now in order: first, notice that if $\esp\leq n$, then $1<\bi_*(\esp-1)<n$ and
\[
[\bi_*(\esp-1)]^*=\frac{n\bi_*(\esp-1)}{n-\bi_*(\esp-1)}>n'(\esp-1)>\xi(\esp-1),
\]
where the first inequality follows from $\bi_*(n+\esp-2)>n$. If $\esp>n$, the exponent $\xi\zeta_{\xi,\esp}'=n'(\esp-1)$ is also admissible in Sobolev's inequality: when $\bi_*(\esp-1)<n$, this follows from the same computation, while for $\bi_*(\esp-1)\geq n$ every finite target exponent is admissible. Thus $\zeta_{\xi,\esp}\in(1,+\infty)$ in both cases. Moreover, since $\xi\zeta_{\xi,\esp}'>\xi(\esp-1)$, we have
\begin{equation}\label{gamma.small.improved}
0\leq\xi\frac{\esp-2}{\esp-1}\zeta_{\xi,\esp}<\xi<\bi_*<\min\{n',q'\}.
\end{equation}
Consequently, all the flux exponents appearing below are admissible; those smaller than one are controlled by H\"older's inequality using the estimate with exponent one. By using H\"older's inequality with conjugate exponents $(\zeta_{\xi,\esp},\zeta_{\xi,\esp}')$, we get
\begin{align*}
V_\esp&\leq  c\,{\bar a}_{\esp}^\xi{\biggl(\mean{B_{\delta R/2}}{|Du-Dv_2|}^{(\esp-2)\xi\zeta_{\smallxi,\esp}}\dx\biggr)}^{\frac{1}{\zeta_{\smallxi,\esp}}}{\biggl(\mean{B_{\delta R/2}}|u-v_2|^{\xi\zeta_{\xi,\esp}'}\dx\biggr)}^\frac1{\zeta_{\xi,\esp}'}\\
&\leq  c\,{\bar a}_{\esp}^\xi R^\xi{\biggl(\mean{B_{\delta R/2}}{|Du-Dv_2|}^{(\esp-2)\xi\zeta_{\smallxi,\esp}}\dx\biggr)}^{\frac{1}{\zeta_{\smallxi,\esp}}}{\biggl(\mean{B_{\delta R/2}}|Du-Dv_2|^{\bi_*(\esp-1)}\dx\biggr)}^{\frac\xi{\bi_*(\esp-1)}}\\
&\leq  c\,R^\xi{\biggl(\mean{B_{\delta R/2}}\bigl[\bar h(Du-Dv_2)\bigr]^{\frac{\esp-2}{\esp-1}\xi\zeta_{\smallxi,\esp}}\dx\biggr)}^{\frac{1}{\zeta_{\smallxi,\esp}}}\biggl(\mean{B_{\delta R/2}}\bigl[\bar h(Du-Dv_2)\bigr]^{\bi_*}\dx\biggr)^{\frac\xi{\bi_*(\esp-1)}}\\
&\leq c\,R^\xi {\Bigl[\frac{|\mu|(B_{\delta R})}{(\delta R)^{n-1}}\Bigr]}^{\xi[\frac{\esp-2}{\esp-1}+\frac1{\esp-1}]}\leq c\,R^\xi d_0^\xi;
\end{align*}
we can indeed use \eqref{Du-Dv_2} with both the values $\bi=\bi_*$ and $\bi=\xi\frac{\esp-2}{\esp-1}\zeta_{\smallxi,\esp}$ thanks to the discussion made above. Similarly, again from \eqref{Du-Dv_2} and \eqref{Dw_1-Dw_2}
\[
VI_\esp\leq  c\,R^\xi\biggl(\mean{B_{R/8}}{\bigl[{\bar a}_{\esp}{|Du-Dw_1|}^{\esp-1}\bigr]}^{\frac{\esp-2}{\esp-1}\xi\zeta_{\smallxi,\esp}}\dx\biggr)^{\frac{1}{\zeta_{\smallxi,\esp}}}\biggl(\mean{B_{\delta R/2}}{\bigl[{\bar a}_{\esp}{|Du-Dv_2|}^{\esp-1}\bigr]}^{\bi_*}\dx\biggr)^{\frac\xi{\bi_*(\esp-1)}}\leq c\,R^\xi d_0^\xi.
\]
For $IV_\esp$, finally, we use \eqref{non.degeneracy} together with Sobolev's inequality (being $\xi\leq n'$) and the following trivial fact: in $B_{\delta R}$ there holds $\bar h(Dw_1) \approx_{K} \lambda,  |Dw_1|\approx_{K}\bar h^{-1}(\lambda)$. We have, also using sub-additivity
\begin{align*}
IV_\esp&\leq c\Big[\frac{\lambda}{\bar h^{-1}(\lambda)}\Big]^\xi\mean{B_{\delta R/2}}|u-v_2|^\xi\dx\\
&\leq c\,R^\xi\Big[\frac{\lambda}{\bar h^{-1}(\lambda)}\Big]^\xi{\biggl(\mean{B_{\delta R/2}}{|Du-Dv_2|}\dx\biggr)}^\xi\\
&\leq c\,R^\xi\Big[\frac{\lambda}{\bar h^{-1}(\lambda)}\Big]^{-\gamma\xi}{\biggl(\mean{B_{\delta R/2}}\Big[\frac{\bar h(Dw_1)}{|Dw_1|}\Big]^{1+\gamma}{|Du-Dv_2|}\dx\biggr)}^\xi\\
&\leq c\,R^\xi\Big[\frac{\lambda}{\bar h^{-1}(\lambda)}\Big]^{-\gamma\xi}{\biggl(\mean{B_{\delta R/2}}\Big[\frac{\bar h(Dv_2)}{|Dv_2|}\Big]^{1+\gamma}{|Du-Dv_2|}\dx\biggr)}^\xi\\
&\qquad+c\,R^\xi\Big[\frac{\lambda}{\bar h^{-1}(\lambda)}\Big]^{-\gamma\xi}{\biggl(\mean{B_{\delta R/2}}\Big[\frac{\bar h(Du-Dw_1)}{|Du-Dw_1|}\Big]^{1+\gamma}{|Du-Dv_2|}\dx\biggr)}^\xi\\
&\qquad\qquad+c\,R^\xi\Big[\frac{\lambda}{\bar h^{-1}(\lambda)}\Big]^{-\gamma\xi}{\biggl(\mean{B_{\delta R/2}}\frac{[\bar h(Du-Dv_2)]^{1+\gamma}}{|Du-Dv_2|^\gamma}\dx\biggr)}^\xi\\
&\leq c\,R^\xi\Big[\frac{\lambda}{\bar h^{-1}(\lambda)}\Big]^{-\gamma\xi}{\biggl(\mean{B_{\delta R/2}}\Big[\frac{\bar h(Dv_2)}{|Dv_2|}\Big]^{1+\gamma}{|Du-Dv_2|}\dx\biggr)}^\xi+c\,R^\xi d_0^\xi;
\end{align*}
we show how to treat the first comparison term in the fifth line of the display; the other terms are treated similarly. We have
\begin{align*}
\mean{B_{\delta R/2}}&\Big[\frac{\bar h(Du-Dw_1)}{|Du-Dw_1|}\Big]^{1+\gamma}{|Du-Dv_2|}\dx\\
&\leq c\sum_{\esp \in\{p,q\}} \mean{B_{\delta R/2}}{\bigl[{\bar a}_{\esp}{|Du-Dw_1|}^{\esp-2}\bigr]}^{1+\gamma}|Du-Dv_2|\dx\\
&\leq c\sum_{\esp \in\{p,q\}} {\bar a}_{\esp}^{\frac{\gamma}{\esp-1}}{\biggl(\mean{B_{R/8}}{\bigl[{\bar a}_{\esp}{|Du-Dw_1|}^{\esp-1}\bigr]}^{1+\gamma}\dx\biggr)}^{\frac{\esp-2}{\esp-1}}{\biggl(\mean{B_{\delta R/2}}{\bar a}_{\esp}{|Du-Dv_2|}^{\esp-1}\dx\biggr)}^{\frac {1}{\esp-1}}\\
&\leq c\sum_{\esp \in\{p,q\}} {\bar a}_{\esp}^{\frac{\gamma}{\esp-1}}d_0^{\frac{\esp-2}{\esp-1}(1+\gamma)}{\biggl[\frac{|\mu|(B_R)}{R^{n-1}}+d_{p,q}\biggr]}^\frac1{\esp-1}\leq c\,d_0\biggl[d_0^\frac{p-2}{p-1}+\bar a_q^{\frac1{q-1}}d_0^{\frac{q-2}{q-1}}\biggr]^{\gamma}\\
&\leq c\,d_0\biggl[\lambda^\frac{p-2}{p-1}+\bar a_q^{\frac1{q-1}}\lambda^{\frac{q-2}{q-1}}\biggr]^{\gamma}\leq c\, d_0\Big[\frac{\lambda}{\bar h^{-1}(\lambda)}\Big]^{\gamma}
\end{align*}
using \eqref{Dw_1-Dw_2}, \eqref{Du-Dv_2} and \eqref{above.below.inverseh}. The estimate for
\[
\mean{B_{\delta R/2}}\frac{[\bar h(Du-Dv_2)]^{1+\gamma}}{|Du-Dv_2|^\gamma}\dx
\]
is exactly the same, or even easier - since H\"older's inequality is not needed. Merging now the estimates for $IV_\esp,V_\esp,VI_\esp$, $\esp \in\{p,q\}$, yields
\[
\bigl[IV_\esp+V_\esp+VI_\esp\bigr]^\frac1\xi\leq c\,R\Big[\frac{\lambda}{\bar h^{-1}(\lambda)}\Big]^{-\gamma}\mean{B_{\delta R/2}}\Big[\frac{\bar h(Dv_2)}{|Dv_2|}\Big]^{1+\gamma}{|Du-Dv_2|}\dx+c\,R d_0;
\]
inserting these into the definition of $\eta_\esp$ in \eqref{def.eta_t} yields
\begin{align*}
\eta_\esp&\leq c\frac{\bar h^{-1}(\lambda)}\lambda\bigl[IV_\esp+V_\esp+VI_\esp\bigr]^\frac1\xi+\varepsilon R\Bigl[\frac{\lambda}{\bar h^{-1}(\lambda)}\Bigr]^{-(1+\gamma)}\\
&\leq c\,R\Big[\frac{\bar h^{-1}(\lambda)}\lambda\Big]^{1+\gamma}\bigl[III_p+III_q\bigr]+c\,R d_0\frac{\bar h^{-1}(\lambda)}\lambda+\varepsilon R\Big[\frac\lambda{\bar h^{-1}(\lambda)}\Big]^{-(1+\gamma)},
\end{align*}
and in turn employing \eqref{III2.pq} and then Young's inequality leads to
\begin{align*}
III_p+III_q&\leq c\bigl[\eta_p+\eta_q\bigr]^\frac12\Big[\frac\lambda{\bar h^{-1}(\lambda)}\Big]^\frac\xi2\biggl[\frac{|\mu|(B_R)}{R^n}\biggr]^{\frac 12}\\
&\leq c\Big[\frac{\bar h^{-1}(\lambda)}\lambda\Big]^\frac{1+\gamma-\xi}2\biggl[\frac{|\mu|(B_R)}{R^{n-1}}\biggr]^{\frac 12}\bigl[III_p+III_q\bigr]^\frac12\\
&\qquad+c\Big[\frac{\bar h^{-1}(\lambda)}\lambda\Big]^\frac{1-\xi}2\biggl[\frac{|\mu|(B_R)}{R^{n-1}}\biggr]^{\frac 12}d_0^\frac12+\varepsilon^\frac12 \Big[\frac{\bar h^{-1}(\lambda)}\lambda\Big]^{-\frac\gamma2}\biggl[\frac{|\mu|(B_R)}{R^{n-1}}\biggr]^{\frac 12}\\
&\leq \frac12\bigl[III_p+III_q\bigr]+c\Big[\frac{\bar h^{-1}(\lambda)}\lambda\Big]^{-\gamma}\frac{|\mu|(B_R)}{R^{n-1}}+c\Big[\frac{\bar h^{-1}(\lambda)}\lambda\Big]^{-\gamma} \biggl[\frac{|\mu|(B_R)}{R^{n-1}}\biggr]^{\frac 12}d_0^\frac12+\varepsilon.
\end{align*}
Absorbing the first term on the right-hand side and using $\frac{|\mu|(B_R)}{R^{n-1}}\leq d_0$ and $\frac{1+\xi}2=1+\gamma$, we obtain
\[
III_p+III_q\leq c\Big[\frac{\bar h^{-1}(\lambda)}\lambda\Big]^{-\gamma} d_0+2\varepsilon
\]
and inserting the latter estimate, together with \eqref{I_pq}-\eqref{II_pq}, into \eqref{est.L1.pq} finally leads to
\begin{equation}\label{2pq.final}
\mean{B_{\delta R/8}}|Du-Dv_2|\dx \leq c\Big[\frac{\bar h^{-1}(\lambda)}\lambda\Big]^{1+\gamma}\biggl[\Big[\frac{\bar h^{-1}(\lambda)}\lambda\Big]^{-\gamma} d_0+\varepsilon\biggr]; 
\end{equation}
this concludes the proof after letting $\varepsilon\searrow0$.

\subsubsection*{The intermediate case $p\leq 2<q$}
In this case we have, employing the monotonicity of $t\mapsto \bigl[\bar h^{-1}(t)\bigr]^{p-2}$ and \eqref{non.degeneracy}
\[
\frac\lambda{\bar h^{-1}(\lambda)}=\frac{\bar h(\bar h^{-1}(\lambda))}{\bar h^{-1}(\lambda)}=\bigl[\bar h^{-1}(\lambda)\bigr]^{p-2}+a_i(B_R)\bigl[\bar h^{-1}(\lambda)\bigr]^{q-2}\leq {\bigl[\bar h^{-1}(\lambda)\bigr]}^{p-2}+ c(q, K)\,a_i(B_R)|Dw_1|^{q-2}
\]
in $B_{\delta R}$. Therefore
\begin{multline*}
\mean{B_{\delta R/8}}|Du-Dv_2|\dx\leq c(n,p,q)\Big[\frac{\bar h^{-1}(\lambda)}\lambda\Big]^{1+\gamma}\,\bigl[\bar h^{-1}(\lambda)\bigr]^{(p-2)(1+\gamma)}\mean{B_{\delta R/8}}|Du-Dv_2|\dx \\
+c(n,q,K)\Big[\frac{\bar h^{-1}(\lambda)}\lambda\Big]^{1+\gamma}\mean{B_{\delta R/8}}{\bigl[a_i(B_R)|Dw_1|^{q-2}\bigr]}^{1+\gamma}|Du-Dv_2|\dx =:\tilde{I}_p+\tilde{I}_q.
\end{multline*}
To estimate the first term, we first observe that if $a_i(B_R)\neq 0$, then for $\esp\in \{p,q\}$
\begin{multline*}
\mean{B_{\delta R/2}}\bigl(|Du|+|Dv_2|\bigr)\dx \leq \mean{B_{\delta R/2}}\bigl(2|Du|+|Du-Dv_2|\bigr)\dx\\
\leq \frac{c}{\bar a_\esp^{\frac1{\esp-1}}}\biggl(\mean{B_{\delta R/2}}\bar a_\esp\Bigl[|Du|^{\max\{1,\esp-1\}}+|Du-Dv_2|^{\max\{1,\esp-1\}}\Bigr]\dx\biggr)^{\frac1{\max\{1,\esp-1\}}};
\end{multline*}
recall that $p-1\leq 1$, $q-1>1$, so we can use H\"older's inequality only when $\esp=q$. Now we can use \eqref{local.energy.SOLA} with $\es=1/(p-1)$ (we already observed that the choice is possible due to our bounds for $p$) together with $d_1\lesssim \lambda$ to estimate the first integral if $\esp=p$, while for $\esp=q$ it is sufficient to invoke \eqref{basic.bounds} together with $a_i(B_R)\leq a(x)$:
\[
{\biggl(\mean{B_{\delta R/2}}\bar a_\esp|Du|^{\max\{1,\esp-1\}}\dx\biggr)}^\frac{1}{\max\{1,\esp-1\}}\leq c\,\lambda^{\frac1{\esp-1}}.
\]
Similarly we can bound the second integral, using \eqref{rough.compa.solopartep} for $\es=1/(p-1)$ if $\esp=p$, \eqref{Du-Dv_2} with $\bi=1$ if $\esp=q$:
\[
{\biggl(\mean{B_{\delta R/2}}\bar a_\esp|Du-Dv_2|^{\max\{1,\esp-1\}}\dx\biggr)}^\frac1{\max\{1,\esp-1\}}\leq c\,\lambda^{\frac1{\esp-1}};
\]
in both cases we also used that we can bound 
\[
\frac{|\mu|(B_{\delta R})}{(\delta R)^{n-1}}+\biggl[\frac{|\mu|(B_R)}{R^{n-1}}\biggr]^{p-1}\biggl(\mean{B_R}h(x,Du)\dx\biggr)^{2-p}\lesssim \lambda
\]
thanks to \eqref{basic.bounds}. Hence we can see, merging the two cases, that 
\begin{equation}\label{two.energies.bound}
\mean{B_{\delta R/2}}\bigl(|Du|+|Dv_2|\bigr)\dx\leq c\min\biggl\{\lambda^{\frac1{p-1}},{\biggl[\frac{\lambda}{a_i(B_R)}\biggr]}^{\frac1{q-1}}\biggr\}\leq c\,\bar h^{-1}(\lambda),   
\end{equation}
thanks to \eqref{above.below.inverseh}, for a constant depending on $\data,\gdata$ and $\delta$. If instead $a_i(B_R)=0$, we still can estimate the sum of the $L^1$ norms by $c\,\bar h^{-1}(\lambda)$, as the minimum can be directly replaced by $\lambda^\frac1{p-1}$ and we use only the estimate with $\esp=p$ above.

\medskip

Applying \eqref{compa.V.component.L1} with $\esp=p$, $v=v_2$ and $\delta R$ replacing $R$, and using $|\mu|(B_{\delta R/2})\leq|\mu|(B_R)$, gives
\[
\mean{B_{\delta R/2}}{\bigl|V_p(Du)-V_p(Dv_2)\bigr|}^{\frac2p}\dx\leq c{\biggl[\frac{|\mu|(B_R)}{R^{n-1}}\mean{B_{\delta R/2}}|Du-Dv_2|\dx\biggr]}^{\frac1p};
\]
if $p<2$, the standard estimate for the nonlinear vector field $V_p$ in \eqref{diffV}, H\"older's inequality and \eqref{two.energies.bound} yield
\begin{align*}
\mean{B_{\delta R/2}}|Du-Dv_2|\dx&\leq c{\biggl(\mean{B_{\delta R/2}}{\bigl|V_p(Du)-V_p(Dv_2)\bigr|}^{\frac2p}\dx\biggr)}^{\frac p2}{\biggl(\mean{B_{\delta R/2}}(|Du|+|Dv_2|)\dx\biggr)}^{\frac{2-p}2}\\
&\leq c{\biggl[\frac{|\mu|(B_R)}{R^{n-1}}\mean{B_{\delta R/2}}|Du-Dv_2|\dx\biggr]}^{\frac12}{\bigl[\bar h^{-1}(\lambda)\bigr]}^{\frac{2-p}2}
\end{align*}
(the same final inequality clearly also holds when $p=2$). Consequently,
\[
\mean{B_{\delta R/2}}|Du-Dv_2|\dx\leq c[\bar h^{-1}(\lambda)]^{2-p}\frac{|\mu|(B_R)}{R^{n-1}}
\]
and since $[\bar h^{-1}(\lambda)]^{p-1}\leq\lambda$, this gives
\begin{align}\label{singular.p.contribution}
\Bigl[\frac{\bar h^{-1}(\lambda)}{\lambda}\Bigr]^{1+\gamma}[\bar h^{-1}(\lambda)]^{(p-2)(1+\gamma)}\mean{B_{\delta R/2}}|Du-Dv_2|\dx&\leq c\,\frac{\bar h^{-1}(\lambda)}{\lambda}\Bigl[\frac{[\bar h^{-1}(\lambda)]^{p-1}}{\lambda}\Bigr]^\gamma\frac{|\mu|(B_R)}{R^{n-1}}\notag\\
&\leq c\,\frac{\bar h^{-1}(\lambda)}{\lambda}\frac{|\mu|(B_R)}{R^{n-1}}.
\end{align}
Enlarging the integration ball in $\tilde I_p$ and applying \eqref{singular.p.contribution}, we therefore find that
\[
\tilde I_p\leq c\,\frac{\bar h^{-1}(\lambda)}{\lambda}\frac{|\mu|(B_R)}{R^{n-1}}.
\]

\medskip

For the estimate of $\tilde{I}_q$, we partially refer to the previous paragraph, since $q\geq2$; in particular, 
\[
\tilde{I}_q\leq c\Big[\frac{\bar h^{-1}(\lambda)}\lambda\Big]^{1+\gamma}[I_q+II_q+III_q] 
\]
as defined in \eqref{I_est} and subsequent lines. To estimate $I_q$ and $II_q$, we retain only the $q$-components of the comparison estimates: if $\bar a_q=0$, these terms vanish, so we assume $\bar a_q>0$. Since $\bar a_q=a_i(B_R)\leq a(x)$ on $B_{R_k/2}$, applying \eqref{compa.q.component.ell} to $u,v_k$, with $R_k$ replacing $R$, gives, for $k=1,2$,
\begin{equation}\label{Du-Dv_2.adhoc}
\biggl(\mean{B_{R_k/2}}\bigl[\bar a_q|Du-Dv_k|^{q-1}\bigr]^\ell\dx\biggr)^{\frac1\ell}\leq c\,\frac{|\mu|(B_R)}{R^{n-1}},\qquad1\leq\ell<\min\{n',q'\}.
\end{equation}
Here we used $R_k\geq\delta R$ and $|\mu|(B_{R_k/2})\leq|\mu|(B_R)$; the constant $c$ depends only on $n,p,q,\nu,\ell$ and $\delta$. We use this estimate in place of \eqref{Du-Dv_2} when estimating $I_q$ and $II_q$. Combining \eqref{Du-Dv_2.adhoc} with \eqref{Dv_2-Dw_2.hbi}, applied at scales $R$ and $\delta R$, we obtain, for every $1\leq\ell<\min\{n',q'\}$,
\begin{equation}\label{q.component.comparison}
\biggl(\mean{B_{R/8}}\bigl[\bar a_q|Du-Dw_1|^{q-1}\bigr]^\ell\dx\biggr)^{\frac1\ell}+\biggl(\mean{B_{\delta R/2}}\bigl[\bar a_q|Du-Dv_2|^{q-1}\bigr]^\ell\dx\biggr)^{\frac1\ell}\leq c\,d_0;
\end{equation}
here we used
\[
\Bigl[\omega(R)+R^{\min\{s_1,\x\}}\Bigr]^{\textnormal{\texttt{E}}}d_1\leq d_0.
\]
Combining these bounds with the comparisons between $v_j$ and $w_j$ in \eqref{Dv_2-Dw_2.hbi}, modified using \eqref{Du-Dv_2.adhoc}, yields
\begin{multline*}
{\biggl(\mean{B_{\delta R/8}}{\bigl[\bar a_q|Dw_1-Dw_2|^{q-1}\bigr]}^{1+\gamma}\dx\biggr)}^{\frac1{1+\gamma}}+{\biggl(\mean{B_{\delta R/8}}{\bigl[\bar a_q|Dv_2-Dw_2|^{q-1}\bigr]}^{1+\gamma}\dx\biggr)}^{\frac1{1+\gamma}}\\
\leq c\Bigl[\omega(R)+R^{\min\{s_1,\x\}}\Bigr]^{\textnormal{\texttt{E}}}d_1+c\,\frac{|\mu|(B_R)}{R^{n-1}}\leq c\,d_0.
\end{multline*}
We recall that $1+\gamma<\min\{n',q'\}$. H\"older's inequality therefore gives
\[
I_q+II_q\leq c\,\bar a_q^{\frac{\gamma}{q-1}}d_0^{(1+\gamma)\frac{q-2}{q-1}}\biggl[\frac{|\mu|(B_R)}{R^{n-1}}\biggr]^{\frac1{q-1}}\leq c\,\bar a_q^{\frac{\gamma}{q-1}}d_0^{1+\gamma\frac{q-2}{q-1}}\leq c\,d_0\Bigl[\bar a_q^{\frac1{q-1}}\lambda^{\frac{q-2}{q-1}}\Bigr]^\gamma\leq c\,d_0{\Bigl[\frac{\lambda}{\bar h^{-1}(\lambda)}\Bigr]}^\gamma,
\]
similarly to \eqref{I_pq} and \eqref{II_pq}, but with \eqref{Du-Dv_2.adhoc} replacing \eqref{Du-Dv_2} - this will avoid the appearance of the $d_{p,q}$ term. To estimate $III_q$, we introduce $\eta_q$ as in \eqref{def.eta_t} and again decompose the integral appearing in its definition into the terms $IV_q$, $V_q$ and $VI_q$, as in the case $p>2$. The additional point concerns $IV_q$: to bound this term through Sobolev's inequality, we first need an auxiliary $L^1$ estimate for $Du-Dv_2$ on $B_{\delta R/2}$. We therefore repeat on this ball the decomposition into the $p$- and $q$-contributions used to define $\tilde I_p$ and $\tilde I_q$. Applying the bound \eqref{singular.p.contribution} already proved for the $p$-contribution, which involves only the measure density, we obtain
\[
\mean{B_{\delta R/2}}|Du-Dv_2|\dx\leq c\,\frac{\bar h^{-1}(\lambda)}{\lambda}\biggl[\frac{|\mu|(B_R)}{R^{n-1}}+c\Bigl[\frac{\bar h^{-1}(\lambda)}{\lambda}\Bigr]^{\gamma}\mean{B_{\delta R/2}}\bigl[a_i(B_R)|Dw_1|^{q-2}\bigr]^{1+\gamma}|Du-Dv_2|\dx\biggr].
\]
Since $q>2$, we can split $Dw_1=Dv_2+(Dw_1-Du)+(Du-Dv_2)$ and use subadditivity to estimate the last integral:
\begin{multline*}
\mean{B_{\delta R/2}}\bigl[a_i(B_R)|Dw_1|^{q-2}\bigr]^{1+\gamma}|Du-Dv_2|\dx\\
\leq c\,III_q+c\mean{B_{\delta R/2}}\bigl[a_i(B_R)|Du-Dw_1|^{q-2}\bigr]^{1+\gamma}|Du-Dv_2|\dx\\
+c\mean{B_{\delta R/2}}\bigl[a_i(B_R)|Du-Dv_2|^{q-2}\bigr]^{1+\gamma}|Du-Dv_2|\dx.
\end{multline*}
For the first comparison term, H\"older's inequality and the preceding comparison estimates yield
\begin{align*}
\mean{B_{\delta R/2}}\bigl[a_i(B_R)&|Du-Dw_1|^{q-2}\bigr]^{1+\gamma}|Du-Dv_2|\dx\\
&\leq c\,a_i(B_R)^{\frac{\gamma}{q-1}}\biggl(\mean{B_{\delta R/2}}{\bigl[\bar a_q|Du-Dw_1|^{q-1}\bigr]}^{1+\gamma}\dx\biggr)^{\frac{q-2}{q-1}}\biggl(\mean{B_{\delta R/2}}\bar a_q|Du-Dv_2|^{q-1}\dx\biggr)^{\frac1{q-1}}\\
&\leq c\,a_i(B_R)^{\frac{\gamma}{q-1}}d_0^{1+\gamma\frac{q-2}{q-1}}\leq c\Bigl[\frac{\lambda}{\bar h^{-1}(\lambda)}\Bigr]^\gamma d_0,
\end{align*}
where we used the bound in \eqref{dzero}-\eqref{dpq} and \eqref{above.below.inverseh}. The second comparison term satisfies the same bound, using the comparison estimate for $Du-Dv_2$ with exponent $1+\gamma(q-2)/(q-1)<1+\gamma$. Thus
\[
\mean{B_{\delta R/2}}|Du-Dv_2|\dx\leq c\,\frac{\bar h^{-1}(\lambda)}{\lambda}d_0+c\,\Bigl[\frac{\bar h^{-1}(\lambda)}{\lambda}\Bigr]^{1+\gamma}III_q.
\]
Finally, since $u-v_2$ has zero trace on $\partial B_{\delta R/2}$, we can use Sobolev's inequality to find 
\[
IV_q\leq c\,R^\xi\Bigl[\frac{\lambda}{\bar h^{-1}(\lambda)}\Bigr]^\xi\biggl(\mean{B_{\delta R/2}}|Du-Dv_2|\dx\biggr)^\xi\leq c\,R^\xi d_0^\xi+c\,R^\xi\Bigl[\frac{\bar h^{-1}(\lambda)}{\lambda}\Bigr]^{\gamma\xi}III_q^\xi;
\]
we additionally relied on the non-degeneracy condition provided by \eqref{non.degeneracy}. The estimates for $V_q$ and $VI_q$ follow through the same chain of H\"older and Sobolev inequalities, keeping exclusively the $q$-components and appealing to \eqref{Du-Dv_2.adhoc} alongside \eqref{q.component.comparison}. In particular, $V_q+VI_q\leq cR^\xi d_0^\xi$ and
\[
\eta_q\leq c\,R\Bigl[\frac{\bar h^{-1}(\lambda)}{\lambda}\Bigr]^{1+\gamma}III_q+c\,R\,\frac{\bar h^{-1}(\lambda)}{\lambda}d_0+\varepsilon R\Bigl[\frac{\bar h^{-1}(\lambda)}{\lambda}\Bigr]^{1+\gamma}.
\]
We can therefore repeat the weighted estimate and the absorption argument for $III_q$ alone that, after letting $\varepsilon\searrow0$, leads to
\[
III_q\leq c\Bigl[\frac{\lambda}{\bar h^{-1}(\lambda)}\Bigr]^\gamma d_0;
\]
merging this with the estimates for $I_q$ and $II_q$, we conclude with
\[
\tilde I_q \leq c\,\frac{\bar h^{-1}(\lambda)}{\lambda}d_0.
\]

\subsubsection*{The case $p\leq q\leq 2$}

In this last case we can directly prove the stronger estimate
\[
\mean{B_{\delta R/2}}|Du-Dv_2|\dx\leq c\,\frac{\bar h^{-1}(\lambda)}{\lambda}\frac{|\mu|(B_R)}{R^{n-1}}.
\]
Applying \eqref{rough.compa.solopartep} with radius $\delta R$ with $\es=1/(p-1)$, which is admissible by \eqref{2-1/n}, and using \eqref{basic.bounds}, we obtain
\[
\mean{B_{\delta R/2}}|Du-Dv_2|\dx\leq c\biggl[{\biggl[\frac{|\mu|(B_R)}{R^{n-1}}\biggr]}^{\frac1{p-1}}+{\frac{|\mu|(B_R)}{R^{n-1}}}\lambda^{\frac{2-p}{p-1}}\biggr]\leq c\,\frac{|\mu|(B_R)}{R^{n-1}}\frac{\lambda^\frac1{p-1}}\lambda,
\]
where the constants also account for the change of scale from $\delta R$ to $R$.

\medskip

In the $p$-phase (first inequality in \eqref{two.different.cases}), Remark~\ref{rem:p.phase.normalization} guarantees $a_i(B_R)\lambda^{\frac{q-p}{p-1}}\leq c$ and $\bar h^{-1}(\lambda)\approx\lambda^{\frac1{p-1}}$; the desired estimate thus follows. In the $(p,q)$-phase (second inequality in \eqref{two.different.cases}), we additionally use the intermediate estimate \eqref{rough.compa.qcomponent} from the proof of Lemma~\ref{lem:rough.compa.pq}, with radius $\delta R$ and exponent $\ti=1/(q-1)$, which is admissible since $2-1/n<p\leq q \leq n$. We apply the estimate at the common center $x_0$: since $a_i(B_R)\leq a(x_0)$ and $h_{x_0}(Du)\leq2h(x,Du)$ in $B_{\delta R}$, we obtain
\[
\mean{B_{\delta R/2}}|Du-Dv_2|\dx\leq c{\bigl[a_i(B_R)\bigr]}^{-\frac1{q-1}}\biggl[{\biggl[\frac{|\mu|(B_R)}{R^{n-1}}\biggr]}^{\frac1{q-1}}+\frac{|\mu|(B_R)}{R^{n-1}}\lambda^{\frac{2-q}{q-1}}\biggr]\leq c\,\frac{|\mu|(B_R)}{R^{n-1}}\,\frac1\lambda{\biggl[\frac{\lambda}{a_i(B_R)}\biggr]}^{\frac1{q-1}}.
\]
Combining the two component estimates, we conclude with \eqref{above.below.inverseh} that
\[
\mean{B_{\delta R/2}}|Du-Dv_2|\dx\leq c\,\frac{1}{\lambda}\min\biggl\{\lambda^{\frac1{p-1}},{\biggl[\frac{\lambda}{a_i(B_R)}\biggr]}^{\frac1{q-1}}\biggr\}\frac{|\mu|(B_R)}{R^{n-1}}\leq c\,\frac{\bar h^{-1}(\lambda)}{\lambda}\frac{|\mu|(B_R)}{R^{n-1}}.
\]

\end{proof}

Merging the two previous results gives the following final $L^1$ linearized comparison Lemma:

\begin{lemma}[A refined comparison Lemma for $Du-Dw_2$]\label{lem:refined.Du-Dw2}
Suppose that \eqref{basic.bounds} is in force for some parameter $\lambda\geq1$ and that \eqref{non.degeneracy} holds for some constant $K\geq1$. Then
\[
\mean{B_{\delta R/8}}|Du-Dw_2|\dx\leq c\,\frac{\bar h^{-1}(\lambda)}\lambda\biggl[\Bigl[\omega(R)+R^{\min\{s_\smallone,\x\}}\Bigr]^{\frac{\textnormal{\texttt{E}}}{\max\{1,q-1\}}}\mean{B_R} h(x,Du)\dx+\frac{|\mu|(B_R)}{R^{n-1}}\biggr]
\]
for a constant depending on $\data,\gdata,\delta,\vartheta$ and $K$. The exponent $\texttt E\equiv\texttt E(\data,\gdata)\in(0,1)$ is as in \eqref{texttt.E}.
\end{lemma}

\subsection{The final comparison Lemma - rough version}

As a first step toward a decay estimate for the flux excess, we compare the flux of $u$ with that of the regular reference map $w_1$. While the resulting lemma renders the comparison error arbitrarily small relative to $\lambda$, it serves only as an intermediate step: the subsequent iteration requires a quantitative bound for this error in terms of the radius and the local measure density, with an error term that is summable over the relevant geometric scales.

\begin{lemma}\label{rough.linearized.comparison}
Let $\epsilon\in(0,1]$ and let $u\in W^{1,H(\cdot)}_{\rm loc}(\Omega)$ be a local energy solution to \eqref{equation} under the structural assumptions \eqref{A.strong}, \eqref{a.holder} to \eqref{bounds}; when $p\leq n$, assume additionally that $\mu\in L^\infty_{\rm loc}(\Omega)$. Then there exist a radius $R_{\epsilon,{\rm rgh}}\in(0,R_{\rm cmp}]$ and a constant $\mathcal B_{\epsilon,{\rm rgh}}\geq1$ such that the following holds. Let $B_R\equiv B_R(x_0)$, with $x_0\in\Omega_0$ and $0<R\leq R_{\epsilon,{\rm rgh}}$, and suppose that
\begin{equation}\label{ass.rough.compa.a}
\mean{B_R}h(x,Du)\dx\leq\lambda,\qquad\frac{|\mu|(B_R)}{R^{n-1}}\leq\frac{\lambda}{\mathcal B_{\epsilon,{\rm rgh}}}
\end{equation}
for some $\lambda\geq1$. Then the $A_{i,1}$-harmonic comparison map $w_1\in v_1+W^{1,H_{i,\smallone}(\cdot)}_0(B_{R/8})$ defined after \eqref{two.different.cases} satisfies
\begin{equation}\label{est.qualitative.exc.small}
\mean{B_{R/16}}\big|A(x,Du)-A(x,Dw_1)\big|\dx\leq\epsilon\lambda.
\end{equation}
The constant $\mathcal B_{\epsilon,{\rm rgh}}$ depends on $\data,\gdata$ and $\epsilon$, while the radius $R_{\epsilon,{\rm rgh}}$ also depends on $\omega(\,\cdot\,)$.
\end{lemma}

\begin{proof}
We recall the definition of $d_0$ in \eqref{dzero} and $d_{p,q}$ in \eqref{d_pq}. Under \eqref{ass.rough.compa.a}, we have
\[
d_0\leq\biggl[\Bigl[\omega(R_{\epsilon,{\rm rgh}})+R_{\epsilon,{\rm rgh}}^{\min\{s_\smallone,\x\}}\Bigr]^{\textnormal{\texttt{E}}}+\frac1{\mathcal B_{\epsilon,{\rm rgh}}}\biggr]\lambda
\]
and 
\[
d_{p,q}\leq \chi_{(-\infty,2)}(p)\biggl[\frac{\lambda}{{\mathcal B}_{\epsilon, {\rm rgh}}}\biggr]^{p-1}\lambda^{2-p}+\chi_{(-\infty,2)}(q)\biggl[\frac{\lambda}{{\mathcal B}_{\epsilon, {\rm rgh}}}\biggr]^{q-1}\lambda^{2-q}\leq \frac2{{\mathcal B}_{\epsilon, {\rm rgh}}^{p-1}}\lambda.
\]
We  want to use Lemma~\ref{lem:a-a}: as we have the pointwise estimate
\begin{multline*}
\big|A(x,Du)-A(x,Dw_1)\big|\leq c(p,q,L)\,h(x,Du-Dw_1)\\
+c(p,q,L)\Big[|Du|^{p-2}\chi_{[2,+\infty)}(p)+a(x)|Du|^{q-2}\chi_{[2,+\infty)}(q)\Big]|Du-Dw_1|,
\end{multline*}
after integration we estimate separately the two possible terms. By \eqref{compa.Du-Dw} in the $p$-phase, and by \eqref{compa.Du-Dw,pq} together with the coefficient comparability in the $(p,q)$-phase, we have
\begin{equation}\label{rough.compa.h.DuDw1}
\mean{B_{R/16}}h(x,Du-Dw_1)\dx\leq c\bigl[d_0+d_{p,q}\bigr], 
\end{equation}
where $c$ depends on $\data$ and $\gdata$. We shall choose $R_{\epsilon,{\rm rgh}}$ and $\mathcal B_{\epsilon,{\rm rgh}}$ at the end of the proof, after estimating the remaining term, that is
\[
II:=\mean{B_{R/16}}\Big[|Du|^{p-2}\chi_{[2,+\infty)}(p)+a(x)|Du|^{q-2}\chi_{[2,+\infty)}(q)\Big]|Du-Dw_1|\dx.
\]
If $q<2$, then $II=0$. Otherwise, we may estimate the nonzero terms directly by H\"older's inequality: for $q>2$, we have
\begin{multline*}
\mean{B_{R/16}}a(x)|Du|^{q-2}|Du-Dw_1|\dx\\
\leq\biggl(\mean{B_{R/16}}a(x)|Du|^{q-1}\dx\biggr)^{\frac{q-2}{q-1}}\biggl(\mean{B_{R/16}}a(x)|Du-Dw_1|^{q-1}\dx\biggr)^{\frac1{q-1}}\leq c\,\lambda^{\frac{q-2}{q-1}}{\bigl[d_0+d_{p,q}\bigr]}^{\frac1{q-1}},
\end{multline*}
where we used \eqref{rough.compa.h.DuDw1} and the bound in \eqref{ass.rough.compa.a}, after enlarging the integral. The case $q=2$ follows directly from \eqref{rough.compa.h.DuDw1}, while the $p$-term is estimated in the same way when $p\geq2$; therefore,
\[
II\leq c\sum_{\esp\in\{p,q\}}\chi_{[2,\infty)}(\esp)\lambda^{\frac{\esp-2}{\esp-1}}{\bigl[d_0+d_{p,q}\bigr]}^{\frac1{\esp-1}}
\]
with $c$ depending on $\data$ and $\gdata$. Combining \eqref{rough.compa.h.DuDw1} with the preceding estimate and Lemma~\ref{lem:a-a}, we obtain, provided $d_0+d_{p,q}\leq5\lambda$,
\[
\mean{B_{R/16}}\bigl|A(x,Du)-A(x,Dw_1)\bigr|\dx\leq \bar c\,\lambda\Bigl(\frac{d_0+d_{p,q}}{\lambda}\Bigr)^{\frac1{\max\{1,q-1\}}},
\]
where $\bar c$ depends only on $\data$ and $\gdata$. On the other hand, the estimates for $d_0$ and $d_{p,q}$ give
\[
\frac{d_0+d_{p,q}}{\lambda}\leq\Bigl[\omega\bigl(R_{\epsilon,{\rm rgh}}\bigr)+R_{\epsilon,{\rm rgh}}^{\min\{s_\smallone,\x\}}\Bigr]^{\textnormal{\texttt{E}}}+\frac1{\mathcal B_{\epsilon,{\rm rgh}}}+\frac2{\mathcal B_{\epsilon,{\rm rgh}}^{p-1}}.
\]
Since $\omega(R)\to0$ as $R\to0$ and $s_1,\x>0$, we can choose $R_{\epsilon,{\rm rgh}}$ sufficiently small and $\mathcal B_{\epsilon,{\rm rgh}}$ sufficiently large so that the right-hand side is bounded by $({\epsilon/\bar c})^{\max\{1,q-1\}}\leq 1$, and this concludes the proof.
\end{proof}

\subsection{The final comparison Lemma - linearized version}
We now combine the preceding estimates to prove the following linearized comparison result.

\begin{proposition}[A refined comparison result]\label{def.comparison.pphase}
Let $u$ and $\mu$ be as in Lemma~\ref{rough.linearized.comparison} and let $\vartheta\in(0,1/2]$ be fixed. There exists a small constant $\delta_{0,\vartheta}\in(0,1/32)$, depending on $\data,\gdata$ and $\vartheta$, such that for every $\delta\in(0,\delta_{0,\vartheta}]$ there exist a large constant $\mathcal B$, depending on $\data,\gdata,\vartheta$ and $\delta$, and a small radius $R_0$, depending on the same quantities and also on $\omega(\,\cdot\,)$, such that if $R\leq R_0$,
\begin{equation}\label{non.degenerate.Du}
\vartheta\lambda\leq\mean{B_{\delta R}}h(x,Du)\dx,\qquad\mean{B_R}h(x,Du)\dx\leq\lambda,\qquad\frac{|\mu|(B_R)}{R^{n-1}}\leq\frac{\lambda}{\mathcal B},
\end{equation}
for some $\lambda\geq1$, then the $A_{i,2}$-harmonic comparison map $w_2\in v_2+W^{1,H_{i,2}(\cdot)}_0(B_{\delta R/8})$ constructed above satisfies
\begin{equation}\label{nice.compa.pphase.final}
\mean{B_{\delta R/16}}\big|A(x,Du)-A(x,Dw_2)\big|\dx\leq c\biggl[\Bigl[\omega(R)+R^{\min\{s_1,\x\}}\Bigr]^{\frac{\textnormal{\texttt{E}}}{\max\{q-1,1\}}}+\frac{|{\mu}|(B_R)}{\lambda\,R^{n-1}}\biggr]\lambda
\end{equation}
where the constant $c$ depends only on $\data,\gdata,\delta$ and $\vartheta$.
\end{proposition}

\begin{proof}
We first estimate the average of the frozen flux, appealing only to the upper bound for the average of $h(x,Du)$ and the measure bound in \eqref{non.degenerate.Du}. The energy estimates for the frozen comparison problem yield
\begin{equation}\label{partial.above.pphase}
\mean{B_{R/8}}\bar h(Dw_1)\dx\leq c\biggl[\mean{B_R}h(x,Du)\dx+\frac{|\mu|(B_R)}{R^{n-1}}\biggr]\leq\bar c\,\lambda.
\end{equation}
In particular, this estimate does not use the lower bound in \eqref{non.degenerate.Du}, and $\bar c$ is independent of $\delta$ and $\vartheta$. In this first part of the proof, the independence of the constants from these two parameters is essential. By \eqref{lip.varphi.harmonic}, applied to the frozen field $A_{i,1}$, we therefore obtain
\begin{equation}\label{pointwise.Dw1}
\mean{B_{R/8}}\bar h(Dw_1)\dx\leq\bar c\,\lambda\qquad\Longrightarrow\qquad\sup_{B_{R/16}}\bar h(Dw_1)\leq c\,\lambda.
\end{equation}
We verify this implication in the two cases of \eqref{two.different.cases}: in the $p$-phase, the frozen energy comparison gives
\[
\mean{B_{R/8}}\bar H(Dw_1)\dx\leq c\,E_1^p.
\]
Using the bound $a_i(B_R)E_1^{q-p}\leq c$ from the normalization of the frozen problem, \eqref{averaged.Psi.inverse} with $\el=1$, and \eqref{en.rh.1}, we obtain
\[
\mean{B_{R/8}}\bar h(Dw_1)\dx\leq c\,\bar h(E_1)\leq c\,E_1^{p-1}\leq c\biggl[\mean{B_R}h(x,Du)\dx+\frac{|\mu|(B_R)}{R^{n-1}}\biggr].
\]
In the $(p,q)$-phase, the same conclusion follows from \eqref{pq-frozen.energy}, \eqref{estimate.E2}, the subsequent bound for the average of $h(x,Dv_1)$, and \eqref{averaged.Psi.inverse}. We now estimate from below, employing \eqref{rough.compa.h.DuDw1}, which follows directly from Proposition~\ref{rough.comparison.DuDw} and does not require the non-degeneracy assumption \eqref{non.degeneracy}. Combining this estimate with \eqref{non.degenerate.Du}, subadditivity and the inclusion $B_{\delta R}\subset B_{R/16}$, valid since $\delta\leq1/16$, we obtain
\begin{align}\label{est.h.below}
\mean{B_{\delta R}}h(x,Dw_1)\dx&\geq\frac1{2^{q-1}}\mean{B_{\delta R}}h(x,Du)\dx-c(n)\delta^{-n}\mean{B_{R/16}}h(x,Du-Dw_1)\dx\notag\\
&\geq\frac{\vartheta}{2^{q-1}}\lambda-\tilde c\bigl[d_0+d_{p,q}\bigr]\geq\frac{\vartheta}{2^{q-1}}\lambda-\tilde c\biggl[\Bigl[\omega(R)+R^{\min\{s_1,\x\}}\Bigr]^{\textnormal{\texttt{E}}}+\frac1{\mathcal B}+\frac2{\mathcal B^{p-1}}\biggr]\lambda,
\end{align}
with $\tilde c=\tilde c(\data,\gdata,\delta,\vartheta)$. We choose $R_0$ sufficiently small and $\mathcal B$ sufficiently large that
\[
\tilde c\Bigl[\omega(R_0)+R_0^{\min\{s_1,\x\}}\Bigr]^{\textnormal{\texttt{E}}}\leq\frac{\vartheta}{2^{q+1}},\qquad\tilde c\biggl[\frac1{\mathcal B}+\frac2{\mathcal B^{p-1}}\biggr]\leq\frac{\vartheta}{2^{q+1}};
\]
consequently, for every $R\leq R_0$,
\[
\mean{B_{\delta R}}h(x,Dw_1)\dx\geq\frac{\vartheta}{2^q}\lambda.
\]
We now transfer the preceding lower bound to the frozen growth function $\bar h$: by the definition of $a_i(B_R)$ and the H\"older continuity of $a$, we have
\[
0\leq h(x,Dw_1)-\bar h(Dw_1)\leq c[a]_{0,\alpha}R^\alpha |Dw_1|^{q-1}\qquad\text{in }B_{R/16}.
\]
If $p\leq n$, estimates \eqref{pointwise.Dw1} and \eqref{R.lambda}, together with the choice of $R_{\rm cmp}$, give
\[
R^\alpha|Dw_1|^{q-p}\leq cR^\alpha\lambda^{\frac{q-p}{p-1}}\leq c\,C_{\rm ref}^{\frac{q-p}{p-1}}R^{\alpha-\frac{(n-1)(q-p)}{p-1}}\leq cR^{s_1},
\]
where $c$ depends only on $\data,\gdata,\delta$ and $\vartheta$. In the last inequality we used $\alpha-(n-1)(q-p)/(p-1)=s_0+\x\geq s_0$ and $C_{\rm ref}^{(q-p)/(p-1)}R^{s_0-s_1}\leq1$, which follows from $R\leq R_{\rm cmp}$ and $s_1=s_0/p$. Since $s_1>0$, we may further reduce $R_0$, depending on $\data,\gdata,\delta,\vartheta$ and $\omega(\,\cdot\,)$, so that
\[
h(x,Dw_1)\leq2\bar h(Dw_1)\qquad\text{in }B_{R/16}.
\]
If $p>n$, we instead use the interior gradient estimate for the frozen comparison map and the energy estimates for $w_1$ and $v_1$ to obtain
\[
\sup_{B_{R/16}}|Dw_1|^p\leq c\mean{B_{R/8}}H_{i,1}(Dw_1)\dx\leq \frac{c}{R^n}\int_{B_{R/2}}H(x,Du)\dx\leq \frac{c\mathcal M}{R^n}.
\]
Consequently, \eqref{rsmaller1} gives $R^\alpha |Dw_1|^{q-p}\leq c(\data,\gdata)$ in $B_{R/16}$ and it follows that
\begin{equation}\label{pointwise.1}
h(\,\cdot\,,Dw_1)\leq \bar h(Dw_1)+c|Dw_1|^{p-1}\leq c\,\bar h(Dw_1)\qquad\text{in }B_{R/16}. 
\end{equation}
In both cases we thus have $h(x,Dw_1)\leq c_f\bar h(Dw_1)$ in $B_{R/16}$, with $c_f\geq2$ depending only on $\data$ and $\gdata$. Combining this comparison with the preceding lower bound gives
\begin{equation}\label{Dw_1.below}
\mean{B_{\delta R}}\bar h(Dw_1)\dx\geq\frac{\vartheta}{2^qc_f}\lambda.
\end{equation}
We now set
\[
\mathcal A_\vartheta:=\max\biggl\{\bar c,\frac{2^qc_f}{\vartheta}\biggr\}
\]
and choose the constant $\delta_{0,\vartheta}$ appearing in the statement of the proposition so that $\delta_{0,\vartheta}\leq\delta_0(\mathcal A_\vartheta)$, where $\delta_0(\mathcal A_\vartheta)$ is given by Lemma~\ref{integral.to.pointwise}. Since \eqref{partial.above.pphase} and \eqref{Dw_1.below} imply
\[
\frac{\lambda}{\mathcal A_\vartheta}\leq\mean{B_{\delta R}}\bar h(Dw_1)\dx,\qquad\mean{B_{R/8}}\bar h(Dw_1)\dx\leq\mathcal A_\vartheta\lambda,
\]
we can apply Lemma~\ref{integral.to.pointwise} with $m=1$ and conclude that
\begin{equation}\label{pointwise.2}
\inf_{B_{\delta R}}\bar h(Dw_1)\geq\frac{\lambda}{2^{q-1}\mathcal A_\vartheta}.
\end{equation}
Let $c_1$ denote the constant in the upper bound \eqref{pointwise.Dw1}, and set $K_\vartheta:=\max\bigl\{c_1,2^{q-1}\mathcal A_\vartheta\bigr\}$. Then \eqref{pointwise.Dw1} and \eqref{pointwise.2} verify \eqref{non.degeneracy} with $K=K_\vartheta$. The scaling properties of $\bar h^{-1}$ give
\[
K_\vartheta^{-\frac1{p-1}}\bar h^{-1}(\lambda)\leq |Dw_1|\leq K_\vartheta^{\frac1{p-1}}\bar h^{-1}(\lambda)\qquad\text{in }B_{\delta R}.
\]
Together with $\bar h(Dw_1)\leq h(x,Dw_1)\leq c_f\bar h(Dw_1)$, this yields
\begin{equation}\label{non.degeneracy.Dw1}
\frac1{c_\vartheta}\frac{\lambda}{\bar h^{-1}(\lambda)}\leq\frac{h(x,Dw_1)}{|Dw_1|}\leq c_\vartheta\frac{\lambda}{\bar h^{-1}(\lambda)}\qquad\text{in }B_{\delta R},
\end{equation}
with $c_\vartheta:=c_fK_\vartheta^{p'}$, depending only on $\data,\gdata$ and $\vartheta$.
\medskip

Now that we have obtained the pointwise estimate that we assumed in \eqref{non.degeneracy}, we consider again three different cases to prove \eqref{nice.compa.pphase.final}:

 \subsubsection*{The case $2\leq p\leq q$}

For $x\in B_{\delta R/16}$, we use Lemma~\ref{lem:a-a}, the triangle inequality several times and Young's inequality in the sense \eqref{conj.inverse} to infer
\begin{align}\label{pgeq2.final.pphase}
\big|A(x,Du)-A(x,Dw_2)\big|&\leq c\,h(x,Du-Dw_2)+c\Big[{|Dw_2|}^{p-2}+a(x){|Dw_2|}^{q-2}\Big]|Du-Dw_2|\notag\\
&\leq c\,h(x,Du-Dw_2)+c\Big[{|Dw_1|}^{p-2}+a(x){|Dw_1|}^{q-2}\Big]|Du-Dw_2|\notag\\
&\hspace{4cm} +c\,\frac{h(x,Dw_1-Dw_2)}{|Dw_1-Dw_2|}|Du-Dw_2|\notag\\
&\leq c\,h(x,Du-Dw_1)+c\,h(x,Du-Dw_2)\notag\\
&\hspace{3cm}+c\Big[|Dw_1|^{p-2}+a(x)|Dw_1|^{q-2}\Big]|Du-Dw_2|,
\end{align}
for $c\equiv c(p,q,L)$; thus, averaging over $B_{\delta R/16}$, using the non-degeneracy in \eqref{non.degeneracy.Dw1}, \eqref{Dw_1-Dw_2.variable} after appropriately enlarging the integrals on the right-hand side and finally the linearized comparison Lemma~\ref{lem:refined.Du-Dw2}, yields
\begin{align}\label{almost.comparison.a}
\mean{B_{\delta R/16}}\big|A(x,Du)-A(x,Dw_2)\big|\dx&\leq c\,d_0+c\,\sup_{B_{\delta R}}\frac{h(x,Dw_1)}{|Dw_1|}\mean{B_{\delta R/8}}|Du-Dw_2|\dx\notag\\
&\leq c\,d_0+c\,\frac{\lambda}{\bar h^{-1}(\lambda)}\mean{B_{\delta R/8}}|Du-Dw_2|\dx\notag\\
&\leq c\biggl[\Bigl[\omega(R)+R^{\min\{s_1,\x\}}\Bigr]^{\frac{\textnormal{\texttt{E}}}{q-1}}\lambda+\frac{|\mu|(B_R)}{R^{n-1}}\biggr],
\end{align}
for $c$ depending on $\data,\gdata,\delta$ and $\vartheta$.

\subsubsection*{The intermediate case $p<2<q$} 
We use again the non-degeneracy bounds \eqref{non.degeneracy.Dw1} established above. By Lemma~\ref{lem:a-a}, in particular \eqref{flux.difference.mixed} for $x\in B_{\delta R/16}$, we have
\begin{align*}
\bigl|A(x,Du)-A(x,Dw_2)\big|&\leq c\,{\bigl(|Du|+|Dw_2|\bigr)}^{p-2}|Du-Dw_2|\\
&\hspace{4cm}+c\,a(x)\Bigl[|Du-Dw_2|^{q-2}+|Dw_2|^{q-2}\Bigr]|Du-Dw_2|\\
&\leq  c|Dw_1|^{p-2}\Bigl[|Du-Dw_1|+|Dw_1-Dw_2|\Bigr]\\
&\hspace{4cm}+c\,a(x)\Bigl[|Du-Dw_2|^{q-2}+|Dw_2|^{q-2}\Bigr]|Du-Dw_2|;
\end{align*}
the last estimate can be proven by distinguishing the cases where $|Du-Dw_1|+|Dw_1-Dw_2|< |Dw_1|/2$ (which guarantees $|Du|+|Dw_2|\geq |Dw_1|$) and the converse, where we have
\[
(|Du|+|Dw_2|)^{p-2}|Du-Dw_2|\leq|Du-Dw_2|^{p-1}\lesssim_p{\bigl[|Du-Dw_1|+|Dw_1-Dw_2|\bigr]}^{p-1}
\]
and we conclude again with the triangle inequality after splitting $p-1=p-2+1$; we repeatedly used that $p\leq 2$.  Since $q>2$, sub-additivity and Young's inequality with conjugate exponents $(q-1,(q-1)/(q-2))$ give
\begin{align*}
a(x){|Dw_2|}^{q-2}|Du-Dw_2|&\leq c\,a(x)|Dw_1|^{q-2}|Du-Dw_2|+c\,a(x){|Dw_1-Dw_2|}^{q-2}|Du-Dw_2|\\
&\leq c\,a(x){|Dw_1|}^{q-2}|Du-Dw_2|+c\,\sum_{k\in\{1,2\}}a(x){|Du-Dw_k|}^{q-1},
\end{align*}
as in \eqref{pgeq2.final.pphase}. Combining the preceding estimates, appealing to \eqref{pointwise.2} and \eqref{non.degeneracy.Dw1}, and averaging over $B_{\delta R/16}$, we obtain 
\begin{multline*} 
\mean{B_{\delta R/16}}\bigl|A(x,Du)-A(x,Dw_2)\bigr|\dx\leq c\,\frac{\lambda}{\bar h^{-1}(\lambda)}\mean{B_{\delta R/16}}|Du-Dw_2|\dx\\
+ c{\bigl[\bar h^{-1}(\lambda)\bigr]}^{p-2}\sum_{k\in\{1,2\}}\mean{B_{\delta R/16}}|Du-Dw_k|\dx+c\sum_{k\in\{1,2\}}\mean{B_{\delta R/16}}a(x)|Du-Dw_k|^{q-1}\dx. 
\end{multline*}
Lemma~\ref{lem:refined.Du-Dw2} allows us to estimate the first term: enlarging the integration ball and using \eqref{non.degenerate.Du}, we directly obtain
\[
\frac{\lambda}{\bar h^{-1}(\lambda)}\mean{B_{\delta R/16}}|Du-Dw_2|\dx \leq c\Bigl[\omega(R)+R^{\min\{s_1,\x\}}\Bigr]^{\frac{\textnormal{\texttt{E}}}{q-1}}\lambda+c\,\frac{|\mu|(B_R)}{R^{n-1}},
\]
where $c$ depends also on $\vartheta$ and $\delta$. 
For $k\in\{1,2\}$, we repeat the argument leading to \eqref{two.energies.bound} and \eqref{singular.p.contribution} on $B_{R_k/2}$, with $R_1=R$ and $R_2=\delta R$. In the $p$-phase we use Remark~\ref{rem:p.phase.normalization}, in the $(p,q)$-phase we take the minimum of the two energy bounds, using $a_i(B_R)\leq a(x)$ on both balls. Estimate \eqref{compa.V.component.L1} with $\esp=p$ and H\"older's inequality then give
\[
\mean{B_{R_k/2}}|Du-Dv_k|\dx\leq c[\bar h^{-1}(\lambda)]^{2-p}\frac{|\mu|(B_R)}{R^{n-1}},\qquad k\in\{1,2\}.
\]
On the other hand, we have for $k\in\{1,2\}$
\[
\mean{B_{R_k/8}}|Dv_k-Dw_k|\dx \leq c\,\bar h^{-1}\biggl(\Bigl[\omega(R)+R^{\min\{s_1,\x\}}\Bigr]^{\textnormal{\texttt{E}}}\lambda\biggr)\leq c\Bigl[\omega(R)+R^{\min\{s_1,\x\}}\Bigr]^{\frac{\textnormal{\texttt{E}}}{q-1}}\bar h^{-1}(\lambda)
\]
making use of the frozen energy comparison \eqref{Dv_2-Dw_2} - coupled with the corresponding argument at scale $R$, Jensen's inequality and the scaling properties of $\bar h^{-1}$. Combining these estimates by the triangle inequality and enlarging the integration balls, we conclude that
\begin{align*}
\bigl[\bar h^{-1}(\lambda)\bigr]^{p-2}\sum_{k\in\{1,2\}}\mean{B_{\delta R/16}}|Du-Dw_k|\dx &\leq c\Bigl[\omega(R)+R^{\min\{s_1,\x\}}\Bigr]^{\frac{\textnormal{\texttt{E}}}{q-1}}\bigl[\bar h^{-1}(\lambda)\bigr]^{p-1}+c\,\frac{|\mu|(B_R)}{R^{n-1}}\\
&\leq c\Bigl[\omega(R)+R^{\min\{s_1,\x\}}\Bigr]^{\frac{\textnormal{\texttt{E}}}{q-1}}\lambda+c\,\frac{|\mu|(B_R)}{R^{n-1}},
\end{align*}
where the constants may also depend on $\delta$, and the last inequality follows from \eqref{above.below.inverseh}. For the third term, we first estimate the contributions of $Du-Dv_k$, and we do it by using \eqref{rough.compa.solopartep} with exponent  $\es=(q-1)/(p-1)$, which is admissible by \eqref{bound.trivial}: this estimate and \eqref{non.degenerate.Du} give 
\[
\mean{B_{R_k/2}}|Du-Dv_k|^{q-1}\dx\leq c\,\lambda^{\frac{q-1}{p-1}}.
\] 
Combining this with \eqref{Du-Dv_2.adhoc}, used with $\ell=1$, and the H\"older continuity of $a$, we obtain, for $k\in\{1,2\}$,
\begin{align*}
\mean{B_{R_k/2}}a(x)|Du-Dv_k|^{q-1}\dx&\leq a_i(B_R)\mean{B_{R_k/2}}|Du-Dv_k|^{q-1}\dx+cR^\alpha\mean{B_{R_k/2}}|Du-Dv_k|^{q-1}\dx\\
&\leq c\,\frac{|\mu|(B_R)}{R^{n-1}}+cR^\alpha\lambda^{\frac{q-1}{p-1}}\\
&\leq c\,\frac{|\mu|(B_R)}{R^{n-1}}+c\,C_{\rm ref}^{\frac{q-p}{p-1}}R^{\alpha-\frac{(n-1)(q-p)}{p-1}}\lambda\\
&\leq c\,\frac{|\mu|(B_R)}{R^{n-1}}+cR^{s_1}\lambda,
\end{align*}
where $c$ depends only on $\data,\gdata,\delta$, and $\vartheta$. After recalling \eqref{R.lambda} and the choice of $R_{\rm cmp}$, we note that for $p<2\leq n$ the exponents satisfy $\alpha-(n-1)(q-p)/(p-1)=s_0+\x\geq s_0$, whereas the bound $R\leq R_{\rm cmp}$ guarantees that $C_{\rm ref}^{(q-p)/(p-1)}R^{s_0-s_1}\leq1$. The contributions of $Dv_k-Dw_k$ are controlled by \eqref{compa.Dw.Lh} in the $p$-phase and by \eqref{compa.Dw.Lh.pq} with $\el=1$, together with coefficient comparability, in the $(p,q)$-phase. Applying these estimates at both scales yields
\[
\sum_{k\in\{1,2\}}\mean{B_{R_k/16}}a(x)|Dv_k-Dw_k|^{q-1}\dx\leq c\Bigl[\omega(R)+R^{\min\{s_1,\x\}}\Bigr]^{\textnormal{\texttt{E}}}\lambda.
\]
Consequently, by the triangle inequality and enlargement of the integration balls,
\begin{align*}
\sum_{k\in\{1,2\}}\mean{B_{\delta R/16}}&a(x)|Du-Dw_k|^{q-1}\dx\\ 
&\leq c\sum_{k\in\{1,2\}}\mean{B_{R_k/2}}a(x)|Du-Dv_k|^{q-1}\dx+c\sum_{k\in\{1,2\}}\mean{B_{R_k/16}}a(x)|Dv_k-Dw_k|^{q-1}\dx\\
&\leq c\,\frac{|\mu|(B_R)}{R^{n-1}}+c\Bigl[\omega(R)+R^{\min\{s_1,\x\}}\Bigr]^{\textnormal{\texttt{E}}}\lambda.
\end{align*}
Combining the estimates for the three terms, we conclude that
\[
\mean{B_{\delta R/16}}\bigl|A(x,Du)-A(x,Dw_2)\bigr|\dx\leq c\Bigl[\omega(R)+R^{\min\{s_1,\x\}}\Bigr]^{\frac{\textnormal{\texttt{E}}}{q-1}}\lambda+c\,\frac{|\mu|(B_R)}{R^{n-1}}.
\]

\subsubsection*{The case $p<2$ and $q\leq2$} In this last case, as we want to obtain a linear dependence on the measure density, we use \eqref{flux.difference.singular} from Lemma~\ref{lem:a-a}, the triangle inequality and \eqref{non.degeneracy.Dw1}: for $x\in B_{\delta R/8}$, we have
\begin{align*}
\bigl|A(x,Du)-A(x,Dw_2)\bigr|&\leq\bigl|A(x,Du)-A(x,Dw_1)\bigr|+\bigl|A(x,Dw_2)-A(x,Dw_1)\bigr|\\
&\leq c(p,q,L)\frac{h(x,Dw_1)}{|Dw_1|}\Bigl[|Du-Dw_1|+|Dw_1-Dw_2|\Bigr]\\
&\leq c(\data,\gdata,\vartheta)\frac{\lambda}{\bar h^{-1}(\lambda)}\Bigl[2|Du-Dw_1|+|Du-Dw_2|\Bigr]
\end{align*}
using also the triangle inequality. Averaging over $B_{\delta R/16}$ gives
\begin{equation}\label{pq<2}
\mean{B_{\delta R/16}}\bigl|A(x,Du)-A(x,Dw_2)\bigr|\dx\leq c\,\frac{\lambda}{\bar h^{-1}(\lambda)}\sum_{k\in\{1,2\}}\mean{B_{\delta R/16}}|Du-Dw_k|\dx.
\end{equation}
Now we show for each $k\in\{1,2\}$ that 
\begin{equation}\label{uv.L1.singular.normalized}
\mean{B_{R_k/2}}|Du-Dv_k|\dx\leq c\,\frac{\bar h^{-1}(\lambda)}{\lambda}\frac{|\mu|(B_R)}{R^{n-1}}
\end{equation}
for constants depending only on $\data,\gdata,\delta$ and $\vartheta$.  We start by applying \eqref{rough.compa.solopartep} at scale $R_k$ with $\es=1/(p-1)$, which is admissible by \eqref{2-1/n}, and using \eqref{non.degenerate.Du} to obtain 
\[
\mean{B_{R_k/2}}|Du-Dv_k|\dx\leq c\biggl[\frac{|\mu|(B_R)}{R^{n-1}}\biggr]^{\frac1{p-1}}+c\,\lambda^{\frac{2-p}{p-1}}\frac{|\mu|(B_R)}{R^{n-1}}\leq c\,\frac{\lambda^{\frac1{p-1}}}{\lambda}\frac{|\mu|(B_R)}{R^{n-1}}.
\]
In the $p$-phase, Remark~\ref{rem:p.phase.normalization} gives $\lambda^{1/(p-1)}\leq c\,\bar h^{-1}(\lambda)$, and \eqref{uv.L1.singular.normalized} is proven; in the $(p,q)$-phase we also apply \eqref{rough.compa.qcomponent} at scale $R_k$ with $\ti=1/(q-1)$, which is admissible since $1/(q-1)\leq1/(p-1)<n/(n-1)$. Using $a_i(B_R)\leq a(x_0)$ and $h_{x_0}(Du)\leq2h(x,Du)$ on $B_{R_k}$, we take the minimum of the two component estimates as in the proof of Lemma~\ref{lem:refined.1}, and \eqref{above.below.inverseh} gives \eqref{uv.L1.singular.normalized}.

\medskip

To prove the comparison between $Dv_k$ and $Dw_k$, we draw upon Jensen's inequality alongside \eqref{Dv_2-Dw_2} and its counterpart at scale $R$. Combined with the bound $d_1\leq2\lambda$, this yields
\begin{align*}
\mean{B_{R_k/8}}|Dv_k-Dw_k|\dx&\leq c\,\bar h^{-1}\biggl(\Bigl[\omega(R)+R^{\min\{s_1,\x\}}\Bigr]^{\textnormal{\texttt{E}}}\lambda\biggr)\\
&\leq c\Bigl[\omega(R)+R^{\min\{s_1,\x\}}\Bigr]^{\frac{\textnormal{\texttt{E}}}{q-1}}\bar h^{-1}(\lambda)\\
&\leq c\Bigl[\omega(R)+R^{\min\{s_1,\x\}}\Bigr]^{\textnormal{\texttt{E}}}\bar h^{-1}(\lambda),
\end{align*}
for $k\in\{1,2\}$, where we used \eqref{constant.hinverse}, $q\leq2$ and the fact $\omega(R)+R^{\min\{s_1,\x\}}\leq1$. Combining this with \eqref{uv.L1.singular.normalized}, using the triangle inequality and enlarging the integration balls, yields
\[
\sum_{k\in\{1,2\}}\mean{B_{\delta R/16}}|Du-Dw_k|\dx\leq c\,\frac{\bar h^{-1}(\lambda)}{\lambda}\biggl[\Bigl[\omega(R)+R^{\min\{s_1,\x\}}\Bigr]^{\textnormal{\texttt{E}}}\lambda+\frac{|\mu|(B_R)}{R^{n-1}}\biggr],
\]
with constants depending only on $\data,\gdata,\delta$ and $\vartheta$. Inserting this estimate into \eqref{pq<2} concludes the proof.

\end{proof}

\subsection{A linearized excess decay with error}

We now turn the comparison estimates above into a decay inequality for the $L^1$-excess of $A(x,Du)$. For $r\in(0,R_0]$, we also define
\begin{equation}\label{small.S}
\mathcal S(r):=\Bigl[\omega(r)+r^{\min\{s_1,\x\}}\Bigr]^{\frac{\textnormal{\texttt{E}}}{\max\{q-1,1\}}};
\end{equation}
from \eqref{Dini.powers} it follows that $\mathcal S$ is non-decreasing and satisfies the Dini condition.

\begin{lemma}[Qualitative excess smallness]\label{qualitative.excess.smallness}
Let $u$ and $\mu$ be as in Lemma~\ref{rough.linearized.comparison} and let $\varepsilon\in(0,1]$.  Then there exists a small constant $\delta_{\varepsilon,{\rm qlt}}\in(0,1/32)$, depending on $\data,\gdata$ and $\varepsilon$, such that for every $0<\delta\leq \delta_{\varepsilon,{\rm qlt}}$ there exist a threshold radius $R_{\varepsilon,{\rm qlt}}\equiv R_{\varepsilon,{\rm qlt}}(\data,\gdata,\omega(\,\cdot\,),\varepsilon,\delta)>0$ and a large constant $\mathcal B_{\varepsilon,{\rm qlt}}\equiv \mathcal B_{\varepsilon,{\rm qlt}}(\data,\gdata,\varepsilon,\delta)\geq1$ such that if $R\leq R_{\varepsilon,{\rm qlt}}$,
\begin{equation}\label{ass.qlt.exc.small}
\mean{B_R}h(x,Du)\dx\leq \lambda,\qquad \frac{|\mu|(B_R)}{R^{n-1}}\leq \frac\lambda{{\mathcal B}_{\varepsilon,{\rm qlt}}}
\end{equation}
for some $\lambda\geq1$, then
\begin{equation}\label{excess.linear1}
\mean{B_{\delta R}}\big|A(x,Du)-{\bigl(A(\,\cdot\,,Du)\bigr)}_{B_{\delta R}}\big|\dx\leq \varepsilon\lambda.
\end{equation}
\end{lemma}

\begin{proof}
Fix $\varepsilon\in(0,1]$ and let $\beta_0,c_{\beta_0}$ be the constants in the autonomous excess decay estimate for the frozen field $A_{i,1}$; the constants $\delta_{\varepsilon,{\rm qlt}}$, $R_{\varepsilon,{\rm qlt}}$ and $\mathcal B_{\varepsilon,{\rm qlt}}$ will be chosen below. For $0<\delta\leq1/32$, let $w_1\in v_1+W^{1,H_{i,\smallone}(\cdot)}_0(B_{R/8})$ be the $A_{i,1}$-harmonic comparison map defined after \eqref{two.different.cases}.

\medskip

We first estimate the error arising from freezing the field. By subadditivity and the preceding scalar comparison estimates, using coefficient comparability in the $(p,q)$-phase, we have
\[
\mean{B_{R/16}}h(x,Dw_1)\dx\leq c\mean{B_{R/16}}h(x,Du)\dx+c\mean{B_{R/16}}h(x,Du-Dw_1)\dx\leq c(\data,\gdata)\lambda:
\]
we used \eqref{ass.qlt.exc.small}, Remark~\ref{rem.obvious.Young} and the fact that $\omega(R)+R^{\min\{s_1,\x\}}\leq1$. By the continuity assumption in \eqref{A.strong} and the inequality $\bar h\leq h(x,\,\cdot\,)$, together with the bound in \eqref{ass.qlt.exc.small} and the concavity of $\omega(\,\cdot\,)$ that implies $\omega(2R)\leq 2\omega(R)$, we now have
\begin{align*}
\mean{B_{R/16}}\bigl|A(x,Dw_1)-A_{i,1}(Dw_1)\bigr|\dx&\leq c\,\omega(2R)\mean{B_{R/16}}h(x,Dw_1)\dx+c\mean{B_{R/16}}[a(x)-a_i(B_R)]|Dw_1|^{q-1}\dx\\
&\leq c\,\omega(R)\lambda+c\mean{B_{R/16}}[a(x)-a_i(B_R)]|Dw_1|^{q-1}\dx. 
\end{align*}
To estimate the second term, in the $p$-phase, we use subadditivity and we enlarge the first integral, then we use the estimate for the $q$-term preceding \eqref{s1} together with \eqref{en.rh.1}, and \eqref{compa.Dw.Lh} to obtain
\begin{align*}
\mean{B_{R/16}}[a(x)-a_i(B_R)]|Dw_1|^{q-1}\dx&\leq c\mean{B_{R/8}}a(x)|Dv_1|^{q-1}\dx+c\mean{B_{R/16}}h(x,Dv_1-Dw_1)\dx\\
&\leq c\Bigl[R^{s_1}+\bigl[\omega(R)+R^{s_1}\bigr]^{\textnormal{\texttt{E}}}\Bigr]\lambda,
\end{align*}
recalling \eqref{ass.qlt.exc.small}. In the $(p,q)$-phase, the relative oscillation estimate for $a(\,\cdot\,)$ \eqref{relative.oscillation.a} instead yields
\[
\mean{B_{R/16}}[a(x)-a_i(B_R)]|Dw_1|^{q-1}\dx\leq cR^\x\mean{B_{R/16}}\bar h(Dw_1)\dx\leq cR^\x\lambda;
\]
consequently, in both phases, 
\begin{equation}\label{freezing.error.w1}
\mean{B_{R/16}}\bigl|A(x,Dw_1)-A_{i,1}(Dw_1)\bigr|\dx\leq c\Bigl[\omega(R)+R^{\min\{s_1,\x\}}\Bigr]^{\textnormal{\texttt{E}}}\lambda,
\end{equation}
provided $R$ is sufficiently small. We now apply the autonomous excess decay estimate to $A_{i,1}(Dw_1)$. By the elementary comparison inequality for mean oscillations, we obtain
\begin{align*}
E(Du,B_{\delta R})&\leq2\mean{B_{\delta R}}\bigl|A(x,Du)-A_{i,1}(Dw_1)\bigr|\dx+\mean{B_{\delta R}}\big|A_{i,1}(Dw_1)-\bigl(A_{i,1}(Dw_1)\bigr)_{B_{\delta R}}\big|\dx\\
&\leq\frac{2}{(16\delta)^n}\mean{B_{R/16}}\bigl|A(x,Du)-A_{i,1}(Dw_1)\bigr|\dx\\
&\qquad\qquad\qquad+c_{\beta_0}(16\delta)^{\beta_0}\mean{B_{R/16}}\big|A_{i,1}(Dw_1)-\bigl(A_{i,1}(Dw_1)\bigr)_{B_{R/16}}\big|\dx\\
&\leq\tilde c\delta^{\beta_0}\lambda+\tilde c\,\delta^{-n}\mean{B_{R/16}}\bigl|A(x,Du)-A(x,Dw_1)\bigr|\dx+\tilde c\,\delta^{-n}\Bigl[\omega(2R)+R^{\min\{s_1,\x\}}\Bigr]^{\textnormal{\texttt{E}}}\lambda,
\end{align*}
where we used \eqref{freezing.error.w1} and the bound $|A_{i,1}(Dw_1)|\leq c\,\bar h(Dw_1)$, and $\tilde c$ depends only on  $\data$ and $\gdata$. Now we choose $\delta_{\varepsilon,{\rm qlt}}\leq1/32$ so that
\[
\tilde c\,\delta_{\varepsilon,{\rm qlt}}^{\beta_0}\leq\frac{\varepsilon}{3}.
\]
For every $0<\delta\leq\delta_{\varepsilon,{\rm qlt}}$, we apply Lemma~\ref{rough.linearized.comparison} with $\epsilon:=\varepsilon\delta^n/[3\tilde c]$ and then we set $\mathcal B_{\varepsilon,{\rm qlt}}:=\mathcal B_{\epsilon,{\rm rgh}}$ and finally choose $R_{\varepsilon,{\rm qlt}}\leq R_{\epsilon,{\rm rgh}}$ sufficiently small that
\[
\tilde c\,\delta^{-n}\Bigl[\omega(2R_{\varepsilon,{\rm qlt}})+R_{\varepsilon,{\rm qlt}}^{\min\{s_1,\x\}}\Bigr]^{\textnormal{\texttt{E}}}\leq\frac{\varepsilon}{3}.
\]
As each of the three terms in the preceding excess estimate is bounded by $\varepsilon\lambda/3$, \eqref{excess.linear1} follows.
\end{proof}

\begin{lemma}[Quantitative excess decay]\label{quantitative.excess.decay}
Let $u$ and $\mu$ be as in Lemma~\ref{rough.linearized.comparison}, and let $\varepsilon\in(0,1]$ and $\vartheta\in(0,1/2]$ be fixed. Then there exists a small constant $\delta_{\varepsilon,\vartheta,{\rm qnt}}$, depending only on $\data,\gdata,\varepsilon$ and $\vartheta$, such that for every $0<\delta\leq\delta_{\varepsilon,\vartheta,{\rm qnt}}$, there exist a threshold radius $R_{\varepsilon,\vartheta,{\rm qnt}}\in(0,1]$ and a constant $\mathcal B_{\varepsilon,\vartheta,{\rm qnt}}\geq1$, depending on $\data,\gdata,\varepsilon,\vartheta$ and $\delta$, with the radius also depending on $\omega(\,\cdot\,)$, such that if $R\leq R_{\varepsilon,\vartheta,{\rm qnt}}$,
\begin{equation}\label{non.degeneracy.quantitativedecay}
\vartheta\lambda\leq\mean{B_{\delta R}}h(x,Du)\dx,\qquad\mean{B_R}h(x,Du)\dx\leq\lambda,\qquad\frac{|\mu|(B_R)}{R^{n-1}}\leq\frac{\lambda}{\mathcal B_{\varepsilon,\vartheta,{\rm qnt}}} 
\end{equation}
for some $\lambda\geq1$, then
\begin{multline}\label{excess.linear3}
\mean{B_{\delta^2R}}\big|A(x,Du)-{\bigl(A(\,\cdot\,,Du)\bigr)}_{B_{\delta^2R}}\big|\dx\\
\leq\varepsilon\mean{B_{\delta R}}\big|A(x,Du)-{\bigl(A(\,\cdot\,,Du)\bigr)}_{B_{\delta R}}\big|\dx+c_{\rm lin}\mathcal S(R)\lambda+c_{\rm lin}\frac{|\mu|(B_R)}{R^{n-1}},
\end{multline}
where $c_{\rm lin}$ depends only on $\data,\gdata,\varepsilon,\vartheta$ and $\delta$.
\end{lemma}

\begin{proof}

We start by taking $\delta_{\varepsilon,\vartheta,{\rm qnt}}\leq\delta_{0,\vartheta}$, where $\delta_{0,\vartheta}$ is the constant appearing in Proposition~\ref{def.comparison.pphase}; for $0<\delta\leq\delta_{\varepsilon,\vartheta,{\rm qnt}}$, we let $w_2\in v_2+W^{1,H_{i,2}(\cdot)}_0(B_{\delta R/8})$ be the $A_{i,2}$-harmonic comparison map defined after \eqref{two.different.cases}. Let $\beta_0$ and $c_{\beta_0}$ be the exponent and the constant in the autonomous excess decay estimate for the frozen field $A_{i,2}$. Using standard properties of the excess, \eqref{nice.compa.pphase.final} and \eqref{decay.varphi.harmonic}, we obtain
\begin{align}\label{exc.delta2}
E(Du,B_{\delta^2R})&\leq2\mean{B_{\delta^2R}}\bigl|A(x,Du)-A_{i,2}(Dw_2)\bigr|\dx+\mean{B_{\delta^2R}}\bigl|A_{i,2}(Dw_2)-\bigl(A_{i,2}(Dw_2)\bigr)_{B_{\delta^2R}}\bigr|\dx\notag\\
&\leq\frac{2}{(16\delta)^n}\mean{B_{\delta R/16}}\bigl|A(x,Du)-A_{i,2}(Dw_2)\bigr|\dx\notag\\
&\hspace{3cm}+c_{\beta_0}(16\delta)^{\beta_0}\mean{B_{\delta R/16}}\bigl|A_{i,2}(Dw_2)-\bigl(A_{i,2}(Dw_2)\bigr)_{B_{\delta R/16}}\bigr|\dx\notag\\
&\leq c_*\delta^{\beta_0}E(Du,B_{\delta R})+c\,\delta^{-n}\mean{B_{\delta R/16}}\bigl|A(x,Du)-A(x,Dw_2)\bigr|\dx\notag\\
&\hspace{4cm}+c\,\delta^{-n}\mean{B_{\delta R/16}}\bigl|A(x,Dw_2)-A_{i,2}(Dw_2)\bigr|\dx
\end{align}
with $c,c_*$ depending on $\data,\gdata$ and independent of $\delta$. To estimate the freezing error, we first observe that
\[
\mean{B_{\delta R}}h(x,Du)\dx+\frac{|\mu|(B_{\delta R})}{(\delta R)^{n-1}}\leq\left(\delta^{-n}+\frac{\delta^{1-n}}{\mathcal B_{\varepsilon,\vartheta,{\rm qnt}}}\right)\lambda\leq 2\delta^{-n}\lambda.
\]
Repeating the argument leading to \eqref{freezing.error.w1} at scale $\delta R$, with the field frozen at $x_{m,2}$, where $a(x_{m,2})=a_i(B_{\delta R})$, therefore gives
\begin{align*}
\mean{B_{\delta R/16}}|A(x,Dw_2)-A_{i,2}(Dw_2)|\dx&\leq c(\data,\gdata,\delta)\Bigl[\omega(2\delta R)+(\delta R)^{\min\{s_1,\x\}}\Bigr]^{\textnormal{\texttt{E}}}\lambda\\
&\leq c(\data,\gdata,\delta)\Bigl[\omega(R)+R^{\min\{s_1,\x\}}\Bigr]^{\textnormal{\texttt{E}}}\lambda,
\end{align*}
where we used $2\delta\leq1$ and the monotonicity of $\omega(\,\cdot\,)$. We now choose
\[
\delta_{\varepsilon,\vartheta,{\rm qnt}}:=\min\biggl\{\delta_{0,\vartheta},\frac1{32},\Bigl(\frac{\varepsilon}{c_*}\Bigr)^{\frac1{\beta_\smallzero}}\biggr\}.
\]
Having fixed any $0<\delta\leq\delta_{\varepsilon,\vartheta,{\rm qnt}}$, we select $R_{\varepsilon,\vartheta,{\rm qnt}}$ sufficiently small and $\mathcal B_{\varepsilon,\vartheta,{\rm qnt}}$ sufficiently large to ensure the validity of Proposition~\ref{def.comparison.pphase} (with parameter $\vartheta$) together with the preceding freezing estimate. A further reduction of $R_{\varepsilon,\vartheta,{\rm qnt}}$ ensures $\omega(R)+R^{\min\{s_1,\x\}}\leq1$ for every $R\leq R_{\varepsilon,\vartheta,{\rm qnt}}$. In view of $\max\{q-1,1\}\geq1$, we then obtain
\[
\Bigl[\omega(R)+R^{\min\{s_1,\x\}}\Bigr]^{\textnormal{\texttt{E}}}\leq\mathcal S(R).
\]
Finally, \eqref{excess.linear3} follows upon inserting \eqref{nice.compa.pphase.final} and the freezing estimate into \eqref{exc.delta2}, absorbing all dependence on $\delta$ and $\vartheta$ into the constant $c_{\rm lin}$.

\end{proof}

\section{Proofs of the main results}\label{proofs}

We begin by showing how the excess estimates obtained in the preceding section extend to local SOLAs. 

\medskip

Let $u$ be a local SOLA  (if $p\leq n$) and let $\{u_j,\mu_j\}$ be an admissible approximating sequence. Lemma~\ref{SOLA.flux.stability} ensures that both $h(\,\cdot\,,Du_j)$ and $A(\,\cdot\,,Du_j)$ converge strongly in $L^1_{\rm loc}(\Omega)$ to $h(\,\cdot\,,Du)$ and $A(\,\cdot\,,Du)$, respectively, which in turn implies the convergence of their averages and flux excesses on every ball compactly contained in $\Omega$. For any given center $x_0\in\Omega_0$, Paragraph~\ref{choice.E} provides a common energy bound for the associated homogeneous comparison maps once finitely many terms are dropped. Exploiting Corollary~\ref{flux.decay} on a fixed reference ball alongside \eqref{SOLA.strong.flux} and \eqref{SOLA.uniform.meas}, we can select the reference radius $\rho_{\rm ref}$, the bound $c_{\rm ref}$, and the comparison radius $R_{\rm cmp}$ uniformly with respect to the sequence. Consequently, $\rho_{\rm ref}$ and $R_{\rm cmp}$ depend on $\data,\gdata$ and $\omega(\,\cdot\,)$, while $c_{\rm ref}$ is governed by $\data,\gdata$ and $\rho_{\rm ref}$.

\medskip

Setting up the quantitative excess estimate requires the parameters $\delta_{\varepsilon,\vartheta,{\rm qnt}}$, $R_{\varepsilon,\vartheta,{\rm qnt}}$, $\mathcal B_{\varepsilon,\vartheta,{\rm qnt}}$, and $c_{\rm lin}$, which we extract from the energy-solution version of Lemma~\ref{quantitative.excess.decay} evaluated for the decay factor $\varepsilon$ and the halved non-degeneracy parameter $\vartheta/2$. Preserving the index $\vartheta$ in our notation ensures these exact constants remain valid for local SOLAs without modifying their underlying dependencies.

\begin{remark}
{\rm 
If $u$ is a local SOLA that satisfies the hypotheses of Lemma~\ref{qualitative.excess.smallness} or Lemma~\ref{quantitative.excess.decay} on $B_R(x_0)$ at level $\lambda\geq1$, then the conclusion of the corresponding lemma (established, at this point, only for energy solutions) still holds as long as we start from a non-degeneracy parameter $\vartheta/2$ replacing $\vartheta$ in the quantitative case. 

If we indeed fix an integer $k\geq3$ and pick radii $\rho<R$ converging to $R$ from below such that $|\mu|(\partial B_\rho(x_0))=0$, the continuity of averages and the continuity from below of $|\mu|$ guarantee that the upper bounds on the average of $h(\,\cdot\,,Du)$ and the density hold with strict inequality on $B_\rho(x_0)$ with level $\lambda_k=(1+1/k)\lambda$, provided $\rho$ is sufficiently close to $R$. Moreover, in the quantitative case, the non-degeneracy assumption at scale $\delta R$ and the continuity of averages imply
\[
\mean{B_{\delta\rho}(x_0)}h(x,Du)\dx\geq\frac34\vartheta\lambda>\frac{\vartheta}{2}\biggl(1+\frac1k\biggr)\lambda=\frac{\vartheta}{2}\lambda_k
\]
again for $\rho$ close to $R$, so that the non-degeneracy condition is preserved at the new level $\lambda_k$, but with parameter $\vartheta/2$. These boundedness and non-degeneracy assumptions thus carry over to $u_j,\mu_j$ for large $j$ thanks to the strong convergence of the flux magnitude combined with \eqref{SOLA.uniform.meas}. Taking limits - first as $j\to\infty$, then as $\rho\nearrow R$, and ultimately as $k\to\infty$ - proves that the conclusions \eqref{excess.linear1} and \eqref{excess.linear3} hold for the local SOLA $u$ with the target decay factor $\varepsilon$ unmodified. The same approximation argument shows that we can employ \eqref{decay.h.Du} and \eqref{local.decay.ultimate} when $u$ is just a local SOLA. We stress that this argument does not require a uniform pointwise bound for the approximating fluxes.
}
\end{remark}
\medskip
Let $B_R(x_0)\subset\Omega$ with $R\leq R_0$. The radius $R_0$ and the contraction factor $\delta\in(0,1/32)$ will be fixed in each proof, with $R_0$ chosen no larger than the threshold radii required by the comparison estimates and Lemmas~\ref{qualitative.excess.smallness} and~\ref{quantitative.excess.decay}, for the corresponding choices of the parameters. We consider the concentric balls
\begin{equation}\label{R_j}
R_j=\delta^jR,\qquad \theta B_j=B_{\theta R_j}(x_0), \qquad j\in\N_0,\quad\theta>0;
\end{equation}
we also define
\begin{equation}\label{def.iteration}
A_j=\Big|\mean{B_j}A(x,Du)\dx\Big|,\quad E_j=\mean{B_j}\big|A(x,Du)-{\bigl(A(\,\cdot\,,Du)\bigr)}_{B_j}\big|\dx,\quad d_j=\frac{|\mu|(B_j)}{R_j^{n-1}}, 
\end{equation}
so that $E_j=E(Du,B_j)$ - we recall that we use the notation $E(\,\cdot\,,B_j)$ to indicate the excess of the flux, see \eqref{def.flux.excess}. Moreover, we stress that, for every fixed $\delta\in(0,1/2]$,
\begin{equation}\label{dyadic.S.smallness}
\sum_{j=0}^\infty\mathcal S(\delta^jR)\leq \mathcal S(R)+\frac1{\log(1/\delta)}\int_0^R\mathcal S(\varrho)\,\frac{d\varrho}{\varrho}\longrightarrow0\qquad\text{as }R\searrow0.
\end{equation} Lemma~\ref{quantitative.excess.decay} can be restated as follows:

\begin{corollary}[Quantitative excess decay]\label{quantitative.excess.decay.E_j}
Assume \eqref{a.holder} to \eqref{bounds} and \eqref{A.strong}, let $u$ be a local solution to \eqref{equation} in the sense of Definition~\ref{def.solution}, and consider the balls in \eqref{R_j} centered at $x_0\in\Omega_0$. Given $\varepsilon\in(0,1]$ and $\vartheta\in(0,1/2]$, choose $0<\delta\leq\delta_{\varepsilon,\vartheta,{\rm qnt}}$ and take the corresponding $R_{\varepsilon,\vartheta,{\rm qnt}}$, $\mathcal B_{\varepsilon,\vartheta,{\rm qnt}}$ and $c_{\rm lin}$, using the constants of Lemma~\ref{quantitative.excess.decay} for $\varepsilon$ and $\vartheta/2$ as specified above. If $R\leq R_{\varepsilon,\vartheta,{\rm qnt}}$ and
\[
\vartheta\lambda\leq\mean{B_{j+1}}h(x,Du)\dx,\qquad\mean{B_j}h(x,Du)\dx\leq\lambda,\qquad d_j\leq\frac{\lambda}{\mathcal B_{\varepsilon,\vartheta,{\rm qnt}}}
\]
for some $\lambda\geq1$ and $j\in\N_0$, then
\begin{equation}\label{excess.decay.E_j}
E_{j+2}\leq\varepsilon E_{j+1}+c_{\rm lin}\mathcal S(R_j)\lambda+c_{\rm lin}d_j,
\end{equation}
where $c_{\rm lin}$ depends only on $\data,\gdata,\varepsilon,\vartheta$ and $\delta$.
\end{corollary}

\subsection{Proof of Theorem~\ref{main.thm} and gradient continuity}\label{main.proof}
We fix a common Lebesgue point $x_0\in\Omega_0$ of $A(\,\cdot\,,Du)$ and $h(\,\cdot\,,Du)$: since both functions are locally integrable, such points form a set of full measure in $\Omega_0$ and, in particular,
\[
\lim_{r\searrow0}\mean{B_r(x_0)}A(x,Du)\dx=A(x_0,Du(x_0)),\qquad\lim_{r\searrow0}\mean{B_r(x_0)}h(x,Du)\dx=h(x_0,Du(x_0)).
\]
Let $R\leq R_0$ be a radius as in the statement of the theorem. If $\bI_1^{|\mu|}(x_0,R)=\infty$, there is nothing to prove. Otherwise, we set $r:=R/2$ and apply the notation in \eqref{R_j} and \eqref{def.iteration} with $r$ in place of $R$, so that $R_j=\delta^jr$ and $B_j=B_{R_j}(x_0)$. We define
\begin{equation}\label{definition.lambda}
\lambda:=H_1\mean{B_r(x_0)}h(x,Du)\dx+H_2\sum_{j=0}^\infty d_j+1,
\end{equation}
where $H_1,H_2\geq1$ and $\delta\in(0,1/32)$ are constants to be fixed below.

\subsubsection{Choice of the constants}

We now fix the constants appearing in Lemma~\ref{qualitative.excess.smallness}, Corollary~\ref{quantitative.excess.decay.E_j} and Lemma~\ref{integral.to.pointwise}. We first set
\[
Q:=\max\left\{1,\frac q\nu\right\},\qquad K:=16Q\Bigl(1+\frac Lp\Bigr)=16\max\Big\{1,\frac q\nu\Big\}\Bigl(1+\frac Lp\Bigr)
\]
and then
\[
\mathcal A_*:=\max\biggl\{\bar c,\frac{12KLq}{p\nu}\biggr\},\qquad \vartheta:=\frac{p\nu}{2^qqL\mathcal A_*},\qquad \delta_{\rm dens}:=\delta_0(\mathcal A_*),
\]
where $\bar c$ is the constant in \eqref{partial.above.pphase}, depending only on $\data$ and $\gdata$ and $\delta_0(\mathcal A_*)$ is given by Lemma~\ref{integral.to.pointwise}. Notice that these constants depend only on $\data,\gdata$. We then fix
\[
\delta=\min\biggl\{\delta_{1/2,\vartheta,{\rm qnt}},\delta_{\rm dens},\frac18\biggr\}.
\]
Here and in the following, $c_{\rm lin}$ denotes the constant in \eqref{excess.decay.E_j} corresponding to $\varepsilon=1/2$, the above choice of $\vartheta$ and $\delta$.

We next set $\varepsilon_*:=\delta^n/[6K]$ and choose an integer $\bar k\geq3$ sufficiently large that $\delta^{\bar k}\leq\delta_{\varepsilon_*,{\rm qlt}}$. Let $R_{\varepsilon_*,{\rm qlt}}$ and $\mathcal B_{\varepsilon_*,{\rm qlt}}$ denote the constants provided by Lemma~\ref{qualitative.excess.smallness} for this choice of $\varepsilon_*$ and contraction factor $\delta^{\bar k}$. Finally, we define
\[
\varepsilon_{\rm cmp}:=\min\biggl\{\frac{\nu}{12Kq}(16\delta^{\bar k})^n,\frac{\nu}{2^qq\mathcal A_*}(16\delta)^n,\frac12\biggr\}.
\]
We can now fix
\begin{equation}\label{Erre.zero}
\mathcal B:=\max\Big\{\mathcal B_{1/2,\vartheta,{\rm qnt}},\mathcal B_{\varepsilon_*,{\rm qlt}},\mathcal B_{\varepsilon_{\rm cmp}/2,{\rm rgh}}\Big\}.
\end{equation}
We further decrease $R_0$, if necessary, so that
\[
R_0\leq\min\Big\{R_{1/2,\vartheta,{\rm qnt}},R_{\varepsilon_*,{\rm qlt}},R_{\varepsilon_{\rm cmp}/2,{\rm rgh}}\Big\},
\]
and
\begin{equation}\label{smallness.S.iteration}
c_{\rm lin}\delta^{-n}\sum_{j=0}^\infty\mathcal S(\delta^jR_0)\leq\frac1{64Q}.
\end{equation}
By further decreasing $R_0$, we also ensure that the right-hand side of \eqref{freezing.error.w1} is bounded by $\varepsilon_{\rm cmp}\lambda/2$ whenever $R\leq R_0$ and 
\[
\mean{B_R}h(x,Du)\dx\leq\lambda,
\]
under the corresponding smallness assumption on the measure. We finally choose
\begin{equation}\label{choice.H1}
H_1:=K\Bigl(3+\frac{6L}{p}\Bigr)\delta^{-4n},\qquad H_2:=\max\Big\{\mathcal B,64Qc_{\rm lin}\delta^{-n}\Big\}.
\end{equation}
In particular, by the definition of $\lambda$, we have $d_j\leq\lambda/\mathcal B$ for every $j\in\N_0$. Lemma~\ref{qualitative.excess.smallness}, applied with contraction factor $\delta^{\bar k}$, therefore gives
\begin{equation}\label{small.excess.k}
\mean{B_i}h(x,Du)\dx\leq\lambda\qquad\Longrightarrow\qquad E_{i+\bar k}\leq\frac{\delta^n}{6K}\lambda.
\end{equation}

\subsubsection{An exit time argument}

Following \cite{KMI1}, for every $j\geq2$ we define the composite quantity
\begin{equation}\label{def.Cj}
C_j:=\sum_{\ell=j-2}^j\mean{B_\ell}h(x,Du)\dx+3\delta^{-n}E_j.
\end{equation}
By \eqref{a-h.weak} and the definition of the excess, we have
\[
E_j\leq\frac{2L}{p}\mean{B_j}h(x,Du)\dx,\qquad \mean{B_j}h(x,Du)\dx\leq\delta^{-jn}\mean{B_0}h(x,Du)\dx;
\]
hence, by \eqref{choice.H1}, we obtain
\[
\max\{C_2,C_3\}\leq\Bigl(3+\frac{6L}{p}\Bigr)\delta^{-4n}\mean{B_0}h(x,Du)\dx\leq\frac{\lambda}K.
\]
We now set $\mathcal J:=\{j\geq2:C_j\leq\lambda/K\}$. If $\mathcal J$ is infinite, there exists an increasing sequence $j_m\to\infty$ such that $C_{j_m}\leq\lambda/K$, and therefore, since $x_0$ is a Lebesgue point of $h(\,\cdot\,,Du)$,
\[
h(x_0,Du(x_0))=\lim_{m\to\infty}\mean{B_{j_m}}h(x,Du)\dx\leq\frac{\lambda}K,
\]
which already gives the desired estimate. We may thus assume that $\mathcal J$ is finite and define $j_e:=\max\mathcal J\geq3$. By construction,
\begin{equation}\label{exit.time.explicit}
C_{j_e}\leq\frac{\lambda}K,\qquad C_j>\frac{\lambda}K\qquad\text{for every }j>j_e,
\end{equation}
and in particular,
\begin{equation}\label{exit.time.initial}
\sum_{\ell=j_e-2}^{j_e}\mean{B_\ell}h(x,Du)\dx+3\delta^{-n}E_{j_e}\leq\frac{\lambda}K.
\end{equation}
Suppose that $i\geq j_e-2$ and 
\begin{equation}\label{nondeg.inductive}
\mean{B_i}h(x,Du)\dx\leq\lambda: 
\end{equation}
since $i+\bar k>j_e$, the exit-time condition and \eqref{small.excess.k} give
\begin{equation}\label{exit.small.balls}
\sum_{\ell=i+\bar k-2}^{i+\bar k}\mean{B_\ell}h(x,Du)\dx=C_{i+\bar k}-3\delta^{-n}E_{i+\bar k}>\frac{\lambda}{K}-\frac{\lambda}{2K}=\frac{\lambda}{2K},
\end{equation}
whence, for some $m\in\{\bar k-2,\bar k-1,\bar k\}$,
\begin{equation}\label{exit.lower.ball}
\mean{B_{i+m}}h(x,Du)\dx>\frac{\lambda}{6K}.
\end{equation}
We first establish the implication from \eqref{nondeg.inductive}, \eqref{exit.lower.ball} and $d_i\leq\lambda/\mathcal B$ to \eqref{exit.non.degeneracy.Du} for an energy solution $u$, with $\mu\in L^\infty_{\rm loc}(\Omega)$ if $p\leq n$. Fix the indices $i$ and $m$ as above and consider the $A_{i,1}$-harmonic comparison map $w_1$ associated with $B_i$. By the choices of $\varepsilon_{\rm cmp}$, $\mathcal B$ and $R_0$ made above, \eqref{est.qualitative.exc.small} and \eqref{freezing.error.w1} yield
\begin{equation}\label{exit.comparison.w1}
\mean{B_{R_i/16}}|A(x,Du)-A_{i,1}(Dw_1)|\dx\leq\varepsilon_{\rm cmp}\lambda.
\end{equation}
Since $\delta^m\leq\delta^{\bar k-2}\leq1/16$ and $\delta^m\geq\delta^{\bar k}$, using \eqref{a-h.weak} and \eqref{exit.lower.ball} we obtain
\begin{align*}
\frac{\nu}{6Kq}\lambda<\mean{B_{i+m}}|A(x,Du)|\dx&\leq\mean{B_{i+m}}|A_{i,1}(Dw_1)|\dx+(16\delta^m)^{-n}\varepsilon_{\rm cmp}\lambda\\
&\leq\frac Lp\mean{B_{i+m}}\bar h(|Dw_1|)\dx+\frac{\nu}{12Kq}\lambda;
\end{align*}
consequently,
\begin{equation}\label{exit.frozen.density}
\mean{B_{\delta^mR_i}}\bar h(|Dw_1|)\dx\geq\frac{p\nu}{12KLq}\lambda\geq\frac{\lambda}{\mathcal A_*}
\end{equation}
by the definition of $\mathcal A_*$. We also have the upper bound
\[
\mean{B_{R_i/8}}\bar h(|Dw_1|)\dx\leq\bar c\lambda\leq\mathcal A_*\lambda;
\]
and this follows from the local estimate \eqref{partial.above.pphase}, whose proof only uses the upper energy bound and the smallness of the measure density, taking into account that $\bar h\leq h(x,\,\cdot\,)$.

\medskip

We are now in a position to apply Lemma~\ref{integral.to.pointwise} to $w_1$ with $\varphi=H_{i,1}$, so that $\varphi'=\bar h$: \eqref{pointwise.bounds} in our case reads as
\begin{equation}\label{exit.pointwise.w1}
\frac{\lambda}{2^{q-1}\mathcal A_*}\leq\inf_{B_{i+1}}\bar h(|Dw_1|)\leq\sup_{B_{R_i/16}}\bar h(|Dw_1|)\leq c_{\sup}\mathcal A_*\lambda
\end{equation}
for constants $\mathcal A_*, c_{\sup}$ that depend only on $\data$ and $\gdata$ and are independent of $i$, $m$, $\bar k$ and $\delta$.

\medskip

Now we estimate from below
\begin{align}\label{exit.non.degeneracy.Du}
\mean{B_{i+1}}h(x,Du)\dx&\geq\frac pL\mean{B_{i+1}}|A(x,Du)|\dx\geq\frac pL\biggl[\mean{B_{i+1}}|A_{i,1}(Dw_1)|\dx-(16\delta)^{-n}\varepsilon_{\rm cmp}\lambda\biggr]\notag\\
&\geq\frac pL\biggl[\frac{\nu}{2^{q-1}q\mathcal A_*}-(16\delta)^{-n}\varepsilon_{\rm cmp}\biggr]\lambda\geq\frac{p\nu}{2^qq L\mathcal A_*}\lambda=:\vartheta\lambda
\end{align}
using \eqref{a-h.weak}, together with the last estimate and \eqref{exit.comparison.w1}; $\vartheta\in(0,1/2)$ depends only on $\data,\gdata$, and the estimate holds independently of $i$, $\delta$ and $\bar k$. This is the integral non-degeneracy bound we need to apply Corollary~\ref{quantitative.excess.decay.E_j}.

The same implication holds when $u$ is a local SOLA: let indeed $\{u_j,\mu_j\}$ be an admissible approximating sequence and fix the indices $i$ and $m$ as above. Since the inequality in \eqref{exit.lower.ball} is strict, we have
\[
\mean{B_{i+m}}h(x,Du)\dx>\frac{(1+1/k)\lambda}{6K}
\]
for every integer $k$ large enough. Choose a radius $\rho\leq R_i$ (that we will ultimately let converge to $R_i$) such that $|\mu|(\partial B_\rho(x_0))=0$: because averages are continuous and $|\mu|$ is continuous from below, \eqref{nondeg.inductive} and $d_i\leq\lambda/\mathcal B$ force the upper bounds to hold strictly at the level $(1+1/k)\lambda$, with the strict lower bound still satisfied on $B_{\delta^m\rho}(x_0)$. By Lemma~\ref{SOLA.flux.stability} and \eqref{SOLA.uniform.meas}, for all sufficiently large $j$ we have
\begin{multline*}
\mean{B_\rho(x_0)}h(x,Du_j)\dx<\Bigl(1+\frac1k\Bigr)\lambda,\qquad \frac{|\mu_j|(B_\rho(x_0))}{\rho^{n-1}}<\frac{(1+1/k)\lambda}{\mathcal B},\\
\mean{B_{\delta^m\rho}(x_0)}h(x,Du_j)\dx>\frac{(1+1/k)\lambda}{6K}
\end{multline*}
and the preceding energy-solution argument applies to $u_j$ on $B_\rho(x_0)$ at the level $(1+1/k)\lambda$, with the same constants, leading to
\[
\mean{B_{\delta\rho}(x_0)}h(x,Du_j)\dx\geq\vartheta\Bigl(1+\frac1k\Bigr)\lambda.
\]
We can pass to the limit as $j\to\infty$, then as $\rho\nearrow R_i$, and finally as $k\to\infty$ to obtain \eqref{exit.non.degeneracy.Du} for $u$.

\medskip

Hence, the induction below applies to both energy solutions and local SOLAs.

\subsubsection{An induction proof}

We now prove that
\[
\mean{B_j}h(x,Du)\dx\leq\lambda\qquad\text{for every }j\geq j_e-1.
\]
By \eqref{exit.time.initial}, this holds for $j=j_e-1,j_e$, since the corresponding averages are both bounded by $\lambda/K$ and $K\geq1$. Moreover, recalling \eqref{a-h.weak}, the same estimate gives
\begin{equation}\label{estimate.je}
A_{j_e}\leq\frac Lp\mean{B_{j_e}}h(x,Du)\dx\leq\frac{L}{pK}\lambda,\qquad E_{j_e}\leq\frac{\delta^n}{3K}\lambda.
\end{equation}
We proceed by induction. Fix $k\geq j_e$ and suppose that
\[
\mean{B_j}h(x,Du)\dx\leq\lambda\qquad\text{for every }j\in\{j_e-1,\ldots,k\}.
\]
By \eqref{exit.non.degeneracy.Du}, for every $i\in\{j_e-1,\ldots,k-1\}$ we also have
\[
\vartheta\lambda\leq\mean{B_{i+1}}h(x,Du)\dx.
\]
Since $d_i\leq\lambda/\mathcal B_{1/2,\vartheta,{\rm qnt}}$ by our choice of $H_2$, all the assumptions of Corollary~\ref{quantitative.excess.decay.E_j} are satisfied with $\varepsilon=1/2$, whence
\[
E_{i+2}\leq\frac12E_{i+1}+c_{\rm lin}\mathcal S(R_i)\lambda+c_{\rm lin}d_i\qquad\text{for }i\in\{j_e-1,\ldots,k-1\}.
\]
Summing these inequalities and absorbing the resulting sum of excesses, we obtain
\begin{equation}\label{sum.excess.iteration}
\sum_{j=j_e+1}^{k+1}E_j\leq\frac12\sum_{j=j_e}^{k}E_j+c_{\rm lin}\lambda\sum_{i=j_e-1}^{k-1}\mathcal S(R_i)+c_{\rm lin}\sum_{i=j_e-1}^{k-1}d_i\leq E_{j_e}+2c_{\rm lin}\lambda\sum_{i=0}^\infty\mathcal S(R_i)+2c_{\rm lin}\sum_{i=0}^\infty d_i.
\end{equation}

On the other hand, by the definition of $A_j$ and the elementary estimate
\[
|A_{j+1}-A_j|\leq\Big|\mean{B_{j+1}}A(x,Du)\dx-\mean{B_j}A(x,Du)\dx\Big|\leq\delta^{-n}E_j,
\]
telescopic summation, \eqref{sum.excess.iteration} and \eqref{estimate.je} give
\begin{align}\label{tel.sum}
A_{k+1}+E_{k+1}&\leq A_{j_e}+\delta^{-n}\sum_{j=j_e}^{k}E_j+E_{k+1}\leq A_{j_e}+\delta^{-n}\sum_{j=j_e}^{k+1}E_j\notag\\
&\leq A_{j_e}+2\delta^{-n}E_{j_e}+2c_{\rm lin}\delta^{-n}\lambda\sum_{i=0}^\infty\mathcal S(R_i)+2c_{\rm lin}\delta^{-n}\sum_{i=0}^\infty d_i\notag\\
&\leq\biggl[\frac{L/p+2/3}{K}+\frac1{16Q}\biggr]\lambda\leq\frac{\lambda}{8Q},
\end{align}
where the last inequality follows from \eqref{smallness.S.iteration}, the choice of $H_2$ and the definition $K=16Q(1+L/p)$. Recalling once again \eqref{a-h.weak}, we conclude that
\[
\mean{B_{k+1}}h(x,Du)\dx\leq\frac q\nu\bigl(A_{k+1}+E_{k+1}\bigr)\leq Q\bigl(A_{k+1}+E_{k+1}\bigr)\leq\frac{\lambda}{8}<\lambda.
\]
This completes the induction. In particular, letting $k\to\infty$ and recalling that $x_0$ is a Lebesgue point of $h(\,\cdot\,,Du)$, we obtain
\[
h(x_0,Du(x_0))\leq\frac{\lambda}{8}.
\]
Finally, using \eqref{dyadic.potential.bound} and comparing the averages over $B_r(x_0)$ and $B_R(x_0)$ in the definition of $\lambda$, we obtain
\[
h(x_0,Du(x_0))\leq\frac{\lambda}{8}\lesssim_n H_1\mean{B_R(x_0)}h(x,Du)\dx+H_2\bI_1^{|\mu|}(x_0,R)+1,
\]
and \eqref{potential.estimate.Du} finally follows, since $H_1$ and $H_2$ depend only on the prescribed data and $r=R/2$.

\medskip

\subsubsection{Continuity proofs}
From the excess estimates of the previous section, we now derive the qualitative oscillation bound in Proposition~\ref{qualitative.flux.oscillation}; Corollary~\ref{cor:flux.continuity} will immediately follow once we combine this with the locally uniform decay of the truncated potentials.

\begin{proof}[Proof of Proposition~\ref{qualitative.flux.oscillation}]
Let us first verify that every comparison and excess decay estimate holds uniformly for centers lying in $B_R(x_0)$. Fixing $R_\varepsilon\leq\min\{R_0/2,d_{\Omega_0}/16\}$, we note that since $x_0\in\Omega_0$ and $R\leq R_\varepsilon$,  $B_r(z)\subset B_{d_{\Omega_0}/8}(x_0)\subset\Omega_1$ whenever $0<r\leq\min\{1,d_{\Omega_0}/16\}$ and $z\in B_R(x_0)$. If $p\leq n$, applying \eqref{local.energy.SOLA} with $\es=b$ on $B_{d_{\Omega_0}/4}(x_0)$, as in the proof of \eqref{global.hi.energy}, gives
\[
\int_{B_r(z)}|Du|^{b(p-1)}\dx\leq\int_{B_{d_{\Omega_0}/8}(x_0)}|Du|^{b(p-1)}\dx\leq c(\data,\gdata),
\]
whereas for $p>n$ we have $\int_{B_r(z)}H(x,Du)\dx\leq\mathcal M$ by \eqref{energy.M}. In both cases, the argument of Paragraph~\ref{choice.E} gives an energy bound $E$ for the homogeneous comparison maps that is independent of $x_0$, $z$ and $r$. For a local SOLA, fix $x_0$ and discard finitely many terms of the approximating sequence as in Paragraph~\ref{choice.E}: the bounds on $B_{d_{\Omega_0}/4}(x_0)$ then give the same $E$ for all the remaining indices and all the balls $B_r(z)$ considered above.

The uniform energy bounds established above allow us to repeat the arguments of Proposition~\ref{decay.measure.data} and Corollary~\ref{flux.decay} with constants independent of the center $z\in B_R(x_0)$: if we (possibly) decrease $\rho_{\rm ref}$ and enlarge $c_{\rm ref}$, we may assume that $\rho_{\rm ref}\leq\min\{1,d_{\Omega_0}/16\}$ and
\begin{equation}\label{uniform.estimate}
t^{n-1}\mean{B_t(z)}h(x,Du)\dx+|\mu|(B_t(z))\leq c_{\rm ref}\qquad\text{for every }0<t\leq\rho_{\rm ref}, 
\end{equation}
where $c_{\rm ref}$ has to be larger than $c(\data,\gdata)\bigl[(1+\mathcal M)/\rho_{\rm ref}+1+|\mu|(\Omega)\bigr]$. The same constants can be used along the approximating sequence, so that the local estimate in \eqref{uniform.estimate} extends to local SOLAs as well. Since the bound in \eqref{uniform.estimate} is uniform in $z$, we can choose a single threshold radius $R_{\rm cmp}\leq\rho_{\rm ref}$ for all $z\in B_R(x_0)$, small enough to satisfy when $p\leq n$ the condition involving $C_{\rm ref}$ described at the beginning of Paragraph~\ref{Other.comparison}. We then take $R_\varepsilon\leq R_{\rm cmp}$.

\medskip

Fix now $\varepsilon\in(0,1]$ and put
\[
\vartheta:=\min\Bigl\{\frac12,\frac{\varepsilon p}{16L}\Bigr\}.
\]

Our aim is to invoke Lemma~\ref{quantitative.excess.decay} with decay factor $1/2$ with this choice of $\vartheta$, and in order to do this we first choose an admissible $\delta\in(0,1/32)$ such that $0<\delta\leq\delta_{1/2,\vartheta,{\rm qnt}}$. We then set $\eta:=\varepsilon\delta^n/256$ since we also want to rely on Lemma~\ref{qualitative.excess.smallness} with smallness parameter $\eta$ (this will compensate for a geometric factor $\delta^{-n}$ later in \eqref{geom.compensation}) and radius ratio $0<\sigma\leq\min\{\delta,\delta_{\eta,{\rm qlt}}\}$. We have at hand now the qualitative estimate \eqref{excess.linear1} (that ensures uniform smallness of the excess) and the quantitative estimate \eqref{excess.decay.E_j}, whose geometric decay will be used whenever the non-degeneracy condition \eqref{non.degeneracy.quantitativedecay} holds. 

Now first we consider the case
\begin{equation}\label{small.potentials.continuity}
\bI_1^{|\mu|}(x,R)+\bI_1^{|\mu|}(y,R)\leq\kappa\lambda,
\end{equation}
where $\kappa>0$ will be chosen sufficiently small, depending only on the prescribed data and $\varepsilon$, and for each $z\in\{x,y\}$, we use the notation introduced in \eqref{R_j} and \eqref{def.iteration} with center $z$ and initial radius $\sigma R/2$, so that for $j\in\N_0$ we have $R_j=\delta^j\sigma R/2$ and $B_j=B_{R_j}(z)$. The elementary potential bound
\begin{equation}\label{density.smallness}
\frac{|\mu|(B_t(z))}{t^{n-1}}\leq \tilde c(n)\bI_1^{|\mu|}(z,R)\qquad\text{for }0<t\leq R/2 
\end{equation}
following from \eqref{dyadic.potential.bound} provides the required smallness of the measure density, provided $\kappa$ is sufficiently small. More precisely, we choose $\kappa$ such that 
\[
\tilde c(n)\kappa\leq\min\biggl\{\frac1{\mathcal B_{\eta,{\rm qlt}}},\frac1{\mathcal B_{1/2,\vartheta,{\rm qnt}}}\biggr\}
\]
and, with $R_{\eta,{\rm qlt}}$ and $R_{1/2,\vartheta,{\rm qnt}}$ denoting the radii provided by Lemmas~\ref{qualitative.excess.smallness} and~\ref{quantitative.excess.decay}, respectively, for the fixed choices of $\eta$, $\delta$ and $\sigma$, we decrease $R_\varepsilon$ so that $R_\varepsilon\leq \min\{R_{\eta,{\rm qlt}},R_{1/2,\vartheta,{\rm qnt}}\}$. For every $0<\rho\leq\sigma R/2$, we apply Lemma~\ref{qualitative.excess.smallness} with outer ball $B_{\rho/\sigma}(z)$ and radius ratio $\sigma$. By \eqref{Linfinity.bound}, the average of $h(\,\cdot\,,Du)$ over this ball is at most $\lambda$ and since $\rho/\sigma\leq R/2$, \eqref{density.smallness} and \eqref{small.potentials.continuity} bound the measure density on the same ball by $\tilde c(n)\kappa\lambda\leq\lambda/\mathcal B_{\eta,{\rm qlt}}$. Consequently, we have \eqref{excess.linear1} uniformly for $0<\rho\leq\sigma R/2$, that is
\begin{equation}\label{continuity.uniform.excess}
\mean{B_\rho(z)}\bigl|A(x,Du)-\bigl(A(\,\cdot\,,Du)\bigr)_{B_\rho(z)}\bigr|\dx\leq\eta\lambda\qquad\text{for every }0<\rho\leq\sigma R/2.
\end{equation}
In particular, taking $\rho=R_j=\delta^j\sigma R/2$, we obtain $E_j\leq\eta\lambda$ for every $j\geq0$. Combining this with \eqref{dyadic.potential.bound} and \eqref{dyadic.S.smallness}, and further decreasing $\kappa$ and $R_\varepsilon$ if necessary, we finally obtain
\begin{equation}\label{continuity.small.excess.and.errors}
E_j\leq\eta\lambda\quad\text{for every }j\geq0,\qquad c_{\rm lin}\sum_{j=0}^\infty\bigl[\mathcal S(R_j)\lambda+d_j\bigr]\leq\eta\lambda.
\end{equation}
Here $c_{\rm lin}$ is the constant supplied by the quantitative excess estimate with the choices made above. All these choices are uniform with respect to $z$ and $\lambda$. If $j\geq1$ and $A_j>\varepsilon\lambda/16$, then \eqref{a-h.weak} gives
\[
\mean{B_j}h(x,Du)\dx\geq\frac pL A_j>\frac{\varepsilon p}{16L}\lambda\geq\vartheta\lambda.
\]
Corollary~\ref{quantitative.excess.decay.E_j}, applied with index $j-1$ replacing $j$, yields
\begin{equation}\label{continuity.nondegenerate.decay}
E_{j+1}\leq\frac12 E_j+c_{\rm lin}\mathcal S(R_{j-1})\lambda+c_{\rm lin}d_{j-1}
\end{equation}
and for every $j\geq0$ we also have
\begin{equation}\label{continuity.consecutive.means}
\bigl|(A(\,\cdot\,,Du))_{B_{j+1}}-(A(\,\cdot\,,Du))_{B_j}\bigr|\leq\delta^{-n}E_j.
\end{equation}
Now we consider two integers $1\leq k<h$: if there are indices in $\mathcal J=\{j\in\N:k\leq j\leq h\}$ for which $A_j\leq\varepsilon\lambda/16$, we take the minimum $j_-$ and the maximum $j_+$ of them and the difference between the corresponding averages of $A(\,\cdot\,,Du)$ is at most $\varepsilon\lambda/8$ by the triangle inequality; if no such indices exist, we apply the forthcoming summation argument to the whole interval $\{k,\ldots,h\}$. We deal with the remaining indices observing that on any interval of indices for which $A_j>\varepsilon\lambda/16$ we can sum up \eqref{continuity.nondegenerate.decay}; reabsorbing and appealing to \eqref{continuity.small.excess.and.errors} bounds the sum of the excesses by $4\eta\lambda$, and in turn \eqref{continuity.consecutive.means} allows us to bound the variation of the averages along the interval by $4\delta^{-n}\eta\lambda$. This will be used to bound the sum of the differences of two consecutive averages between the indices $k$ and $j_-$, and $j_+$ and $h$, if such sets are not empty. If we split $\mathcal J$ into these three, possibly empty, intervals, merging these two arguments and taking into account that the single step leaving the last small average costs at most $\delta^{-n}\eta\lambda$, we have proven that
\begin{equation}\label{geom.compensation}
\bigl|(A(\,\cdot\,,Du))_{B_h}-(A(\,\cdot\,,Du))_{B_k}\bigr|\leq\frac{\varepsilon\lambda}{8}+(2\cdot4+1)\delta^{-n}\eta\lambda\leq\frac{\varepsilon\lambda}{4}. 
\end{equation}
Taking $k=1$ and letting $h\to\infty$, since $z$ is a Lebesgue point, we obtain
\begin{equation}\label{continuity.point.to.mean}
\bigl|A(z,Du(z))-(A(\,\cdot\,,Du))_{B_{\delta\sigma R/2}(z)}\bigr|\leq\frac{\varepsilon\lambda}{4},\qquad z\in\{x,y\}.
\end{equation}
Finally, the bound $|A(\,\cdot\,,Du)|\leq L \lambda/p$ and the elementary estimate for the measure of the symmetric difference of two balls give
\[
\bigl|(A(\,\cdot\,,Du))_{B_{\delta\sigma R/2}(x)}-(A(\,\cdot\,,Du))_{B_{\delta\sigma R/2}(y)}\bigr|\leq\frac{c(n)L}{p\delta\sigma}\frac{|x-y|}{R}\lambda\leq\frac{\varepsilon\lambda}{2},
\]
provided $\tau_\varepsilon$ is chosen sufficiently small. Together with \eqref{continuity.point.to.mean}, this proves that the left-hand side of \eqref{qualitative.flux.oscillation.estimate} is at most $\varepsilon\lambda$ under \eqref{small.potentials.continuity}. If \eqref{small.potentials.continuity} fails, then we can simply estimate
\[
\bigl|A(x,Du(x))-A(y,Du(y))\bigr|\leq\frac{2L}{p}\lambda\leq\frac{2L}{p\kappa}\bigl[\bI_1^{|\mu|}(x,R)+\bI_1^{|\mu|}(y,R)\bigr]
\]
and this completes the proof for an appropriate choice of $C_\varepsilon$.
\end{proof}

\begin{proof}[Proof of Corollary~\ref{cor:flux.continuity}]
Since $\Omega_0\Subset\Omega$ is arbitrary, Theorem~\ref{main.thm} and our assumption \eqref{uniform.vanishing.potential} imply that $h(\,\cdot\,,Du)$ is locally bounded in $\Omega$. Fix a compact set $K\Subset\Omega$ and two intermediate open sets $U,\Omega_0$ with $K\Subset U\Subset\Omega_0\Subset\Omega$, and take $\lambda=\|h(\,\cdot\,,Du)\|_{L^\infty(U)}+1$. Given any positive $\tilde \varepsilon$, first choose $\varepsilon$ so that $\varepsilon\lambda$ is smaller than $\tilde \varepsilon/2$, and this fixes $\tau_\varepsilon, C_\varepsilon$ and $R_\varepsilon$ in Proposition~\ref{qualitative.flux.oscillation}. Then we choose $R\leq R_\varepsilon$ so that the potential $\mathbf I_1^{|\mu|}(\cdot,R)$ is uniformly small on $K$, smaller than $\tilde\varepsilon/[4C_\varepsilon]$, and this is possible thanks to \eqref{uniform.vanishing.potential}; we also require $2R<{\rm dist}(K,\partial U)$, so that Proposition~\ref{qualitative.flux.oscillation} can be applied on $B_{2R}(x)\subset U$ for every $x\in K$. Consequently, $|A(x,Du(x))-A(y,Du(y))|$ is smaller than $\tilde \varepsilon$ whenever $x,y\in K$ are Lebesgue points sufficiently close (i.e., satisfying $|x-y|\leq\tau_\varepsilon R$). We have proven that the flux is locally uniformly continuous on its set of Lebesgue points and, as such, it admits a unique continuous representative in $\Omega$.

The continuity and coercivity of $A$, together with its strict monotonicity in the gradient variable, imply that the inverse of $z\mapsto A(x,z)$ depends continuously on both $x$ and the flux value. Hence $Du$ also admits a continuous representative in $\Omega$, and $u$ admits a representative in $C^1(\Omega)$. Finally, if $\mu=f\dx$ with $f\in L^{n,1}_{\rm loc}(\Omega)$, \eqref{uniform.vanishing.potential} is ensured by the endpoint Riesz potential estimate and the absolute continuity of the Lorentz norm, see for instance \cite{KMGuide}.

\end{proof}

\subsection{Proof of Theorem~\ref{main.thm.2}}
To prove Theorem~\ref{main.thm.2} we need to modify some of the comparison estimates obtained in the previous sections; the first estimate we need is the following one, an appropriate variation of Lemma~\ref{lem:rough.compa.solopartep}. As it involves the comparison map $v$, we state it as an {\em a priori} estimate for energy solutions $u\in W^{1,H(\cdot)}_{\rm loc}(\Omega)$, therefore assuming that $\mu\in L^\infty_{\rm loc}(\Omega)$ when $p\leq n$.

\begin{lemma}\label{u.comparison.L1}
Let $u$ and $\mu$ be as in Lemma~\ref{comparison.uv.weighted} and let $v\in u+W^{1,H(\cdot)}_0(B_\rho)$ be the homogeneous comparison map defined by \eqref{def.v} with $R=2\rho$, for some ball $B_\rho=B_\rho(x_0)$ with $x_0\in\Omega_0$ and $0<2\rho\leq\min\{1,d_{\Omega_0}/8\}$. If $p<2$, then for every $\varepsilon\in(0,1]$ we have
\begin{equation}\label{u.comparison.gradient.epsilon}
\mean{B_\rho}|Du-Dv|\dx\leq\varepsilon\mean{B_\rho}|Du|\dx+c\,\varepsilon^{-\frac{2-p}{p-1}}h_{x_0}^{-1}\biggl(\frac{|\mu|(B_\rho)}{\rho^{n-1}}\biggr);
\end{equation}
if $p\geq2$, then
\begin{equation}\label{u.comparison.gradient.nonsingular}
\mean{B_\rho}|Du-Dv|\dx\leq c\,h_{x_0}^{-1}\biggl(\frac{|\mu|(B_\rho)}{\rho^{n-1}}\biggr).
\end{equation}
In both cases, the constants depend only on $\data$ and $\gdata$.
\end{lemma}

\begin{proof}
Consider the phase conditions \eqref{p.phase} and \eqref{pq.phase} with $\x=0$ and $M=2$, with $\esp\in\{p,q\}$ ($\esp=q$ only when the $(p,q)$-phase holds), $a_p=1$, $a_q=a(x_0)$. Applying \eqref{compa.q.component.ell} with $R/2=\rho$, $\ti=1/(\esp-1)$ (this exponent is admissible also when $\esp<2$, thanks to \eqref{2-1/n}), $\bar a_p=1$, and $\bar a_q=a_i(B_\rho)$ and dividing through by $\bar a_\esp$, raising to the power $1/(\esp-1)$, and invoking \eqref{comparability.a} when $\esp=q$, we obtain for $\esp\geq2$ that
\begin{equation}\label{u.comparison.component.degenerate}
\mean{B_\rho}|Du-Dv|\dx\leq c\biggl[\frac{|\mu|(B_\rho)}{a_\esp\rho^{n-1}}\biggr]^{\frac1{\esp-1}}
\end{equation}
while when $\esp<2$, the same line of reasoning gives
\begin{align}\label{u.comparison.component.singular}
\mean{B_\rho}|Du-Dv|\dx&\leq c\biggl[\frac{|\mu|(B_\rho)}{a_\esp\rho^{n-1}}\biggr]^{\frac1{\esp-1}}+c\,\frac{|\mu|(B_\rho)}{a_\esp\rho^{n-1}}\biggl(\mean{B_\rho}|Du|\dx\biggr)^{2-\esp}\notag\\ 
&\leq\varepsilon\mean{B_\rho}|Du|\dx+c\varepsilon^{-\frac{2-\esp}{\esp-1}}\biggl[\frac{|\mu|(B_\rho)}{a_\esp\rho^{n-1}}\biggr]^{\frac1{\esp-1}}
\end{align}
for every $\varepsilon\in(0,1]$, also employing Young's inequality with a suitable parameter choice. 

\medskip

We split the analysis now into the two phases. If \eqref{p.phase} holds, then \eqref{sup-inf_withx} with $\x=0$ implies $a(x_0)\leq c\rho^\alpha$; hence, since $|\mu|(B_\rho)\leq|\mu|(\Omega)$,
\[
a(x_0)\biggl[\frac{|\mu|(B_\rho)}{\rho^{n-1}}\biggr]^{\frac{q-p}{p-1}}\leq c|\mu|(\Omega)^{\frac{q-p}{p-1}}\rho^{\alpha-\frac{(n-1)(q-p)}{p-1}}\leq c(\data,\gdata),
\]
where the final inequality follows from \eqref{bounds} since $\rho\leq 1$. The conclusion then follows from \eqref{u.comparison.component.degenerate}-\eqref{u.comparison.component.singular} with $\esp=p$, because the bound in the display above precisely ensures that
\[
h_{x_0}^{-1}\biggl(\frac{|\mu|(B_\rho)}{\rho^{n-1}}\biggr)\approx\biggl[\frac{|\mu|(B_\rho)}{\rho^{n-1}}\biggr]^{\frac1{p-1}}.
\]
Assume instead that \eqref{pq.phase} holds: the comparability of $a(x)$ and $a(x_0)$ allows us to draw upon both the $p$- and the $q$-component estimates. Recalling \eqref{above.below.inverseh}, we can select the smaller of the two bounds: for $p\geq2$, both components satisfy \eqref{u.comparison.component.degenerate}, yielding \eqref{u.comparison.gradient.nonsingular}; for $p<2$, we apply \eqref{u.comparison.component.singular} to the $p$-component and either \eqref{u.comparison.component.singular} or \eqref{u.comparison.component.degenerate} to the $q$-component, depending on whether $q<2$ or $q\geq2$. Because 
\[
\frac{1}{q-1}-1=\frac{2-q}{q-1}\leq\frac{2-p}{p-1}=\frac{1}{p-1}-1
\] 
whenever $p\leq q<2$ (with $\varepsilon\leq1$), the dependence on $\varepsilon$ across both cases remains controlled by $\varepsilon^{-(2-p)/(p-1)}$. Taking the smaller estimate proves \eqref{u.comparison.gradient.epsilon}, and the proof is complete even in this case.
\end{proof}

In order to derive an oscillation comparison estimate we now rely on Lemma~\ref{u.comparison.L1}. If $p<2$ a preliminary iteration is needed to absorb the gradient term, and for this reason we first focus on the singular case. We fix two intermediate radii $R/2\leq r<s\leq R$ and we let $v_s\in u+W^{1,H(\cdot)}_0(B_s)$ be the homogeneous comparison map defined in $B_s\equiv B_s(x_0)$. We recall that, by the argument in Paragraph~\ref{choice.E} (applied here with outer radius $2s$) the energy bound \eqref{bound.energy.mdata} holds for $v_s$ with a constant $E\equiv E(\data,\gdata)$ independent of $s\in[R/2,R]$. $R$ here must be sufficiently small so that $B_{2R}(x_0)$ is admissible, and we can therefore apply Corollary~\ref{Caccio.L1} uniformly in $s$ to the map $v_s$ with $k=(u)_{B_R}$ to find
\[
\int_{B_r}|Dv_s|\dx\leq\frac c{s-r}\int_{B_s}|v_s-(u)_{B_R}|\dx.
\]
Since $u-v_s\in W^{1,H(\cdot)}_0(B_s)$, Poincar\'e's inequality and the triangle inequality give
\[
\int_{B_r}|Du|\dx\leq\frac{cR}{s-r}\int_{B_s}|Du-Dv_s|\dx+\frac{c}{s-r}\int_{B_s}|u-(u)_{B_R}|\dx.
\]
Recalling that $R/2\leq r<s\leq R$, we infer
\[
\mean{B_r}|Du|\dx\leq \frac{c}{s-r}\mean{B_R}|u-(u)_{B_R}|\dx+c\frac{R}{s-r}\mean{B_s}|Du-Dv_s|\dx.
\]
We apply \eqref{u.comparison.gradient.epsilon} on $B_s$, choosing its parameter $\varepsilon$ as a sufficiently small multiple of $(s-r)/[2R]$. Since $s\geq R/2$, the monotonicity of $|\mu|(B_s)$ and the scaling properties of $h_{x_0}^{-1}$ give
\begin{equation}\label{u.gradient.intermediate}
\mean{B_r}|Du|\dx\leq\frac12 \mean{B_s}|Du|\dx+\frac{c}{s-r}\mean{B_R}|u-(u)_{B_R}|\dx+c\Bigl(\frac{R}{s-r}\Bigr)^{\frac1{p-1}} h_{x_0}^{-1}\biggl(\frac{|\mu|(B_R)}{R^{n-1}}\biggr)
\end{equation}
and we can, at this point, apply Giusti's Lemma~\ref{Giusti.fundamental} to the bounded function $\phi(r):=R\mean{B_r}|Du|\dx$ on $[R/2,R]$, with exponent $\beta=1/(p-1)$, obtaining
\begin{equation}\label{u.gradient.absorbed}
R\mean{B_{R/2}}|Du|\dx\leq c\mean{B_R}|u-(u)_{B_R}|\dx+cR\, h_{x_0}^{-1}\biggl(\frac{|\mu|(B_R)}{R^{n-1}}\biggr).
\end{equation}
We now let $v\in u+W^{1,H(\cdot)}_0(B_{R/2})$ be the homogeneous comparison map in $B_{R/2}$. By Poincar\'e's inequality and \eqref{u.comparison.gradient.epsilon}, for every $\eta\in(0,1]$ we have
\[
\mean{B_{R/2}}|u-v|\dx\leq cR\mean{B_{R/2}}|Du-Dv|\dx\leq c\,\eta R\mean{B_{R/2}}|Du|\dx+c\,\eta^{-\frac{2-p}{p-1}}R h_{x_0}^{-1}\biggl(\frac{|\mu|(B_R)}{R^{n-1}}\biggr).
\]
Combining this inequality with \eqref{u.gradient.absorbed} and choosing $\eta$ sufficiently small in terms of $\varepsilon\in(0,1)$, we conclude that
\begin{equation}\label{u.comparison.oscillation}
\mean{B_{R/2}}|u-v|\dx\leq\varepsilon\mean{B_R}|u-(u)_{B_R}|\dx+c_\varepsilon R\, h_{x_0}^{-1}\biggl(\frac{|\mu|(B_R)}{R^{n-1}}\biggr).
\end{equation}
When $p\geq2$, the same conclusion follows directly from Poincar\'e's inequality and \eqref{u.comparison.gradient.nonsingular}, without any preliminary iteration; thus \eqref{u.comparison.oscillation} holds throughout the range of the theorem.

\medskip

We now derive a decay estimate for the mean oscillation of $u$. Fix $\theta\in(0,1/32)$. By the triangle inequality and the oscillation decay for the homogeneous comparison map $v$ in $B_{R/2}$, we have
\begin{align*}
\mean{B_{\theta R}}|u-(u)_{B_{\theta R}}|\dx&\leq2\mean{B_{\theta R}}|u-v|\dx+\osc_{B_{\theta R}}v\\
&\leq c\,\theta^{-n}\mean{B_{R/2}}|u-v|\dx+c\,\theta^{\alpha_0}\mean{B_{R/4}}|v-(v)_{B_{R/4}}|\dx\\
&\leq c\,\theta^{\alpha_0}\mean{B_R}|u-(u)_{B_R}|\dx+c\,\theta^{-n}\mean{B_{R/2}}|u-v|\dx,
\end{align*}
where the last inequality follows from the triangle inequality and the inclusion chain $B_{R/4}\subset B_{R/2}\subset B_R$. Moreover, combining \eqref{oscillation.v} with the local boundedness estimate \eqref{sup.L1.preliminary} applied to $v-(v)_{B_{R/4}}$, we obtain
\[
\osc_{B_{\theta R}}v\leq c\,\theta^{\alpha_0}\osc_{B_{R/8}}v\leq c\,\theta^{\alpha_0}\mean{B_{R/4}}|v-(v)_{B_{R/4}}|\dx,\qquad 0<\theta\leq1/32.
\]
Using \eqref{u.comparison.oscillation} to bound the remaining average of $|u-v|$, we arrive at
\[
\mean{B_{\theta R}}|u-(u)_{B_{\theta R}}|\dx\leq c\bigl(\theta^{\alpha_0}+\theta^{-n}\varepsilon\bigr)\mean{B_R}|u-(u)_{B_R}|\dx+c_\varepsilon\theta^{-n}R\, h_{x_0}^{-1}\biggl(\frac{|\mu|(B_R)}{R^{n-1}}\biggr).
\]
We first choose $\theta$ sufficiently small and then $\varepsilon$ sufficiently small, both depending only on $\data,\gdata$, so that $c(\theta^{\alpha_0}+\theta^{-n}\varepsilon)\leq1/2$. We thus obtain
\begin{equation}\label{u.oscillation.decay}
\mean{B_{\theta R}}|u-(u)_{B_{\theta R}}|\dx\leq\frac12\mean{B_R}|u-(u)_{B_R}|\dx+c\,R\, h_{x_0}^{-1}\biggl(\frac{|\mu|(B_R)}{R^{n-1}}\biggr).
\end{equation}
Having established this estimate for energy solutions, if $u$ is a local SOLA we can apply the preceding argument only to an approximating sequence $u_j$ given by its definition. For each fixed $x_0\in\Omega_0$, Paragraph~\ref{choice.E} provides (after discarding finitely many terms) a bound $E\equiv E(\data,\gdata)$ for the $L^{b(p-1)}$ integrals of the homogeneous comparison maps $v_j$ associated with $u_j$, together with the required uniform local control of $|\mu_j|$. The conclusions of Paragraph~\ref{choice.E} apply to the comparison maps used above in $B_s(x_0)$, whose outer radii satisfy $2s\leq2R_0$, if we take $R_0\leq\min\{1/2,d_{\Omega_0}/16\}$. Consequently, the constants, the choices of $\alpha_0,\theta,\varepsilon$, and the radius $R_0$ itself can be chosen to depend only on $\data$ and $\gdata$ and \eqref{u.oscillation.decay} holds with $u_j,\mu_j$ in place of $u,\mu$ for every $0<R\leq R_0$, again with constants independent of $j$ and $x_0$. Accounting for the normalization by $|B_R|$ via the scaling properties of $h_{x_0}^{-1}$, Lemma~\ref{SOLA.oscillation.stability} yields \eqref{u.oscillation.decay} for the local SOLA - up to modifying $c$ by a factor depending only on $n,p$, and $q$ - and this enables a direct iteration on $u$.

\medskip

Setting $r_j=\theta^jR/2$ for $j\in\N_0$, \eqref{u.oscillation.decay} gives
\[
\mean{B_{r_{j+1}}}|u-(u)_{B_{r_{j+1}}}|\dx\leq\frac12\mean{B_{r_j}}|u-(u)_{B_{r_j}}|\dx+cr_j h_{x_0}^{-1}\biggl(\frac{|\mu|(B_{r_j})}{r_j^{n-1}}\biggr).
\]
Summing over $j=0,\ldots,m-1$ and absorbing the resulting half-sum on the right-hand side, we obtain
\begin{align*}
\sum_{j=0}^m\mean{B_{r_j}}|u-(u)_{B_{r_j}}|\dx&\leq2\mean{B_{R/2}}|u-(u)_{B_{R/2}}|\dx+c\sum_{j=0}^{m-1}r_j h_{x_0}^{-1}\biggl(\frac{|\mu|(B_{r_j})}{r_j^{n-1}}\biggr)\\
&\leq2\mean{B_{R/2}}|u-(u)_{B_{R/2}}|\dx+c\sum_{j=0}^{m-1}\int_{r_j}^{2r_j}h_{x_0}^{-1}\biggl(\frac{|\mu|(B_\rho)}{\rho^{n-1}}\biggr)\,d\rho 
\end{align*}
using the monotonicity of $|\mu|(B_\rho)$ and the scaling properties of $h_{x_0}^{-1}$. Now, since $\theta\leq1/32$ the intervals $[r_j,2r_j]$ are pairwise disjoint and contained in $(0,R]$, and this gives
\begin{equation}\label{u.oscillation.summation}
\sum_{j=0}^\infty\mean{B_{r_j}}|u-(u)_{B_{r_j}}|\dx\leq c\,\mean{B_R}|u-(u)_{B_R}|\dx+c\,\int_0^R h_{x_0}^{-1}\biggl(\frac{|\mu|(B_\rho)}{\rho^{n-1}}\biggr)\,d\rho
\end{equation}
enlarging the mean oscillation on $B_{R/2}$ and using the property of the excess in \eqref{exc.elementary}, by monotone convergence; we have assumed that the potential is finite, since otherwise the desired estimate is immediate. Moreover, by the inclusion $B_{r_{j+1}}\subset B_{r_j}$,
\[
\bigl|(u)_{B_{r_{j+1}}}-(u)_{B_{r_j}}\bigr|\leq\theta^{-n}\mean{B_{r_j}}|u-(u)_{B_{r_j}}|\dx.
\]
Since $x_0$ is a Lebesgue point of $u$, the averages $(u)_{B_{r_j}}$ converge to $u(x_0)$. Telescoping the preceding inequalities, using \eqref{u.oscillation.summation} and estimating the initial difference between $(u)_{B_{R/2}}$ and $(u)_{B_R}$, we finally obtain
\[
|u(x_0)-(u)_{B_R}|\leq c\mean{B_R}|u-(u)_{B_R}|\dx+c\int_0^R h_{x_0}^{-1}\biggl(\frac{|\mu|(B_\rho)}{\rho^{n-1}}\biggr)\,d\rho
\]
and the use of the triangle inequality now completes the proof of \eqref{potential.estimate.u}.

\medskip

\begin{proof}[Proof of Remark~\ref{rem:potential.u.sup}]
Since $a(\,\cdot\,)$ is  H\"older continuous, we have
\[
\frac{|\mu|(B_\rho(x_0))}{\rho^{n-1}}\leq h_\rho\biggl(h_{x_\smallzero}^{-1}\biggl(\frac{|\mu|(B_\rho(x_0))}{\rho^{n-1}}\biggr)\biggr)\leq\biggl[1+c\rho^{\alpha-\frac{(n-1)(q-p)}{p-1}}\biggr]\frac{|\mu|(B_\rho(x_0))}{\rho^{n-1}}
\]
where we used the bounds $h_{x_\smallzero}^{-1}(t)\leq t^{\frac1{p-1}}$ and $|\mu|(B_\rho(x_0))\leq |\mu|(\Omega)$. The exponent of $\rho$ is positive by \eqref{bounds} so the factor in brackets is uniformly bounded for $0<\rho\leq1$ and the equivalence follows by applying $h_\rho^{-1}$, with constants independent of $\rho$.
\end{proof}

\appendix\section{Technical results for frozen double phase functions}

In this appendix we collect three properties of the frozen double phase functions, used repeatedly in the comparison and decay estimates. We fix $1<p\leq q$ and $a_0\geq0$, and write
\[
H_0(t):=\frac{t^p}{p}+a_0\frac{t^q}{q},\qquad h_0(t):=t^{p-1}+a_0t^{q-1},\qquad t\geq0.
\]
We stress that all the conclusions below hold uniformly with respect to $a_0$, without any restriction on $q-p$.

\subsection{Convexity of powers of the frozen function}

\begin{lemma}\label{convexity.h.power}
Let $b\geq 1/(p-1)>0$; then the map $t\mapsto[h_0(t)]^b$ is convex on $[0,\infty)$. Consequently,
\begin{equation}\label{jensen.h.power}
h_0\biggl(\mean{B}f\dx\biggr)\leq\biggl(\mean{B}[h_0(f)]^b\dx\biggr)^{\frac1b}
\end{equation}
for every measurable set $B\subset\mathbb R^n$ with $0<|B|<\infty$ and every non-negative measurable function $f$ such that $[h_0(f)]^b\in L^1(B)$.
\end{lemma}

\begin{proof}
We first prove the convexity of
\[
F(t):={\bigl[h_0(t)\bigr]}^{\frac1{p-1}}=t{\bigl(1+a_0t^{q-p}\bigr)}^{\frac1{p-1}}.
\]
If $a_0=0$ or $p=q$, then $F$ is linear; otherwise, we differentiate $F$ twice, for $t>0$, to obtain
\[
F''(t)=\frac{a_0(q-p)}{p-1}\,t^{q-p-1}{\bigl(1+a_0t^{q-p}\bigr)}^{\frac1{p-1}-2}\biggl[q-p+1+\frac{q-1}{p-1}a_0t^{q-p}\biggr]\geq0.
\]
Since $F$ is continuous at zero, it is convex on $[0,\infty)$; moreover, also the function  $[h_0(t)]^b=[F(t)]^{b(p-1)}$ is convex, since $b(p-1)\geq1$.  Applying Jensen's inequality to the convex function $[h_0]^b$ (this is possible since $[h_0(f)]^b\geq f^{b(p-1)}$ and this implies $f\in L^1(B)$) and taking the power $1/b$ proves \eqref{jensen.h.power}.
\end{proof}

\subsection{Concavity of small powers of the inverse composition}

\begin{lemma}
The map
\begin{equation}\label{Psi.small.concave}
t\longmapsto \bigl[(H_0\circ h_0^{-1})(t)\bigr]^\sigma
\end{equation}
is concave on $[0,\infty)$ for every $0<\sigma\leq1/p'$.
\end{lemma}

\begin{proof}
Set $\tau:=h_0^{-1}(t)$ for $t>0$; the function
\[
\tau\mapsto\frac{\tau h_0'(\tau)}{h_0(\tau)}=p-1+(q-p)\frac{a_0\tau^{q-p}}{1+a_0\tau^{q-p}}
\]
is positive and non-decreasing, so its reciprocal has non-positive derivative. Since the first derivative of $t\mapsto[(H_0\circ h_0^{-1})(t)]^\sigma$ can be written as $\sigma\,\tau H_0(\tau)^{\sigma-1}\cdot h_0(\tau)/[\tau h_0'(\tau)]$, we obtain
\[
\frac{d^2}{dt^2}\bigl[(H_0\circ h_0^{-1})(t)\bigr]^\sigma=\frac{\sigma}{h_0'(\tau)}\biggl\{H_0(\tau)^{\sigma-2}\bigl[H_0(\tau)-(1-\sigma)\tau h_0(\tau)\bigr]\frac{h_0(\tau)}{\tau h_0'(\tau)}+\tau H_0(\tau)^{\sigma-1}\frac{d}{d\tau}\biggl[\frac{h_0(\tau)}{\tau h_0'(\tau)}\biggr]\biggr\}
\]
and this quantity is non-positive. Indeed, $1-\sigma\geq1/p$ and $\tau h_0(\tau)\geq pH_0(\tau)$ give $(1-\sigma)\tau h_0(\tau)\geq H_0(\tau)$, so the first term between parentheses is non-positive; the second is non-positive too for the consideration made above. This proves concavity on $(0,\infty)$ and hence, by continuity, on $[0,\infty)$.
\end{proof}

\subsection{An averaged estimate for the frozen flux}

\begin{lemma}
Let $1\leq\el\leq q'$, let $B\subset\mathbb R^n$ be measurable with $0<|B|<\infty$, and let $G:B\to\mathbb R^n$ be measurable with $H_0(G)\in L^1(B)$. Then $[h_0(G)]^\el\in L^1(B)$ and
\begin{equation}\label{averaged.Psi.inverse}
\biggl(\mean{B}[h_0(G)]^\el\dx\biggr)^{\frac1\el}\leq c\bigl(h_0\circ H_0^{-1}\bigr)\biggl(\mean{B}H_0(G)\dx\biggr),
\end{equation}
where $c=c(p,q)$ is independent of $a_0$ and $\el$.
\end{lemma}

\begin{proof}
We borrow an idea from \cite{OkUhlenbeck}: set
\[
\mathfrak E:=H_0^{-1}\biggl(\mean{B}H_0(G)\dx\biggr).
\]
The energy assumption ensures that $\mathfrak E<\infty$. If $\mathfrak E=0$, then $G=0$ almost everywhere in $B$ and the conclusion is immediate; we therefore assume $\mathfrak E>0$. Since $1\leq\el\leq q'$, the scaling properties of $h_0$ and the inequalities in \eqref{H.h.comparison} give, for every $t\geq\mathfrak E$,
\[
\frac{[h_0(t)]^\el}{H_0(t)}\leq q\,\frac{[h_0(t)]^{\el-1}}t\leq q\,\frac{[h_0(\mathfrak E)]^{\el-1}}{\mathfrak E}{\Bigl(\frac{t}{\mathfrak E}\Bigr)}^{(q-1)(\el-1)-1}\leq\frac qp\,\frac{[h_0(\mathfrak E)]^\el}{H_0(\mathfrak E)},
\]
where we used $(q-1)(\el-1)\leq1$. For $0\leq t\leq\mathfrak E$, the monotonicity of $h_0$ instead gives $h_0(t)\leq h_0(\mathfrak E)$. Combining the two cases, we obtain
\[
{\biggl[\frac{h_0(t)}{h_0(\mathfrak E)}\biggr]}^\el\leq1+\frac qp\,\frac{H_0(t)}{H_0(\mathfrak E)}\qquad\text{for every }t\geq0.
\]
Taking $t=|G|$, averaging over $B$ and using the definition of $\mathfrak E$, we conclude that
\[
\mean{B}[h_0(G)]^\el\dx\leq\Bigl(1+\frac qp\Bigr)[h_0(\mathfrak E)]^\el.
\]
This proves the asserted integrability and, after taking the power $1/\el$, \eqref{averaged.Psi.inverse}, with a constant depending only on $p$ and $q$.
\end{proof}

\vspace{5mm}

\noindent\textbf{Statement on the use of AI tools.}
During the finalization of this manuscript, the author used ChatGPT 5.6 Sol and 6 Astra to assist with language editing and \LaTeX{} formatting; to check notational clarity, internal consistency, and computations; and to suggest ways of avoiding repetition and simplifying the presentation, making this version of the paper significantly cleaner and less repetitive than any original draft. The vast majority of the manuscript had already been worked out and completed before LLM models became available; the mathematical arguments and results were developed and independently verified by the author, who bears sole responsibility for the content of the article, its conclusions, and any remaining errors.

\medskip

{\bf Acknowledgments}. The author is grateful to M. Colombo and G. Mingione for several useful discussions on the topic of this paper. He thanks G. Mingione for suggesting the investigation of the radial setting, M. Colombo for her contribution to the derivation of the corresponding radial fundamental solution, and both for their insights into the importance of the reverse-H\"older inequality \eqref{nice.reverse.2}. The author is supported by the European Research Council through the ERC Starting Grant project NEW, no.~101220121.

\end{document}